\documentclass[oneside,german,english,reqno]{amsart}
\usepackage{palatino}
\usepackage[T1]{fontenc}
\usepackage[latin1]{inputenc}
\usepackage{graphicx}
\usepackage{amssymb}
\IfFileExists{url.sty}{\usepackage{url}}
                      {\newcommand{\url}{\texttt}}

\makeatletter

\providecommand{\tabularnewline}{\\}

 \theoremstyle{plain}
 \theoremstyle{plain}    
 \newtheorem{thm}{Theorem} 
 \theoremstyle{plain}    
 \newtheorem{lem}{Lemma} 
 \theoremstyle{definition}
 \newtheorem*{defn*}{Definition}
 \theoremstyle{remark}
 \newtheorem*{rem*}{Remark}
 \theoremstyle{plain}    
 \newtheorem*{cor*}{Corollary}
 \theoremstyle{plain}    
 \newtheorem*{conjecture*}{Conjecture} 
 \newenvironment{lyxlist}[1]
   {\begin{list}{}
     {\settowidth{\labelwidth}{#1}
      \setlength{\leftmargin}{\labelwidth}
      \addtolength{\leftmargin}{\labelsep}
      }}
   {\end{list}}
 \theoremstyle{remark}    
 \newtheorem*{conclusion*}{Conclusion} 

\DeclareMathOperator{\diam}{diam}
\DeclareMathOperator{\supp}{supp}
\DeclareMathOperator{\len}{len}

\DeclareMathOperator{\LE}{LE}
\DeclareMathOperator{\QL}{QL}
\DeclareMathOperator{\cut}{cut}
\DeclareMathOperator{\Haus}{Haus}
\DeclareMathOperator{\sh}{\subset_{\Haus}}
\DeclareMathOperator{\capa}{Cap}
\renewcommand{\theenumi}{(\roman{enumi})}

\numberwithin{lem}{section}
\usepackage{graphicx}
\usepackage{color}
\usepackage{caption2}
\allowdisplaybreaks[1]

\ifx\hyperlink\undefined
\newcommand{\refs}[1]{\ref{#1}}
\newcommand{\toolong}[2]{\url{#1}}
\newcommand{\hypertarget}[1]{}
\newcommand{\hrlb}{}
\else
\newcommand{\refs}[1]{\ref*{#1}}
\newcommand{\toolong}[2]{\href{#1}{\nolinkurl{#2} (too long)}}
\newcommand{\hrlb}{\\}
\fi

\usepackage{babel}
\makeatother
\begin{document}
\theoremstyle{remark}
\newtheorem*{discussion}{Discussion}
\theoremstyle{plain}
\newtheorem{sublem}{Sublemma}[lem] 
\newtheorem*{dsh}{Theorem (Delmotte; Holopainen and Soardi)}
\newtheorem*{dmt}{Theorem (Delmotte)}
\newtheorem*{hsc}{Theorem (Hebisch and Saloff-Coste)}
\newtheorem*{varop}{Theorem (Varopoulos)}
\newtheorem*{bpp}{Theorem (Benjamini, Pemantle and Peres)}
\newcounter{const}
\newcounter{Const}
\def\newc#1{
\refstepcounter{const}
\label{#1}
}
\def\newC#1{
\refstepcounter{Const}
\label{#1}
}

\title{The scaling limit of loop-erased random walk in three dimensions}

\author{Gady Kozma}

\begin{abstract}
We show that the scaling limit exists and is invariant to dilations
and rotations. We give some tools that might be useful to show universality.
\end{abstract}
\maketitle
\tableofcontents{}

\section{Introduction}

Loop-erased random walk (LERW) is a model for a random simple path,
created by taking a simple random walk and, whenever the random walk
hits its path, removing the resulting loop and continuing. See section
\ref{sub:Loop-erasure} for a precise definition. It is strongly related
to the uniform spanning tree (UST), a random spanning tree of a graph
$G$ selected uniformly between all spanning trees of $G$: the path
in the UST between two points is distributed like a LERW between them
\cite{P91}, and further, the entire UST can be generated using repeated
use of LERW by Wilson's algorithm \cite{W96}. Both models, and the
connections between them are interesting on a general graph, but we
shall be most interested in lattices on $\mathbb{R}^{d}$ and open
subsets thereof, in which case these models arise naturally in statistical
mechanics in conjunction with the Potts model.

Of all the non-Gaussian models in statistical mechanics, LERW is probably
the most tractable. Above five dimensions, it can be analyzed using
the non-intersections of simple random walk directly \cite[chapter 7]{L91}
giving an easy proof that the scaling limit is Brownian motion. In
4 dimensions a logarithmic correction is required \cite{L95}, and
that too has been proved with no use of the difficult technique of
lace expansion (see \cite{BS85,HvdHS03} for lace expansion). Borrowing
a term from physics we might say that the \emph{upper critical dimension}
for this model is $4$. In $2$ dimensions, LERW is conformally invariant
in the limit as the lattice becomes finer and finer. This allowed
physicists to make precise conjectures about fractal dimensions, critical
exponents and winding numbers \cite{D92,M92}. Rigorously, $3$ different
approaches proved fruitful: the connection to random domino tilings
\cite{K00a,K00b}, the connection to SLE \cite{LSW04}, and the approach
we will pursue in this paper, \cite{K}. In fact, SLE was discovered
\cite{S00} in the context of LERW. 

Attempts to understand LERW in dimension $3$ focused mainly on the
number of steps it takes to reach the distance $r$. Physicists conjecture
that it is $\approx r^{\xi}$ and did numerical experiments to show
that $\xi=1.62\pm0.01$ \cite{GB90}. Rigorously the existence of
$\xi$ is not proved (so we must talk about an upper and lower exponents
$\underline{\xi}\leq\overline{\xi}$), and the best estimates known
are $1<\underline{\xi}\leq\overline{\xi}\leq5/3$ \cite{L99}. LERW
has no natural continuum equivalent in dimensions smaller than $4$
--- Brownian motion has a dense set of loops and therefore it is not
clear how to remove them in chronological order. In two dimensions
the scaling limit is radial SLE $2$, but it is not clear if this
can be interpreted as a {}``Brownian motion with loops removed''.
For example, take a coupling of Brownian motion and SLE 2 which is
the scaling limit of the couple $(R,\LE(R))$ --- it is not proved
that this limit exists, but for the purpose of the discussion we may
assume it does or alternatively take a subsequential limit. It is
not known whether in that coupling the SLE $2$ path is a function
of the Brownian path (I was informed of this question by O. Schramm).

In this paper we shall show that LERW has a scaling limit in three
dimensions. More precisely we shall show the following theorem:

\begin{thm}
\label{thm:Z3scal}Let $\mathcal{D}\subset\mathbb{R}^{3}$ be a polyhedron
and let $a\in\mathcal{D}$. Let $\mathbb{P}_{n}$ be the distribution
of the loop-erasure of a random walk on $\mathcal{D}\cap2^{-n}\mathbb{Z}^{3}$
starting from $a$ and stopped when hitting $\partial\mathcal{D}$.
Then $\mathbb{P}_{n}$ converge in the space $\mathcal{M}(\mathcal{H}(\overline{\mathcal{D}}))$.
\end{thm}
Here $\mathcal{H}(\mathcal{X})$ is the space of compact subsets of
$\mathcal{X}$ with the Hausdorff metric, and $\mathcal{M}(\mathcal{X})$
is the space of measures on $\mathcal{X}$ with the topology of weak
convergence (these, and a couple of other notations are explained
in section \ref{sub:Preliminaries}). In general, the choice of topologies
above is not canonical. For example, \cite{LSW04} shows the existence
of a scaling limit for LERW in two dimensions replacing $\mathcal{H}$
above with the somewhat stronger topology of {}``minimal distance
after optimal change of variables''. However, for our techniques
the Hausdorff metric is the natural choice. I believe that the tools
that will be developed here can be used for a number of convergence
questions for LERW (e.g.~the existence of $\xi$, the existence of
the scaling limit on more general domains, and in stronger topologies,
universality and so on). However, as this paper is long as it is,
I chose to show only the simplest consequences: that the limit exists
and is invariant to dilations and rotations.%
\begin{figure}
\includegraphics[%
  scale=0.6]{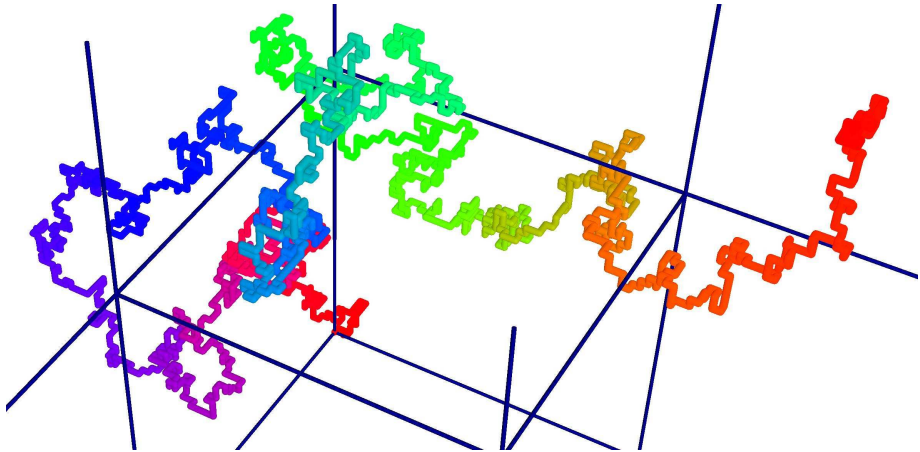}

\caption{A loop-erased random walk in three dimensions. Color indicates position
on the curve (red to blue to green and back to red).}
\end{figure}

Since we are interested in scaling limits it might be useful to review
quickly known results of this type. The archetypical example is of
course the Donsker invariance principle \cite[page 16]{RW94} stating
that the scaling limit of simple random walk is Brownian motion in
any dimension. As already remarked, in two dimensions the scaling
limit of LERW is radial SLE $2$, and a good deal of other discrete
models have been shown to converge to SLE: critical percolation on
the triangular lattice converges to chordal SLE $6$ \cite{S01,SW01},
the Peano curve of the UST converges to chordal SLE $8$ \cite{LSW04},
and the harmonic explorer converges to SLE $4$ \cite{SS}. The case
of the self-avoiding walk demonstrates the difficulties involved nicely:
it has been proved that if the limit exists and is conformally invariant,
it would be chordal SLE $8/3$ \cite{LSW}, but the existence of the
limit is still open. In high dimensions lace expansion has been used
to show that the scaling limit of the self avoiding walk is Brownian
motion \cite{BS85,HS92}, and that the scaling limit of percolation,
oriented percolation and lattice trees is integrated super-Brownian
excursion \cite{NY95,DS98,S99,vdHS03}. In intermediate dimensions
much less is known. The discrete Gaussian free field converges to
the continuum Gaussian free field and the Richardson model was shown
to have a limit shape from subadditivity arguments in any dimension
\cite{R73}, and I cannot resist citing the beautiful work on branched
polymers in dimension $3$ \cite{BI03a}. But these examples are the
exception, not the rule.

\subsection{Sketch of the proof}

The core of the argument is very similar to that of \cite{K}, so
let us recall the argumentation there. Let $R$ be a random walk on
a {}``three dimensional graph'' starting from $a$ and stopped on
the boundary of some domain $\mathcal{D}$. Let $B\subset\mathcal{D}$
be some (small) ball. We write \[
\LE(R)=\gamma_{1}\cup\gamma_{2}\cup\gamma_{3}\]
 where $\gamma_{1}$ is the portion of $\LE(R)$ until the first time
when $\LE(R)$ hits $B$. Notice that this is \textbf{not} the same
as the loop-erasure of a random walk stopped on $\partial B$! $\gamma_{2}$
is the portion of $\LE(R)$ until the last time when $\LE(R)$ is
inside $B$, and $\gamma_{3}$ is the reminder (the precise form of
this division is in the proof of lemma \ref{lem:0_to_Y}, page \pageref{lem:0_to_Y}).
Tracing the process of loop-erasure in $\mathcal{D}$ one sees that
$\gamma_{1}$ does not depend on anything that happens inside $B$:
when one knows all entry and exit points of $R$ from $B$, and all
the trajectories that $R$ does outside $B$, one can calculate $\gamma_{1}$.
In particular, if we compare random walks $R^{1}$ and $R^{2}$ on
graphs $G^{1}$ and $G^{2}$, where $G^{1}\setminus B=G^{2}\setminus B$
and inside $B$ we have some estimate of the sort\begin{equation}
p^{1}(v)\simeq p^{2}(v)\label{eq:pisketch}\end{equation}
where $p^{i}(v)$ is the probability of a random walk on $G^{i}$
to exit $B$ in a particular vertex $v$, then we should have that
\[
\gamma_{1}^{1}\simeq\gamma_{1}^{2}.\]
This argument and the precise meaning of {}``$\simeq$'' are contained
in lemma \ref{lem:xizLE}. To make this argument work for $\gamma_{3}$,
we have to use the symmetry of loop-erased random walk.

Now, if $\gamma_{2}$ is large then we are in the situation that was
coined in \cite{S00} a {}``quasi-loop'', namely two close points
on the path with a long way between them. In two dimensions it is
possible to show that $\LE(R)$ has few quasi-loops by tracing the
process that creates them and seeing that it necessitates that a random
walk that starts quite close to a loop-erased random walk will avoid
hitting it for a long while. This, however, contradicts the discrete
Beurling projection principle (see \cite{K87}) that states that a
random walk starting near any path has a high probability to intersect
it%
\footnote{Kesten's theorem is stronger, and claims that the minimum probability
is achieved, up to a constant, for a straight line, in which case
it can be estimated directly. However we will not need this level
of accuracy.%
}. See \cite[lemma 2.1]{S00} or \cite[lemma 18]{K}. Unfortunately
this argument no longer holds in three dimensions. Random walk starting,
say, in a distance $n$ from a straight line, has a reasonable probability
to intersect it only after extending to a distance of $n^{2}$, and
even after $n^{2}$ the probability to not intersect the line only
decreases logarithmically. In other words, in three dimensions not
all paths are {}``hittable'', and we have to show specifically that
loop-erased random walk is, using facts about its structure. This
will be done in chapter \ref{sec:Quasi-loops} and we shall show that
the probability that a random walk starting near a loop-erased random
walk will avoid hitting it decreases like a power law. This will allow
to repeat the above argumentation in three dimensions, show that there
are no quasi-loops and hence that LERW is similar on our $G^{1}$
and $G^{2}$.

The proof that loop-erased random walk is hittable is based on searching
for (local) cut times. A cut time for a random walk $R$ is a time
$t$ such that $R[0,t]\cap R[t+1,\infty[\;=\emptyset$. The number
of cut times is connected to the non-intersection exponent: by considering
the parts of $R$ up to $t$ and from $t$ on as two random walks,
and reversing the first part, we see that it is important to estimate
the probability that two random walks of length $t$ will not intersect.
This is $\approx t^{-2\xi}$ \cite{L96b} where $\xi$ is the famous
non-intersection exponent of Brownian motion. See section \ref{sub:Background-exponent}
for a description of this topic. Heuristically speaking, a set is
{}``hittable'' if its Hausdorff dimension is $>1$, and the set
of cut times has dimension $2-\xi$ so the argument terminates by
the well known fact that $\xi<1$ in three dimensions \cite{BL90b}.

Once the fact that $\LE(R)$ has no quasi-loops on either $G^{1}$
or $G^{2}$ is established, we get that it is similar on these two
graphs. Hence we can show that $\LE(R)$ is similar on $\mathbb{Z}^{3}$
and $2\mathbb{Z}^{3}$ by interpolating between them: dissecting into
a grid of cubes of intermediate size, and at each step change one
cube from $\mathbb{Z}^{3}$ to $2\mathbb{Z}^{3}$. Hence all of the
above discussion actually referred to graphs formed by cutting and
gluing together cubes of these two graphs.

Here are some corresponding reading recommendations:

\begin{itemize}
\item If you are familiar with \cite{K}, the part most interesting for
you would probably be the proof that there are no quasi loops. Read
the definition of a Euclidean net (section \ref{sub:DefinitionNet},
page \pageref{sub:DefinitionNet}), the definition of an isotropic
graph (section \ref{sub:DefinitionIsotropic}, page \pageref{sub:DefinitionIsotropic})
and the statements of theorems \ref{lem:escape} and \ref{thm:kof}
(pages \pageref{lem:escape} and \pageref{thm:kof} respectively)
and jump directly to chapter \ref{sec:Quasi-loops} (page \pageref{sec:Quasi-loops}).
\item If you are unfamiliar with \cite{K}, the part most interesting for
you would probably be the proof core sketched above. Read the definitions
of a Euclidean net and an isotropic graph as above; and the definition
of a quasi loop and the statement of theorem \ref{thm:QL} (page \pageref{thm:QL}).
Then jump directly to chapter \ref{sec:Isotropic-gluing} (page \pageref{sec:Isotropic-gluing})
or even to section \ref{sub:Definition-II} (page \pageref{sub:Definition-II}).
\item If you are the kind of person who prefers explicit examples to generalizations,
start with chapter \ref{sec:Examples} (page \pageref{sec:Examples})
and read a few examples of isotropic graphs and isotropic interpolations.
Then you can read the rest of the paper keeping in mind that an {}``isotropic
graph'' is really cubes of $\mathbb{Z}^{3}$ and $2\mathbb{Z}^{3}$
(and other variations) cut and sewn together so that random walk would
behave like Brownian motion. Section \ref{sub:Invariance} also contains
the proof of the invariance of the scaling limit to dilations and
rotations.
\item Chapters \ref{sec:euclidean-nets} and \ref{sec:Brownian-graphs}
are mostly recommended for students and non-special\-ists. Chapter
\ref{sec:euclidean-nets} consists mostly of citing well known connections
between rough isometries, heat kernel decay, Harnack inequality and
similar topics. In chapter \ref{sec:Brownian-graphs} we are forced
to replicate the results of Lawler \cite{L96b} in our settings. Roughly
we show that a relatively simple estimate of hitting probabilities
allow to couple random walk and Brownian motion and then Lawler's
argument goes through almost unchanged, giving that on the graphs
that interest us random walk has many cut times.
\end{itemize}
Finally I wish to point out how this paper improves over \cite{K}
in the two dimensional case. The use of a computer to calculate precise
estimates for the harmonic potential on $\mathbb{Z}^{2}$ \cite{KS04}
and on {}``hybrid graphs'' has been made completely unnecessary
by the use of electrical conductance techniques (see lemma \ref{lem:a}).
The use of {}``nice rectangles'' to ensure that hitting probabilities
are comparable was replaced by a multi-scale application of Harnack's
inequality together with a coupling argument (and in particular we
use spheres throughout rather than rectangles). See lemmas \ref{lem:sball}--\ref{lem:striag}
and \ref{lem:nodir}--\ref{lem:iiexit2}. Here the representation
is of comparable length, but is possibly less cumbersome. Finally,
the proof here of the final limit process is much shorter and simpler.

\subsection{About the settings}

The usual settings for these problems is that of a lattice in $\mathbb{R}^{d}$.
However, as explained above, the proof has to cut and saw together
different graphs, and even if these graphs were to be grids, the \emph{intermediate
objects} we must handle would not be. Hence we need to understand
random walks on graphs which are {}``similar'' to $\mathbb{Z}^{d}$.
It turns out that there are two important levels of similarity, which
correspond to {}``metric'' and {}``conformal'' properties. 

Much effort has gone into understanding what properties of random
walk are related to the metric structure only, or, more formally,
are satisfied by any graph \emph{roughly isometric} to $\mathbb{Z}^{d}$
(see definition and background on page \pageref{sub:Background-on-RI}).
However, one cannot expect LERW on a graph roughly isometric to $\mathbb{Z}^{d}$
to converge to a limit independent of the graph. Indeed, the scaling
limit of the random walk on the graph $\mathbb{Z}^{2}\times2\mathbb{Z}$
is not Brownian motion but a stretched version of it. These are not
identical --- indeed, even their hitting distribution on, say, a sphere,
differ, which implies that the scaling limit of LERW on $\mathbb{Z}^{3}$
and $\mathbb{Z}^{2}\times2\mathbb{Z}$ also differ. Hence we need
some condition to ensure that locally the graph is not stretched in
any direction. In other words, we need to preserve the conformal structure
of $\mathbb{R}^{d}$.

Properties related to the conformal structure are less well understood.
The {}``invariance principle'', that is the fact that random walk
converges to Brownian motion, which is a conformal property, has been
researched intensively, but it seems in different contexts than here.
Hence we will use a definition of \emph{isotropic graph} which is,
to the best of my knowledge, new. These graphs will satisfy the invariance
principle (this is more or less a tautology) and they preserve many
properties which are not preserved by the metric structure alone,
such as escape probabilities from a line, the non-intersection exponent,
and so on. 

Our definition of an isotropic graph (see chapter \ref{sec:Brownian-graphs})
is definitely not the most general imaginable. There are at least
two important examples which fall out of its scope. The first is a
conformal map of a grid --- for example the graph $\mathbb{Z}^{2}/\{(a,b)\sim(-a,-b)\}$
embedded into $\mathbb{C}$ via the map $(a,b)\mapsto(a+ib)^{2}$.
The second is random graphs, such as the Delaunay triangulation of
a Poisson process or the infinite cluster of super-critical percolation.
These graphs are not even roughly isometric to $\mathbb{R}^{2}$ and
yet are {}``isotropic'' in some heuristic sense. For example, the
percolation cluster is isotropic in the sense that it satisfies the
invariance principle, see \cite{DFGW89,SS04,BB,MP}. I conjecture
that the results here extend to these graphs, but will not complicate
the paper by considering them.

We shall prove theorem \ref{thm:Z3scal} (and other results) in both
the two and three dimensional cases. While the three dimensional case
is the more interesting one, the two dimensional proof is not quite
a subset of known results: it is proved for multiply connected domains
and for graphs more general than grids. However, at points the presentation
of specifically two dimensional issues will be sketchy.

\subsection{Acknowledgements}

Enormous thanks go to Itai Benjamini for many useful discussions,
encouragements, and for pointing out to me the relevance of rough
isometries and of the non-intersection exponent to this project. Many
thanks go to Gidi Amir and Omer Angel for useful discussions, in particular
with respect to counterexamples around lemmas \ref{lem:roughxi} and
\ref{lem:plane}. Lemma \ref{lem:omer} was discovered together with
Omer Angel.

This project was carried out while I was enjoying the hospitality
of, in chronological order, Université Bordeaux I, The Weizmann Institute
of Science (Charles Clore fund), Tel Aviv University and the Institute
of Advanced Study in Princeton (Oswald Veblen fund). I wish to thank
all these institutions, and especially A. Olevski\u\i{} from Tel
Aviv University who went to great efforts for me at unusual times.

\subsection{\label{sub:Preliminaries}Preliminaries}

A weighted graph is a couple $(G,\omega)$ where $G$ is a set and
$\omega:G\times G\to\left[0,\infty\right[$ such that $\omega(v,w)=\omega(w,v)$.
We shall often call $(G,\omega)$ simply $G$ and use $\omega$ only
in the places it is needed. For $v\in G$ the neighbors of $v$ are
the vertices $w$ such that $\omega(v,w)>0$. We denote by $v\sim w$
the neighborhood relation. We shall assume always that the number
of neighbors of every vertex is bounded and that the graph has \emph{bounded
weights} i.e.\[
\sup_{v,w}\omega(v,w)<\infty\quad\inf_{v\sim w}\omega(v,w)>0.\]
We do not assume $\omega(v,v)=0$ i.e.~we allow self loops.

A directed graph is a graph where $\omega(v,w)$ might be different
from $\omega(w,v)$. We will only use directed graphs once, in section
\ref{sub:Beurling}. Unless specifically marked {}``directed'' everything
below should be assumed to hold for undirected graphs only.

For a subset $X\subset G$ we denote by $\partial X$ the external
boundary, namely all vertices of $G\setminus X$ with a neighbor in
$X$. When this is not clear from the context, we shall denote $\partial_{G}X$
for the graph boundary and $\partial_{\textrm{cont}}$ for the boundary
of subsets of $\mathbb{R}^{d}$ in the usual sense. We write $\overline{X}=X\cup\partial X$.

A path in a graph is a function $\gamma$ from $\{1,\dotsc,n\}$ to
$G$ such that $\gamma(n)$ and $\gamma(n+1)$ are neighbors. $n$
is the length of the path, denoted by $\len\gamma$. If $\gamma_{1}$
and $\gamma_{2}$ are two paths and $\gamma_{1}(\len\gamma_{1})$
is a neighbor of $\gamma_{2}(1)$ (in which case we call $\gamma_{1}$
and $\gamma_{2}$ {}``concatenatable'') we shall define $\gamma_{1}\cup\gamma_{2}$
to be the path of length $\len\gamma_{1}+\len\gamma_{2}$ obtained
by concatenating them. It will be convenient to regard $\emptyset$
as a path of length $0$ and define $\gamma\cup\emptyset=\emptyset\cup\gamma=\gamma$.
The notation $\gamma[a,b]$ will be a short for the path of length
$b-a+1$ defined by $\gamma'(i)=\gamma(a+i-1)$, and also for the
set $\{\gamma(t):t\in[a,b]\}$ (there will rarely be a need to differentiate
between a path and its image). The same holds for other types of segments
(open, half-open). We say that $\gamma$ is {}``between'' $\gamma(1)$
and $\gamma(\len\gamma)$ and call $G$ connected if there exists
a path between any two vertices.

A random walk on a weighted graph is a process $R$ in discrete time
such that $R(t+1)$ depends only on $R(t)$ and \begin{equation}
\mathbb{P}(R(t+1)=w\,|\, R(t)=v)=\frac{\omega(v,w)}{\omega(v)}\quad\omega(v):=\sum_{w}\omega(v,w).\label{eq:omega}\end{equation}
$\mathbb{P}$ denoting the probability. We shall denote by $\mathbb{E}$
the expectation. When we shall need to specify the starting point,
we shall do so using $\mathbb{P}^{v}$ for the probability when $R(0)=v$,
and similarly $\mathbb{E}^{v}$. When we shall need to specify the
graph, we shall do so using $\mathbb{P}_{G}^{v}$ etc. Occasionally
(as in the statement of theorem \ref{thm:Z3scal}) we will have a
graph with an embedding in $\mathbb{R}^{d}$ and the {}``starting
point'' would be an $a\in\mathbb{R}^{d}$. In this case we mean by
$\mathbb{P}_{G}^{a}$ a random walk on $G$ starting from the point
of $G$ closest to $a$ (if more than one such point exist, choose
one, say by lexicographic order).

For a subset $X\subset G$ and a random walk $R$ we denote by $T(X)$
the hitting time of $X$ i.e.\[
T(X):=\min\{ t\geq1:R(t)\in X\}\]
or $\infty$ if the set is never hit. If $X=\{ x\}$ we shall write
$T(x)$ as short for $T(\{ x\})$. If $d$ is some metric on $G$
and if $X=\partial B(v,r)$ where $B(v,r)$ is a ball around $v$
of radius $r$ in the metric $d$, we will denote for short $T_{v,r}:=T(\partial B(v,r))$
(and assume the metric is clear from the context). Note that even
if we start from $v$, $T(v)$ is non-trivial since hitting is defined
only for $t\geq1$.

Sometimes we will have a few independent walks denoted by $R^{1},R^{2},\dotsc$.
In this case the corresponding stopping times will be denoted by $T^{i}(X)$
and $T_{v,r}^{i}$. Similarly we shall denote $\mathbb{P}^{1,v,2,w}$
when we want to denote that $R^{1}$ started from $v$ while $R^{2}$
started from $w$.

The strong Markov property says that for any stopping time $T$ the
random walk after $T$ behaves like a regular random walk. We shall
often use it, say for an event $E$ that depends only on what happened
after $T$, in the form $\mathbb{P}(E)=\mathbb{E}\mathbb{P}^{R(T)}(E)$.
Here $\mathbb{E}$ denotes expectation over the value of $R(T)$. 

Random walk is symmetric in the sense that the probabilities to traverse
a given path in one direction and in the opposite direction are equal
up to the ratio of $\omega$ at the beginning and end. In particular
we can sum over all paths of length $t$ and get\begin{equation}
\omega(v)\mathbb{P}^{v}(R(t)=w)=\omega(w)\mathbb{P}^{w}(R(t)=v)\quad\forall t,v,w.\label{eq:sym}\end{equation}
A similar argument shows that if $v,w\in A\subset G$ then\begin{equation}
\omega(v)\mathbb{P}^{v}(T(A)<\infty,R(T(A))=w)=\omega(w)\mathbb{P}^{w}(T(A)<\infty,R(T(A))=v).\label{eq:symT}\end{equation}

For a function $f:G\to\mathbb{R}$ (or to any linear space over $\mathbb{R}$)
we define the (discrete) Laplacian of $f$, $\Delta f$ by \[
(\Delta f)(v)=-f(v)+\sum_{w\sim v}\frac{\omega(v,w)}{\omega(v)}f(w)\]
A function $f$ such that $\Delta f$ is zero will be called (discretely)
harmonic. If $\Delta f$ is zero on a set $A\subset G$ we shall call
$f$ {}``harmonic on $A$''. Harmonic functions satisfy the maximum
principle, i.e.~a function harmonic on $A$ attains its maximum in
$\overline{A}$ on the boundary $\partial A$. Harmonic functions
are related to random walks by the following simple and well known
fact: if $f$ is harmonic on $A$ and $v\in A$ then\begin{equation}
f(v)=\mathbb{E}^{v}(f(R(T(\partial A)))).\label{eq:fEfRT}\end{equation}

For two sets $A$ and $B$ in a metric space $(X,d)$, we define their
distance by \[
d(A,B):=\inf_{a\in A,\, b\in B}d(a,b).\]
If $x\in X$ we write $d(x,A)$ as a short for $d(\{ x\},A)$. The
Hausdorff distance between $A$ and $B$ is defined by \[
d_{\Haus}(A,B):=\max(\sup_{a\in A}d(a,B),\sup_{b\in B}d(b,A)).\]
The diameter of a set is defined by\[
\diam A:=\sup_{a,b\in A}d(a,b).\]
If the metric space has an addition structure, we will use the notation
$A+B$ for the Minkowski sum of the sets $A$ and $B$ i.e.\[
A+B:=\{ a+b:a\in A,b\in B\}.\]
In particular, if $B$ is a ball centered at $0$ then \[
A+B(0,r)=\{ x:d(x,A)<r\}.\]
We will sometimes abuse notations by denoting the right hand side
by $A+B(r)$ even when the metric space has no addition structure.

A domain is a non-empty bounded open connected subset of $\mathbb{R}^{d}$.
A polyhedron is a domain whose boundary is composed of a collection
of non-degenerate linear polyhedra of dimension $d-1$. In particular,
we do not require that the boundary of the polyhedron is connected.
For simplicity, however, we will not allow slits.

By $C$ and $c$ we denote absolute constants which may be different
from place to place. $C$ will usually denote constants which are
{}``large enough'' and $c$ {}``small enough''. We shall number
($c_{1},c_{2},\dotsc$) only constants to which we will need to refer
to again. Sometimes we shall also write $C(\cdot)$ and $c(\cdot)$
for a constant which is not properly absolute --- it depends on some
parameters --- but is best thought of as absolute. This notation implicitly
means I cannot think of any applications where the parameters are
not themselves constants. Again, $C(G)$ could change from place to
place, and we shall number only those that we shall need to refer
to in the future. If, say, $C_{187}$ depends on some parameters we
shall only note this once and from that place on refer to it as simply
$C_{187}$, not $C_{187}(\alpha,\tau,\mathcal{H})$.

As usual we denote by $\left\lfloor x\right\rfloor $ the largest
integer $\leq x$ and by $\left\lceil x\right\rceil $ the smallest
integer $\geq x$.

\subsection{\label{sub:Loop-erasure}Loop erasure}

For a finite path $\gamma:\{1,\dotsc,n\}\rightarrow G$ we define
its loop erasure, $\LE(\gamma)$, which is a simple path in $G$,
by the consecutive removal of loops from $\gamma$. Formally, \begin{eqnarray}
\LE(\gamma)(1) & := & \gamma(1)\nonumber \\
\LE(\gamma)(i+1) & := & \gamma(j_{i}+1)\quad j_{i}:=\max\{ j:\gamma(j)=\LE(\gamma)(i)\}\quad.\label{eq:defLE}\end{eqnarray}
Which is defined for all $i$ such that $j_{i-1}<n$. 

\begin{lem}
\label{lem:condLE_sym}Let $b^{0},b^{1}\in B\subset G$. Let $R^{i}$
be a random walk starting at $b^{i}$, stopped at $B$ and conditioned
to hit $b^{1-i}$. Then $\LE(R^{0})$ has the same distribution as
the reversal of $\LE(R^{1})$.
\end{lem}
\noindent This is well known. See e.g.~\cite[lemma 2]{K}.

The following lemma was discovered with Omer Angel. To the best of
our knowledge, it has never been published before.

\begin{lem}
\label{lem:omer}Let $G$ be a weighted graph and let $v,w\in G$.
Let $R$ be a random walk on $G$ starting from $v$ and let $T_{n}$
be the $n$-th time $R$ is at $w$. Then\[
\LE(R[0,T_{1}])\sim\LE(R[0,T_{n}])\quad\forall n=2,3,\dotsc\]

\end{lem}
The notation $\sim$ here stands, as usual, for {}``having the same
distribution''.

\begin{proof}
Let $\gamma$ be any path starting from $v$ not containing $w$ and
let $k=1,\dotsc,n$. Denote $l=\len\gamma$. Define $X_{\gamma,k}$
to be the event that $\gamma=\LE(R[0,T_{n}])[1,l]$ and that $j_{l}\in[T_{n-k},T_{n-k+1}]$
where $j_{l}$ is from the definition of $\LE$ above and where we
consider $T_{0}$ to be $0$. Let $x$ be any neighbor of $\gamma(l)$.
We have that $x$ is the next element of $\LE(R[0,T_{n}])$ if and
only if $R(j_{l}+1)=x$. Denote this event by $N_{x}$.

Conditioning by $X_{\gamma,k}$ we get that $R[j_{l},T_{n}]$ is a
random walk on $G$ starting from $\gamma(l)$ and conditioned not
to hit $\gamma$ before hitting $w$ for $k$ times. Therefore \[
\mathbb{P}(N_{x}\,|\, X_{\gamma,k})=\mathbb{P}^{\gamma(l)}(R(1)=x\,|\, T_{k}<T(\gamma)).\]
The point of the lemma is that the right hand side is independent
of $k$ --- after $T_{1}$ it is no longer possible to know anything
about the value of $R(1)$. Therefore\[
\mathbb{P}(N_{x}\,|\, X_{\gamma,k})=\mathbb{P}(N_{x}\,|\, X_{\gamma,1})\]
and summing over $k$ we get (denoting $X_{\gamma}=\bigcup X_{\gamma,k}$)\[
\mathbb{P}(N_{x}\,|\, X_{\gamma})=\sum_{k=1}^{n}\frac{\mathbb{P}(N_{x}\cap X_{\gamma,k})}{\mathbb{P}(X_{\gamma})}=\sum_{k=1}^{n}\frac{\mathbb{P}(N_{x}|X_{\gamma,1})\mathbb{P}(X_{\gamma,k})}{\mathbb{P}(X_{\gamma})}=\mathbb{P}(N_{x}\,|\, X_{\gamma,1})\]
and this last term is equal to the probability that $\LE(R[0,T_{1}])$
will, if conditioned to start from $\gamma$, will have as its next
vertex $x$. Indeed, this is the well known {}``Laplacian random
walk'' representation of loop-erased random walk, see \cite{L87}.
\end{proof}
\begin{lem}
\label{lem:omerB}Lemma \ref{lem:omer} holds also when the random
walk is conditioned not to hit a given $B\subset G$, $w\not\in B$.
In a formula,\[
\LE(R[0,T_{1}])\,|\, R[1,T_{1}]\cap B=\emptyset\sim\LE(R[0,T_{n}])\,|\, R[1,T_{n}]\cap B=\emptyset\quad\forall n=2,3,\dotsc\]

\end{lem}
\noindent The proof is identical to that of the previous lemma, except
the random walk is conditioned to not hit $B\cup\gamma$ instead of
just $\gamma$.

\section{\label{sec:euclidean-nets}Euclidean nets}

\subsection{\label{sub:Background-on-RI}Background on rough isometries}

Let $X$ and $Y$ be two metric spaces. A function $f:X\to Y$ is
called a \textbf{rough morphism} if \[
d(f(x),f(y))\leq Cd(x,y)+C\]
for some $C$ which depends on $f$. A \textbf{rough identity} is
a function $f:X\to X$ satisfying \[
d(f(x),x)\leq C\]
for some $C$ which depends on $f$. Notice that $f$ need be neither
one-to-one nor onto! $X$ and $Y$ would be called \textbf{roughly
isometric} if there exist rough morphisms $f:X\to Y$ and $g:Y\to X$
such that both $f\circ g$ and $g\circ f$ are rough identities. In
this case we call both $f$ and $g$ \textbf{rough isometries}. The
term was introduced by Kanai in \cite{K85}, though in more restricted
settings it already appeared in \cite{G81}. There are various equivalent
definitions in the literature, but I prefer the above {}``categorical''
one. A rough isometry completely ignores all local structure, and
in fact $\mathbb{R}^{d}$ is roughly isometric to $\mathbb{Z}^{d}$
and more generally, any manifold is roughly isometric to any net inside
it.

To talk about rough isometry of graphs, we need to introduce a metric.
Let therefore $G$ be a weighted connected graph. Define \[
\delta(v,w):=\min_{\gamma:v\to w}\len\gamma\]
where the minimum is taken over all paths $\gamma$ from $v$ to $w$.
Clearly this makes $G$ a metric space.

The {}``Euclidean nets'' we are going to define in the next section
are graphs roughly isometric to $\mathbb{R}^{d}$. Whether properties
of random walks are preserved under rough isometries is in general
not obvious. In some cases (e.g.~transience) this requires an equivalent
representation as a geometric property. In others (e.g.~Harnack's
inequality) it is actually unknown. Let us therefore state here some
of connections between random walks and the geometry of the graph
that we will use.

\begin{defn*}
We say that a graph $G$ satisfies the volume doubling property if
there exists a constant $C$ such that for any $v\in G$ and any $r\geq1$,
$\omega(B(v,2r))\leq C\omega(B(v,\linebreak[1]r))$ where $\omega(A):=\sum_{x\in A}\omega(x)$,
$\omega$ from (\ref{eq:omega}).
\end{defn*}

\begin{defn*}
We say that a graph $G$ satisfies the weak Poincaré inequality if
there exists a constant $C$ such that for any function $f:G\to\mathbb{R}$,
any $v\in G$ and any integer $r$,\begin{gather*}
\sum_{w\in B(v,r)}\omega(w)\left|f(w)-\overline{f}\right|^{2}\leq Cr^{2}\sum_{w\sim x\in B(v,2r)}\omega(w,x)|f(w)-f(x)|^{2}\\
\overline{f}:=\frac{1}{\omega(B(v,r))}\sum_{x\in B(v,r)}\omega(x)f(x).\end{gather*}

\end{defn*}
The inequality is called {}``weak'' because the sum on the right
hand side is over a ball of radius $2r$. The regular Poincaré inequality
is defined with the sum over the ball of radius $r$. However, under
the assumption of the volume doubling property, these properties are
equivalent, see \cite[\S 5]{J86} (the settings there are a little
different but the proof carries through literally the same). In fact,
the equivalence (under volume doubling) of the weak Poincaré inequality
under different constants $>1$ replacing the {}``$2$'' in the
radius of the ball is much easier and the only thing we will use:
this easily implies that the combination of volume doubling and weak
Poincaré inequality is invariant to rough isometries.

Another common variation on this inequality is an $L^{1}$ version
i.e.~$\sum\omega|f(w)-\overline{f}|\leq Cr\sum\omega|f(w)-f(x)|.$
The $L^{1}$ version is stronger --- indeed, since the $L^{\infty}$
version $|f(w)-\overline{f}|\leq2r\max|f(w)-f(x)|$ is obviously always
true, the $L^{2}$ version follows from the $L^{1}$ version by interpolation.

\begin{defn*}
We say that a graph $G$ satisfies the elliptic Harnack inequality
if there exists a constant $C$ such that for any $v\in G$ and $r\geq1$
and any function $f$ harmonic and positive on $B(v,2r)$ one has\begin{equation}
\max\{ f(x):x\in B(v,r)\}\leq C\min\{ f(x):x\in B(v,r)\}.\label{eq:ellipt-Harnack}\end{equation}

\end{defn*}

\begin{defn*}
We say that a graph $G$ satisfies the parabolic Harnack inequality
if there exists a constant $C$ such that for any $v\in G$ and $r\geq1$
and any positive function $f$ on $B(v,2r)\times[0,4r^{2}]$ satisfying
\begin{equation}
f(\cdot,t+1)-f(\cdot,t)=\Delta f(\cdot,t)\label{eq:parab-Harnack}\end{equation}
one has\begin{multline}
\max\{ f(x,t):(x,t)\in B(v,r)\times[r^{2},2r^{2}]\}\leq\\
C\min\{ f(x,t):(x,t)\in B(v,r)\times[3r^{2},4r^{2}]\}.\label{eq:parab-res}\end{multline}

\end{defn*}
Clearly, the parabolic Harnack inequality is stronger than the elliptic
one. A difficulty in applying this fact is as follows: if the graph
is bipartite (say $\mathbb{Z}^{d}$) then the parabolic Harnack inequality
cannot hold --- for example, $f(x,t)=\mathbb{P}(R(t)=x)$ satisfies
(\ref{eq:parab-Harnack}) but the right hand side of (\ref{eq:parab-res})
is $0$ for any $r>1$, since $f(x,t)=0$ whenever $t+\sum x_{i}$
is odd. However, adding self loops will allow the graph to satisfy
the parabolic Harnack inequality without changing the set of harmonic
functions at all. After adding self-loops it is not at all easy to
construct examples of graphs satisfying the elliptic Harnack inequality
without satisfying the parabolic one. See \cite{BB99,GSC05} for some
constructions (see also \cite{HS01}).

\begin{dmt}Let $G$ be an infinite connected graph and assume $\omega(v,v)>c$
for all $v\in G$. Then the following are equivalent

\begin{enumerate}
\item $G$ satisfies the volume doubling property and the weak Poincaré
inequality.
\item $G$ satisfies the parabolic Harnack inequality.
\item The random walk on $G$ satisfies upper and lower Gaussian estimates,
namely\[
\frac{c}{\omega(B(v,\sqrt{t}))}e^{-C\delta(v,w)^{2}/t}\leq\mathbb{P}^{v}(R(t)=w)\leq\frac{C}{\omega(B(v,\sqrt{t}))}e^{-c\delta(v,w)^{2}/t}\]
for all $\delta(v,w)\leq t$.
\end{enumerate}
Further, for any two clauses, all constants in the first depend only
on the constants in the second.\end{dmt}

See \cite{D99}. One of the important consequences of this theorem
is that the parabolic Harnack inequality is invariant to rough isometries:
as already remarked, the combination of the volume doubling property
and the Poincaré inequality is invariant to rough isometries. For
the elliptic Harnack inequality, the question of its invariance is
still open.

\subsection{\label{sub:DefinitionNet}Definition}

A $d$-dimensional \textbf{Euclidean net} is a graph $(G,\omega)$
such that $G\subset\mathbb{R}^{d}$ and

\begin{enumerate}
\item $G$ has bounded weight;
\item The inclusion $i:G\to\mathbb{R}^{d}$ is a rough isomorphism between
$(G,\delta)$ and $\mathbb{R}^{d}$;
\item $\inf\{|v-w|:v\neq w\in G\}>0$.
\end{enumerate}
We shall mostly be interested in the $\mathbb{R}^{d}$ distance on
$G$, which we will denote by $d(v,w)$ or $|v-w|$. Likewise, the
notation $B(v,r)$ will relate to a ball in the $\mathbb{R}^{d}$
distance, while a ball in the metric $\delta$ will be denoted by
$B_{\delta}$.

\subsection{Harnack's inequality}

\begin{lem}
\label{lem:Harnack}A Euclidean net satisfies the elliptic Harnack
inequality (\ref{eq:ellipt-Harnack}) for $r$ sufficiently large.
\end{lem}
Two comments should be made on the formulation of the lemma. First
is that in the definition of the elliptic Harnack inequality (\ref{eq:ellipt-Harnack})
we mean balls in the $\mathbb{R}^{d}$ metric and not in the graph
metric. The second is about the constant in (\ref{eq:ellipt-Harnack}).
We implicitly assume that the constant $C=C(G)$ depends only on the
following parameters:\label{page:eucCstruct}

\begin{enumerate}
\item The bounds for $\omega$;
\item The constants of the rough isometries $i$ and $f$ between $G$ and
$\mathbb{R}^{d}$; 
\item The constants of the rough identities $i\circ f$ and $f\circ i$;
\item The lower bound for $|v-w|$;
\item $d$.
\end{enumerate}
We call the aggregation of these parameters the \textbf{Euclidean
net structure constants}. Whenever we use the notation $C(G)$ we
mean a constant depending only on these parameters. Similarly, constants
implicit in the notations $o,O$ and $\approx$ notation are not universal
but may depend on the structure constants of the Euclidean net. The
phrase {}``sufficiently large'' means {}``larger than a constant
depending on the Euclidean structure constants only''.

In chapter \ref{sec:Isotropic-gluing} we shall apply results obtained
up to that point to families of graphs with uniformly bounded structure
constants, hence it is important that $C$ does not depend on other
properties of $G$.

\begin{proof}
Construct an auxiliary graph $G^{*}$ with the same vertex set as
$G$ and with the weights defined by $\omega_{G^{*}}(v,v)=\omega_{G}(v,v)+\omega_{G}(v)$.
In other words, the random walk on $G^{*}$ is a random walk on $G$
with a probability of $\frac{1}{2}$ to stay at the same spot at each
step (additional to any such probability already existing for $G$).
The random walk on $G^{*}$ is sometimes called the lazy walk on $G$.
Clearly $\mathbb{Z}^{d}$ satisfies the volume doubling property and
it is easy to see that $\mathbb{Z}^{d}$ satisfies the weak Poincaré
inequality --- every group does, see e.g.~\cite[4.1.1]{PSC99} (the
other conditions of Delmotte's theorem are also easy to verify, if
you prefer). Since volume doubling and the weak Poincaré inequality
are preserved by rough isometries, $G^{*}$ satisfies them. Hence
by Delmotte's theorem it satisfies the parabolic Harnack inequality
with respect to the graph metric $\delta$. Hence it satisfies the
elliptic Harnack inequality, and since $G$ and $G^{*}$ have the
same harmonic functions, $G$ also satisfies the elliptic Harnack
inequality with respect to $\delta$.

Now, to prove Harnack's inequality for the $\mathbb{R}^{d}$ metric,
cover $B(v,r)$ by a constant number of balls $B_{\delta}(w_{i},cr)$
such that $B_{\delta}(w_{i},4cr)\subset B(v,2r)$. It is easy to see
that this can be done with the number of balls uniformly bounded.
From the above discussion we have, for every $f$ harmonic and positive
on $B(v,2r)$ that \[
\max\{ f:B_{\delta}(w_{i},2cr)\}\leq C\min\{ f:B_{\delta}(w_{i},2cr)\}\]
for every $i$. Therefore if we have that $B_{\delta}(w_{j},2cr)\cap B_{\delta}(w_{j+1},2cr)\neq\emptyset$
for $j=1,\dotsc,k$ then we get\[
\max\Big\{ f:\bigcup_{j}B_{\delta}(w_{j},2cr)\Big\}\leq C^{k}\min\Big\{ f:\bigcup_{j}B_{\delta}(w_{j},2cr)\Big\}.\]
Now, if $r$ is sufficiently large then the balls $B_{\delta}(w_{i},2cr)$
form a connected graph with respect to intersection, so we get the
required result.
\end{proof}
\begin{lem}
\label{lem:Harnack-general}Let $\mathcal{E},\mathcal{D}$ be domains
in $\mathbb{R}^{d}$ such that $\overline{\mathcal{E}}\subset\mathcal{D}$.
Then Harnack's inequality holds with respect to $\mathcal{E}$ and
$\mathcal{D}$ i.e. for any $v\in\mathbb{R}^{d}$, any $r>C(\mathcal{E},\mathcal{D},G)$
and any function $f$ positive and harmonic on $(r\mathcal{D}+v)\cap G$,
\[
\max\{ f:(r\mathcal{E}+v)\cap G\}\leq C(\mathcal{E},\mathcal{D},G)\min\{ f:(r\mathcal{E}+v)\cap G\}.\]
Further, if $\mathcal{K}$ be a family of $(\mathcal{E},\mathcal{D})$
with $\diam\mathcal{E}$ bounded above and $d(\mathcal{E},\partial\mathcal{D})$
bounded below, then all $C(\mathcal{E},\mathcal{D},G)$ are bounded
by some constant $C(\mathcal{K},G)$.
\end{lem}
\begin{proof}
Use the same covering trick as above.
\end{proof}
\begin{rem*}
The fact that a graph roughly isometric to $\mathbb{R}^{d}$ satisfies
the elliptic Harnack inequality was known before \cite{D99}. In the
setting of graphs, it was proved concurrently by Delmotte \cite{D97}
and Holopainen-Soardi \cite{HS97} (who proved it for the $p$-Laplacian
for any $p$). In the setting of manifolds this goes back to Kanai
\cite{K85}, who proved that a manifold roughly isometric to $\mathbb{R}^{d}$
satisfies the (continuous) Harnack inequality by showing that it follows
from a $d$-dimensional isoperimetric inequality.
\end{rem*}

\subsection{Green's function}

Let $H$ be any graph (possibly directed). Let $v,w\in H$ and $S\subset H$.
Then Green's function with respect to $S$ is defined by\[
G(v,w;S)=\sum_{t=0}^{\infty}\mathbb{P}^{v}(R(t)=w,\, R[0,t]\subset S)\]
or in other words, the expected number of visits to $w$ before leaving
$S$. If $S=H$ we shall omit it in the notation and write $G(v,w)$.
In general there is nothing forcing $G$ to be finite.

If $G$ is finite then it is zero outside $S$ and inside $S$ satisfies\begin{equation}
\Delta G(\cdot,w;S)=-\delta_{w}\label{eq:GreenLap}\end{equation}
i.e.~$G$ is harmonic on $S\setminus\{ w\}$ and $\Delta G(w,w)=-1$.
These conditions uniquely determine $G(v,w;S)$. The symmetry of random
walk (\ref{eq:sym}) translates to a symmetry of $G$ in the form\[
\omega(v)G(v,w;S)=\omega(w)G(w,v;S).\]

\begin{lem}
\label{lem:a}Let $H$ be a $d$-dimensional Euclidean net. Then
\begin{enumerate}
\item If $d=2$ then $H$ is recurrent and Green's function satisfies\begin{equation}
G(v,w;S)\leq C(H)\log\diam S.\label{eq:arecur}\end{equation}

\item If $d\geq3$ then $H$ is transient and Green's function satisfies\begin{equation}
\begin{split}G(v,w;S) & \leq C(H)|v-w|^{2-d}\quad\forall v\neq w.\\
G(v,v;S) & \leq C(H).\end{split}
\label{eq:atrans}\end{equation}
If $B(x,2r)\subset S$ and $r$ is sufficiently large then inside
$B(x,r)$ a lower bound also holds, \begin{equation}
\begin{split}G(v,w;S) & \approx|v-w|^{2-d}\quad\forall v\neq w\in B(x,r)\\
G(v,v;S) & \approx1.\end{split}
\label{eq:GvwS-lower-bound}\end{equation}
 
\end{enumerate}
\end{lem}
\begin{rem*}
For $d=2$ we have (from recurrence) that $G(v,w)=\infty$ for all
$v$ and $w$. The natural analog of $G(v,w)$ in this case is the
harmonic potential of $H$ defined by $a(v,w)=\lim_{r\to\infty}G(v,v;B(v,r))-G(v,w;B(v,r))$.
It is possible to show that for any two dimensional Euclidean net
$a(v,w)$ is well defined and $a(v,w)\approx\log|v-w|$, but we will
have no use for this fact.
\end{rem*}
\begin{proof}
We start with the case of $d\geq3$. Let $H^{*}$ be the lazy version
of $H$ as in the proof of lemma \ref{lem:Harnack}. By Delmotte's
theorem,\[
ct^{-d/2}e^{-C|v-w|^{2}/t}\leq\mathbb{P}_{H^{*}}^{v}(R(t)=w)\leq Ct^{-d/2}e^{-c|v-w|^{2}/t}.\]
Summing this over all $t$ we get \[
G_{H^{*}}(v,w)\approx|v-w|^{2-d}.\]
Now, $G_{H}=\frac{1}{2}G_{H^{*}}$ because one may couple the walks
on $H$ and $H^{*}$ so that each step of $H$ the walker on $H^{*}$
walks the same step and then waits for an expected time of $1$. Therefore
we get $G_{H}(v,w)\approx|v-w|^{2-d}$. We will not need the graph
$H^{*}$ again, so all Green functions henceforth are with respect
to $H$.

Now, $G(v,w;S)\leq G(v,w)$ gives us (\ref{eq:atrans}). To get the
lower bound under the assumption $B(x,2r)\subset S$ take $f=f_{S,w}$
to be a harmonic function on $S$ with $f(v)=G(v,w)$ for all $v\in\partial S$.
$G(\cdot,w)-f$ will satisfy (\ref{eq:GreenLap}) which defines $G(\cdot,w;S)$
so they are equal. By the maximum principle we get\[
f(v)\leq\max_{y\in\partial S}f(y)\leq Cr^{2-d}\quad\forall v\in S\]
so we get $G(v,w)\approx|v-w|^{2-d}$ inside a ball $B(v,\lambda r)$
for some constant $\lambda=\lambda(G)$ sufficiently small. Using
Harnack's inequality (lemma \ref{lem:Harnack-general}) for the domains
$B(0,1)\setminus B(0,\lambda)\subset B(0,2)\setminus B(0,\frac{1}{2}\lambda)$
proves (\ref{eq:GvwS-lower-bound}).

The two dimensional case follows from electrical resistance arguments.
See \cite{S94} for background on this topic. The maximum principle
shows that $G(v,w)\leq G(v,v)$ and the latter is equal to the resistance
between $v$ and $\partial S$. The electrical resistance is preserved
(up to a constant) by rough isometries, and so we get \[
G(v,v;S)\leq G(v,v;B(v,2\diam S))\approx G_{\mathbb{Z}^{2}}(0,0;B(0,2\diam S))\approx\log\diam S.\qedhere\]

\end{proof}
\begin{rem*}
The use of Delmotte's theorem here is somewhat of an overkill. The
fact that the a graph roughly isometric to $\mathbb{R}^{d}$ has a
$d$-dimensional heat kernel decay follows essentially from Varopoulos
\cite{V85}. To get an estimate for the Green function one can apply
e.g.~Hebisch and Saloff-Coste \cite{HS93} which gives a square exponential
decay upper bound. 
\end{rem*}
\begin{lem}
\label{lem:r2}Let $G$ be a Euclidean net and let $v\in G$ and $r>1$.
Then for some constant $\lambda(G)$ sufficiently large,\[
\mathbb{P}^{w}(R(\left\lfloor \lambda r^{2}\right\rfloor )\in B(v,r))\leq\frac{1}{2}\quad\forall w.\]

\end{lem}
\begin{proof}
Let $G^{*}$ be the lazy version of $G$, as in lemma \ref{lem:Harnack}.
Again we use Delmotte's theorem to show that for any $\lambda$, \[
\mathbb{P}_{G^{*}}^{w}(R(n)=x)\leq C\left\lfloor \lambda r^{2}\right\rfloor ^{-d/2}\quad\forall w,\forall n\geq\left\lfloor \lambda r^{2}\right\rfloor .\]
 We have that $\# B(v,r)\leq C(G)r^{d}$. Hence summing gives that
for $\lambda$ a constant sufficiently large \[
\mathbb{P}_{G^{*}}^{w}(R(n)\in B(v,r))\leq\frac{1}{2}\quad\forall w,\forall n\geq\left\lfloor \lambda r^{2}\right\rfloor .\]
But the coupling between the walk and the lazy walk shows that after
the walk did $n$ steps the lazy walk did at least $n$ steps. Hence\begin{align*}
\lefteqn{\mathbb{P}_{G}^{w}(R(\left\lfloor \lambda r^{2}\right\rfloor )\in B(v,r))=}\qquad\\
 & =\sum_{n\geq\left\lfloor \lambda r^{2}\right\rfloor }\mathbb{P}(\textrm{the lazy walk did }n\textrm{ steps and is in }B(v,r))=\\
 & =\sum_{n\geq\left\lfloor \lambda r^{2}\right\rfloor }\mathbb{P}(\textrm{the lazy walk did }n\textrm{ steps})\mathbb{P}_{G^{*}}(R(n)\in B(v,r))\leq\\
 & \leq\frac{1}{2}\sum_{n\geq\left\lfloor \lambda r^{2}\right\rfloor }\mathbb{P}(\textrm{the lazy walk did }n\textrm{ steps})=\frac{1}{2}.\qedhere\end{align*}

\end{proof}
Here too Delmotte's theorem can be replaced by Varopoulos \cite{V85}.

The following lemma basically states that a random walk has positive
probability to hit large objects. We will only use the lemma for simple
domains with piecewise smooth boundary, so the requirements of clause
\ref{enu:lemDD-K} will always be satisfied.

\begin{lem}
\label{lem:DD}Let $H$ be a $d$-dimensional Euclidean net.\newC{C:DD}\newc{c:DDi}\newc{c:DDo}
\begin{enumerate}
\item \label{enu:lemDD-hit}Let $\mathcal{D}$, $\mathcal{S}$ (start) and
$\mathcal{H}$ (hit) be domains in $\mathbb{R}^{d}$ with $\overline{\mathcal{S}},\mathcal{H}\subset\mathcal{D}$,
$\mathcal{H}\neq\mathcal{D}$. Then there exists a $C_{\ref{C:DD}}(\mathcal{D},\mathcal{S},\mathcal{H},H)$
such that for all $r>C_{\ref{C:DD}}$; all $v\in H$ and all $w\in\left(v+r\mathcal{S}\right)\cap H$,
if $R$ is a random walk starting from $w$ then\begin{equation}
\mathbb{P}(T(\partial(v+r\mathcal{H}))<T(\partial(v+r\mathcal{D})))>c_{\ref{c:DDi}}(\mathcal{D},\mathcal{S},\mathcal{H},H).\label{eq:hitinside}\end{equation}

\item \label{enu:lemDD-avoid}If $\overline{\mathcal{S}}\cap\overline{\mathcal{H}}=\emptyset$
and $\mathcal{S}$ is a subset of the unbounded component of $\mathbb{R}^{d}\setminus\overline{\mathcal{H}}$
then in addition\begin{equation}
\mathbb{P}(T(\partial(v+r\mathcal{H}))\geq T(\partial(v+r\mathcal{D})))>c_{\ref{c:DDo}}(\mathcal{D},\mathcal{S},\mathcal{H},H).\label{eq:hitoutside}\end{equation}

\item \label{enu:lemDD-K}Let $\mathcal{K}$ is a family of triplets $(\mathcal{D},\mathcal{S},\mathcal{H})$
such that for every $x\in\mathcal{S}$ there exists a path $\gamma$
with $\len\gamma$ bounded above leading from $x$ to $\partial_{\textrm{cont}}\mathcal{H}$
with $d(\gamma,\partial\mathcal{D})$ bounded below (case \ref{enu:lemDD-hit})
or from $x$ to $\partial\mathcal{D}$ with $d(\gamma,\partial_{\textrm{cont}}\mathcal{H})$
bounded below (case \ref{enu:lemDD-avoid}). Then $C_{\ref{C:DD}},c_{\ref{c:DDi}}$
and $c_{\ref{c:DDo}}$ are bounded on $\mathcal{K}$.
\end{enumerate}
\end{lem}
\begin{proof}
Let us start with (\ref{eq:hitinside}). Assume first that $\mathcal{D}=B(0,1)$,
$\mathcal{H}=B(0,\frac{1}{22})$ and $\mathcal{S}=B(0,\frac{1}{2})\setminus B(0,\frac{2}{22})$,
and that $d\geq3$. Examine the Green function $G(w)=G(v,w;v+r\mathcal{D})$.
Let $w\in(v+r\mathcal{S})\cap H$, let $R$ be a random walk starting
from $w$ and let $T$ be its stopping time on $\partial(v+r\mathcal{H})\cup\partial(v+r\mathcal{D})$.
Then, since $G$ is harmonic on $(v+r\mathcal{D})\setminus(v+r\mathcal{H})$
we get (from (\ref{eq:fEfRT})) that\[
G(w)=\mathbb{E}G(R(T)).\]
Denote $p=\mathbb{P}(R(T)\in\partial(v+r\mathcal{H}))$ which is the
probability we want to estimate. Now for every $x\in\partial(v+r\mathcal{D})$,
$G(x)=0$ while for $x\in\partial(v+r\mathcal{H})$ we have from (\ref{eq:atrans})
that $G(x)\leq C(H)r^{2-d}$. At $w$ itself we have from (\ref{eq:GvwS-lower-bound})
that $G(w)\geq cr^{2-d}$ for some $c(H)$ and $r$ sufficiently large.
We get\begin{equation}
cr^{2-d}\leq G(w)=\mathbb{E}G(R(T))\leq p\cdot Cr^{2-d}\label{eq:1213to113}\end{equation}
so $p\geq c$ for $r$ sufficiently large.

To see that the same holds for $d\leq2$, construct an auxiliary graph
$H'=H\times\mathbb{Z}^{3-d}$ weighted so that the projection of the
random walk on $H'$ on $H$ is (a time change of) the random walk
on $H$. Then the fact that there is a positive probability to hit
$\partial B(v,\frac{1}{22}r)\subset H'$ before hitting $\partial B(v,r)\subset H'$
immediately implies the same for $H$. 

We now consider general domains. Let $x\in\mathcal{S}\setminus\mathcal{H}$
and choose an $\epsilon>0$ and a sequence of points $\left\{ x_{i}\right\} _{i=0}^{n}$
with the following properties:
\begin{enumerate}
\item $x_{0}=x$ and $B(x_{n},\epsilon)\subset\mathcal{H}$.
\item $B(x_{i},2\epsilon)\subset\mathcal{D}$.
\item $|x_{i+1}-x_{i}|\in\left]\frac{6}{22}\epsilon,\frac{7}{22}\epsilon\right[$.
\end{enumerate}
It is a simple exercise to show that $\epsilon$ and $x_{i}$ can
always be chosen, and furthermore that both $\epsilon$ and $n$ may
be bounded on all $\mathcal{S}\setminus\mathcal{H}$ and in case \ref{enu:lemDD-K}
on all $\mathcal{K}$ i.e.~for any $x$ in any $\mathcal{S}\setminus\mathcal{H}$
such that $(\mathcal{D},\mathcal{S},\mathcal{H})\in\mathcal{K}$ we
have $\epsilon(x,\mathcal{H},\mathcal{D})>c(\mathcal{K})$ and $n(x,\mathcal{H},\mathcal{D})<C(\mathcal{K})$.
Assume now that $r$ is sufficiently large such that the following
are satisfies:
\begin{enumerate}
\item \label{enu:every-ball-point}every ball of radius $\frac{1}{22}\epsilon r$
contains at least one point of $H$;
\item $\partial B(w,s)\subset B(w,s+\frac{1}{22}\epsilon r)$ for any $s>0$
and $w\in H$;
\item $\epsilon r$ is sufficiently large so as to satisfy (\ref{eq:1213to113}).
\end{enumerate}
Clearly this is an assumption of the type $r>C(\mathcal{D},\mathcal{S},\mathcal{H},H)$.
Condition \ref{enu:every-ball-point} allows to choose a point $w_{i}\in H\cap B(v+rx_{i},\frac{1}{22}\epsilon r)$
for every $i$. Let $T_{0}=0$ and \[
T_{i}=\min\{ t>T_{i-1}:R(t)\in\partial B(w_{i},{\textstyle \frac{1}{22}}\epsilon r)\}.\]
We wish to use the case already established for the portion of the
random walk after $T_{i-1}$, with the radius being $\epsilon r$
instead of $r$ and the center of the balls being $w_{i}$ instead
of $v$. We may do this because $R(T_{i-1})\in\partial B(w_{i-1},\frac{1}{22}\epsilon r)\subset B(w_{i-1},\frac{2}{22}\epsilon r)\subset B(v+rx_{i-1},\frac{3}{22}\epsilon r)\subset B(v+rx_{i},\frac{10}{22}\epsilon r)\setminus B(v+rx_{i},\frac{3}{22}\epsilon r)\subset B(w_{i},\frac{1}{2}\epsilon r)\setminus B(w_{i},\frac{2}{22}\epsilon r)$.
However, since $B(w_{i},\epsilon r)\subset B(v+rx_{i},2\epsilon r)\subset v+r\mathcal{D}$
then not hitting $\partial B(w_{i},\epsilon r)$ means staying inside
$v+r\mathcal{D}$ and so we get\[
\mathbb{P}\big(T_{i}<T(\partial(v+r\mathcal{D}))\,|\, T_{i-1}<T(\partial(v+r\mathcal{D}))\big)>c(H).\]
This immediately gives\[
p\stackrel{(*)}{\geq}\mathbb{P}(T_{n}<T(\partial(v+r\mathcal{D})))>c^{n}\]
where $(*)$ comes from the fact that $\partial B(w_{n},\frac{1}{22}\epsilon r)\subset B(v+rx_{n},\frac{3}{22}\epsilon r)\subset v+r\mathcal{H}$.
This proves the direction (\ref{eq:hitinside}) for $x\not\in\mathcal{H}$.
The case $x\in\mathcal{H}$ is proved identically but taking the $x_{i}$-s
from $x$ to $\mathcal{D}\setminus\mathcal{H}$. For the direction
(\ref{eq:hitoutside}) take the $x_{i}$-s outside, i.e.~with $B(x_{n},\epsilon)\cap\mathcal{D}=\emptyset$
and $B(x_{i},2\epsilon)\cap\mathcal{H}=\emptyset$.
\end{proof}
\begin{lem}
\label{lem:DDunbounded}With the notations of the previous lemma,
if $d\geq3$ then (\ref{eq:hitinside}) and (\ref{eq:hitoutside})
hold even if $\mathcal{D}$ is allowed to be unbounded. If $d\leq2$
then only (\ref{eq:hitinside}) holds for unbounded $\mathcal{D}$.
\end{lem}
\begin{proof}
The only part not following directly from lemma \ref{lem:DD} is the
proof of (\ref{eq:hitoutside}) when $d\geq3$. We start with the
case that $\mathcal{D}=\mathbb{R}^{d}$ (so $T(\partial(v+r\mathcal{D}))$
is always $\infty$), $\mathcal{H}=B(0,1)$ and $\mathcal{S}\cap B(0,\lambda)=\emptyset$
where $\lambda(G)>2$ is some constant sufficiently large that will
be fixed later. Let $s>\lambda r$ and denote \[
p(s)=\mathbb{P}^{w}(T_{v,s}<T_{v,r})\quad w\in v+r\mathcal{S}.\]
 Let $G(x)=G(v,x;B(v,s))$ i.e.~Green's function. From (\ref{eq:GvwS-lower-bound})
we see that for $x\in\partial B(v,r)$ we have $G(x)\geq c(H)r{}^{2-d}$
if $r$ is sufficiently large. At $w$ itself we have $G(w)\leq C(H)(\lambda r)^{2-d}$
and by definition $G|_{\partial B(v,s)}\equiv0$ so\[
C(\lambda r)^{2-d}\geq G(w)=\mathbb{E}G(R(T))\geq(1-p)c(H)r^{2-d}\]
so $p\geq1-C\lambda^{2-d}$ and for $\lambda$ sufficiently large
this would be $\geq\frac{1}{2}$. Fix $\lambda$ to be some such constant.
We get that $\lim_{s\to\infty}p\geq\frac{1}{2}$. Hence\begin{equation}
\mathbb{P}(T(\partial(v+r\mathcal{H}))=\infty)=\lim_{s\to\infty}\mathbb{P}(T(\partial(v+r\mathcal{H}))>T_{v,s})\geq{\textstyle \frac{1}{2}}\label{eq:infcG}\end{equation}
and this case is finished. For general $\mathcal{D}$, $\mathcal{S}$
and $\mathcal{H}$, let $B(0,\rho)$ ($\rho>1$) be a ball sufficiently
large such that \emph{$\mathcal{D}\cap B(0,\rho)$} contains $\mathcal{S}$,
$\mathcal{H}$ and a path between them. Let $\mathcal{D}_{2}$ be
the component of $\mathcal{D}\cap B(0,\lambda\rho)$ containing $\mathcal{S}$
and $\mathcal{H}$. Then lemma \ref{lem:DD} shows that for $r$ sufficiently
large we have \begin{equation}
\mathbb{P}(T(\partial(v+r\mathcal{H}))\geq T(\partial(v+r\mathcal{D}_{2})))>c(\mathcal{D}_{2},\mathcal{S},\mathcal{H},G).\label{eq:HD2c}\end{equation}
If $R(T(\partial(v+r\mathcal{D}_{2})))\not\in\partial(v+r\mathcal{D})$
then it is $\in\partial B(v,\lambda\rho r)$ and then the previous
case (rescaled by $\rho$) with the strong Markov property shows that
\begin{eqnarray*}
\lefteqn{\mathbb{P}(T(\partial(v+r\mathcal{H}))\geq T(\partial(v+r\mathcal{D}))\,|\, T(\partial(v+r\mathcal{H})\geq T(\partial(v+r\mathcal{D}_{2})))=}\\
 &  & =\mathbb{P}(T(\partial(v+r\mathcal{H}))\geq T(\partial(v+r\mathcal{D})),\, T(\partial(v+r\mathcal{D}_{2}))=T_{v,\lambda\rho r}\,|\,\\
 &  & \qquad\qquad T(\partial(v+r\mathcal{H})\geq T(\partial(v+r\mathcal{D}_{2})))=\\
 &  & =\mathbb{E}\mathbb{P}^{R(T_{v,\lambda\rho r})}(T(\partial(v+r\mathcal{H}))\geq T(\partial(v+r\mathcal{D})))\geq\\
 &  & \geq\mathbb{E}\mathbb{P}^{R(T_{v,\lambda\rho r})}(T_{v,\rho r}=\infty)\stackrel{(\ref{eq:infcG})}{\geq}\mathbb{E}c(G)=c(G)\end{eqnarray*}
and together with (\ref{eq:HD2c}) we get (\ref{eq:hitoutside}).\newC{C:DDp}
\end{proof}
\begin{lem}
\label{lem:DDp}Let $\mathcal{D}$, $\mathcal{S}$ and $\mathcal{H}$
be domains in $\mathbb{R}^{d}$ with $\overline{\mathcal{S}},\overline{\mathcal{H}}\subset\mathcal{D}$
and $\overline{\mathcal{S}}\cap\overline{\mathcal{H}}=\emptyset$.
Then there exists a $C_{\ref{C:DDp}}(\mathcal{D},\mathcal{S},\mathcal{H},G)$
such that for all $r>C_{\ref{C:DDp}}$; all $v\in G$; all $w\in\left(v+r\mathcal{S}\right)\cap G$
and all $x\in(v+r\mathcal{H})\cap G$, if $R$ is a random walk starting
from $w$ then\[
\mathbb{P}(T(\{ x\})<T(\partial(v+r\mathcal{D})))\approx\begin{cases}
r^{2-d} & d\geq3\\
1/\log r & d=2\end{cases}\]
where the constants implicit in the $\approx$ may depend on $\mathcal{D},\mathcal{S},\mathcal{H}$
and $G$. Further, if $\mathcal{K}$ is family of $(\mathcal{D},\mathcal{H},\mathcal{S})$
triplets satisfying the conditions of lemma \ref{lem:DD}, clause
\ref{enu:lemDD-K} then $C_{\ref{C:DDp}}$ and the implicit constants
are bounded on $\mathcal{K}$.
\end{lem}
\begin{proof}
Let $\epsilon=\epsilon(\mathcal{K},G)$ be sufficiently small. Use
lemma \ref{lem:DD} to show that the probability to hit a ball of
radius $r\epsilon$ around $x$ is $\approx1$ and then the same Green's
function calculations as in that lemma to show that the probability
to hit a point before exiting from a ball containing $\mathcal{D}$
are $\approx r^{2-d}$ for $d\geq3$ and $\approx1/\log r$ for $d=2$.
\end{proof}
\begin{lem}
\label{lem:hitbsame}Let $G$ be a $d$-dimensional Euclidean net,
$d\geq3$. Let $v\in G$, $s>4r>C(G)$, let $A\subset B(v,r)$ and
$w\in\partial B(v,2r)$. Then\[
\mathbb{P}^{w}(T(A)<T_{v,4r})\approx\mathbb{P}^{w}(T(A)<T_{v,s}).\]
Further, if $B\subset A$ then\[
\mathbb{P}^{w}(R(T(A\cup\partial B(v,4r)))\in B)\approx\mathbb{P}^{w}(R(T(A\cup\partial B(v,s)))\in B).\]

\end{lem}
\begin{proof}
We shall only show the first estimate, the second one is proved identically.
Clearly $\mathbb{P}^{w}(T(A)<T_{v,4r})\leq\mathbb{P}^{w}(T(A)<T_{v,s})$,
so we need to show the other direction. Define stopping times $T_{0}=0$
and\begin{align*}
T_{2i+1} & :=\min\{ t>T_{2i}:R(t)\in\partial B(v,4r)\cup A)\}\\
T_{2i} & :=\min\{ t>T_{2i-1}:R(t)\in\partial B(v,2r)\cup\partial B(v,s)\}.\end{align*}
Let $I$ be the first time when $R(T_{I})\in A$ (for $i$ odd) or
$R(T_{i})\in\partial B(v,s)$ (for $i$ even). We consider the process
stopped at $I$. From lemma \ref{lem:DDunbounded} we see that \begin{equation}
\mathbb{P}(R(T_{2i})\in\partial B(v,s)\,|\, I>2i-1)\geq\mathbb{EP}^{R(T_{2i-1})}(T_{v,2r}=\infty)\geq c(G).\label{eq:RT2iexp}\end{equation}
From Harnack's inequality we get that\begin{multline*}
\mathbb{P}^{w}(R(T_{2i+1})\in A\,|\, I>2i)=\mathbb{E}\mathbb{P}^{R(T_{2i})}(T(A)<T_{v,4r})\approx\\
\approx\min_{x\in\partial B(v,2r)}\mathbb{P}^{x}(T(A)<T_{v,4r})\end{multline*}
and hence this is (up to a constant) independent of $i$. Hence we
get\begin{align*}
\mathbb{P}^{w}(T(A)<T_{v,s}) & =\sum_{i=0}^{\infty}\mathbb{P}^{w}(I>2i,\, R(T_{2i+1})\in A)\leq\\
 & \leq C(G)\sum_{i=0}^{\infty}\mathbb{P}^{w}(R(T_{1})\in A)\mathbb{P}(I>2i)\leq\\
 & \!\!\stackrel{(\ref{eq:RT2iexp})}{\leq}C(G)\sum_{i=0}^{\infty}\mathbb{P}^{w}(R(T_{1})\in A)(1-c(G))^{i}=\\
 & =C(G)\mathbb{P}^{w}(T(A)<T_{v,4r}).\qedhere\end{align*}

\end{proof}

\subsection{\label{sub:Beurling}The discrete Beurling projection in three dimensions}

The Beurling projection theorem says that the probability of a two
dimensional Brownian motion starting at $0$ to hit a given set $K\subset B(0,1)$
before hitting $\partial B(0,1)$ is larger than the probability to
hits its \emph{angular projection}, namely the set $\{|z|:z\in K\}$.
In particular, if $K$ is connected and intersects both $\partial B(0,\epsilon)$
and $\partial B(0,1)$ then the probability to avoid it is maximal
when $K=[\epsilon,1]$, and in this case it may be calculated explicitly
from the conformal invariance of Brownian motion and is $\approx\sqrt{\epsilon}$.
A discrete version of this result (up to constants) was achieved by
Kesten \cite{K87}. In this section we shall prove a three dimensional
variation on this result, namely the following lemma:

\begin{lem}
\label{lem:Beurling}Let $H$ be a three dimensional Euclidean net.
Then there exists a constant $C(H)$ such that for all $v\in H$,
for all $r>C(H)$ and for all connected sets $A\subset H$ that intersect
both $B(v,r)$ and $H\setminus B(v,2r)$ one have\[
\mathbb{P}^{v}(R[0,T_{v,4r}]\cap A\neq\emptyset)\geq\frac{c}{\log r}.\]

\end{lem}
While Kesten's version of Beurling's arguments may be applied to three
dimensions without much change, in the setting of lemma \ref{lem:Beurling}
the notion of capacity, particularly of Martin capacity, can be used
to shorten the argument significantly. We shall first give the relevant
definitions, and the proof will follow after.

\subsubsection{Martin capacity}

\begin{defn*}
Let $H$ be a countable set and let $K(v,w)$ be some function ({}``the
kernel''). The capacity of a set $S\subset H$ with respect to $K$
is defined by \[
\capa_{K}(S):=\Big(\inf_{\mu(S)=1}\int_{S}\int_{S}K(v,w)\, d\mu(v)\, d\mu(w)\Big)^{-1}.\]
The infimum here is over all probability measures $\mu$ supported
on $S$.
\end{defn*}

\begin{defn*}
Let $H$ be a directed graph and let $v\in H$. Then the \textbf{Martin
capacity} of the graph with respect to $v$ is the capacity with respect
to the Martin kernel, defined by\[
K(w,x):=\frac{G(w,x)}{G(v,x)}.\]
where $G$ is Green's function.
\end{defn*}
\begin{bpp}Let $H$ be a directed graph, let $v\in H$ and let $S\subset H$
satisfy that Green's function $G(w,x)$ is finite for all $w,x\in S$.
Let $\capa$ be the Martin capacity with respect to $v$. Then for
any $S$ we have\[
{\textstyle \frac{1}{2}}\capa(S)\leq\mathbb{P}^{v}(R\left[0,\infty\right[\cap S\neq\emptyset)\leq\capa(S).\]
\end{bpp}The nice and simple proof may be found in \cite{BPP95},
theorem 2.2.

\begin{proof}
[Proof of lemma \ref{lem:Beurling}]Let $\lambda=\lambda(H)$ be some
constant such that every edge of $H$ has length $\leq\lambda$. Then
$A$ intersects every spherical shell $B(v,s+\lambda)\setminus B(v+s)$,
$r\leq s\leq2r$. Let $a_{i}\in A\cap(B(v,r+(i+1)\lambda)\setminus B(v,r+i\lambda))$
for $i=0,\dotsc,\left\lfloor r/\lambda\right\rfloor $. Let $A^{*}=\{ a_{i}\}$.
The lemma will be proved if we show\[
\mathbb{P}^{v}(T(A^{*})<T_{v,4r})\geq\frac{c}{\log r}.\]
Let $H'$ be the directed graph given by taking $\overline{H\cap B(v,4r)}$
and making each point of $\partial B(v,4r)$ a {}``sink'' i.e. a
point with the only exit being a self loop. By definition, $G_{H'}(w,x)=G_{H}(w,x;B(v,4r))$.
We get the equivalent formulation $\mathbb{P}_{H'}^{v}(T(A^{*})<\infty)\geq c/\log r$.
We now use Benjamini-Pemantle-Peres on $H'$. We get that it is enough
to estimate $\capa(A^{*})\geq c/\log r$. By the definition of capacity
we that we need to show that there exists a $\mu$ on $A^{*}$ such
that \begin{equation}
\int_{A^{*}}\int_{A^{*}}\frac{G(w,x)}{G(v,x)}d\mu(w)d\mu(x)\leq C\log r.\label{eq:req-mu}\end{equation}
Let $\mu$ be the uniform measure on $A^{*}$. Then by (\ref{eq:GvwS-lower-bound})
we have that \[
K(a_{i},a_{j})=\frac{G(a_{i},a_{j})}{G(v,a_{j})}\leq C\frac{|a_{i}-a_{j}|^{-1}}{r^{-1}}\leq Cr|i-j|^{-1}.\]
Summing gives (\ref{eq:req-mu}) and the lemma.
\end{proof}

\subsection{Intersection probabilities}

\begin{lem}
\label{lem:Rintrsct}Let $H$ be a Euclidean net of dimension $\leq3$
and let $\epsilon\in\left]0,1\right[$. Let $v\in H$ and $r>0$ and
let $R^{1}$ and $R^{2}$ be random walks starting from vertices $v^{1}$
and $v^{2}$, $v^{i}\in B(v,(1-\epsilon)r)$. Then\newC{C:minr}\newc{c:intrsct}\[
\mathbb{P}(R^{1}[0,T_{v^{1},\epsilon r}^{1}]\cap R^{2}[0,T_{v,r}^{2}]\neq\emptyset)>c_{\ref{c:intrsct}}(\epsilon,H)\]
if only $r>C_{\ref{C:minr}}(\epsilon,H)$.
\end{lem}
\begin{proof}
We shall only show the case $d=3$ --- the case $d=2$ is identical
and will be left to the reader. Let $\lambda=\lambda(\epsilon,H)$
be some constant whose value will be fixed later. Let $s\geq\epsilon\lambda$
be some number, $w\in H$ and $\delta>0$. For any $x\in B(w,(1-\delta)s$),
let $\Gamma(y)$ be the probability that a random walk starting from
$y$ will hit $x$ before hitting $\partial B(w,s)$, and define $\Gamma(x):=1$.
$\Gamma$ is harmonic on $B(w,s)\setminus\{ x\}$ and $0$ on $\partial B(w,s)$.
Hence $\Gamma(y)=G(y,x;B(w,s))/G(x,x;B(w,s))$ where $G$ is Green's
function. Using (\ref{eq:atrans}), (\ref{eq:GvwS-lower-bound}) and
Harnack's inequality (lemma \ref{lem:Harnack-general}) we get\begin{align*}
\Gamma(y) & \geq c(\delta,H)/|y-x| & y & \in B(w,(1-\delta)s)\setminus\{ x\}\\
\Gamma(y) & \leq C(\delta,H)/|y-x| & y & \in B(w,s)\setminus\{ x\}.\end{align*}
for $s$ sufficiently large.

Define now $A:=B(v^{1},\frac{1}{2}\epsilon r)\setminus B(v^{1},\frac{1}{4}\epsilon r)$
and $X:=A\cap R^{1}[0,T_{v^{1},\epsilon r}^{1}]\cap R^{2}[0,T_{v,r}^{2}]$.
For any $x\in A$, the preceding calculation (used once for $w=v^{1}$,
$s=\epsilon r$ and $\delta=\frac{1}{2}$ and a second time for $w=v^{2}$,
$s=r$ and $\delta=\frac{1}{2}\epsilon$) shows that $\mathbb{P}(x\in X)\approx r^{-2}$.
The $\approx$ sign here and below may depend on $\epsilon$ and on
the Euclidean net structure constants of $H$. Rough isometry preserves
(up to a constant) volumes of balls and shells, hence if $\epsilon r$
is sufficiently large we get $\# A\approx r^{3}$ and hence $\mathbb{E}\# X\approx r$.
Next we want to calculate $\mathbb{E}(\# X)^{2}$. For any $x\neq y\in A$
we have\begin{equation}
\mathbb{P}(x,y\in R^{i}[0,T^{i}])\approx r^{-1}|x-y|^{-1}.\label{eq:probxy}\end{equation}
Indeed, this probability is $\geq$ the probability to hit $x$ first
(which is $\approx1/r$) and then to hit $y$ (which is $\approx1/|x-y|$).
On the other hand it is $\leq$ the sum of this probability and its
symmetric image. So (\ref{eq:probxy}) is explained. This shows that
$\mathbb{P}(x,y\in X)\approx r^{-2}|x-y|^{-2}$ and summing over $y$
we get\begin{align*}
\sum_{y:y\neq x}\mathbb{P}(x,y\in X) & \stackrel{(*)}{=}\sum_{n=-C(H)}^{\left\lfloor \log_{2}r\right\rfloor }\sum_{y\in B(x,2^{n+1})\setminus B(x,2^{n})}\mathbb{P}(x,y\in X)\\
 & \leq C(\epsilon,H)r^{-2}\sum_{n=-C(H)}^{\left\lfloor \log_{2}r\right\rfloor }4^{-n}\#(B(x,2^{n+1})\setminus B(x,2^{n}))\\
 & \!\!\stackrel{(**)}{\leq}\!\! C(\epsilon,H)r^{-2}\sum_{n=-C(H)}^{\left\lfloor \log_{2}r\right\rfloor }2^{n}\leq C(\epsilon,H)r^{-1}\end{align*}
where $(*)$ comes from the fact that $H$ is separated in $\mathbb{R}^{3}$
hence $|x-y|\geq c(H)$ and $(**)$ uses again the fact that rough
isometry preserves volumes. Summing over $x$ we get \[
\mathbb{E}(\# X)^{2}\leq\# A\cdot C(\epsilon,H)r^{-1}+\sum_{x\in A}\mathbb{P}(x\in X)\leq C(\epsilon,H)r^{2}.\]
 Hence the well known inequality $\mathbb{P}(\# X>0)\geq(\mathbb{E}\# X)^{2}/\mathbb{E}(\# X)^{2}$
finishes the lemma.
\end{proof}
\begin{lem}
\label{lem:roughxi}Let $G$ be a Euclidean net of dimension $d\leq3$
and let $v^{1}$ and $v^{2}\in G$. Let $R^{i}$ be random walks starting
from $v^{i}$ and stopped on $\partial B(v^{1},r)$. Then\[
\mathbb{P}(R^{1}\cap R^{2}=\emptyset)\leq C(G)\left(\frac{|v^{1}-v^{2}|}{r}\right)^{c(G)}.\]

\end{lem}
\begin{proof}
Let $a_{j}:=2^{j}|v^{1}-v^{2}|$ and assume without loss of generality
that $r=a_{n}$ for some integer $n$. Let $T_{j}^{i}$ be the stopping
time of $R^{i}$ on the shell $\partial B(v^{1},a_{j})$. Examine
the events \[
\mathcal{E}_{j}:=\{ R^{1}[T_{j}^{1},T_{j+1}^{1}]\cap R^{2}[T_{j}^{2},T_{j+1}^{2}]\neq\emptyset\}.\]
For $j>C(G)$ we have that $\partial B(v^{1},a_{j})\subset B(v^{1},\frac{2}{3}a_{j+1})$
and we may use lemma \ref{lem:Rintrsct} with $\epsilon=\frac{1}{3}$
and the strong Markov property to get\[
\mathbb{P}(\mathcal{E}_{j}\,|\, R^{1}(T_{j}^{1}),R^{2}(T_{j}^{2}))\geq c_{\ref{c:intrsct}}({\textstyle \frac{1}{3}},G)\quad\forall j>C(G).\]
However, $\mathcal{E}_{j}$ may depend on $\mathcal{E}_{0},\dotsc,\mathcal{E}_{j-1}$
only through $R^{i}(T_{j}^{i})$ so we get \[
\mathbb{P}(\mathcal{E}_{j}\,|\,\mathcal{E}_{0},\dotsc,\mathcal{E}_{j-1})\geq c(G).\]
And hence\[
\mathbb{P}(R^{1}\cap R^{2}=\emptyset)\leq\mathbb{P}\Big(\bigcap_{j=C}^{n-1}\neg\mathcal{E}_{j}\Big)\leq(1-c(G))^{n-C}\leq C\left(\frac{|v^{1}-v^{2}|}{r}\right)^{c(G)}\]
and the lemma is proved.
\end{proof}
This basic proof method is known as the {}``Wiener shell test''. 

\begin{rem*}
Given theorem \ref{thm:kof} below (page \pageref{thm:kof}) it might
be tempting to conjecture that $c(G)$ is in effect $\xi$, the non-intersection
exponent of $d$-dimensional Brownian motion. However, this is not
true. Indeed, the intersection exponent is a {}``conformally invariant''
property rather than a metric property. Unfortunately, I don't know
any example sufficiently simple to explain here.
\end{rem*}
\begin{lem}
\label{lem:plane}Let $G$ be a Euclidean net and let $H\subset\mathbb{R}^{d}$
be a half-space and let $v\in G\setminus\partial H$. Let $R$ be
a random walk starting from $v$. Then\[
\mathbb{P}(T_{v,r}<T(\partial H))\leq C(G)\left(\frac{d(v,\partial H)}{r}\right)^{c(G)}\quad\forall r>2d(v,\partial H).\]

\end{lem}
\begin{proof}
Let $s>2d(v,\partial H)$ and examine a random walk $R$ starting
from any point in $\partial B(v,s)$ and stopped on $\partial B(v,2s)$.
We use lemma \ref{lem:DD} with $\mathcal{D}=B(0,2)$, $\mathcal{S}=B(0,\frac{3}{2})$
and $\mathcal{H}=((H-v)/s)\cap\mathcal{D}$. If $s>C(G)$ then $\partial B(v,s)\subset v+s\mathcal{S}$
and $\mathcal{H}$ is non empty so the lemma applies to our $R$.
We get that the probability of $R$ to hit $\partial(v+s\mathcal{H})\subset\partial H\cap B(v,2s)$
before $\partial B(v,2s)$ is $\geq c_{\ref{c:DDi}}(\mathcal{D},\mathcal{S},\mathcal{H},G)$
if only $s>C_{\ref{C:DD}}(\mathcal{D},\mathcal{S},\mathcal{H},G)$.
Further, if $s$ is sufficiently large then the family of possible
$\mathcal{H}$-s satisfies the requirements of clause \ref{enu:lemDD-K}
of lemma \ref{lem:DD} and these $c$-s and $C$-s are bounded. The
Wiener shell test now gives the lemma.
\end{proof}
Again, it is not necessarily true that $c=1$ as in $\mathbb{Z}^{d}$.
This is only true with additional assumptions, such as isotropicity,
see theorem \ref{lem:escape} below (page \pageref{lem:escape}).
\textbf{}A counterexample may be constructed as follows: Let $\phi:\overline{\mathbb{H}}\to\mathbb{C}$
be defined by $\phi(re^{i\theta})=re^{2i\theta}$ where $\mathbb{H}$
is as usual the upper half space $\{\mathrm{Im\,}z>0\}$. Let $G:=\phi(\mathbb{Z}^{2}\cap\overline{\mathbb{H}})$
and identify $\phi(n)$ and $\phi(-n)$ so that $G$ contains edges
from $n$ to both $\sqrt{n^{2}+1}\, e^{\pm i2\arctan1/n}$. Then it
is easy to see that $G$ is a Euclidean net while the escape probability
from, say, $v=i\sqrt{2}$ to $\partial B(v,r)$ without hitting $\partial\mathbb{H}$
are the same as the escape probabilities of a random walk on $\mathbb{Z}^{2}$
from a corner, which are well known to be $\approx r^{-2}$.

An argument identical to that of lemma \ref{lem:plane} works for
any polyhedron:

\begin{lem}
\label{lem:polyh}Let $G$ be a Euclidean net and let $Q\subset\mathbb{R}^{d}$
be a polyhedron. Let $1<r_{1}<r_{2}<s$ be some numbers and let $v\in G$
satisfy that $d(v,sQ)\leq r_{1}$. Let $R$ be a random walk starting
from $v$. Then \[
\mathbb{P}(T_{v,r_{2}}<T(sQ))\leq C(Q,G)\left(\frac{r_{1}}{r_{2}}\right)^{c(Q,G)}.\]

\end{lem}
We omit the proof.

The next lemma is technical and is only here for completeness. In
fact all the graphs we will consider in this paper have no {}``dangling
ends'' and it is straightforward to see that for every $v^{1}\neq v^{2}$
one can construct disjoint paths going in opposite directions (so
clause \ref{enu:gamisep} of the lemma is satisfied).

\begin{lem}
\label{lem:0n0}Let $G$ be an Euclidean net of dimension $\geq2$
and let $0<\epsilon<\frac{1}{2}$, $s>0$. Then there exists a $\kappa=\kappa(\epsilon,s,G)$
such that for any $v,v^{1}\neq v^{2}\in G$ with $|v-v^{i}|\leq s$
one of the following holds:
\begin{enumerate}
\item There are no two disjoint paths starting from $v^{i}$ and ending
outside $B(v,\kappa)$.
\item \label{enu:gamisep}For any $r\geq\kappa$ there exists two disjoint
simple paths $\gamma^{i}\subset B(v,r)$ starting from $v^{i}$ and
ending in some $w^{i}$ such that \begin{equation}
B(w^{i},\left(1-\epsilon\right)r)\cap\gamma^{3-i}=\emptyset,\quad w^{i}\in B(v,r)\setminus\overline{B(v,(1-\epsilon)r)}\label{eq:lem0n0}\end{equation}

\end{enumerate}
\end{lem}
\begin{proof}
Let $\lambda=\lambda(G)$ satisfy that any ball in $\mathbb{R}^{d}$
with radius $\geq\lambda$ contains at least one point of $G$, and
such that no edge of $G$ has length $>\lambda$. Let $\mu$ satisfy
that for any $x$, $y$ in $G$ such that $|x-y|\leq4\lambda$ there
is a path $\gamma$ from $x$ to $y$, $\gamma\subset B(x,\mu)$.
Now take any point $x\in G$ and any direction $\theta$ and construct
an infinitely long {}``ray'' $x\in R\subset G$ by taking a sequence
of tangent balls on the half line in direction $\theta$, taking a
point of $G$ in every ball and connecting them by short paths as
above. The result is that $R$ is contained in the open infinite cylinder
whose basis is a $d-1$ dimensional ball of radius $\mu+\lambda$
orthogonal to $\theta$, centered at $x$. Actually, $R$ is contained
only in the half cylinder starting $\mu$ before $x$.

It is now clear that if $|v^{1}-v^{2}|>2\mu$ then we may simply extend
such rays in the directions $\pm(v^{1}-v^{2})$ and they will not
intersect. This allows to prove the case $s>2\mu$ given the case
$s=2\mu$ --- this is a simple geometric exercise (I got that it is
enough to define $\kappa(s,\epsilon)=(s+C(G)+\kappa(\frac{1}{10},2\mu))/\epsilon$).
Hence we will assume $s=2\mu$.

Define now $\nu:=2\mu+2\lambda$ and  impose the condition $\kappa\geq16\nu$.
Let now $\delta^{1}$ and $\delta^{2}$ be two disjoint paths starting
from $v^{i}$ and going to a distance of $\kappa$. The lemma will
be proved once we construct $\gamma^{i}$ satisfying \ref{enu:gamisep}.
Let $x^{i}=\delta^{i}(j^{i})$ be the first point of $\delta^{i}$
outside $B(v,16\nu)$. A simple exercise in plane geometry shows that,
if $|v-x^{i}|\geq\alpha$ then one can find an $\eta$ satisfying
that $\left|\langle\eta,v\rangle-\langle\eta,x^{i}\rangle\right|\geq\frac{1}{2}\alpha$
for $i=1,2$. Applying this in our case gives $\left|\langle\eta,v\rangle-\langle\eta,x^{i}\rangle\right|\geq8\nu$.
Define now\[
I^{i}=\left[\min_{j\leq j^{i}}\left\langle \eta,\delta^{i}(j)\right\rangle ,\max_{j\leq j^{i}}\left\langle \eta,\delta^{i}(j)\right\rangle \right].\]

In the sequel, we will say about a (half-)cylinder $\mathcal{C}$
in a direction orthogonal to $\eta$ that it is {}``in elevation
$e$'' if \[
\min_{v\in\mathcal{C}}\left\langle \eta,v\right\rangle =\left\langle \eta,x\right\rangle +e.\]

\noindent \emph{Case 1}. Consider the case that $|I^{1}\cap I^{2}|<6\nu$.
In this case one of the $I^{i}$-s --- without loss of generality,
we may assume $I^{1}$ --- contains $\left\langle \eta,v\right\rangle +[\nu,8\nu]$
and $I^{2}$ contains $\left\langle \eta,v\right\rangle +[-8\nu,-\nu]$
(otherwise replace $\eta$ with $-\eta$). Therefore $\left\langle \eta,v\right\rangle \pm[7\nu,8\nu]$
is contained in $I^{1}\setminus I^{2}$ and $I^{2}\setminus I^{1}$
respectively. Let $y^{i}=\delta^{i}(k^{i})$ satisfy that \[
\Big|\left|\left\langle \eta,y^{i}\right\rangle -\left\langle \eta,v\right\rangle \right|-{\textstyle \frac{15}{2}}\nu\Big|\leq{\textstyle \frac{1}{2}}\lambda.\]
Such $y^{i}$ always exist since every edge in $G$ has length $\leq\lambda$.
Let $\theta$ be some vector orthogonal to $\eta$. Let $R^{1}$ (respectively
$R^{2}$) be an infinite path starting from $w^{1}$ (respectively
$w^{2}$) and contained in the half cylinder of radius $\frac{1}{2}\nu$
in the direction $\theta$ (respectively $-\theta$) and in elevation
$7\nu$ (respectively $-8\nu$) . Since the cylinders are disjoint
so are the $R^{i}$-s. See figure \ref{cap:I1uI2d}, left. %
\begin{figure}
\input{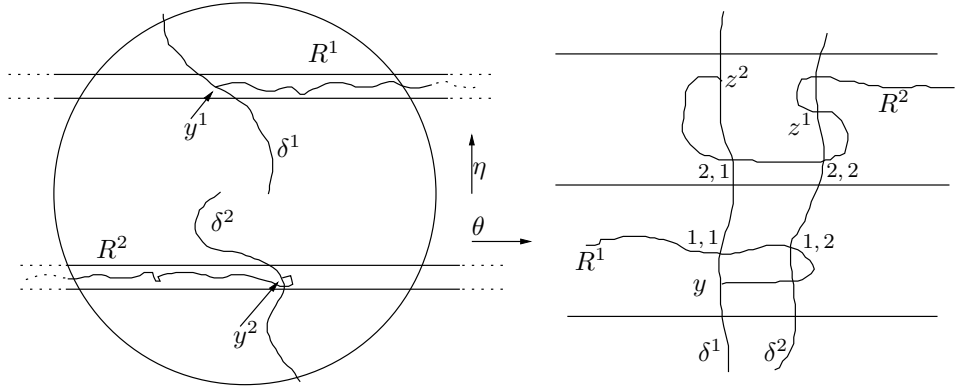}

\setcaptionmargin{0.5in}

\caption{\label{cap:I1uI2d}On the left, case 1, the case of $|I^{1}\cap I^{2}|$
small. On the right, case 2. The points denoted by $i,j$ are $\delta^{j}(k^{i,j})=R^{i}(m^{i,j})$.}
\end{figure}
Since $I^{1}\cap\left\langle \eta,v\right\rangle +[-8\nu,-7\nu]=\emptyset$
we get that $\delta^{1}\cap R^{2}=\emptyset$ and symmetrically we
have $\delta^{2}\cap R^{1}=\emptyset$. Let $l^{i}$ be the first
time such that $R^{i}(l^{i})\not\in B(x,r$). Then we define\[
\gamma^{i}=\LE(\delta^{i}[0,k^{i}]\cup R^{i}[0,l^{i}-1]).\]
Clearly $\delta^{i}\cup R^{i}$ are disjoint and the operation of
taking $\LE$ conserves this, hence the $\gamma^{i}$-s are simple
and disjoint. (\ref{eq:lem0n0}) now follows from simple plane geometry
(recall that $\LE$ conserves the end points) if only $\kappa$ is
sufficiently large, so this case is finished.

\medskip{}
\noindent \emph{Case 2}. We now assume that $|I^{1}\cap I^{2}|\geq6\nu$.
Without loss of generality we may assume that for some $a\geq\left\langle \eta,v\right\rangle +\nu$
we have $[a,a+2\nu]\subset I^{1}\cap I^{2}$ (if not, replace $\eta$
with $-\eta$). Let $y\in\delta^{1}$ be the first that satisfies
that $|\left\langle \eta,y\right\rangle -(a+\frac{1}{2}\nu)|\leq\frac{1}{2}\lambda$
and let $z^{i}\in\delta^{i}$ satisfy that $|\left\langle \eta,z^{i}\right\rangle -(a+\frac{3}{2}\nu)|\leq\frac{1}{2}\lambda$.

This time define $\theta$ be the projection of $z^{1}-z^{2}$ into
the (hyper) plane orthogonal to $\eta$ (if $z^{1}-z^{2}$ is collinear
with $\eta$, just pick $\theta$ an arbitrary vector orthogonal to
$\eta$). Let $\mathcal{C}^{2}$ be a cylinder of radius $\frac{1}{2}\nu$
in the direction $\theta$ with elevation $a+\nu$ such that both
$z^{i}$ are in the middle cylinder of side length $\lambda$. Let
$R^{2}$ be an infinite path in a half of $\mathcal{C}^{2}$ in the
direction $\theta$ starting from $z^{2}$ and containing $z^{1}$.
Let $\mathcal{C}^{1}$ be a similar cylinder in elevation $a$ containing
$y$ in its middle and let $R^{1}$ be an infinite path in the half
of $\mathcal{C}^{1}$ in the direction $-\theta$ starting from $y$.
Also let $R^{1}$ be simple, which can be done, say by taking $\LE$.
See figure \ref{cap:I1uI2d}, right. Let $k^{2,i}$ be the first time
when $\delta^{i}(k^{2,i})\in R^{2}$ and let $m^{2,i}$ be such that
$R^{2}(m^{2,i})=\delta^{i}(k^{2,i})$. Let \[
m^{1,i}:=\max\left\{ m:R^{1}(m)\in\delta^{i}\left[0,k^{2,i}\right[\right\} \quad k^{1,i}:=\min\{ k:\delta^{i}(k)=R^{1}(m^{1,i})\}.\]
 By definition $R^{1}$ always intersects $\delta^{1}[0,k^{2,1}[$
hence $m^{1,1}$ and $k^{1,1}$ are well defined. If $R^{1}$ does
not intersect $\delta^{2}[0,k^{2,2}[$ we consider $m^{1,2}$ to be
$-\infty$ and $k^{1,2}$ to be undefined. As before define $l^{i}$
to be the first time when $R^{i}(l^{i})\not\in B(v,r)$.

We can now define $\gamma^{i}$ by connecting a $\delta^{i}$ to an
$R^{i}$ according to the relation between $m^{1,1}$ and $m^{1,2}$.
In formulas, if $m^{1,1}<m^{1,2}$ define\[
\gamma^{1}:=\LE(\delta^{1}[0,k^{2,1}[\;\cup\; R^{2}[m^{2,1},l^{2}-1])\quad\gamma^{2}:=\LE(\delta^{2}[0,k^{1,2}[\;\cup\; R^{1}[m^{1,2},l^{1}-1])\]
while if $m^{1,1}>m^{1,2}$ (which includes $m^{1,2}=-\infty$) define\[
\gamma^{1}:=\LE(\delta^{1}[0,k^{1,1}[\;\cup\; R^{1}[m^{1,1},l^{1}-1])\quad\gamma^{2}:=\LE(\delta^{2}[0,k^{2,2}[\;\cup\; R^{2}[m^{2,2},l^{2}-1]).\]
In both cases it is easy to verify that $\gamma^{1}\cap\gamma^{2}=\emptyset$
and that \ref{enu:gamisep} holds, if $\kappa$ is sufficiently large,
just like in case 1. Hence the lemma is concluded.
\end{proof}

\section{\label{sec:Brownian-graphs}Isotropic graphs}

\subsection{\label{sub:Brownian}Preliminaries}

In this chapter we will need to compare random walk and Brownian motion.
For definition and basic properties of Brownian motion see any standard
text book, e.g.~\cite{RW94,B95}. 

To avoid confusion with the use of the letter $B$ for a ball in a
metric space, we will denote Brownian motion by $W$, giving homage
to Wiener, even though he seems to have only been interested in the
one dimensional case.

The equivalent of the stopping times $T(X)$ and $T_{v,r}$ will be
denoted by $S$ i.e. $S(X)$ is the time when the Brownian motion
hits $X$ for the first time and $S_{v,r}=S(\partial_{\textrm{cont}}B(v,r))$.
Similarly we shall use $S^{i}(X)$ and $S_{v,r}^{i}$ when we have
more than one Brownian motion involved.

We will also need in a few places the following {}``Hausdorff distance
from a subset'', defined by\[
\subset_{\Haus}(A,B):=\sup_{a\in A}d(a,B)=\inf_{C\subset B}d_{\Haus}(A,C).\]
$\subset_{\Haus}$ is monotone in the sense that if $A_{1}\subset A_{2}$
and $B_{1}\subset B_{2}$ then \[
\subset_{\Haus}(A_{1},B_{2})\leq\;\subset_{\Haus}(A_{2},B_{1}).\]

\subsection{\label{sub:Background-exponent}Background on the non-intersection
exponent}

We shall need some known results about the non-intersection exponent
$\xi_{d}$, so let us start with a quick survey of this topic (mostly
developed by Lawler and coauthors).

\subsubsection*{1}

Let $x$ and y be two points on $\partial B(0,1)$ and let $W^{x}$
and $W^{y}$ be independent Brownian motions in $\mathbb{R}^{d}$
starting from $x$ and $y$ respectively and stopped when hitting
$\partial B(0,r)$. Define the non-intersection probability by \[
p(x,y,r):=\mathbb{P}(W^{x}\cap W^{y}=\emptyset)\quad P(r):=\max_{x,y}p(x,y,r).\]
 The scaling invariance of Brownian motion and the strong Markov property
easily give that $P$ is submultiplicative in $r$ (i.e.~$P(rs)\leq P(r)P(s)$,
or $\log P$ is subadditive in $\log r$) and we get \begin{equation}
P(r)=r^{-\xi_{d}+o(1)}\textrm{ as }r\rightarrow\infty.\label{eq:defxi}\end{equation}
$\xi_{d}$ is the well known non-intersection exponent. The invariance
of Brownian motion to rotations, scaling and translations allows to
map $\vec{1}:=(1,0\dotsc,0)$ and $-\vec{1}$ to the $x$ and $y$
where the maximum occurs and conclude that $p(\vec{1},-\vec{1},r)=r^{-\xi_{d}+o(1)}$.

We will also need a generalization of this quantity: let $W^{x_{i}}$
($i=1,\dotsc,k)$ and $W^{y_{i}}$ ($i=1,\dotsc,l$) be independent
Brownian motions in $\mathbb{R}^{d}$ starting from $x_{1},\dotsc,x_{k}$
and $y_{1},\dotsc,y_{l}$ in $\partial B(0,1)$ respectively and stopped
on $\partial B(0,r)$. Define equivalently \begin{align}
p(x_{1},\dotsc,x_{k},y_{1},\dotsc,y_{l},r) & :=\mathbb{P}\Big(\Big(\bigcup_{i=1}^{k}W^{x_{i}}\Big)\cap\Big(\bigcup_{i=1}^{l}W^{y_{i}}\Big)=\emptyset\Big)\nonumber \\
P(k,l,r) & :=\max_{x_{1},\dotsc,x_{k},y_{1},\dotsc,y_{l}}p(x_{1},\dotsc,x_{k},y_{1},\dotsc,y_{l},r).\label{eq:defPmax}\end{align}
Again, submultiplicativity shows that $P=r^{-\xi_{d}(k,l)+o(1)}$
and $\xi_{d}(k,l)$ is called the $k,l$-nonintersection exponent.
With this notation $\xi_{d}=\xi_{d}(1,1)$. A simple {}``choosing
the best point'' argument shows that the maximum in (\ref{eq:defPmax})
is achieved, up to a factor $\leq kl$ when all the $x_{i}$-s are
the same and all the $y_{i}$-s are the same. And invariance again
shows us that\[
p(\underbrace{\vec{1},\dotsc,\vec{1}}_{k},\underbrace{-\vec{1},\dotsc,-\vec{1}}_{l},r)=r^{-\xi_{d}(k,l)+o(1)}.\]

These $\xi$-s are non-trivial only in dimensions $2$ and $3$. In
dimension $1$ it follows from the {}``gambler ruin problem'' that
$\xi_{1}(k,l)=k+l$ while in dimensions $\geq4$ Brownian motions
never intersect \cite{DEK50} so $\xi\equiv0$. See \cite{L91} for
a more detailed explanation of these facts. Hence from now on we will
only relate to dimensions $2$ and $3$.

\subsubsection*{2}

In \cite{BL90a} it was shows that the same $\xi(k,l)$ hold for the
equivalent problem for random walks (see also \cite{CM91}). Other
relevant variations consider using Brownian motions (or random walks)
$W_{i}^{*}$ with fixed length $t$, or with the length an exponential
variable with expectation $t$ ({}``a random walk with killing rate
$1/t$''). In either case,\begin{equation}
\mathbb{P}\Big(\Big(\bigcup_{i=1}^{k}W_{i}^{x}\Big)\cap\Big(\bigcup_{i=1}^{l}W_{i}^{y}\Big)=\emptyset\Big)=t^{-\xi_{d}(k,l)/2+o(1)}.\label{eq:2zeta}\end{equation}
For example, notice that, if $\tau_{R}$ is the stopping time when
$W$ exits $\partial B(0,r)$, then the probability that either $\tau_{R}>Cr^{2}\log r$
or $\tau_{R}<cr^{2}/\log r$ are negligible, which explains (\ref{eq:2zeta}).

\subsubsection*{3}

In \cite{L89} it was shown that $\xi_{d}(2,1)=4-d$. Very roughly,
the proof uses the fact that two random walks starting from the same
point can be thought of as one bi-directional walk, which allows to
{}``reduce one parameter'' and get an estimate for the probability.
We remark that a similar technique was used in \cite[section 12.5]{L99}
to calculate some intersection exponents for combinations of random
walks and loop-erased random walks.

\subsubsection*{4}

In \cite{BL90b} it was shown that the $\xi(k,l)$ are strictly increasing,
and in particular that $\xi_{3}(1,1)<1$. The proof uses the Wiener
shell test, somewhat like the techniques we will use in chapter \ref{sec:Quasi-loops}.

\subsubsection*{5}

In \cite{L96a} the estimate (\ref{eq:defxi}) was improved to\begin{equation}
P(x,y,r)\approx r^{-\xi(1,1)}\label{eq:xiConly}\end{equation}
i.e.~the error was shown to be in a constant only (for better comparison,
write $P(x,y,r)=r^{-\xi(1,1)+O(1/\log r)}$). Roughly, this follows
by proving {}``supermultiplicativity'' in the sense that $P(rs)\geq cP(r)P(s)$.
This, in turn, follows after proving that two Brownian motions conditioned
not to intersect will also be quite far along the path and in their
end points. In \cite{L96b} this result was extended to simple random
walk via the so-called Skorokhod embedding, a coupling of Brownian
motion and random walk on the same probability space so as to be quite
close.

The analog of \cite{L96a} for general $\xi(k,l)$ was proved in \cite{L98}
while the analog of \cite{L96b} is \cite{LP00}. See also \cite{LSW02b}.

\subsubsection*{6}

Both \cite{L96a} and \cite{L96b} used the estimate (\ref{eq:xiConly})
to prove the existence of many cut times or cut points for random
walk, using relatively straightforward second moment methods. Since
$\xi_{3}(1,1)<1$ we get that the Hausdorff dimension of the cut points
is strictly bigger than $1$, which implies that the set of cut points
of Brownian motion is hittable by a second Brownian motion meaning
that the hitting probability is positive. See \cite[section 12.4]{L99}
for the corresponding calculation for random walk.

\subsubsection*{7}

While we will not use it, it is impossible not to mention that in
dimension $2$ there is a precise formula for $\xi_{2}(n,k)$, conjectured
by Duplantier and Kwon \cite{DK88} and proved by Lawler, Schramm
and Werner \cite{LSW02a}. Both the heuristic arguments and the final
proof depend crucially on the Riemann conformal mapping theorem and
are therefore specifically two dimensional.

\subsection{\label{sub:DefinitionIsotropic}Definition}

Let $G$ be a $d$-dimensional Euclidean net. Let $v\in G$ and $r>0$.
Let $A$ be a $d-1$ dimensional spherical simplex (since we are only
interested in $d=2,3$ we have in effect an arc or a spherical triangle)
in $\partial_{\textrm{cont}}B(v,r)$. Let $|A|$ be $(d-1)$-volume
of $A$ normalized so that $|\partial_{\textrm{cont}}B(v,r)|=1$.
We wish to define discrete versions of $A$. For this purpose, identify
each edge $(v,w)$ of $G$ with the linear segment in $\mathbb{R}^{d}$
between the two vertices, and say that $w\in A^{-}$ if $w\in\partial B(v,r)$
and all edges $(v,w)$, $v\in B(v,r)$ intersect $A$. Say that $w\in A^{+}$
if $w\in\partial B(v,r)$ and some edge $(v,w)$ intersects $A$.
Any set $A^{*}$ between $A^{-}$ and $A^{+}$ will be called a discrete
version of $A$. Denote by $p_{A^{*}}=\mathbb{P}^{v}(R(T_{v,r})\in A^{*})$.
We call $G$ \textbf{isotropic} if \begin{equation}
\big|p_{A^{*}}-|A|\big|\leq Kr^{-\alpha}\quad\forall v,r,A,A^{*}.\label{eq:defbrown}\end{equation}
$K>0$ and $\alpha>0$ are parameters of $G$, so it would be more
precise to call $G$ $(d,\alpha,K)$-isotropic. We will rarely need
to do so, though. As in the previous chapter, when we write $C(G)$
we mean a constant that depends only on the isotropicity parameters
$d,\alpha,K$ and the Euclidean net structure constants (see page
\pageref{page:eucCstruct}), but not on other properties of $G$.
Together we call these numbers the \textbf{isotropicity structure
constants}.

We haven't defined whether we are talking about an open, closed or
other simplex because by expanding or contracting slightly it is obvious
that if (\ref{eq:defbrown}) holds for one than it holds for any and
all. We also remark that by examining triangles intersecting no edge
of $G$, it is obvious that $\alpha\leq d-1$, and if $G$ is a grid
then $\alpha\leq1$. This last inequality tight: it is possible to
show that the grid $\mathbb{Z}^{d}$ is isotropic with $\alpha=1$,
though we will have no use for this fact. It would be interesting
to construct an example in $d\geq3$ of an isotropic graph with $\alpha>1$,
even if one weakens the definition to require that (\ref{eq:defbrown})
holds only for a specific choice of discrete version of $A$.

\subsection{Coupling with Brownian motion}

In this section with shall show how to couple random walk on $G$
with Brownian motion on $\mathbb{R}^{d}$. This will be the main tool
for using isotropic graphs and indeed, it is probably possible to
define isotropic graphs via the coupling. However, we will need some
specific properties of the coupling (see below) that are cumbersome
to formulate.

We will construct the coupled walk and motion by considering a random
walk $R$ on $G$ and constructing an appropriate Brownian motion
$W$. Let therefore $R$ be given and define inductively a sequence
of stopping times, $\tau_{0}=0$ and\begin{equation}
\tau_{i}:=\min\{ t>\tau_{i-1}:|R(t)-R(\tau_{i-1})|>r_{i}\}\quad r_{i}:=i^{4/\alpha}.\label{eq:defri}\end{equation}
The reason behind the choice of $r_{i}$ will become evident later
on, during the proof of lemma \ref{lem:RWclose} --- we remark only
that the connection between $4$ and the dimension $d$ is $4\geq2(d-1)$.
Construct now fixed divisions $\Delta_{i}$ of the sphere $\mathbb{S}^{d-1}$
into $D_{i}:=\left\lfloor r_{i}^{\alpha/2}\right\rfloor +4$ disjoint
spherical simplices of $(d-1)$-normalized volume $\approx1/D_{i}$
and diameter $\approx D_{i}^{-1/(d-1)}$ (associate the boundaries
of the simplices to them as you please --- this is not important).
In two dimensions one may just take $\Delta_{i}$ to be a collection
of (half-closed half-open) arcs of length $1/D_{i}$. In three dimensions
it is an easy geometric exercise to show that such a {}``triangulation''
exists, knowing only that $D_{i}$ is $\geq4$. For every $\delta\in\Delta_{i}$
define $\delta^{*}$ which, unlike $\delta$, may depend on $v$ and
on the walk up to $R(\tau_{i-1})$, to be a discrete version of $\delta$
such that the $\delta^{*}$ cover $\partial B(R(\tau_{i-1}),r_{i})$
and are disjoint%
\footnote{This definition is not unique, but everything will do will not depend
on the choice of which {}``boundary vertex'' to associate with which
$\delta^{*}$. If one prefers a uniquely defined coupling, just order
$\Delta_{i}$ and then associate each boundary vertex to the $\delta^{*}$
first in this order.%
}. Define, \begin{equation}
p_{i,\delta}:=\mathbb{P}(R(\tau_{i})\in\delta^{*}\,|\, R(\tau_{i-1})).\label{eq:defexit}\end{equation}
We get from (\ref{eq:defbrown}) that\begin{equation}
\big|p_{i,\delta}-|\delta|\big|\leq Kr_{i}^{-\alpha}.\label{eq:pidel_areadel}\end{equation}
Define therefore $\eta_{i}:=\min_{\delta}p_{i,\delta}/|\delta|$ and
get \[
1\geq\eta_{i}\geq1-CKr_{i}^{-\alpha/2}.\]
We can now construct $W$, and we shall do so in parts, in parallel
with times $\sigma_{i}$ which would be the analogs of $\tau_{i}$.
Define $W(0):=v$ and $\sigma_{0}:=0$. Assume $R(\tau_{i})\in\delta^{*}$.
Throw a random independent coin $X_{i}$ with probability $\eta_{i}|\delta|/p_{i,\delta}$
for $1$. The definition of $\eta_{i}$ ensures that this number is
$\in[0,1]$. If $X_{i}=1$, define $W'_{i}$ to be a Brownian motion
starting from $0$ and conditioned to exit $B(0,r_{i})$ at $r_{i}\delta$.
If $X_{i}=0$, let $W'_{i}$ be an unconditioned Brownian motion.
In both cases define $\sigma_{i}'$ to be the time when $W_{i}'$
exits $B(0,r_{i})$. Finally define $\sigma_{i}=\sigma_{i-1}+\sigma_{i}'$
and $W$ on the interval $\left]\sigma_{i-1},\sigma_{i}\right]$ by
$W(t):=W(\sigma_{i-1})+W_{i}'(t-\sigma_{i-1})$.

\begin{lem}
The $W$ constructed above is regular Brownian motion.
\end{lem}
\begin{proof}
Since $\mathbb{E}\sigma_{i}'>c>0$ and they are independent we get
that almost surely $\sum\sigma_{i}=\infty$ and hence $W$ is an almost
surely well defined function $[0,\infty[\to\mathbb{R}^{d}$. Now compare
$W$ to a regular Brownian motion $W^{*}$. Let $\sigma_{i}^{*}$
be stopping times defined by\[
\sigma_{i}^{*}=\inf\{ t>\sigma_{i-1}^{*}:W^{*}(t)\not\in B(W^{*}(\sigma_{i-1}^{*},r_{i}))\}.\]
Using the strong Markov property \cite[page 21]{RW94} inductively
gives that $W^{*}(t-\sigma_{i-1}^{*})-W^{*}(\sigma_{i-1}^{*})$ is
distributed like Brownian motion starting from $0$ and stopped when
exiting $B(0,r_{i})$. On the other hand, it follows from the definition
that each $W_{i}'$ has probability $\eta_{i}|\delta|$ to be a Brownian
motion conditioned to hit $r_{i}\delta$ (for every $\delta\in\Delta_{i}$)
and probability $1-\eta_{i}(\sum p_{i})$ to be unconditioned, hence
$W_{i}'$ is also a regular Brownian motion starting from $0$ and
stopped on $\partial B(0,r_{i})$. Hence $W[0,\sigma_{i}]\sim W^{*}[0,\sigma_{i}^{*}]$
for all $i$. Taking limit as $i\to\infty$ shows that $W\sim W^{*}$.
\end{proof}
\begin{lem}
\label{lem:RWclose}Let $G$ be an isotropic graph and $v\in G$.
Let $R$ and $W$ be the coupled walk and motion starting from $v$.
Let $r_{i}$, $\tau_{i}$ and $\sigma_{i}$ be as in the definition
of the coupling (\ref{eq:defri}). Then\[
\mathbb{P}(\exists j\leq i:|R(\tau_{j})-W(\sigma_{j})|\geq\lambda r_{i})\leq C(G)\exp(-c(G)\lambda)\]
for any $\lambda>0$.
\end{lem}
\begin{proof}
We use the notation $X_{i}$ from the definition of the coupling.
If $X_{j}=1$ for some $j$ then we get\[
R(\tau_{j-1})-R(\tau_{j})\in r_{j}\delta+B(0,C(G))\quad W(\sigma_{j-1})-W(\sigma_{j})\in r_{j}\delta\]
Hence\[
|R(\tau_{j-1})-R(\tau_{j})-W(\sigma_{j-1})+W(\sigma_{j})|\leq Cr_{j}^{1-\alpha/2(d-1)}+C(G),\]
and since $d\leq3$ we may simply write $\leq Cr_{j}^{1-\alpha/4}+C$.
Summing we get that if $X_{j}=1$ for all $j\leq i$ then\[
|R(\tau_{i})-W(\sigma_{i})|\leq C\sum_{j=1}^{i}j^{4/\alpha-1}+C\leq C(G)i^{4/\alpha}=C(G)r_{i}.\]
Hence we need to estimate \[
\Sigma:=\sum_{j:X_{j}=0}|R(\tau_{j-1})-R(\tau_{j})-W(\sigma_{j-1})+W(\sigma_{j})|\leq\sum_{j:X_{j}=0}2r_{j}+C.\]
Divide this sum into blocks \[
\Sigma_{k}:=\!\!\!\sum_{\substack{j:X_{j}=0\\
2^{k-1}<r_{j}\leq2^{k}}
}\!\!2r_{j}+C.\]

Now, each $\Sigma_{k}$ contains $\leq C(G)2^{k\alpha/4}$ summands,
and each summand is zero with probability $\geq1-C(G)2^{-k\alpha/2}$
independently so a very rough estimate gives\[
\mathbb{P}(\Sigma_{k}>\lambda2^{k})\leq C(G)\exp(-c(G)\lambda).\]
Define $l:=\left\lceil \log_{2}r_{j}\right\rceil $ and sum over $k$
from $0$ to $l$ to get \begin{align*}
\mathbb{P}\left(\Sigma>\lambda r_{i}\right) & \stackrel{(*)}{\leq}\mathbb{P}(\exists k:\Sigma_{k}>c\lambda2^{(l+k)/2})\leq C(G)\sum_{k}\exp(-c(G)\lambda2^{(l-k)/2})\\
 & \leq C(G)\exp(-c(G)\lambda).\end{align*}
$(*)$ comes from the fact that if $\Sigma_{k}\leq\lambda2^{(k+l)/2}$
for all $k$ then $\Sigma=\sum\Sigma_{k}\leq(2+\sqrt{2})\lambda2^{l}\leq2(2+\sqrt{2})\lambda r_{i}$
so one may take $c=1/2(2+\sqrt{2})$ on the right hand side of $(*)$.
Since $\Sigma\leq\lambda r_{i}$ implies $|R(\tau_{j})-W(\sigma_{j})|\leq(\lambda+C(G))r_{i}$
for all $j\leq i$, the lemma is proved.
\end{proof}
Lemma \ref{lem:RWclose} is not really convenient to use as is, because
one needs to relate $i$ to more natural events. Here is one such
useful relation:

\begin{lem}
\label{lem:rismall}Let $G$ be an isotropic graph and let $v\in G$.
Let $R$ and $W$ be the coupled walk and motion starting from $v$.
Let $r>1$ and let $T=T_{v,r}$ and $S=S_{v,r}$. Let $r_{i}$, $\tau_{i}$
and $\sigma_{i}$ be as in the definition of the coupling. Let\newc{c:RWclose}
\[
I:=\min\{ i:\tau_{i}\geq T,\,\sigma_{i}\geq S\}.\]
 Then for some constant $c_{\ref{c:RWclose}}(G)$,\[
\mathbb{P}(r_{I}>\lambda r^{1-c_{\ref{c:RWclose}}(G)})\leq C(G)\exp(-\lambda^{c(G)}).\]

\end{lem}
\begin{proof}
Since $W(\sigma_{i})-W(\sigma_{i-1})$ are independent variables with
mean zero and variance $cr_{i}$ we get (say by second moment methods)
that for some $c>0$, \[
\mathbb{P}\Big(|W(\sigma_{j})-W(\sigma_{i})|>c\Big(\sum_{k=i}^{j}r_{k}^{2}\Big)^{1/2}\Big)>c\quad\forall j>i.\]
Lemma \ref{lem:RWclose} allows us to replace $W$ with $R$: we get
that for some constants $\mu=\mu(G)$ and $\nu=\nu(G)$,\[
\mathbb{P}\Big(|R(\tau_{j})-R(\tau_{i})|>\mu\Big(\sum_{k=i}^{j}r_{k}^{2}\Big)^{1/2}\Big)>c\quad\forall j>i+\nu.\]
In particular, if $r_{i}\sqrt{j-i}>(2/\mu)r$ then $\mathbb{P}(R(\tau_{j})\not\in B(v,r))>c$
for any $R(\tau_{i})\in B(v,r)$. Hence if we define $J:=C(G)r^{1/(4/\alpha+1/2)}$
for some $C$ sufficiently large we get both $r_{J}\sqrt{J}>(2/\mu)r$
as well as $J>\nu$. Hence \[
\mathbb{P}(R(\tau_{(n+1)J})\not\in B(v,r)\,|\, R[0,\tau_{nJ}],\, R(\tau_{nj})\in B(v,r))>c\quad\forall n>1\]
and hence $\mathbb{P}(T>\tau_{nJ})\leq Ce^{-cn}$. An identical calculation
shows that $\mathbb{P}(S>\sigma_{nJ})\leq Ce^{-cn}$. Hence $\mathbb{P}(r_{I}>r_{nJ})\leq Ce^{-cn}$
and since\[
r_{nJ}=(nJ)^{4/\alpha}=n^{C(G)}r^{1-c(G)}\]
the lemma is proved. 
\end{proof}
\begin{cor*}
With the notations of lemma \ref{lem:rismall}, \[
\mathbb{P}(\exists j\leq I:d_{\Haus}(R[0,\tau_{j}],W[0,\sigma_{j}])>\lambda r^{1-c_{\ref{c:RWclose}}})\leq C(G)\exp(-\lambda^{c(G)}).\]

\end{cor*}
\begin{proof}
Clearly we may assume $\lambda>5$. Let $i_{*}$ by the maximal $i$
such that $r_{i_{*}}\leq\sqrt{\lambda}r^{1-c_{\ref{c:RWclose}}}$.
Then lemma \ref{lem:RWclose} shows that \[
\mathbb{P}(\exists j\leq i_{*}:|R(\tau_{j})-W(\sigma_{j})|\geq(\sqrt{\lambda}-2)r_{i_{*}})\leq C(G)\exp(-c(G)(\sqrt{\lambda}-2)).\]
Now, the point $R(\tau_{j})$ are an approximation (in the Hausdorff
distance) of the entire path, i.e. \[
d_{\Haus}(R[0,\tau_{j}],\{ R(0),\dotsc,R(\tau_{j})\})\leq r_{j}\]
and similarly for $W$. Thus we get\[
\mathbb{P}(\exists j\leq i_{*}:d_{\Haus}(R[0,\tau_{j}],W[0,\sigma_{j}])\geq\sqrt{\lambda}r_{i_{*}})\leq C(G)\exp(-c(G)(\sqrt{\lambda}-2))\]
and from the definition of $i_{*}$, \[
\mathbb{P}(\exists j\leq i_{*}:d_{\Haus}(R[0,\tau_{j}],W[0,\sigma_{j}])\geq\lambda r^{1-c_{\ref{c:RWclose}}})\leq C(G)\exp(-c(G)(\sqrt{\lambda}-2))\]
Estimating the probability that $I>i_{*}$ using lemma \ref{lem:rismall}
proves the corollary.
\end{proof}
\begin{lem}
\label{lem:subhaus}Let $G$ be an isotropic graph, let $\nu>1$ and
let $v\in G$. Let $R$ and $W$ be the coupled walk and motion starting
from $v$. Let $T_{r}=T_{v,r}$ and $S_{r}=S_{v,r}$. Then for all
$r>\max1,s$, and for all $\lambda>0$,\begin{align}
\mathbb{P}(\sh(R[T_{s},T_{r}],W[S_{s/\nu},S_{\nu r}])>\lambda r^{1-c_{\ref{c:RWclose}}(G)}) & \leq C(\nu,G)\exp(-\lambda^{c(G)})\label{eq:RinW2}\\
\mathbb{P}(\sh(W[S_{s},S_{r}],R[T_{s/\nu},T_{\nu r}])>\lambda r^{1-c_{\ref{c:RWclose}}(G)}) & \leq C(\nu,G)\exp(-\lambda^{c(G)}).\label{eq:WinR2}\end{align}

\end{lem}
We explicitly include the case $s=0$ in which case we define $T_{0}=S_{0}=0$.

\begin{proof}
Let us prove (\ref{eq:RinW2}). Let $r_{i}$, $\tau_{i}$ and $\sigma_{i}$
be as in the definition of the coupling (\ref{eq:defri}). Define\[
I_{1}:=\max\{ i:\tau_{i}<T_{s}\},\quad I_{2}:=\min\{ i:\tau_{i}\geq T_{r}\}.\]
 Then lemma \ref{lem:rismall} gives that\begin{align}
\mathbb{P}(r_{I_{1}}>\lambda s^{1-c_{\ref{c:RWclose}}}) & \leq C(G)\exp(-\lambda^{c(G)}),\label{eq:flem22s}\\
\mathbb{P}(r_{I_{2}}>\lambda r^{1-c_{\ref{c:RWclose}}}) & \leq C(G)\exp(-\lambda^{c(G)}).\label{eq:flem22r}\end{align}
Denote by $i_{1}$ the last $i$ such that $r_{i}\leq\lambda s^{1-c_{\ref{c:RWclose}}}$
and by $i_{2}$ the last $i$ such that $r_{i}\leq\lambda r^{1-c_{\ref{c:RWclose}}}$.
Lemma \ref{lem:RWclose} shows that \begin{align}
\mathbb{P}(\exists j\leq i_{k}:|R(\tau_{j})-W(\sigma_{j})|\geq\lambda r_{i_{k}}) & \leq C(G)\exp(-c(G)\lambda)\quad k=1,2.\label{eq:flem23rs}\end{align}
Together with the estimates of $r_{I_{k}}$ this gives the following
corollary:\begin{align}
\mathbb{P}(d_{\Haus}(R[\tau_{I_{1}},\tau_{I_{2}}],W[\sigma_{I_{1}},\sigma_{I_{2}}])\geq\lambda^{2}r^{1-c_{\ref{c:RWclose}}}) & \leq C(G)\exp(-\lambda^{c(G)}).\label{eq:RtauI1I2W}\end{align}
We can replace $\lambda^{2}$ with $\lambda$ on the left hand side
paying only in the constant inside the exponent on the right hand
side.

\medskip{}
\noindent \textbf{case 1:} If $s<r^{1-c_{\ref{c:RWclose}}}$ then
the fact that $|W(S_{s/\nu})-R(T_{s})|\leq C(G)+2s$ shows that \[
\sh(R[T_{s},T_{r}],W[0,S_{\nu r}])\geq\sh(R[T_{s},T_{r}],W[S_{s/\nu},S_{\nu r}])-(C(G)+2s)\]
which allows to estimate\begin{align*}
\lefteqn{\mathbb{P}(\sh(R[T_{s},T_{r}],W[S_{s/\nu},S_{\nu r}])>\lambda r^{1-c_{\ref{c:RWclose}}})\leq}\quad\\
 & \leq\mathbb{P}(\sh(R[T_{s},T_{r}],W[0,S_{\nu r}])>(\lambda-C(G))r^{1-c_{\ref{c:RWclose}}})\leq\\
 & \leq\mathbb{P}(\sh(R[\tau_{I_{1}},\tau_{I_{2}}],W[0,\sigma_{I_{2}}])>(\lambda-C(G))r^{1-c_{\ref{c:RWclose}}})+\mathbb{P}(\sigma_{I_{2}}>S_{\nu r})\leq\\
 & \!\!\stackrel{(\ref{eq:RtauI1I2W})}{\leq}\! C(G)\exp(-(\lambda-C(G))^{c(G)})+\mathbb{P}(\sigma_{I_{2}}>S_{\nu r}).\end{align*}
Now, to estimate $\mathbb{P}(\sigma_{I_{2}}>S_{\nu r})$ we use the
fact that $|R(\tau_{i})-v|\leq r$ for any $i<I_{2}$ and get \begin{alignat*}{2}
\mathbb{P}(\sigma_{I_{2}}>S_{\nu r}) &  & \leq\;\;\: & \mathbb{P}(\exists j\leq i_{2}:|W(\sigma_{j})-R(\tau_{j})|\geq(\nu-1)r-\lambda r^{1-c_{\ref{c:RWclose}}})+\mathbb{P}(I_{2}>i_{2})\\
 &  & \stackrel{(\ref{eq:flem22r},\ref{eq:flem23rs})}{\leq} & C(G)\exp(-c(G)((\nu-1)\lambda^{-1}r^{c_{\ref{c:RWclose}}}-1))+C(G)\exp(-\lambda^{c(G)}).\end{alignat*}
Hence, if $\lambda<r^{c(G)}$, (\ref{eq:RinW2}) is proved.

\medskip{}
\noindent \textbf{case 2:} If $s\geq r^{1-c_{\ref{c:RWclose}}}$ then
we estimate\begin{align*}
\lefteqn{\mathbb{P}(\sh(R[T_{s},T_{r}],W[S_{s/\nu},S_{\nu r}])>\lambda r^{1-c_{\ref{c:RWclose}}})\leq}\\
 & \quad\leq\mathbb{P}(\sh(R[\tau_{I_{1}},\tau_{I_{2}}],W[\sigma_{I_{1}},\sigma_{I_{2}}])>\lambda r^{1-c_{\ref{c:RWclose}}})+\mathbb{P}(\sigma_{I_{1}}<S_{s/\nu})+\mathbb{P}(\sigma_{I_{2}}>S_{\nu r})\\
 & \quad\!\!\stackrel{(\ref{eq:RtauI1I2W})}{\leq}\! C(G)\exp(-\lambda^{c(G)})+\mathbb{P}(\sigma_{I_{1}}<S_{s/\nu})+\mathbb{P}(\sigma_{I_{2}}>S_{\nu r}).\end{align*}
Now, the estimate of $\mathbb{P}(\sigma_{I_{2}}>S_{\nu r})$ is as
in case 1. The estimate of $\mathbb{P}(\sigma_{I_{1}}<S_{s/\nu})$
is similar. Using $|R(\tau_{I_{1}})-v|\geq s-r_{I_{1}}-C(G)$, we
get \begin{align*}
\lefteqn{\mathbb{P}(\sigma_{I_{1}}<S_{s/\nu})\leq\mathbb{P}(|W(\sigma_{I_{1}})-v|<s/\nu)\leq}\\
 & \quad & \leq\;\;\: & \mathbb{P}(\exists j\leq i_{1}:|W(\sigma_{j})-R(\tau_{j})|\geq s(1-1/\nu)-\lambda s^{1-c_{\ref{c:RWclose}}}-C(G))+\\
 &  &  & \qquad+\mathbb{P}(I_{1}>i_{1})\leq\\
 &  & \stackrel{(\ref{eq:flem22s},\ref{eq:flem23rs})}{\leq} & C(G)\exp(-c(G)((1-1/\nu)\lambda^{-1}s^{c_{\ref{c:RWclose}}}-C(G))+C(G)\exp(-\lambda^{c(G)})\end{align*}
and again, if $\lambda\leq s^{c(G)}$, (\ref{eq:RinW2}) is proved.
Since in this case $s>r^{c(G)}$ it is enough to assume $\lambda\leq r^{c(G)}$
in order to get $\lambda\leq s^{c(G)}$ and consequently (\ref{eq:RinW2}).

\medskip{}
\noindent \textbf{case 3:}\newc{c:lamr}The previous calculations
proved the case $\lambda\leq r^{c_{\ref{c:lamr}}}$ for some $c_{\ref{c:lamr}}(G)$.
However, this implies that for any $\lambda\leq2r^{c_{\ref{c:RWclose}}}$,\[
\mathbb{P}(\sh(R[T_{s},T_{r}],W[S_{s/\nu},S_{\nu r}])>\lambda r^{1-c_{\ref{c:RWclose}}(G)})\leq C(\nu,G)\exp\left(-\left({\textstyle \frac{1}{2}}\lambda\right)^{c(G)c_{\ref{c:lamr}}/c_{\ref{c:RWclose}}}\right)\]
or in other words, (\ref{eq:RinW2}) holds with different constants
on the right hand side. However, for $\lambda>2r^{c_{\ref{c:RWclose}}}$,
(\ref{eq:RinW2}) holds trivially because \[
\sh(R[T_{s},T_{r}],W[S_{s/\nu},S_{\nu r}])\leq d_{\Haus}(R[T_{s},T_{r}],W[S_{s/\nu},S_{r}])\leq2r\]
so the probability in (\ref{eq:RinW2}) is zero. This finishes the
proof of (\ref{eq:RinW2}). The proof of (\ref{eq:WinR2}) is identical.
\end{proof}
\begin{cor*}
With the notations of lemma \ref{lem:subhaus}, if $R$ and $W$ start
from a $w$, $|v-w|\leq\frac{1}{4}s(1-1/\nu)$, then (\ref{eq:RinW2})
and (\ref{eq:WinR2}) still hold, possibly with different constants.
Further, this holds if $s=0$ and $|v-w|\leq\frac{1}{4}r(\nu-1)$.
\end{cor*}
This follows from lemma \ref{lem:subhaus} and the monotonicity of
$\sh$.

\subsection{Hitting of small balls}

From now on we will prove {}``natural'' facts about walk on isotropic
graphs, natural in the sense that they don't need the coupling (or
other special notations) to be stated. In this section we shall prove
two lemmas about the hitting probability of {}``intermediate scale''
objects, i.e. of the size $r^{1-c}$ (both will be used in chapter
\ref{sec:Isotropic-gluing}). In the next sections we shall focus
on more delicate facts.

\begin{lem}
\label{lem:sball}Let $G$ be an isotropic graph, and let $\epsilon>0$.
Let $r>C(\epsilon,G)$ and let $v,w\in B(x,r(1-\epsilon))$, $|v-w|>\frac{1}{2}\epsilon r$.
Let \[
p:=\mathbb{P}^{v}(T_{w,s}<T_{x,r})\quad q:=\mathbb{P}^{v}(S_{w,s}<S_{x,r}).\]
Where $r^{1-c_{\ref{c:RWclose}}/2}\leq s\leq|v-w|-\frac{1}{2}\epsilon r$,
$c_{\ref{c:RWclose}}(G)$ from lemma \ref{lem:rismall}. Then\[
\left|p-q\right|\leq C(\epsilon,G)r^{-c(G)}\max p,q.\]

\end{lem}
\begin{proof}
The first step is to get a simple lower bound for $q$. In the case
$|v-w|\leq\frac{3}{4}d(w,\partial_{\textrm{cont}}B(x,r))$ then a
calculation using the continuous analog of (\ref{eq:fEfRT}) with
the Newtonian potential \cite[chapter II 3]{B95} around $w$ shows
that\begin{equation}
q\geq c(\epsilon,d)\begin{cases}
s/r & d=3\\
1/\log(r/s) & d=2\end{cases}\label{eq:lowq}\end{equation}
($d$ being the dimension). Removing the condition $|v-w|\leq\frac{3}{4}d(w,\partial B(0,r))$,
we still have (\ref{eq:lowq}), perhaps with a different constant.
Indeed, using the continuous Harnack inequality \cite[chapter II 1]{B95}
for the domain $B(x,r)\setminus B(w,s+\frac{1}{4}\epsilon r)$  shows
(\ref{eq:lowq}) for all $w$ and $v$. 

Next define \[
q^{\pm}=\mathbb{P}^{v}(S_{w,s^{\pm}}<S_{x,r})\quad s^{\pm}=s\pm r^{1-(3/4)c_{\ref{c:RWclose}}}\]
so $q^{-}\leq q\leq q^{+}$. Now, the strong Markov property gives
us that\[
q-q^{-}\leq\mathbb{P}^{v}\left(S_{w,s}<S_{x,r}\right)\cdot\max_{y\in\partial_{\textrm{cont}}B(w,s)}\mathbb{P}^{y}\left(W[0,S_{x,r}]\cap B(w,s^{-})=\emptyset\right)\]
and again, similar calculations with the continuous Newtonian potential
gives, for any $y\in\partial_{\textrm{cont}}B(w,s)$, \[
\mathbb{P}^{y}\left(W[0,S_{x,r}]\cap B(w,s^{-})=\emptyset\right)\leq C(\epsilon,d)\frac{r^{1-(3/4)c_{\ref{c:RWclose}}}}{s}\leq Cr^{-c_{\ref{c:RWclose}}/4}\]
So we get $q-q^{-}\leq Cqr^{-c}$. A similar calculation shows that
$q^{+}-q\leq Cq^{+}r^{-c}$ and for $r$ sufficiently large we may
write $q^{+}-q\leq Cqr^{-c}$. 

To extract from these inequalities conclusions about $|q-p|$, couple
$R$ and $W$ as above. Let $r_{i}$, $\tau_{i}$ and $\sigma_{i}$
be as in the definition of the coupling (\ref{eq:defri}). Let\[
I:=\min\{ i:\tau_{i}\geq T_{v,r},\,\sigma_{i}\geq S_{v,r}\}.\]
 When comparing $p$ to $q^{+}$ we need to consider two cases: the
first that $R[0,\tau_{I}]$ and $W[0,\sigma_{I}]$ are not very close;
and the second is as in figure \ref{cap:RBnWB}. %
\begin{figure}
\input{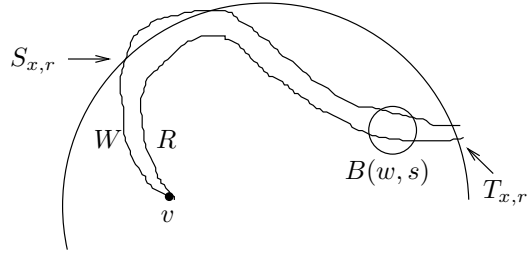}

\caption{$R[0,T_{x,r}]$ intersects the ball $B(w,s)$ but $W[0,S_{x,r}]$
is quite far from it.\label{cap:RBnWB}}
\end{figure}
In a formula: if $\lambda$ is the length of the longest edge in $G$
then\begin{align}
p-q^{+} & \leq\mathbb{P}\left(\left\{ R[0,T_{x,r}]\cap\partial B(w,s)\neq\emptyset\right\} \cap\left\{ W[0,S_{x,r}]\cap\partial B(w,s^{+})=\emptyset\right\} \right)\leq\nonumber \\
 & \leq\mathbb{P}(d_{\Haus}(R[0,\tau_{I}],W[0,\sigma_{I}])\geq r^{1-(3/4)c_{\ref{c:RWclose}}}-\lambda)\;+\nonumber \\
 & \qquad\mathbb{P}(W[S_{x,r},S_{x,r+r^{1-(3/4)c_{\ref{c:RWclose}}}}]\cap B(w,s^{+})\neq\emptyset).\label{eq:WSAnn}\end{align}
Now, the corollary to lemma \ref{lem:rismall} for the ball $B(v,2r)$
gives, since $B(x,r)\subset B(v,2r)$, that\[
\mathbb{P}(d_{\Haus}(R[0,\tau_{I}],W[0,\sigma_{I}])\geq r^{1-(3/4)c_{\ref{c:RWclose}}}-\lambda)\leq C(G)\exp(-r^{c(G)}).\]
while for the second summand we have from the strong Markov property,\begin{eqnarray*}
\lefteqn{\mathbb{P}(W[S_{x,r},S_{x,r+r^{1-(3/4)c_{\ref{c:RWclose}}}}]\cap B(w,s^{+})\neq\emptyset)\leq}\\
 & \qquad & \leq\max_{y\in\partial_{\textrm{cont}}B(x,r)}\mathbb{P}^{y}(W[0,S_{x,r+r^{1-(3/4)c_{\ref{c:RWclose}}}}]\cap B(x,r(1-\epsilon/4))\neq\emptyset)\cdot\\
 &  & \qquad\max_{y\in\partial_{\textrm{cont}}B(x,r(1-\epsilon/4))}\mathbb{P}^{y}(W[0,S_{x,r+r^{1-(3/4)c_{\ref{c:RWclose}}}}]\cap B(w,s^{+})\neq\emptyset)\leq\\
 &  & \stackrel{(*)}{\leq}Cr^{-(3/4)c_{\ref{c:RWclose}}}\begin{cases}
s/r & d=3\\
1/\log(r/s) & d=2\end{cases}\leq Cr^{-c}q.\end{eqnarray*}
Both estimates of $(*)$ follow from the continuous Newtonian potential.
Note that we assumed here that $\frac{1}{8}\epsilon r>r^{1-(3/4)c_{\ref{c:RWclose}}}$
which we may, if $r$ is sufficiently large. This finishes the proof
that $p\leq q+Cr^{-c}q$.

The proof of the other direction is similar. We have\begin{align*}
q^{-}-p & \leq\mathbb{P}\left(\left\{ W[0,S_{x,r}]\cap\partial B(w,s^{-})\neq\emptyset\right\} \cap\left\{ R[0,T_{x,r}]\cap\partial B(w,s)=\emptyset\right\} \right)\leq\\
 & \leq\mathbb{P}(d_{\Haus}(R[0,\tau_{I}],W[0,\sigma_{I}])\geq r^{1-(3/4)c_{\ref{c:RWclose}}})\;+\\
 & \qquad\mathbb{P}(W[S_{x,r-r^{1-(3/4)c_{\ref{c:RWclose}}}},S_{x,r}]\cap B(w,s^{-})\neq\emptyset)\end{align*}
and an identical calculation finishes this case, and the lemma.
\end{proof}
\begin{lem}
\label{lem:striag}Let $G$ be an isotropic graph and let $\epsilon>0$.
Let $r>C(\epsilon,G)$ and let $v\in B(x,(1-\epsilon)r)$. Let $\Delta$
be a spherical triangle on $\partial_{\textrm{cont}}B(x,r)$ and let
$\Delta^{*}$ be a discrete version of it. Let $q:=\mathbb{P}^{v}(W(S_{x,r})\in\Delta)$
and $p:=\mathbb{P}^{v}(R(T_{x,r})\in\Delta^{*})$. Then\newc{c:striag}\[
\left|p-q\right|\leq C(G,\epsilon)r^{-c_{\ref{c:striag}}(G)}.\]

\end{lem}
\noindent (the only difference between lemma \ref{lem:striag} and
the definition of an isotropic graph is that here the starting point
of the walk $v$ might be different from the center of the stopping
ball $x$). The proof is very similar to the proof of the previous
lemma, so we indicate only the differences. We define $q^{\pm}$ as
the probabilities of $W$ to hit $\partial B(x,r)$ at $\Delta^{\pm}$
where\begin{align*}
\Delta^{+} & :=\left(\Delta+B(0,r^{1-c_{\ref{c:RWclose}}/2})\right)\cap\partial_{\textrm{cont}}B(x,r)\\
\Delta^{-} & :=\left\{ x\in\Delta:B(x,r^{1-c_{\ref{c:RWclose}}/2})\cap\partial_{\textrm{cont}}B(x,r)\subset\Delta\right\} .\end{align*}
 The proof that $q^{+}-q$, $q-q^{-}\leq Cr^{-c}$ is direct calculation
for $v=x$, and for general $v$ follows from the continuous Harnack
inequality. The proof that $p-q^{+}\leq Cr^{-c}$ is similar, except
the last term on the right should be replaced, for example, with\begin{multline*}
\mathbb{P}\left(\left\{ W[0,S_{x,r}]\cap\left(\Delta+B(0,r^{1-(3/4)c_{\ref{c:RWclose}}})\right)\neq\emptyset\right\} \cap\left\{ W(S_{x,r})\not\in\Delta^{+}\right\} \right)\leq\\
\leq\max_{d(y,\Delta)\leq r^{1-(3/4)c_{\ref{c:RWclose}}}}\mathbb{P}^{y}(W(S_{x,r})\not\in\Delta^{+})\leq Cr^{-c}.\end{multline*}
The proof that $q^{-}-p\leq Cr^{-c}$ and the rest of the lemma are
similar.\qed

\subsection{Escape probabilities}

In the section we move from the {}``intermediate scale'' objects
of the previous section to single points. This is more delicate, and
we shall employ techniques similar to those of Lawler \cite{L96b}.
Our main goal is theorem \ref{lem:escape}, but first we need to state
and prove two simple claims.

Henceforth $R$ and $W$ will always be a random walk and a Brownian
motion coupled as above.

\begin{lem}
\label{lem:Brownescape}Let $H\subset\mathbb{R}^{d}$ be a half space.
Let $W$ be a Brownian motion on $\mathbb{R}^{d}$ starting from some
vertex $v\not\in\partial_{\textrm{cont}}H$. Let $r>2d(v,\partial_{\textrm{cont}}H)$.
Then\[
p:=\mathbb{P}(S_{v,r}<S(\partial_{\textrm{cont}}H))\approx\frac{d(v,\partial_{\textrm{cont}}H)}{r}.\]

\end{lem}
\begin{proof}
By translation, scaling and rotation invariance we may assume that
$r=1$, that $v=\epsilon e_{1}$ for some $\frac{1}{2}>\epsilon>0$
($e_{1}$ being the first basis element) and that $H=\left\{ x:\left\langle x,e_{1}\right\rangle >0\right\} $.
Examine two positive harmonic functions: $f(x)=\mathbb{P}^{x}(S_{v,r}<S(\partial H_{i}))$
and $g(x)=\left\langle x,e_{1}\right\rangle $. Both $f$ and $g$
are zero on $\partial_{\textrm{cont}}H$ so the boundary Harnack principle
for Lipschitz domains \cite[theorem III.1.2, page 178]{B95} shows
that $g(v)/f(v)\approx g(x_{0})/f(x_{0})$ where $x_{0}$ is any reference
point. But $g(x_{0})/f(x_{0})$ is just a number, and the lemma is
proved.
\end{proof}
\begin{lem}
\label{lem:Ers}Let $G\subset\mathbb{R}^{d}$ be an isotropic graph.
Let $H\subset\mathbb{R}^{d}$ be a half space. Let $v\in G\setminus\partial_{\textrm{cont}}H$.
Let $\frac{1}{2}\rho>s>r>2d(v,\partial_{\textrm{cont}}H)$. Let $R$
and $W$ be coupled walk and motion starting from $v$. Let $\mathcal{E}$
be an event depending on $R[0,T_{v,r}]$ and $W[0,S_{v,r}]$ only.
Then\[
\mathbb{P}(\mathcal{E}\cap\{ W[S_{v,r},S_{v,s}]\cap\partial_{\textrm{cont}}H=\emptyset\})\leq C\frac{r}{s}(\mathbb{P}(\mathcal{E})+C(G)\exp(-r^{c(G)})).\]

\end{lem}
\begin{proof}
We may assume $s>4r$ and (using lemma \ref{lem:Brownescape}) $r>1$.
Let $r_{i}$, $\tau_{i}$ and $\sigma_{i}$ be as in the definition
of the coupling. Define $I:=\min\{ i:\tau_{i}\geq T_{v,r},\,\sigma_{i}\geq S_{v,r}\}$
and examine $W$ after $\sigma_{I}$. This is a regular Brownian motion,
starting from $W(\sigma_{I})$ and \emph{independent of both $W[0,\sigma_{I}]$
and $R[0,\tau_{I}]$.} Therefore, if we denote by $\mathcal{F}$ the
event $W[\sigma_{I},S_{v,s}]\cap\partial_{\textrm{cont}}H=\emptyset$
we get from lemma \ref{lem:Brownescape}, \begin{equation}
\mathbb{P}(\mathcal{F}\,|\, W(\sigma_{I}))\leq C\frac{d(W(\sigma_{I}),\partial_{\textrm{cont}}H)}{d(W(\sigma_{I}),\partial_{\textrm{cont}}B(v,s))}.\label{eq:S1SHWsigI}\end{equation}
Actually, this holds only if $W(\sigma_{I})\in B(v,s)$ --- in the
other case we will simply estimate $\mathbb{P}\leq1$. This shows
that, \begin{equation}
\mathbb{P}(\mathcal{E}\cap\mathcal{F}\cap\{|W(\sigma_{I})-v|\leq2r\})\leq C\mathbb{P}(\mathcal{E})\frac{3r}{s-2r}\leq C\mathbb{P}(\mathcal{E})\frac{r}{s}.\label{eq:smallWv}\end{equation}
 For the case that $|W(\sigma_{I})-v|>2r$ we use the corollary to
lemma \ref{lem:rismall} to get $\mathbb{P}(|R(\tau_{I})-W(\sigma_{I})|>\lambda r)\leq C\exp(-c(\lambda r)^{c}).$
By the definition of $I$, either $|W(\sigma_{I})-v|=r$ or $|R(\tau_{I})-v|\leq r+C(G)$
so we get \begin{equation}
\mathbb{P}(|W(\sigma_{I})-v|>\lambda r)\leq C(G)\exp(-(\lambda r)^{c(G)}).\label{eq:WsigIvlog}\end{equation}
Hence we get that \begin{eqnarray}
\lefteqn{\mathbb{P}(\{|W(\sigma_{I})-v|>2r\}\cap\mathcal{F})=}\nonumber \\
 & \qquad & =\sum_{\lambda=3}^{\infty}\mathbb{P}(\{|W(\sigma_{I})-v|\in\left](\lambda-1)r,\lambda r\right]\}\cap\mathcal{F})\nonumber \\
 &  & \stackrel{(*)}{\leq}C(G)\sum_{\lambda=3}^{\infty}\exp(-(\lambda r)^{c(G)})\cdot\begin{cases}
2(\lambda+1)r/s & \lambda r<\frac{1}{2}s\\
1 & \textrm{otherwise}\end{cases}\nonumber \\
 &  & \leq C(G)\exp(-r^{c(G)})\frac{r}{s}\label{eq:largeWv}\end{eqnarray}
where the inequality $(*)$ comes from (\ref{eq:WsigIvlog}) for the
left multiplicand and (\ref{eq:S1SHWsigI}) for the right multiplicand.
Combining (\ref{eq:smallWv}) and (\ref{eq:largeWv}) ends the lemma.
\end{proof}
\begin{thm}
\label{lem:escape}\newC{C:mindvH}Let $G\subset\mathbb{R}^{d}$
be an isotropic graph. Let $H\subset\mathbb{R}^{d}$ be a half space.
Then there exists a constant $C_{\ref{C:mindvH}}(G)$ such that for
any such $H$, any $v$ with $d(v,\partial_{\textrm{cont}}H)>C_{\ref{C:mindvH}}$
and any $r>2d(v,\partial_{\textrm{cont}}H)$ one has that the probability
$p$ that a random walk $R$ on $G$ starting from $v$ will hit $\partial B(v,r)$
before hitting $\partial H$ satisfies \[
p\approx\frac{d(v,\partial_{\textrm{cont}}H)}{r}.\]

\end{thm}
Remember that the constants implicit in the $\approx$ above, like
$C_{\ref{C:mindvH}}$, may depend on the Brownian structure constants
of $G$. Actually, the proof below shows that the lower bound does
not depend on $G$ at all, and the upper bound can be shown to do
the same easily. However, we will have no use for these facts. 

\begin{proof}
Denote $\mu=d(v,\partial_{\textrm{cont}}H)/r$. We may assume w.l.o.g.~$\mu=2^{-M}$
for some integer $M$. Before starting with estimates for $R$ we
need to know a fact about Brownian motion, roughly speaking that Brownian
motion conditioned to have $\{ S_{v,4r}<S(\partial_{\textrm{cont}}H)\}$
also avoids $\partial_{\textrm{cont}}H$ along (most of) its path.
To formulate precisely, let $\beta<1$ be some parameter which will
be fixed later, and denote \[
a_{k}=d(v,\partial_{\textrm{cont}}H)2^{k}\quad b_{k}=a_{k}^{\beta}.\]
So that $r=a_{M}$. Define stopping times $S_{k}:=S_{v,a_{k}}$ for
$k\in\left]-\infty,M+2\right]$ and let $\mathcal{E}:=\bigcup_{k=-\infty}^{M+1}\mathcal{E}_{k}$
be defined by \[
\mathcal{E}_{k}:=\{\exists S_{k-1}<t\leq S_{k}:d(W(t),\partial_{\textrm{cont}}H)\leq b_{k}\}\cap\{ S_{M+2}<S(\partial_{\textrm{cont}}H))\}.\]
Note that because we end at $\mathcal{E}_{M+1}$ we actually ignore
the event of getting close to $\partial_{\textrm{cont}}H$ on the
last stretch of the Brownian motion, namely $\left]S_{M+1},S_{M+2}\right]$
(this makes the proof a little simpler). We assume $C_{\ref{C:mindvH}}\geq1$
and then $\mathcal{E}_{k}$ is empty for $k\leq-1$. For $0\leq k\leq M+1$
we shall use lemma \ref{lem:Brownescape} in the form $\mathbb{P}(S_{k-1}<S(\partial_{\textrm{cont}}H))\leq C2^{-k}$
and again (together with the strong Markov property) to get\[
\mathbb{P}(S_{M+2}<S(\partial_{\textrm{cont}}H)\,|\, W[0,t^{*}])\le C\frac{b_{k}}{r}\quad t^{*}:=\inf_{t>S_{k-1}}d(W(t),\partial_{\textrm{cont}}H)\leq b_{k}.\]
(clearly $t^{*}$ is a stopping time hence we may use the strong Markov
property). Hence we get\[
\mathbb{P}(\mathcal{E}_{k})\leq C2^{-k}\frac{b_{k}}{r}=C\mu a_{k}^{\beta-1}\]
and summing (remember that $\beta<1$)\begin{equation}
\mathbb{P}(\mathcal{E})\leq C\mu\frac{d(v,\partial_{\textrm{cont}}H)^{\beta-1}}{1-2^{\beta-1}}.\label{eq:Fsmall}\end{equation}

We now move to examine the random walk $R$. Couple $W$ with $R$
as above. Let $T_{k}=T_{v,a_{k}}$ that is the $R$-equivalents of
the $S_{k}$. Assume $C_{\ref{C:mindvH}}$ is sufficiently large so
that $R[0,T_{-1}]\cap\partial H=\emptyset$ always. We start with
a lower bound for $p$. A little set calculus gives \begin{align}
\mathbb{P}(T_{M}<T(\partial H)) & \geq\mathbb{P}(S_{M+2}<S(\partial_{\textrm{cont}}H))-\mathbb{P}(\mathcal{E})-\nonumber \\
 & \quad-\mathbb{P}(\{ S_{M+2}<S(\partial_{\textrm{cont}}H)\}\setminus(\{ T_{M}<T(\partial H)\}\cup\mathcal{E})).\label{eq:setcalc}\end{align}
Now, if $T_{M}\geq T(\partial H)$ then we have that $R(n)\in\partial H$
for some $n\in[T_{k-1},T_{k}]$. If $S_{M+2}<S(\partial_{\textrm{cont}}H)$
but $\mathcal{E}$ did not happen then we must have that $d(W[S_{k-2},S_{k+1}],\linebreak[4]\partial_{\textrm{cont}}H)>b_{k-1}$
and hence $\sh(R[T_{k-1},T_{k}],W[S_{k-2},S_{k+1}])>b_{k-1}-C$. Thus
we arrive at\begin{eqnarray}
\lefteqn{\mathbb{P}(\{ S_{M+2}<S(\partial_{\textrm{cont}}H)\}\setminus(\{ T_{M}<T(\partial H)\}\cup\mathcal{E}))\le}\label{eq:RScapT}\\
 &  & \le\sum_{k=0}^{M}\mathbb{P}(\{\sh(R[T_{k-1},T_{k}],W[S_{k-2},S_{k+1}])>b_{k-1}-C\}\cap\nonumber \\
 &  & \qquad\qquad\cap\;\{ S_{M+2}<S(\partial_{\textrm{cont}}H)\}).\nonumber \end{eqnarray}
The estimate of (\ref{eq:RScapT}) follows from lemma \ref{lem:subhaus}
but first we need to chose $\beta$ and we choose $\beta=1-\frac{1}{2}c_{\ref{c:RWclose}}(G)$
where $c_{\ref{c:RWclose}}(G)$ comes from lemma \ref{lem:subhaus}.
The lemma then claims \begin{equation}
\mathbb{P}(\sh(R[T_{k-1},T_{k}],W[S_{k-2},S_{k+1}])>b_{k-1}-C)\leq C(G)\exp(-a_{k}^{c(G)}).\label{eq:noSMSH}\end{equation}
The condition $S_{M+2}<S(\partial_{\textrm{cont}}H)$ can be added
via lemma \ref{lem:Ers}, and we get \begin{multline*}
\mathbb{P}\Big(\{\sh(R[T_{k-1},T_{k}],W[S_{k-2},S_{k+1}])>b_{k-1}-C\}\:\cap\\
\cap\{ S_{M+2}<S(\partial_{\textrm{cont}}H)\}\Big)\leq C(G)\exp(-a_{k}^{c(G)})\mu2^{k}.\end{multline*}
 Plugging this into (\ref{eq:RScapT}) and summing we get \[
\mathbb{P}(\{ S_{M+2}<S(\partial_{\textrm{cont}}H)\}\setminus(\{ T_{M}<T(\partial H)\}\cup\mathcal{E}))\le C(G)\exp(-d(v,\partial_{\textrm{cont}}H)^{c(G)})\mu.\]
This we may plug into (\ref{eq:setcalc}) together with (\ref{eq:Fsmall})
and lemma \ref{lem:Brownescape} and get\[
\mathbb{P}(T_{M}<T(H))\geq\mu\left(c-C\frac{d(v,\partial_{\textrm{cont}}H)^{-c_{\ref{c:RWclose}}(G)/2}}{1-2^{-c_{\ref{c:RWclose}}(G)/2}}-C(G)\exp(-d(v,\partial_{\textrm{cont}}H)^{c(G)})\right)\]
and it is now clear that if $C_{\ref{C:mindvH}}$ is chosen sufficiently
large, then $d(v,\partial_{\textrm{cont}}H)>C_{\ref{C:mindvH}}$ would
give that everything inside the parenthesis is $>c$ and the direction
$p\geq c\mu$ is proved.

The proof that $p\leq C\mu$ is, generally speaking, a mirror image
exchanging the roles of $R$ and $W$ in the proof of $p\geq c\mu$.
Since our a-priori knowledge about the random walk is smaller (it
is, essentially, lemma \ref{lem:plane}), the proof is somewhat rearranged.
Here are the details: Define \[
g_{i}:=2^{i}\mathbb{P}(T_{i}<T(\partial H)).\]
Fix one $i\geq4$. For every $j\in[2,i-2]$ examine the event\[
\mathcal{F}_{j}:=\{ W[S_{j-1},S_{j}]\cap\partial_{\textrm{cont}}H\neq\emptyset\}\cap\{ W\left]S_{j},S_{i-2}\right]\cap\partial_{\textrm{cont}}H=\emptyset\}.\]
The event $\{ T_{i}<T(\partial H)\}\cap\mathcal{F}_{j}$ has a number
of consequences:

\smallskip{}
\noindent \textbf{(i)} $R[0,T_{j-2}]\cap\partial H=\emptyset$. By
definition this event has probability $2^{2-j}g_{j-2}$.

\smallskip{}
\noindent \textbf{(ii)} $R[T_{j-2},T_{j+1}]\cap\partial H=\emptyset$.
Lemma \ref{lem:subhaus} shows that \[
\mathbb{P}(\sh(W[S_{j-1},S_{j}],R[T_{j-2},T_{j+1}])>a_{j}^{1-c_{\ref{c:RWclose}}/2})\leq C(G)\exp(-a_{j}^{c(G)})\]
and since $W[S_{j-1},S_{j}]\cap\partial_{\textrm{cont}}H\neq\emptyset$
we get that with probability $1-C\exp(-a_{j}^{c})$, $d(R[T_{j-2},T_{j+1}],\partial H)\leq a_{j}^{1-c_{\ref{c:RWclose}}/2}+C(G)$.
Denote this event by $\mathcal{G}$. Lemma \ref{lem:plane} and the
strong Markov property now show that\[
\mathbb{P}(\mathcal{G}\cap\{ T_{j+2}>T(\partial H)\}\,|\, R[0,T_{j-2}])\leq C(G)\left(\frac{a_{j}^{1-c_{\ref{c:RWclose}}/2}}{a_{j}}\right)^{c(G)}.\]
Together with clause (i) we get\begin{eqnarray*}
\lefteqn{\mathbb{P}(\{ R[T_{0},T_{j+2}]\cap\partial H=\emptyset\}\cap\{ W[S_{j-1},S_{j}]\cap\partial_{\textrm{cont}}H\neq\emptyset\})\leq}\\
 & \qquad & \leq C(G)2^{2-j}g_{j-2}a_{j}^{-c(G)}+C(G)\exp(-a_{j}^{c(G)})\stackrel{(*)}{\leq}C(G)2^{2-j}g_{j-2}a_{j}^{-c(G)}\end{eqnarray*}
where in $(*)$ we used the lower bound $g_{j}\geq c$ already established.

\smallskip{}
\noindent \textbf{(iii)} $W[S_{j+2},S_{i}]\cap\partial_{\textrm{cont}}H=\emptyset$.
Here we employ lemma \ref{lem:Ers} and get \begin{align*}
\mathbb{P}(\{ T_{i}<T(\partial H)\}\cap\mathcal{F}_{j}) & \leq C2^{j+2-i}\big(C(G)2^{2-j}g_{j-2}a_{j}^{-c(G)}+C(G)\exp(-a_{j}^{c(G)})\big)\\
 & \leq C(G)2^{-i}g_{j-2}a_{j}^{-c(G)}.\end{align*}

With the estimate of $\mathbb{P}(\{ T_{i}<T(\partial H)\}\cap\mathcal{F}_{j})$
complete we need only sum on $j$ and use lemma \ref{lem:Brownescape}
to get\[
\mathbb{P}(\{ T_{i}<T(\partial H)\}\setminus\bigcup_{j=2}^{i-2}\mathcal{F}_{j})\leq\mathbb{P}(W[S_{1},S_{i-2}]\cap\partial_{\textrm{cont}}H=\emptyset)\leq C2^{-i}\leq C2^{-i}g_{0}\]
so we get\[
g_{i}\leq C(G)\sum_{j=0}^{i-4}g_{j}\left(2^{-c(G)}\right)^{j}.\]
By lemma 4.5 of \cite{L96b}, $g_{i}$ are bounded and the bound depends
only on the isotropic structure constants of $G$ (we use here that
$g_{0}\leq1$) and the theorem is proved. \newC{C:arc}
\end{proof}
\begin{lem}
\label{lem:arc}With the notations of theorem \ref{lem:escape} (but
$d(v,\partial_{\textrm{cont}}H)\geq C_{\ref{C:arc}}$), let $p$ be
the probability that a random walk $R$ on $G$ starting from $v$
will hit $\partial B(v,r)\cup\partial H$ in the arc\[
\alpha:=\partial B(v,r)\cap\{ x:d(x,\partial_{\textrm{cont}}H)>{\textstyle \frac{1}{2}}r\}.\]
Then $p\approx d(v,\partial_{\textrm{cont}}H)/r$.
\end{lem}
\begin{proof}
$p\leq C(G)d(v,\partial_{\textrm{cont}}H)/r$ is an immediate consequence
of theorem \ref{lem:escape}. For the other direction, first assume
w.l.o.g.~that $r>4d(v,\partial_{\textrm{cont}}H)$, which can be
done by lemma \ref{lem:DD}. Let $\lambda>0$ be some parameter, and
let $\beta$ be the two arcs $\partial B(v,r/2)\cap\{ x:d(x,\partial_{\textrm{cont}}H)<\lambda r\}$.
Examine the event\[
\mathcal{E}:=\{ R(T(\partial B(v,r/2)\cup\partial H))\in\beta\}\cap\{ T_{v,r}<T(\partial H)\}.\]
Theorem \ref{lem:escape} shows that \[
\mathbb{P}(R(T(\partial B(v,r/2)\cup\partial H))\in\beta)\leq\mathbb{P}^{v}\left(T_{v,r/2}<T(\partial H)\right)\leq C(G)d(v,\partial_{\textrm{cont}}H)/r.\]
For any $x\in\beta$ we have, again from theorem \ref{lem:escape}
\[
\mathbb{P}^{x}(T_{v,r}<T(\partial H))\leq\mathbb{P}^{x}(T_{x,r/4}<T(\partial H))\leq C(G)\lambda\]
if only $\lambda<\frac{1}{8}$ and $\lambda r>C(G)$. Hence we get\[
\mathbb{P}(\mathcal{E})\leq C(G)\lambda d(v,\partial_{\textrm{cont}}H)/r.\]
Combining this with the lower bound of theorem \ref{lem:escape} we
get \begin{eqnarray*}
\lefteqn{c(G)\frac{d(v,\partial_{\textrm{cont}}H)}{r}\leq\mathbb{P}\left(T_{v,r/2}<T(\partial H)\right)\leq}\\
 & \qquad & \leq\mathbb{P}(R(T(B(v,r/2)\cup\partial H))\in\partial B(v,r/2)\setminus\beta)+\lambda C(G)\frac{d(v,\partial_{\textrm{cont}}H)}{r}.\end{eqnarray*}
Choose $\lambda=\lambda(G)$ some constant sufficiently small and
get that \[
\mathbb{P}(R(T(B(v,r/2)\cup\partial H))\in\partial B(v,r/2)\setminus\beta)>c(G)\frac{d(v,\partial_{\textrm{cont}}H)}{r}.\]
Lemma \ref{lem:DD} now shows that there is a probability $>c(G)$
to hit $\partial B(v,r)\cup\partial H$ at $\alpha$ if you start
from any point of $\partial B(v,r/2)\setminus\beta$, and we are done.
\end{proof}

\subsection{Lower bound for the non-intersection probability}

The proof of theorem \ref{lem:escape} in the previous section was
modeled roughly on Lawler \cite{L96b}. In contrast, theorem \ref{thm:kof},
which will be proved in this section and the next, is a completely
straightforward generalization of \cite{L96b}. 

\begin{thm}
\label{thm:kof}Let $G$ be an isotropic graph of dimension $2$
or $3$. Then for any $v^{1},v^{2}\in G$ with $|v^{1}-v^{2}|>C(G)$,
If $R^{1}$ and $R^{2}$ are two walks with $R^{i}$ starting from
$v^{i}$ and stopped on $\partial B(v^{1},r)$, $r>2|v^{1}-v^{2}|$,
then \begin{equation}
c(G)\left(\frac{|v_{1}-v_{2}|}{r}\right)^{-\xi}\leq\mathbb{P}(R^{1}\cap R^{2}=\emptyset)\leq C(G)\left(\frac{|v_{1}-v_{2}|}{r}\right)^{-\xi}\label{eq:thmkof}\end{equation}
where $\xi=\xi_{d}(1,1)$ is the intersection exponent from (\ref{eq:defxi}).
\end{thm}
We shall not repeat the argumentation of \cite{L96b}, we shall only
note the pieces that require changes. Hence the rest of the chapter
should be read side by side with \cite{L96b}. Chapter 2 of \cite{L96b}
requires almost no changes: the following lemma, which is a replacement
for (7) in lemma 2.5 is perhaps worth proving here.\newC{C:min7}

\begin{lem}
\label{lem:lwlr25}Let $\epsilon>0$ and let $G$ be an isotropic
graph. Then there exists a $\delta=\delta(\epsilon,G)$ such that
for any $r>C_{\ref{C:min7}}(\epsilon,G)$; any $v\in G$ and any $w\in G$
with $|v-w|\leq r$ one has\[
\mathbb{P}^{1,w}\left(\inf_{|z-v|<r}\mathbb{P}^{2,z}\big(R^{2}[0,T_{v,2r}^{2}]\cap R^{1}[T_{v,r}^{1},T_{v,2r}^{1}]\neq\emptyset\,|\, R^{1}[T_{v,r}^{1},T_{v,2r}^{1}]\big)<\delta\right)<\epsilon\]
where $R^{i}$ are two independent random walks.
\end{lem}
In words, if we consider a path to be {}``$\delta$-hittable from
$z$'' if the probability of a random walk ($R^{2}$) starting from
$z$ to hit it is $\geq\delta$, then what we prove here is that random
walk ($R^{1}$) is, with probability $1-\epsilon$, $\delta$-hittable
from any $z\in B(v,r)$. (to understand the formula formally, remember
that the conditional probability $\mathbb{P}(\cdot|*)$ is a function
of $*$ and note that the $\inf$ relates to a pointwise infimum of
these functions). 

\begin{proof}
Let $x\in G$ and let $s>\lambda$ where $\lambda=\lambda(G)$ will
be some constant sufficiently large that will be fixed later. Assume
for now that $\lambda>C_{\ref{C:minr}}(\frac{1}{8},G)$ where $C_{\ref{C:minr}}$
comes from lemma \ref{lem:Rintrsct}. Hence we use lemma \ref{lem:Rintrsct}
and get for any $y\in B(x,\frac{7}{4}s)$, \begin{equation}
\mathbb{P}^{1,x,2,y}(R^{1}[0,T_{x,s/4}^{1}]\cap R^{2}[0,T_{x,2s}^{2}]\neq\emptyset)>c(G).\label{eq:P12}\end{equation}
For any path $\gamma$ from $x$ to $\partial B(x,s/4)$ define $Y(y,\gamma):=\mathbb{P}^{2,y}(\gamma\cap R^{2}[0,T_{x,2s}^{2}]\neq\emptyset)$.
Then (\ref{eq:P12}) implies that\[
\mathbb{P}^{1,x}(Y(y,R^{1}[0,T_{x,s/4}^{1}])>c(G))>c(G)\quad\forall y\in B(x,{\textstyle \frac{7}{4}}s).\]
If $\lambda$ is sufficiently large then we have $\partial B(x,\frac{1}{4}s)\subset B(x,\frac{1}{2}s)$
and then $Y(\cdot,\gamma)$ is harmonic on $\{ y:\frac{1}{2}s\leq|x-y|\leq\frac{3}{2}s\}$,
and we may use Harnack's inequality (lemma \ref{lem:Harnack-general})
to show that for some constant $\mu=\mu(G)$, \begin{equation}
\mathbb{P}^{1,x}\Big(\inf_{\frac{3}{4}s\leq|x-y|\leq\frac{5}{4}s}Y(y,R^{1}[0,T_{x,s/2}^{1}])>\mu\Big)>c(G).\label{eq:defZ}\end{equation}
Denote by $Z(\gamma)$ the event $\{\inf_{3s/4\leq|x-y|\leq5s/4}Y(y,\gamma)>\mu\}$
where $x$ is the beginning of the path $\gamma$.

Next, let $N=N(\epsilon,G)$ be an integer parameter which will be
fixed later. For $i=1,\dotsc,N-1$ define $T_{i}:=T_{v,(1+i/N)r}^{1}$.
Let $x_{i}=R^{1}(T_{i})$. Let $s=r/4N$. Define $U_{i}$ to be the
stopping times\[
U_{i}:=\min\{ t>T_{i}:R^{1}(t)\in\partial B(x_{i},{\textstyle \frac{1}{4}}s)\}.\]
Finally define $Z_{i}=Z(R^{1}[T_{i},U_{i}])$. Then (\ref{eq:defZ})
says that $\mathbb{P}(Z_{i}\,|\, x_{i})>c(G)$. Since the only effect
of $Z_{1},\dotsc,Z_{i-1}$ on $Z_{i}$ is through $x_{i}$ we get
in fact that\[
\mathbb{P}(Z_{i}\,|\, Z_{1},\dotsc,Z_{i-1})>c(G)\]
and hence\[
\mathbb{P}\Big(\bigcap_{i=1}^{N-1}\neg Z_{i}\Big)<(1-c(G))^{N-2}.\]
Denote $\mathcal{Z}:=\cup_{i=1}^{N-1}Z_{i}$ and choose our parameter
$N$ such that $\mathbb{P}(\mathcal{Z})>1-\epsilon$. Lemma \ref{lem:DD}
shows that for $r$ bigger than some constant $\nu(N,G)$ we have
that the probability of $R^{2}$ to hit $\partial B(x_{i},s)$ for
any $i$ and for any starting point $z$ of $R^{2}$ is $\geq c(N,G)$.
If $\lambda$ is sufficiently large then $\partial B(x_{i},s)\subset B(x_{i},\frac{5}{4}s)$.
Hence we get for any $i\in\{1,\dotsc,N\}$,\begin{eqnarray*}
\lefteqn{\mathbb{P}(R^{2}[0,T_{v,2r}^{2}]\cap R^{1}[T_{v,r}^{1},T_{v,2r}^{1}]\neq\emptyset\,|\, Z_{i})\geq}\\
 & \qquad & \geq\mathbb{P}(R^{2}[0,T_{v,2s}^{2}]\cap R^{1}[T_{i},U_{i}]\neq\emptyset\,|\, Z_{i})\geq\\
 &  & \stackrel{(*)}{\geq}\mathbb{P}(T_{x_{i},s}^{2}<T_{v,2r}^{2})\cdot\mathbb{E}(Y(R^{2}(T_{x_{i},s}^{2}),R^{1}(T_{i},U_{i}))\,|\, Z_{i})\geq\\
 &  & \geq\mu\mathbb{P}(T_{x_{i},s}^{2}<T_{v,2s}^{2})\geq\mu c(N,G)\end{eqnarray*}
where $(*)$ comes from the strong Markov property at the stopping
time $T_{x_{i},s}^{2}$. Hence we get \[
\mathbb{P}(R^{2}[0,T_{v,2s}^{2}]\cap R^{1}[T_{v,s}^{1},T_{v,2s}^{1}]\,|\,\mathcal{Z})\geq\mu c(N,G).\]
This finishes the lemma: we fix $\lambda$ and define $\delta(\epsilon,G):=\mu c(N,G)$
and $C_{\ref{C:min7}}(\epsilon,G):=\max(8N\lambda,\nu(N,G))$ and
we are done.
\end{proof}
\begin{lem}
\label{lem:lwlr26}\newC{C:MKepsG}Let $G$ be an isotropic graph
and let $M$, $K$ and $\epsilon$ be some parameters. Then there
exists a $\delta(M,K,\epsilon,G)$ and a $C_{\ref{C:MKepsG}}(M,K,\epsilon,G)$
such that for all $v,w\in G$, and all $r>|v-w|$,\[
\mathbb{P}^{1,w}\left(\inf_{(*)}\mathbb{P}^{2,z}\big(R^{2}[0,T_{v,2r}^{2}]\cap R^{1}[0,T_{v,2r}^{1}]\neq\emptyset\,|\, R^{1}[0,T_{v,2r}^{1}]\big)\geq r^{-\delta}\right)<C_{\ref{C:MKepsG}}r^{-M}\]
where the $(*)$ stands for all the $z\in B(v,r)$ such that $d(z,R^{1}[0,T_{v,2r}^{1}])\leq Kr^{1-\epsilon}$.
\end{lem}
The proof is identical to that of lemma 2.6 from \cite{L96b} and
we shall omit it. Very roughly, it uses the previous lemma and the
Wiener shell test.

Chapter 3 of \cite{L96b} has no real equivalence here. The Skorokhod
embedding used in \cite{L96b} has the convenient property that the
random walk and the Brownian motion have comparable times, that is
$|R(t)-W(t)|\ll t^{1/4+\epsilon}$ (after linear calibration). This
is just not true in our case, or anyway would require non-linear adaptive
calibration which is not worth messing with --- measuring the Hausdorff
distance between $R$ and $W$ is a completely adequate replacement.
Hence we shall make no effort to give analogs of the results of chapter
3 of \cite{L96b} and continue immediately to chapter 4. Lemma \ref{lem:Lwlr41}
is a replacement for Lawler's lemma 4.1, but first an auxiliary result:

\begin{lem}
\label{lem:Erskof}Let $G\subset\mathbb{R}^{d}$ be an isotropic graph.
Let $R^{1},W^{1}$ and $R^{2},W^{2}$ be two independent pairs of
coupled random walk and Brownian motion on $G$ starting from $v^{1}$
and $v^{2}$ respectively. Let $s>r>2|v^{1}-v^{2}|$. Let $\mathcal{E}$
be an event depending on $R^{i}[0,T_{v,r}^{i}]$ and $W^{i}[0,S_{v,r}^{i}]$
only. Then\[
\mathbb{P}(\mathcal{E}\cap\{ W^{1}[S_{v,r}^{1},S_{v,s}^{1}]\cap W^{2}[S_{v,r}^{2},S_{v,s}^{2}]=\emptyset\})\leq C\left(\frac{r}{s}\right)^{\xi}(\mathbb{P}(\mathcal{E})+C(G)\exp(-r^{c(G)})).\]
where $\xi$ is from (\ref{eq:defxi}).
\end{lem}
The proof is identical to that of lemma \ref{lem:Ers}, with the use
of lemma \ref{lem:Brownescape} replaced by estimates for the non-intersection
probability of two Brownian motion, see \cite[(2)]{L96a}. We omit
the details.

\begin{lem}
\label{lem:Lwlr41}Let $G$ be an isotropic graph, let $v\in G$ and
let $R^{1},W^{1}$ and $R^{2},W^{2}$ be two independent pairs of
coupled random walk and Brownian motion on $G$ starting from $v^{1}$
and $v^{2}$ respectively, $|v-v^{i}|\leq2^{m}$ and stopped on $\partial B(v,2^{n})$.
Define $T_{j}^{i}:=T_{v,2^{j}}^{i}$ and $S_{j}^{i}:=S_{v,2^{j}}^{i}$
and\begin{alignat*}{2}
Q_{j}^{i} & := & \; & \{\sh(R^{i}[T_{j-1}^{i},T_{j}^{i}],W^{i}[S_{j-2}^{i},S_{j+1}^{i}])\geq2^{j(1-c_{\ref{c:RWclose}}/2)}\}\;\cup\\
 &  &  & \{\sh(W^{i}[S_{j-1}^{i},S_{j}^{i}],R^{i}[T_{j-2}^{i},T_{j+1}^{i}])\geq2^{j(1-c_{\ref{c:RWclose}}/2)}\},\\
Q_{*}^{i} & := &  & \{\sh(R^{i}[0,T_{m+2}^{i}],W^{i}[0,S_{m+3}^{i}])\geq2^{(m+2)(1-c_{\ref{c:RWclose}}/2)}\}\;\cup\\
 &  &  & \{\sh(W^{i}[0,S_{m+2}^{i}],R^{i}[0,T_{m+3}^{i}])\geq2^{(m+2)(1-c_{\ref{c:RWclose}}/2)}\}.\\
\mathcal{Q} & := &  & Q_{*}^{1}\cup Q_{*}^{2}\cup\bigcup_{j=m+3}^{n-1}Q_{j}^{1}\cup Q_{j}^{2}.\end{alignat*}
Then\[
\mathbb{P}(\mathcal{Q}\cap\{ W^{1}[0,S_{n+1}^{1}]\cap W^{2}[0,S_{n+1}^{2}]=\emptyset\})\leq C(G)\exp(-2^{mc(G)})2^{-(n-m)\xi}.\]

\end{lem}
\begin{proof}
The corollary to lemma \ref{lem:subhaus} with $\nu=2$ shows that\[
\mathbb{P}(Q_{j}^{i})\leq C(G)\exp(-2^{jc(G)}).\]
Next, lemma \ref{lem:Erskof} shows that\[
\mathbb{P}(Q_{j}^{i}\cap\{ W^{1}[0,S_{n+1}^{1}]\cap W^{2}[0,S_{n+1}^{2}]=\emptyset\})\leq C(G)\exp(-2^{jc(G)})2^{-(n-j)\xi}.\]
 $Q_{*}^{i}$ have a similar estimate. Summing on $i$ and $j$ we
get the lemma.
\end{proof}
We now prove a lemma, the equivalent of corollary 4.2 of \cite{L96b},
somewhat stronger than the direction $\mathbb{P}(R^{1}\cap R^{2}=\emptyset)\geq c\left(|v_{1}-v_{2}|/r\right)^{-\xi}$
of (\ref{eq:thmkof}). We will need the strengthening in the next
chapter.\newC{C:cor42}

\begin{lem}
\label{lem:lwlr42}Let $G$ be an isotropic graph of dimension $2$
or $3$ and let $v\in G$ and $s\geq C_{\ref{C:cor42}}(G)$. Let $v^{1},v^{2}\in G\cap\big(B(v,s)\setminus B(v,\frac{7}{8}s)\big)$,
$|v^{1}-v^{2}|>\frac{1}{4}s$, let $r>4s$ and let $\eta$ be a unit
vector in $\mathbb{R}^{d}$ and define two subsets of $G$,\begin{align}
 & U^{1}:=\big(B(v,r/2)\setminus B(v,s)\big)\cup B(v^{1},{\textstyle \frac{1}{4}}s)\cup\left(\overline{B(v,r)}\cap\left\{ w:\left\langle w,\eta\right\rangle \geq{\textstyle \frac{1}{4}}r\right\} \right),\label{eq:defMi}\\
 & U^{2}:=\big(B(v,r/2)\setminus B(v,s)\big)\cup B(v^{2},{\textstyle \frac{1}{4}}s)\cup\left(\overline{B(v,r)}\cap\left\{ w:\left\langle w,\eta\right\rangle \leq-{\textstyle \frac{1}{4}}r\right\} \right).\nonumber \end{align}
Let $R^{1}$ and $R^{2}$ are two walks with $R^{i}$ starting from
$v^{i}$ and stopped on $\partial B(v,r)$. Then\begin{equation}
\mathbb{P}(\{ R^{1}\cap R^{2}=\emptyset\}\cap\{ R^{i}\subset U^{i},i=1,2\})\geq c(G)\left(\frac{s}{r}\right)^{\xi}.\label{eq:mshrm}\end{equation}

\end{lem}
\begin{proof}
This is now immediate. Indeed, consider slightly smaller domains (but
extended outward), \begin{align*}
 & V^{1}:=\big(B(v,{\textstyle \frac{5}{12}}r)\setminus B(v,{\textstyle \frac{25}{24}}s)\big)\cup B(v^{1},{\textstyle \frac{5}{24}}s)\cup\left\{ w\in B(v,2r):\left\langle w,\eta\right\rangle \geq{\textstyle \frac{1}{3}}r\right\} ,\\
 & V^{2}:=\big(B(v,{\textstyle \frac{5}{12}}r)\setminus B(v,{\textstyle \frac{25}{24}}s)\big)\cup B(v^{2},{\textstyle \frac{5}{24}}s)\cup\left\{ w\in B(v,2r):\left\langle w,\eta\right\rangle \leq-{\textstyle \frac{1}{3}}r\right\} .\end{align*}
Consider also the event $\mathcal{F}$ that $W^{1}$ and $W^{2}$
are reasonably far apart along their paths, namely\begin{align*}
\mathcal{F}_{j}^{i} & :=\{ d(W^{i}[S_{v,2^{j-3}}^{i},S_{v,2^{j}}^{i}],W^{3-i}[0,S_{v,2^{j}}^{3-i}])>2^{j(1-c_{\ref{c:RWclose}}/4)}\}\\
\mathcal{F} & :=\bigcap_{i=1,2}\bigcap_{j=\left\lfloor \log_{2}s\right\rfloor }^{\left\lceil \log_{2}2r\right\rceil }\mathcal{F}_{j}^{i}\end{align*}
 Then it follows using techniques similar to \cite{L96a}, see corollaries
3.9, 3.11 and 3.12 ibid. and lemma 2.8 of \cite{L96b} that for $s>C(G)$,
\[
\mathbb{P}(\mathcal{F}\cap\{ W^{i}[0,2r]\subset V^{i},\, i=1,2\})\geq c\left(\frac{s}{r}\right)^{\xi}.\]
We couple $W^{i}$ to $R^{i}$ such that $(R^{1},W^{1})$ is independent
from $(R^{2},W^{2})$, and consider the event $\mathcal{Q}$ from
lemma \ref{lem:Lwlr41}. If $\mathcal{Q}$ did not occur, then $R^{i}$
is sufficiently close to $W^{i}$ such that $W^{i}[0,S_{v,2r}^{i}]\subset V^{i}$
implies $R^{i}[0,T_{v,r}^{i}]\subset U^{i}$, if $s>C(G)$. Further,
$\mathcal{F}\setminus\mathcal{Q}$ also implies that $R^{1}[0,T_{v,r}^{1}]\cap R^{2}[0,T_{v,r}^{2}]=\emptyset$.
Finally, lemma \ref{lem:Lwlr41} shows that if $s>C(G)$ then $\mathbb{P}(\mathcal{Q})\leq\frac{1}{2}c(s/r)^{\xi}$
which finishes the lemma.
\end{proof}
\begin{cor*}
Let $G$ be an isotropic graph, let $v^{1},v^{2}\in G$ and let $r>4|v^{1}-v^{2}|$.
Let $R^{1}$ and $R^{2}$ be two walks starting from $v^{i}$ and
stopped on $\partial B(v^{1},r)$. Then \[
\mathbb{P}(R^{1}\cap R^{2}=\emptyset)\begin{cases}
>c(G)(|v^{1}-v^{2}|/r)^{\xi}\\
=0\end{cases}.\]

\end{cor*}
\begin{proof}
Using lemma \ref{lem:lwlr42} we can fix a constant $\lambda=\lambda(G)$
such that for all $|v^{1}-v^{2}|>\lambda$ the first choice happens.
Hence assume $|v^{1}-v^{2}|\leq\lambda$ and use lemma \ref{lem:0n0}
with $\epsilon=\frac{1}{8}$ and $s=\lambda$ and get that either
\begin{enumerate}
\item There are no two disjoint paths going from $v^{1}$ and $v^{2}$ to
the exterior of $B(v^{1},\kappa(\frac{1}{8},\lambda,G))$, $\kappa$
from lemma \ref{lem:0n0}. In this case the probability is $0$ for
every $r>\kappa$.
\item For $\mu=\max\kappa(\frac{1}{8},\lambda,G),C_{\ref{C:cor42}}$ there
are two disjoint simple paths $\gamma^{i}$ starting from $v^{i}$
and ending at $w^{i}$ satisfying $w^{i}\in B(v^{1},\mu)\setminus\overline{B(v^{1},\frac{7}{8}\mu)}$
and $B(w^{i},\frac{7}{8}\mu)\cap\gamma^{3-i}=\emptyset$.
\end{enumerate}
In the second case we use lemma \ref{lem:lwlr42} and get\[
\mathbb{P}^{1,w^{1},2,w^{2}}((\gamma^{1}\cup R^{1}[0,T_{v,r}^{1}])\cap(\gamma^{2}\cup R^{2}[0,T_{v,r}^{2}])=\emptyset)>c(G)(\mu/r)^{\xi}\]
where $\gamma^{i}\cap R^{3-i}=\emptyset$ is satisfied because $B(w^{i},\frac{1}{4}\mu)\cap\gamma^{3-i}=\emptyset$,
$\gamma^{i}\subset B(v^{1},\mu)$ and the event of lemma \ref{lem:lwlr42}
includes that $R^{i}\cap B(v^{1},\mu)\subset B(w^{i},\frac{1}{4}\mu)$.
In the case that the $R^{i}$ start from $v^{i}$, the probability
that both follow $\gamma^{i}$ until its end is $>c(G)$, which proves
the corollary for $r>4\mu$. For $r<4\mu$ the lemma will hold automatically
for a sufficiently small constant in its definition.
\end{proof}

\subsection{The upper bound}

Having settled the lower bound in theorem \ref{thm:kof}, we need
only the following lemma, which is slightly stronger than the upper
bound (again, we will need the stronger version in the next chapter).

\begin{lem}
\label{lem:kofup}Let $G$ be an isotropic graph of dimension $2$
or $3$. Then for any $v^{1},v^{2}\in G$, if $R^{1}$ and $R^{2}$
are two walks with $R^{i}$ starting from $v^{i}$ and stopped on
$\partial B(v^{1},r)$ then \[
\mathbb{P}(R^{1}\cap R^{2}=\emptyset)\leq C(G)\left(\frac{|v_{1}-v_{2}|}{r}\right)^{-\xi}\]

\end{lem}
\begin{proof}
Assume that $\mathbb{P}(R^{1}\cap R^{2}=\emptyset)>0$ (in particular
that $v^{1}\neq v^{2}$). Also assume w.l.o.g.~that $r>4|v^{1}-v^{2}|$.
Let $a_{j}=2^{j}|v^{1}-v^{2}|$, $b_{j}=a_{j}^{1-c_{\ref{c:RWclose}}/2}$
and $T_{j}^{i}:=T_{v^{1},a_{j}}^{i}$. Define\[
g_{j}:=2^{j\xi}\mathbb{P}(R^{1}[0,T_{j}^{1}]\cap R^{2}[0,T_{j}^{2}]=\emptyset).\]
The corollary to lemma \ref{lem:lwlr42} shows that $g_{j}>c(G)$.
We need to show that $g_{j}\leq C(G)$. Let $W^{1}$ and $W^{2}$
be Brownian motions coupled to $R^{1}$ and $R^{2}$ respectively,
i.e.~the couples $(R^{1},W^{1})$ and $(R^{2},W^{2})$ are independent.
Let $S_{j}^{i}=S_{v^{1},a_{j}}^{i}$. Examine the event $\mathcal{F}_{j}$
that $j$ is the last step where the $W^{i}$ intersect, namely,\[
\mathcal{F}_{j}^{1}=\{ W^{1}\left]S_{j}^{1},S_{n}^{1}\right]\cap W^{2}\left]S_{j}^{2},S_{n}^{2}\right]=\emptyset\}\cap\{ W^{1}[S_{j-1}^{1},S_{j}^{1}]\cap W^{2}[0,S_{j}^{2}]\neq\emptyset\}\]
Define $\mathcal{F}_{j}^{2}$ replacing the roles of $W^{1}$ and
$W^{2}$. $\{ R^{1}[0,T_{n}^{1}]\cap R^{2}[0,T_{n}^{2}]=\emptyset\}\cap\mathcal{F}_{j}^{1}$
has a number of consequences:

\smallskip{}
\noindent \textbf{(i)} $R^{1}[0,T_{j-2}^{1}]\cap R^{2}[0,T_{j-2}^{2}]=\emptyset$.
By definition this event has probability $2^{(2-j)\xi}g_{j-2}$.

\smallskip{}
\noindent \textbf{(ii)} Next we use the fact that the $W^{i}$ intersect
while the $R^{i}$ don't. The corollary to lemma \ref{lem:subhaus}
shows that \begin{align*}
\mathbb{P}(\sh(W^{1}[S_{j-1}^{1},S_{j}^{1}],R^{1}[T_{j-2}^{1},T_{j+1}^{1})>b_{j}) & \leq C(G)\exp(-a_{j}^{c(G)})\\
\mathbb{P}(\sh(W^{2}[0,S_{j}^{2}],R^{2}[0,T_{j+1}^{2}])>b_{j}) & \leq C(G)\exp(-a_{j}^{c(G)})\end{align*}
and hence if we define $\mathcal{A}:=\{ d(R^{1}[T_{j-2}^{1},T_{j+1}^{1}],R^{2}[0,T_{j+1}^{2}])\leq2b_{j}\}$
($\mathcal{A}$ standing for {}``almost intersecting'') we get\begin{equation}
\mathbb{P}\left(\left\{ W^{1}[S_{j-1}^{1},S_{j}^{1}]\cap W^{2}[0,S_{j}^{2}]\neq\emptyset\right\} \setminus\mathcal{A}\right)\leq C(G)\exp(-a_{j}^{c(G)}).\label{eq:WIRA}\end{equation}
Next define the event $\mathcal{N}:=\left\{ R^{1}[T_{j-2}^{1},T_{j+2}^{1}]\cap R^{2}[0,T_{j+2}^{2}]=\emptyset\right\} $
($\mathcal{N}$ standing for {}``not intersecting''). Lemma \ref{lem:lwlr26}
allows as to estimate $\mathbb{P}(\mathcal{A}\cap\mathcal{N})$: we
use it with the parameters $v=v^{1}$, $w=R^{1}(T_{j-2}^{1})$, $r=a_{j+1}$,
$M=2\xi$, $K=2$ and $\epsilon=c_{\ref{c:RWclose}}/2$. We get that
with probability $\geq1-C_{\ref{C:MKepsG}}(2\xi,2,c_{\ref{c:RWclose}}/2,G)a_{j+1}^{-2\xi}$
in $R^{2}$, \begin{equation}
\mathbb{P}^{1}\left(\mathcal{A}\cap\mathcal{N}\Bigm|R^{2}\left[0,T_{j+2}^{1}\right],R^{1}(T_{j-2}^{1})\right)\leq a_{j}^{-\delta(2\xi,2,c_{\ref{c:RWclose}}/2,G)}.\label{eq:ANsmall}\end{equation}
Notice that we used the strong Markov property from the stopping time
$\min\big\{ t>T_{j-2}^{1}:R^{1}(t)\in B(v^{1},a_{j+1}),\, d(R^{1}(t),R^{2}[0,T_{j+2}^{1}])\leq2b_{j}\big\}$.
Since the events that $R^{i}$ do not intersect up to $T_{j-2}^{i}$
and everything that happens after the $T_{j-2}^{i}$ are dependant
only through $R^{i}(T_{j-2}^{i})$, and since (\ref{eq:ANsmall})
holds for any values of $R^{1}(T_{j-2}^{1})$ we get\begin{align*}
\mathbb{P}(\mathcal{A}\cap\{ R^{1}[0,T_{j+2}^{1}]\cap R^{2}[0,T_{j+2}^{2}]=\emptyset\}) & \leq C(G)a_{j+1}^{-2\xi}+2^{(2-j)\xi}g_{j-2}a_{j}^{-c(G)}\leq\\
 & \stackrel{(*)}{\leq}C(G)2^{-j\xi}g_{j-2}a_{j}^{-c(G)}\end{align*}
where in $(*)$ we used the lower bound. Adding (\ref{eq:WIRA}) we
get\textbf{}\begin{multline}
\mathbb{P}(\left\{ W^{1}[S_{j-1}^{1},S_{j}^{1}]\cap W^{2}[0,S_{j}^{2}]\neq\emptyset\right\} \cap\{ R^{1}[0,T_{j+1}^{1}]\cap R^{2}[0,T_{j+1}^{2}]=\emptyset\})\leq\\
\leq C(G)2^{-j\xi}g_{j-2}a_{j}^{-c(G)}+C(G)\exp(-a_{j}^{c(G)})\leq C(G)2^{-j\xi}g_{j-2}a_{j}^{-c(G)}.\label{eq:kjsmall}\end{multline}

\smallskip{}
\noindent \textbf{(iii)} Finally we use the condition $W^{1}\left]S_{j+2}^{1},S_{n}^{1}\right]\cap W^{2}\left]S_{j+2}^{2},S_{n}^{2}\right]=\emptyset$.
Here we employ lemma \ref{lem:Erskof} and together with (\ref{eq:kjsmall})
we get\begin{equation}
\mathbb{P}\left(\left\{ R^{1}[0,T_{n}^{1}]\cap R^{2}[0,T_{n}^{2}]\neq\emptyset\right\} \cap\mathcal{F}_{j}^{1}\right)\leq C(G)2^{-n\xi}g_{j-2}a_{j}^{-c(G)}.\label{eq:Fnk1}\end{equation}

An identical calculation holds for $\mathcal{F}_{j}^{2}$. We are
almost done! We need only remark that\begin{eqnarray*}
\lefteqn{\mathbb{P}\bigg(\left\{ R^{1}[0,T_{n}^{1}]\cap R^{2}[0,T_{n}^{2}]=\emptyset\right\} \setminus\bigcup_{\substack{2\leq j\leq n-2\\
i=1,2}
}\mathcal{F}_{j}^{i}\bigg)\leq}\\
 & \qquad\qquad & \leq\mathbb{P}\left(W^{1}[S_{1}^{1},S_{n-2}^{1}]\cap W^{2}[S_{1}^{2},S_{n-2}^{2}]=\emptyset\right)\leq C2^{-n\xi}\leq C2^{-n\xi}g_{0}\end{eqnarray*}
and we get\[
g_{n}\leq C(G)\sum_{j=0}^{n-4}g_{j}\left(2^{-c(G)}\right)^{j}.\]
By lemma 4.5 of \cite{L96b}, $g_{i}$ are bounded and the bound depends
only on the isotropic structure constants of $G$ (we use here that
$g_{0}\leq1$) so the lemma and theorem \ref{thm:kof} are proved.
\end{proof}

\section{\label{sec:Quasi-loops}Quasi-loops}

Let $\gamma$ be a path in a $d$-Euclidean net and let $v\in\mathbb{R}^{d}$.
We say that $\gamma$ has an $(s,r)$-quasi-loop near $v$ if there
exists a couple of points $\gamma(i),\gamma(j)\in B(v,s)$ such that
$\diam\gamma[i,j]\geq r$. In this case we write $v\in\mathcal{QL}(s,r,\gamma)$.
We take the $v$-s in a grid such that the balls $B(v,s)$ cover $\mathbb{R}^{d}$
and define\[
\QL(s,r,\gamma):=\#(\mathcal{QL}(s,r,\gamma)\cap{\textstyle \frac{1}{d}}s\mathbb{Z}^{d}).\]

Our purpose in this chapter is to prove that loop-erased random walk
has no quasi loops in the following sense:

\begin{thm}
\label{thm:QL}Let $G$ be an isotropic graph of dimension two or
three, and let $0<\epsilon<1$. Then there exists a $\delta=\delta(\epsilon,G)>0$
such that for all $v\in G$, all $r>C(\epsilon,G)$ and any subset
$v\in\mathcal{D}\subset B(v,r)$, \[
\mathbb{E}^{v}\QL(r^{1-\epsilon},r^{1-\delta},\LE(R[0,T(\partial\mathcal{D})]))\leq C(\epsilon,G)r^{-\delta}.\]

\end{thm}
Dimensions two and three are very different. The proof for dimension
two was done in the case of $\mathbb{Z}^{2}$ by Schramm \cite[lemma 3.4]{S00}
and is practically the same in our more general settings (\cite[lemma 18]{K}
is another variation on Schramm's argument). Therefore we shall only
sketch the required elements in the end of the chapter. We shall concentrate
on dimension three. It turns out that the techniques we use will rely
heavily on the non-intersection exponent and therefore work only for
isotropic graphs. Hence an interesting conjecture appears

\begin{conjecture*}
Theorem \ref{thm:QL} holds for any Euclidean net.
\end{conjecture*}
Again, this is true in dimension two, hence the interesting case is
dimension three.

It will be convenient in many places to consider discontinuous paths.
Therefore, if $\gamma:\{1,\dotsc,n\}\to G$ is some function (without
the restriction that $\gamma(i)$ and $\gamma(i+1)$ are neighbors),
$\LE(\gamma)$ will be defined using the formula (\ref{eq:defLE})
literally, and is a simple discontinuous path. Likewise we will define
$\gamma_{1}\cup\gamma_{2}$ even if $\gamma_{1}(\len\gamma_{1})$
is not a neighbor of $\gamma_{2}(1)$. If $\gamma$ is a (possibly
discontinuous) path and $A$ is some set, then $\gamma\cap A$ would
stand for the discontinuous path created in the natural way from the
parts of $\gamma$ inside $A$, in order.

Here and below when we say {}``$\gamma$ is a discontinuous path'',
we do not exclude the possibility that it is in effect continuous.

\subsection{Cut times}

For any path $\gamma$ we define\[
\cut(\gamma):=\{\gamma(i):\gamma[0,i]\cap\gamma[i+1,\len\gamma]=\emptyset\}.\]
$i$-s satisfying the condition will be called \textbf{cut times}
and the $\gamma(i)$-s will be called \textbf{cut points}. It is clear
that $\cut\gamma\subset\LE(\gamma)$, indeed $\cut\gamma$ is contained
in any connected subset of $\gamma$. It will also be convenient to
define \[
\cut(\gamma;t):=\{\gamma(i):i<t,\,\gamma[0,i]\cap\gamma[i+1,\len\gamma]=\emptyset\}.\]
It has the useful property that $\cut(\gamma;t)$ is increasing in
$t$ and decreasing as $\gamma$ is extended. 

For a random walk $R$, $\cut R$ is intimately related with the non-intersection
exponent $\xi$ via time symmetry. Lemma \ref{lem:CARb0} below has
the details, but first we need some simple preparations.

\begin{lem}
\label{lem:condnomtr}Let $G$ be a three dimensional Euclidean net,
let $v\in\mathbb{R}^{3}$ and let $r>C(G)$. Let $R^{i}$ be random
walks on $G$ starting from points in $B(v,r)$. Let $\mathcal{E}$
be an event with depends only on $R^{i}[0,T_{v,r}^{i}]$ and let $\mathcal{F}$
be an event that depends only on $R^{i}[T_{v,2r}^{i},\infty[$. Then\[
\mathbb{P}(\mathcal{E}\cap\mathcal{F})\approx\mathbb{P}(\mathcal{E})\mathbb{P}(\mathcal{F}).\]
The constant implicit in the $\approx$ notation may depend on the
number of walks, and on the isotropic structure constants.
\end{lem}
\begin{proof}
For every $x\in\partial B(v,r)$ and $y\in\partial B(v,2r)$ let $\pi_{x,y}$
be the probability that a random walk starting from $x$ will hit
$\partial B(v,2r)$ in $y$, and let $\pi_{y}$ be the probability
that $R(T_{v,2r})=y$. By Harnack's inequality (lemma \ref{lem:Harnack})
we have that $\pi_{x,y}\approx\pi_{x',y}$ for any $x,x'\in\partial B(v,r)$.
Hence \[
\pi_{y}\approx\pi_{x,y}.\]
This gives\begin{align*}
\mathbb{P}(\mathcal{E}\cap\mathcal{F}) & =\sum_{x^{i},y^{i}}\mathbb{P}\left(\mathcal{E}\cap\left\{ R^{i}(T_{v,r}^{i})=x^{i}\right\} _{i}\right)\prod_{i}\pi_{x^{i},y^{i}}\mathbb{P}(\mathcal{F}\,|\, R^{i}(T_{v,2r}^{i})=y^{i}\:\forall i)\\
 & \approx\sum_{x^{i},y^{i}}\mathbb{P}\left(\mathcal{E}\cap\left\{ R^{i}(T_{v,r}^{i})=x^{i}\right\} _{i}\right)\prod_{i}\pi_{y^{i}}\mathbb{P}(\mathcal{F}\,|\, R^{i}(T_{v,2r}^{i})=y^{i}\:\forall i)\\
 & =\mathbb{P}(\mathcal{E})\mathbb{P}(\mathcal{F})\qedhere\end{align*}

\end{proof}
\begin{lem}
\label{lem:CARb0}Let $G$ be a three dimensional isotropic graph.
Let $v\in G$ and $r>C(G)$. Define the annulus $A:=B(v,2r)\setminus\overline{B(v,r)}$.
Let $w\in B(v,\frac{1}{2}r)$ and let $R^{1}$ be a random walk starting
from $w$. Let \[
\mathcal{C}:=\cut\left(R^{1}[0,\infty[;T_{v,4r}^{1}\right).\]
Let $z\in B(v,\frac{1}{2}r)$ and let $R^{2}$ be a random walk starting
from $z$ and stopped on $\partial B(v,4r)$. Then\newc{c:CARb0}\[
\mathbb{P}(\mathcal{C}\cap R^{2}[0,T_{v,4r}^{2}]\cap A\neq\emptyset)>c_{\ref{c:CARb0}}(G).\]

\end{lem}
The proof is a relatively straightforward application of second moment
methods, but is quite long. Hence we shall divide it into several
shorter claims.

\begin{sublem}\label{sublem:0n0}\newC{C:sub0n0}There exists a $C_{\ref{C:sub0n0}}(G)$
such that for any $x\in G$ one of the following holds:

\begin{enumerate}
\item There are no two disjoint paths leading from $x$ to $\partial B(x,C_{\ref{C:sub0n0}})$.
\item For any $r\geq C_{\ref{C:sub0n0}}$ there exists two disjoint simple
paths $\gamma^{i}\subset B(x,r)$ that satisfy that if $y^{i}$ is
the end point of $\gamma^{i}$ then \begin{equation}
B(y^{i},{\textstyle \frac{1}{4}}r)\cap\gamma^{3-i}=\emptyset,\quad y^{i}\in B(x,r)\setminus\overline{B(x,{\textstyle \frac{7}{8}}r)}\label{eq:gamiB3i}\end{equation}

\end{enumerate}
\end{sublem}

{}``disjoint paths'' here mean except the point $x$ common to both

\begin{proof}
[Subproof]Let $\lambda=\lambda(G)$ satisfy that any edge in $G$
has length $\leq\lambda$. Then lemma \ref{lem:0n0} for $\epsilon=\frac{1}{8}$,
$s=\lambda$ and all neighbors of $x$ gives the result with $C_{\ref{C:sub0n0}}=\kappa(\frac{1}{8},\lambda,G)$.
\end{proof}
Points $x$ for which there exist two disjoint paths leading outside
$B(x,C_{\ref{C:sub0n0}})$ will be called $\mathcal{C}$\textbf{-capable}.

\begin{sublem}\label{sublem:anyball}\newC{C:hasx}Any ball $B$
of radius $C_{\ref{C:hasx}}(G)$ contains at least one $\mathcal{C}$-capable
point $x$.

\end{sublem}

\begin{proof}
[Subproof]Let $\gamma$ be a path in $B\cap G$ that the distance
between its two ends $y^{1},y^{2}$ is $\geq2C_{\ref{C:hasx}}-C(G)$.
We may assume $\gamma$ is simple (say by taking its loop-erasure).
Let $x$ be the point of $\gamma$ closest to the plane exactly between
$y^{1}$ and $y^{2}$. Then clearly the portions of $\gamma$ up to
$x$ and from $x$ on are disjoint paths that lead to distance at
least $C_{\ref{C:hasx}}-C(G)$, which proves the sublemma, if $C_{\ref{C:hasx}}$
is sufficiently large.
\end{proof}
\begin{sublem}

\label{sublem:NERxi}Let $x\in G$ and $\rho>C(G)$. Let $R^{1}$
and $R^{2}$ be two random walks starting from $x$. Define subsets
similar to (\ref{eq:defMi}) as follows:\begin{align}
 & V^{1}:=B(x,\rho/2)\cup\left(\overline{B(x,\rho)}\cap\left\{ y:\left\langle y,(1,0,0)\right\rangle \geq{\textstyle \frac{1}{4}}\rho\right\} \right),\nonumber \\
 & V^{2}:=B(x,\rho/2)\cup\left(\overline{B(x,\rho)}\cap\left\{ y:\left\langle y,(1,0,0)\right\rangle \leq-{\textstyle \frac{1}{4}}\rho\right\} \right).\label{eq:defMisimp}\end{align}
Notice that $\overline{B(x,\rho)}$ above refers to closure in $G$.
Define further events $\mathcal{V}:=\{ R^{i}[0,T_{x,\rho}^{i}]\subset V^{i}\}_{i=1,2}$
and $\mathcal{N}:=\{ R^{1}[0,T_{x,\rho}^{1}]\cap R^{2}[1,T_{x,\rho}^{2}]=\emptyset\}$.
Then\[
\mathbb{P}^{x}(\mathcal{N})\approx\mathbb{P}^{x}\big(\mathcal{N}\cap\mathcal{V})\begin{cases}
\approx\rho^{-\xi} & x\textrm{ is }\mathcal{C}\textrm{-capable}\\
=0 & \textrm{otherwise}\end{cases}.\]

\end{sublem}

\begin{proof}
[Subproof] The case that $x$ is not $\mathcal{C}$-capable is obvious
if $\rho>C_{\ref{C:sub0n0}}$. In the second case, use sublemma \ref{sublem:0n0}
with its $r$ equal to $\sigma:=\max\{ C_{\ref{C:sub0n0}},C_{\ref{C:cor42}}\}$
($C_{\ref{C:cor42}}$ from lemma \ref{lem:lwlr42}) and get two disjoint
paths $\gamma^{i}$ ending in $y^{i}$ satisfying (\ref{eq:gamiB3i}).
This allows to use lemma \ref{lem:lwlr42} with walks starting from
the $y^{i}$, the $v$, $s$ and $r$ of lemma \ref{lem:lwlr42} equal
to $x$, $\sigma$ and $\frac{1}{2}\rho$ respectively, and with $\eta=(1,0,0)$.
We get\begin{equation}
\mathbb{P}^{1,y_{1},2,y_{2}}\left(\mathcal{N}\cap\left\{ R^{i}[0,T_{x,\rho}^{i}]\subset U^{i}\right\} _{i=1,2}\right)\approx\rho^{-\xi}\label{eq:rhoxiM}\end{equation}
where $U^{i}$ is defined in (\ref{eq:defMi}). In particular, $R^{i}\subset U^{i}$
shows that $R^{i}\cap B(x,\sigma)\subset B(y^{i},\frac{1}{4}\sigma)$
and hence from (\ref{eq:gamiB3i}) $R^{i}\cap\gamma^{3-i}=\emptyset$.
Further, $R^{i}\subset U^{i}$ implies $\gamma^{i}\cup R^{i}\subset V^{i}$.
Finally, since the probability that the $R^{i}$-s starting from $x$
follow $\gamma^{i}$ until $y^{i}$ is $\approx1$ we get $\mathbb{P}(\mathcal{N}\cap\mathcal{V})\approx\rho^{-\xi}.$
To finish the sublemma, notice that $\mathbb{P}(\mathcal{N})\leq C(G)\rho^{-\xi}$
follows from lemma \ref{lem:kofup}.
\end{proof}
\begin{sublem}\label{sublem:NtHRt}

Let $x\in A$ and let $R^{1}$ and $R^{2}$ be two random walks starting
from $x$. Define $\mathcal{H}:=T^{1}(w)<T_{v,4r}^{1}$ and $\mathcal{N}':=\{ R^{1}[0,T^{1}(w)]\cap R^{2}[1,\infty]=\emptyset\}$.
Then\[
\mathbb{P}(\mathcal{N}'\cap\mathcal{H})\begin{cases}
\geq c(G)r^{-1-\xi} & x\textrm{ is }\mathcal{C}\textrm{-capable}\\
=0 & \textrm{otherwise}.\end{cases}\]

\end{sublem}

\begin{proof}
[Subproof]We use sublemma \ref{sublem:NERxi} with $\rho=\frac{1}{4}r$
and get that (assuming $x$ is $\mathcal{C}$-capable), that \begin{equation}
\mathbb{P}(\mathcal{N}\cap\mathcal{V})\approx r^{-\xi}.\label{eq:NVRnz}\end{equation}
Examining the structure of the $V^{i}$-s it is not difficult to see
that one may construct six domains $\mathcal{S}^{i},\mathcal{H}^{i},\mathcal{D}^{i}\subset B(0,5)$
with the following properties (see figure \ref{cap:SDH})%
\begin{figure}
\input{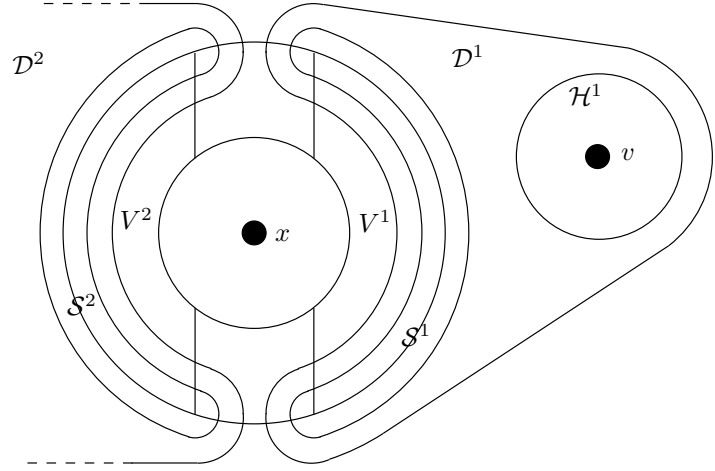}

\caption{\label{cap:SDH}The $\mathcal{S}^{i}$, $\mathcal{H}^{i}$ and $\mathcal{D}^{i}$
in the figure are actually $v+r\mathcal{S}^{i}$, $v+r\mathcal{H}^{i}$
and $v+r\mathcal{D}^{i}$ respectively. $\mathcal{H}^{2}$ is not
shown, imagine it {}``far away'' inside $\mathcal{D}^{2}$.}
\end{figure}

\begin{enumerate}
\item $\overline{\mathcal{S}^{i}},\overline{\mathcal{H}^{i}}\subset\mathcal{D}^{i}$
and $\overline{\mathcal{S}^{i}}\cap\overline{\mathcal{H}^{i}}=\emptyset$. 
\item \label{enu:SiSi}If $r>C(G)$ then $\partial V^{i}\cap\partial B(x,\frac{1}{4}r)\subset v+r\mathcal{S}^{i}$. 
\item \label{enu:MiD3i}$\mathcal{D}^{1}\cap\mathcal{D}^{2}=\emptyset$
and, if $r>C(G)$ then $V^{i}\cap(v+r\mathcal{D}^{3-i})=\emptyset$.
\item $\mathcal{D}^{1}\subset B(0,3)$\emph{,} $B(0,\frac{1}{2})\subset\mathcal{H}^{1}$,
and $B(0,4)\cap\mathcal{H}^{2}=\emptyset$.
\item The collection of $\mathcal{S}^{i},\mathcal{H}^{i},\mathcal{D}^{i}$
for all $x,v$ and $r$ satisfies the conditions of lemma \ref{lem:DD},
\ref{enu:lemDD-K}.
\end{enumerate}
Condition \ref{enu:SiSi} ensures that under the event $\mathcal{V}$
we have $R^{i}(T_{x,\rho}^{i})\in v+r\mathcal{S}^{i}$. Hence we can
apply lemma \ref{lem:DDp} with $\mathcal{D}^{1},\mathcal{S}^{1},\mathcal{H}^{1}$and
lemma \ref{lem:DD} with $\mathcal{D}^{2},\mathcal{S}^{2},\mathcal{H}^{2}$
for the continuation of $R^{i}$ after $T_{x,\rho}^{i}$. We get that\begin{align*}
\mathbb{P}(T^{1}(w)<T^{1}(\partial(v+r\mathcal{D}^{1}))\,|\, R^{1}[0,T_{x,\rho}^{1}]\subset V^{1}) & >c(G)/r.\\
\mathbb{P}(T_{v,4r}^{2}<T^{2}(\partial(v+r\mathcal{D}^{2}))\,|\, R^{2}[0,T_{x,\rho}^{2}]\subset V^{2}) & >c(G)\end{align*}
On the other hand, condition \ref{enu:MiD3i} ensures that if $R^{i}[T_{x,\rho}^{i},T(v+r\mathcal{H}^{i})]\subset v+r\mathcal{D}^{i}$
then $R^{1}[0,T(w)]\cap R^{2}[T_{x,\rho}^{2},T(v+r\mathcal{H}^{2})]=\emptyset$
and vice versa. Together with (\ref{eq:NVRnz}) we get\begin{gather*}
\mathbb{P}\left(\mathcal{N}''\cap\left\{ R^{1}[0,T^{1}(w)]\subset B(v,3r)\right\} \right)>c(G)r^{-1-\xi}.\\
\mathcal{N}'':=\left\{ R^{1}[0,T^{1}(\{ w\})]\cap R^{2}[1,T_{v,4r}^{2}]=\emptyset\right\} \end{gather*}
This ends the sublemma since lemma \ref{lem:DDunbounded} with $\mathcal{D}=\mathbb{R}^{3}$
and $\mathcal{H}=B(0,3)$ shows that the probability of $R^{1}$ to
never hit $B(v,3r)$ after hitting $\partial B(v,4r)$ is $\geq c(G)$.
\end{proof}
Let $\mathcal{E}$ be an event on a space of curves. We say that $\mathcal{E}$
is \textbf{loop-monotone} if $\mathcal{E}(\gamma)\Rightarrow\mathcal{E}(\gamma')$
whenever $\gamma$ is $\gamma'$ with some loops added. In other words,
adding loops can only hurt $\mathcal{E}$. A typical example of a
loop-monotone event is $\{ x\in\mathcal{C}\}$ for some $x$ ($\mathcal{C}$
from the statement of lemma \ref{lem:CARb0}). We shall use loop-monotonicity
to encapsulate the idea of time reversal in a convenient way in the
following sublemma:

\begin{sublem}\label{sublem:timerevers}

Let $v,w,x\in G$ and let $\mathcal{E}\subset\{ T(x)<T_{v,r}\}$ be
a loop-monotone event on the space of curves on $G$ starting from
$w$. Then\[
\mathbb{P}^{w}(\mathcal{E}(R\left[0,\infty\right[))\approx\mathbb{P}^{1,x,2,x}\left(\mathcal{E}\left(\overleftarrow{R^{1}}\cup R^{2}\right)\cap\left\{ T^{1}(w)<T_{v,r}^{1}\right\} \right)\]
where the notation $\overleftarrow{R^{1}}\cup R^{2}$ means taking
$R^{1}[0,T^{1}(w)]$, reversing it (so that it starts from $w$ and
ends at $x$), and concatenating $R^{2}[1,\infty]$ at its end $x$.

\end{sublem}

\begin{proof}
[Subproof]Denote $\mathcal{W}=\#\{ t\geq1:R(t)=w\}$. The loop monotonicity
of $\mathcal{E}$ gives \begin{equation}
\mathbb{P}^{w}\left(\mathcal{E}\,|\,\mathcal{W}=k\right)\leq\mathbb{P}\left(\mathcal{E}\,|\,\mathcal{W}=0\right)\quad\forall k>0\label{eq:FWksW0}\end{equation}
since conditioning by $\mathcal{W}=k$ is equivalent to adding $k$
closed paths from $w$ to itself and then starting a walk conditioned
to have $\mathcal{W}=0$. Hence we get\begin{align}
\mathbb{P}(\mathcal{E}) & =\sum_{k\geq0}\mathbb{P}(\mathcal{E}\,|\,\mathcal{W}=k)\mathbb{P}(\mathcal{W}=k)\stackrel{(\ref{eq:FWksW0})}{\leq}\mathbb{P}(\mathcal{E}\,|\,\mathcal{W}=0)\leq\nonumber \\
 & \stackrel{(*)}{\leq}C(G)\mathbb{P}(\mathcal{E}\cap\{\mathcal{W}=0\}).\label{eq:noWW2}\end{align}
where $(*)$ comes from the transience of $G$. 

Next we use the time-symmetry of random walk in the form (\ref{eq:symT})
for the portion of the walk between $w$ and $x$. We get\begin{equation}
\mathbb{P}^{w}(\mathcal{E}\cap\{\mathcal{W}=0\})\approx\mathbb{P}^{1,x,2,x}\left(\mathcal{E}\left(\overleftarrow{R^{1}}\cup R^{2}\right)\cap\left\{ T^{1}(w)<T^{1}(\{ x\}\cup\partial B(v,r))\right\} \right)\label{eq:timesym}\end{equation}
where the $\approx$ sign hides the bounded quantity $\omega(w)/\omega(x)$. 

The last step is defining $\mathcal{G}:=\mathcal{E}\left(\overleftarrow{R^{1}}\cup R^{2}\right)\cap\left\{ T^{1}(w)<T_{v,r}^{1}\right\} $
and $\mathcal{X}:=\#\{ t\in[1,T^{1}(w)]:R^{1}(t)=x\}$. The loop-monotonicity
of $\mathcal{E}$ gives\[
\mathbb{P}^{1,x,2,x}\left(\mathcal{G}\,|\,\mathcal{X}=k,\, R^{2}=\gamma\right)\leq\mathbb{P}\left(\mathcal{G}\,|\,\mathcal{X}=0,\, R^{2}=\gamma\right)\quad\forall k>0,\,\gamma\]
which gives, like (\ref{eq:noWW2}), \[
\mathbb{P}(\mathcal{G}\,|\, R^{2}=\gamma)\approx\mathbb{P}(\mathcal{G}\cap\{\mathcal{X}=0\}\,|\, R^{2}=\gamma)\quad\forall\gamma\]
and summing over all paths $\gamma$ starting from $x$ we get\begin{equation}
\mathbb{P}(\mathcal{G})\approx\mathbb{P}(\mathcal{G}\cap\{\mathcal{X}=0\}).\label{eq:returnloops}\end{equation}
(\ref{eq:noWW2}), (\ref{eq:timesym}) and (\ref{eq:returnloops})
together finish the proof.
\end{proof}
\begin{sublem}\label{sublem:condintrst}For any $x\in A$,\[
\mathbb{P}(x\in\mathcal{C})\begin{cases}
>c(G)r^{-1-\xi} & x\textrm{ is }\mathcal{C}\textrm{-capable}\\
=0 & \textrm{otherwise.}\end{cases}\]
\end{sublem}

\begin{proof}
[Subproof]This is an immediate consequence of sublemmas \ref{sublem:timerevers}
and \ref{sublem:NtHRt}.
\end{proof}
This completes what we would need for the estimate of the first moment,
and we move to the second moment, which is not really all that more
complicated --- the complication from the fact that it is second moment
are partially compensated by the fact that we need an upper bound
rather than a lower.

\begin{sublem}\label{sublem:1mlow}

For any $x\in A$,\[
\mathbb{P}(x\in\mathcal{C})\leq C(G)r^{-1-\xi}.\]

\end{sublem}

\begin{proof}
[Subproof]Sublemma \ref{sublem:timerevers} shows that\begin{align*}
\mathbb{P}^{w}(x\in\mathcal{C}) & \approx\mathbb{P}^{1,x,2,x}\left(\left\{ T^{1}(w)<T_{v,4r}^{1}\right\} \cap\left\{ R^{1}[0,T^{1}(w)]\cap R^{2}\left[1,\infty\right[=\emptyset\right\} \right)\leq\\
 & \leq\mathbb{P}^{1,x,2,x}\left(\left\{ T^{1}(w)<T_{v,4r}^{1}\right\} \cap\left\{ R^{1}[0,T_{x,r/4}^{1}]\cap R^{2}[1,T_{x,r/4}^{2}]=\emptyset\right\} \right)\\
\intertext{\textrm{and lemma \ref{lem:condnomtr} shows that}} & \approx\mathbb{P}^{1,x}\left(T^{1}(w)<T_{v,4r}^{1}\right)\mathbb{P}^{1,x,2,x}\left(R^{1}[0,T_{x,r/4}^{1}]\cap R^{2}[1,T_{x,r/4}^{2}]=\emptyset\right).\end{align*}
 Sublemma \ref{sublem:NERxi} shows that the term on the right is
$\leq C(G)r^{-\xi}$ while (\ref{eq:GvwS-lower-bound}) shows that
the term on the left is $\approx r^{-1}$.
\end{proof}
\begin{sublem}\label{sublem:2m}

For any $x_{1},x_{2}\in A$,\[
\mathbb{P}(x_{1},x_{2}\in\mathcal{C})\leq C(G)(r|x_{1}-x_{2}|)^{-1-\xi}.\]

\end{sublem}

\begin{proof}
[Subproof]First let us note that it is possible to assume $|x_{1}-x_{2}|>C(G)$
since otherwise $\mathbb{P}(x_{1},x_{2}\in\mathcal{C})\leq\mathbb{P}(x_{1}\in\mathcal{C})$
and then sublemma \ref{sublem:1mlow} applies. Moreover, it is enough
to prove that\begin{gather*}
\mathbb{P}\left(\left\{ x_{1},x_{2}\in\mathcal{C}\right\} \cap\mathcal{O}\right)\leq C(G)(r|x_{1}-x_{2}|)^{-1-\xi}\end{gather*}
where $\mathcal{O}:=\{ T_{1}<T_{2}\}$ and $T_{i}$ is the last time
$R$ is in $x_{i}$. The other case is just a renaming of $x_{1}$
and $x_{2}$. 

Define now $\rho:=\frac{1}{16}|x_{1}-x_{2}|$ and $x=(x_{1}+x_{2})/2$.
Denote $\mathcal{X}=\{ x_{1},x_{2}\in\mathcal{C}\}\cap\mathcal{O}$.
$\mathcal{X}$ is loop-monotone, hence we may use sublemma \ref{sublem:timerevers}
for $x_{1}$ and get\begin{align*}
\mathbb{P}^{w}(\mathcal{X}) & \approx\mathbb{P}^{1,x_{1},2,x_{1}}\Big(\left\{ T^{1}(w)<T_{v,4r}^{1}\right\} \cap\left\{ R^{1}\cap R^{2}=\emptyset\right\} \cap\\
 & \quad\Big\{\exists t<T_{v,4r}^{2}:\left\{ R^{2}(t)=x_{2}\right\} \cap\left\{ \left(R^{1}\cup R^{2}[0,t])\right)\cap R^{2}\left[t+1,\infty\right[=\emptyset\right\} \Big\}\Big).\end{align*}
 where the $R^{i}$-s in the expression $R^{1}\cap R^{2}=\emptyset$
stand for the walks until their natural ending, namely $R^{1}[0,T^{1}(w)]$
and $R^{2}[1,\infty[$ respectively. Denote $R^{2}(T_{x_{1},2\rho}^{2})$
by $y$ and {}``stop'' $R^{2}$ there, and consider the rest of
$R^{2}$ as a new random walk $R^{3}$ starting from $y$. We get\begin{align}
\mathbb{P}^{w}(\mathcal{X}) & \approx\sum_{y\in\partial B(x_{1},2\rho)}\mathbb{P}^{1,x_{1},2,x_{1},3,y}\Big(\left\{ T^{1}(w)<T_{v,4r}^{1}\right\} \cap\left\{ R^{2}(T_{x_{1},2\rho}^{2})=y\right\} \cap\nonumber \\
 & \qquad\cap\left\{ R^{1}\cap(R^{2}\cup R^{3})=\emptyset\right\} \cap\Big\{\exists t<T_{v,4r}^{3}:\left\{ R^{3}(t)=x_{2}\right\} \cap\nonumber \\
 & \qquad\cap\left\{ \left(R^{1}\cup R^{2}\cup R^{3}[0,t]\right)\cap R^{3}\left[t+1,\infty\right[=\emptyset\right\} \Big\}\Big)\label{eq:onesym}\end{align}
where $R^{2}$ stands for $R^{2}[1,T_{x_{1},2\rho}^{2}]$ and $R^{3}$
stands for $R^{3}[0,\infty[$. We use sublemma \ref{sublem:timerevers}
again, this time for the random walk $R^{3}$ and the point $x_{2}$
(it is easy to see that the corresponding event is loop-monotone for
any value of $R^{1}$ and $R^{2}$). We get\begin{align}
\mathbb{P}^{w}(\mathcal{X}) & \approx\sum_{y}\mathbb{P}^{1,x_{1},2,x_{1},3,x_{2},4,x_{2}}\Big(\left\{ T^{1}(w)<T_{v,4r}^{1}\right\} \cap\left\{ R^{2}(T_{x_{1},2\rho}^{2})=y\right\} \cap\nonumber \\
 & \cap\left\{ T^{3}(y)<T_{v,4r}^{3}\right\} \cap\left\{ R^{1}\cap(R^{2}\cup R^{3}\cup R^{4})=\emptyset\right\} \cap\nonumber \\
 & \cap\left\{ (R^{1}\cup R^{2}\cup R^{3})\cap R^{4}=\emptyset\right\} \Big)\label{eq:R1234}\end{align}
where $R^{3}$ stands for $R^{3}[0,T^{3}(y)]$ and $R^{4}$ stands
for $R^{4}[1,\infty[$. Reducing slightly the non-intersecting sections
we may write\begin{align*}
\mathbb{P}^{w}(\mathcal{X}) & \leq C(G)\sum_{y}\mathbb{P}\Big(\textrm{same first three conditions }\cap\\
 & \cap\left\{ R^{1}[0,T_{x_{1},\rho}^{1}]\cap R^{2}[1,T_{x_{1},\rho}^{2}]=\emptyset\right\} \cap\left\{ R^{3}[0,T_{x_{2},\rho}^{3}]\cap R^{4}[1,T_{x_{2},\rho}^{4}]=\emptyset\right\} \cap\\
 & \cap\left\{ R^{1}[T_{x,18\rho}^{1},T_{x,r/4}^{1}]\cap R^{4}[T_{x,18\rho}^{4},T_{x,r/4}^{4}]=\emptyset\right\} \Big).\end{align*}

Denote the three non-intersection events above by $\mathcal{N}_{1}$,
$\mathcal{N}_{2}$ and $\mathcal{N}_{3}$ by order. We understand
that if $18\rho>r/4$ then $\mathcal{N}_{3}$ is considered to always
be satisfied. Now, sublemma \ref{sublem:NERxi} shows that \[
\mathbb{P}(\mathcal{N}_{1})\leq C(G)\rho^{-\xi}\quad\mathbb{P}(\mathcal{N}_{2})\leq C(G)\rho^{-\xi}\]
and since these events are independent the probability of their intersection
is $\leq C(G)\rho^{-2\xi}$. Assume for a moment that $18\rho\leq r/4$.
Then we use lemma \ref{lem:condnomtr} for the ball $B(x,9\rho$)
and get\[
\mathbb{P}(\mathcal{N}_{1}\cap\mathcal{N}_{2}\cap\mathcal{N}_{3})\approx\mathbb{P}(\mathcal{N}_{1}\cap\mathcal{N}_{2})\mathbb{P}(\mathcal{N}_{3})\]
and theorem \ref{thm:kof} shows that $\mathbb{P}(\mathcal{N}_{3})\leq C(G)(\rho/r)^{\xi}$,
so in total \begin{equation}
\mathbb{P}(\mathcal{N}_{1}\cap\mathcal{N}_{2}\cap\mathcal{N}_{3})\leq C(G)(r\rho)^{-\xi}.\label{eq:N123xi}\end{equation}
If $18\rho>r/4$ then $(r\rho)^{-\xi}\approx\rho^{-2\xi}$ and (\ref{eq:N123xi})
is again satisfied, so we can continue without the assumption $18\rho\leq r/4$.

Finally we need to accommodate the various hitting and exit conditions
in (\ref{eq:R1234}). Let $\mathcal{E}$ be the end points of the
portions of the $R^{i}$-s needed for the $\mathcal{N}_{i}$-s, namely\[
\mathcal{E}:=(R^{1}(T_{x,r/4}^{1}),R_{x_{1},\rho}^{2},R_{x_{2},\rho}^{3}).\]
 For the condition $T^{3}(y)<T_{v,4r}^{3}$ we use the fact that for
any point $z$ where $R^{3}$ exits $B(x_{2},\rho)$ we have $|z-y|\geq\rho$
and therefore the estimate of the harmonic potential (\ref{eq:atrans})
gives\begin{equation}
\mathbb{P}(T^{3}(y)<T_{v,4r}^{3}\,|\,\mathcal{E})\leq\mathbb{P}(T^{3}(y)<\infty\,|\,\mathcal{E})\leq C(G)\rho^{-1}.\label{eq:EE1}\end{equation}
A similar argument for $R^{1}$ gives\begin{equation}
\mathbb{P}(T^{1}(w)<T_{v,4r}^{1}\,|\,\mathcal{E})\leq C(G)r^{-1}.\label{eq:EE2}\end{equation}
Conditioning over $\mathcal{E}$ the events of (\ref{eq:EE1}) and
(\ref{eq:EE2}); $\mathbb{P}(R^{2}(T_{x_{1},2\rho}^{2})=y)$ and $\mathcal{N}_{1}\cap\mathcal{N}_{2}\cap\mathcal{N}_{3}$
are all independent. Hence we get \begin{align*}
\mathbb{P}^{w}(\mathcal{X}) & \stackrel{(*)}{\leq}\sum_{E}C(G)\rho^{-1}r^{-1}\mathbb{P}(\mathcal{N}_{1}\cap\mathcal{N}_{2}\cap\mathcal{N}_{3}\cap\{\mathcal{E}=E\})\:\cdot\\
 & \qquad\qquad\cdot\:\sum_{y\in\partial B(x_{1},2\rho)}\mathbb{P}(R^{2}(T_{x_{1},2\rho}^{2})=y\,|\,\mathcal{E}=E)=\\
 & =C(G)\rho^{-1}r^{-1}\mathbb{P}(\mathcal{N}_{1}\cap\mathcal{N}_{2}\cap\mathcal{N}_{3})\stackrel{(\ref{eq:N123xi})}{\leq}C(G)(\rho r)^{-1-\xi}.\end{align*}
where in $(*)$ we used (\ref{eq:EE1}), (\ref{eq:EE2}) and independence.
\end{proof}

\begin{proof}
[Proof of lemma \ref{lem:CARb0}]Let \[
\mathcal{X}=\#\{\mathcal{C}\cap R^{2}[0,T_{v,4r}^{2}]\cap A\}.\]
Sublemma \ref{sublem:condintrst} shows that\begin{align*}
\mathbb{E}\mathcal{X} & =\sum_{x\in A}\mathbb{P}\left(x\in\mathcal{C}\right)\mathbb{P}\left(x\in R^{2}[0,T_{v,4r}^{2}]\right)\geq\\
 & \stackrel{(*)}{\geq}c(G)r^{-2-\xi}\#\{ x\in A:x\textrm{ is }\mathcal{C}\textrm{ capable}\}\stackrel{(**)}{\geq}c(G)r^{1-\xi}\end{align*}
where in $(*)$ we used sublemma \ref{sublem:condintrst} to estimate
$\mathbb{P}(x\in\mathcal{C})$ and (\ref{eq:GvwS-lower-bound}) to
estimate $\mathbb{P}\left(x\in R^{2}[0,T_{v,4r}^{2}]\right)$; and
$(**)$ follows from sublemma \ref{sublem:anyball}. Correspondingly
we have\begin{align*}
\mathbb{E}\mathcal{X}^{2} & =\sum_{x_{1},x_{2}\in A}\mathbb{P}(x_{1},x_{2}\in\mathcal{C})\mathbb{P}\left(x_{1},x_{2}\in R^{2}[0,T_{v,4r}^{2}]\right)\leq\\
 & \stackrel{(*)}{\leq}C(G)\sum_{x_{1},x_{2}\in A}(r|x_{1}-x_{2}|)^{-2-\xi}\stackrel{(**)}{\leq}C(G)r^{-2-\xi}\sum_{x_{1}\in A}\sum_{n=1}^{\log_{2}r}2^{n(1-\xi)}\leq\\
 & \stackrel{(\dagger)}{\leq}C(G)r^{2-2\xi}\\
\end{align*}
where $(*)$ follows from sublemma \ref{sublem:2m} for $\mathbb{P}(x_{1},x_{2}\in\mathcal{C})$
and (\ref{eq:atrans}) for $\mathbb{P}(x_{1},x_{2}\in R^{2}[0,T_{v,4r}^{2}])$;
where $(**)$ comes from the volume estimate $\#\{ x_{2}:|x_{1}-x_{2}|\in[2^{n},2^{n+1}[\}\approx2^{3n}$
for $n>C(G)$ since our graph $G$ is roughly isometric to $\mathbb{R}^{3}$;
and where $(\dagger)$ comes from the same volume estimate since $r>C(G)$,
and (finally!) from $\xi<1$. The well known inequality $\mathbb{P}(\mathcal{X}>0)\geq(\mathbb{E}\mathcal{X})^{2}/\mathbb{E}\mathcal{X}^{2}$
now finishes the lemma.
\end{proof}
\begin{cor*}
Under the assumptions of lemma \ref{lem:CARb0}, \[
\mathbb{P}^{1,w}(\forall z\in B(v,{\textstyle \frac{1}{2}}r),\,\mathbb{P}^{2,z}(\mathcal{C}\cap R^{2}[0,T_{v,4r}^{2}]\cap A\neq\emptyset)>c(G))>c(G).\]

\end{cor*}
\begin{proof}
Denote the event inside the inner $\mathbb{P}$ symbol by $\mathcal{E}$.
Then lemma \ref{lem:CARb0} shows that $\mathbb{P}^{1,w,2,v}(\mathcal{E})\geq c(G)$.
This shows that\[
\mathbb{P}^{1,w}(\mathbb{P}^{2,v}(\mathcal{E})>c(G))>c(G).\]
Now, for any infinite path $\gamma$ starting from $w$, the function
\[
f(z)=\mathbb{P}^{2,z}(\mathcal{E}\,|\, R^{1}[0,\infty[=\gamma)\]
is harmonic outside $A$ and in particular in $B(v,\frac{1}{2}r)$.
Hence Harnack's inequality (lemma \ref{lem:Harnack}) shows that $\min_{z\in B(v,\frac{1}{2}r)}f(z)\geq cf(y)$
which proves the corollary.
\end{proof}

\subsection{Conditioned random walks}

\begin{lem}
\label{lem:RAH}Let $G$ be a three dimensional isotropic graph. Let
$v\in G$ and let $H\subset\mathbb{R}^{3}$ be a closed half space
with $v\in\partial_{\textrm{cont}}H$. Let $r>C(G)$ and let $\Gamma\subset B(v,r)$,
$d(\Gamma,H)>C_{\ref{C:arc}}$. Let $R$ be a random walk starting
from $v$. Then\newc{c:RAH}\[
\mathbb{P}(R(T_{v,r})\in H\,|\, R[0,T_{v,r}]\cap\Gamma=\emptyset)\geq c_{\ref{c:RAH}}(G).\]

\end{lem}
($C_{\ref{C:arc}}$ is from lemma \ref{lem:arc}, page \pageref{lem:arc}.
In particular $d(\Gamma,H)>C_{\ref{C:arc}}$ implies that the set
of paths from $v$ to $\partial B(v,r)$ not intersecting $\Gamma$
is non-empty --- use the lemma for a translation of $H$ by $C_{\ref{C:arc}}$)

\begin{proof}
The equivalent question for a Brownian motion can be solved by reflecting
through $\partial_{\textrm{cont}}H$ the last section of the motion
not intersecting $\partial_{\textrm{cont}}H$, with the result that
the corresponding probability is $\geq\frac{1}{2}$. Our proof is
a discrete version of this idea. Formally, denote by $\mathcal{H}$
(respectively $\mathcal{H}^{-}$) the space of all paths from $v$
to $\partial B(v,r)\cap H$ (respectively $\partial B(v,r)\setminus H$)
not intersecting $\Gamma$. We shall dissect $\mathcal{H}^{-}$ to
disjoint sets $N_{\gamma}^{-}$ indexed by $\mathcal{G}$:\[
\mathcal{H}^{-}=\bigcup_{\gamma\in\mathcal{G}}N_{\gamma}^{-}\]
and map each $N_{\gamma}^{-}$ into a set $N_{\gamma}\subset\mathcal{H}$
such that the following holds:
\begin{enumerate}
\item \label{enu:prob}$\mathbb{P}(N_{\gamma})\geq c(G)\mathbb{P}(N_{\gamma}^{-})$.
\item \label{enu:cover}Every path $h\in\mathcal{H}$ is contained in at
most $C(G)$ different $N_{\gamma}$-s.
\end{enumerate}
Together these properties show that $\mathbb{P}(\mathcal{H})\geq c(G)\mathbb{P}(\mathcal{H}^{-})$
or equivalently\[
\mathbb{P}(R(T_{v,r})\in H\,|\, R[0,T_{v,r}]\cap\Gamma=\emptyset)\geq c(G)\mathbb{P}(R(T_{v,r})\not\in H\,|\, R[0,T_{v,r}]\cap\Gamma=\emptyset)\]
which would conclude the lemma.

The set $\mathcal{G}$ is the set of all paths $\gamma$ in $B(v,r)$
avoiding $\Gamma$ such that $x:=\gamma(\len\gamma)\not\in H$ but
its neighbor $x':=\gamma(\len\gamma-1)\in H$. For each $\gamma\in\mathcal{G}$
the set $N_{\gamma}^{-}$ is the set of all paths that follow $\gamma$
until its end and then avoid hitting $H\cup\Gamma$ until hitting
$\partial B(v,r)$. It is clear that $\left\{ N_{\gamma}^{-}\right\} _{\gamma\in\mathcal{G}}$
are disjoint sets covering $\mathcal{H}^{-}$. Take one $\gamma\in\mathcal{G}$,
let $x$ be its end and denote $\rho:=d(x,\partial B(v,r))$. Clearly
$\mathbb{P}(N_{\gamma}^{-})$ is the probability that $R$ follows
$\gamma$ (denote it by $p_{\gamma}$) multiplied by the escape probability\[
\mathbb{P}^{x}(T_{v,r}<T(\Gamma\cup H))\leq\mathbb{P}^{x}(T_{v,r}<T(H))\leq\mathbb{P}^{x}(T_{x,\rho}<T(H))\stackrel{(*)}{\leq}C(G)/\rho\]
where $(*)$ comes from theorem \ref{lem:escape}.

We shall now construct $N_{\gamma}$ under the assumption that $\rho$
is bigger than some constant $\rho_{0}(G)$. The value of $\rho_{0}$
will be fixed later on, but for now we need $\rho_{0}>4C_{\ref{C:arc}}$.
We use lemma \ref{lem:arc} with the point $x'$, the radius $4C_{\ref{C:arc}}$
and with the half-space $H'=H+B(0,C_{\ref{C:arc}})$ and we get that
there exists a simple path $\delta'\subset\overline{B(x',4C_{\ref{C:arc}})}\setminus H'$
from $x'$ to $\partial B(x',4C_{\ref{C:arc}})\cap\{ y\in H:d(y,\partial_{\textrm{cont}}H>C_{\ref{C:arc}}\}$.
Let $\delta=\gamma\cup(x,x')\cup\delta'$. Note that $\Gamma\cap H'=\emptyset$
and therefore also $\Gamma\cap\delta=\emptyset$. Let $N_{\gamma}$
be the family of all paths that follow $\delta$ until its end and
then stay inside $H$ until they exit $B(v,r)$. If $\rho\leq\rho_{0}$
simply let $\delta'$ be the shortest path from $x$ to $\partial B(v,r)\setminus H$
not intersecting $\Gamma$ and let $N_{\gamma}$ to contain only the
path $\gamma\cup\delta'$.

The lemma will be concluded once we show \ref{enu:prob} and \ref{enu:cover}.
To see \ref{enu:prob}, first note that the case when $\rho\leq\rho_{0}$
is obvious since then $\mathbb{P}(N_{\gamma}^{-})\approx p_{\gamma}\approx\mathbb{P}(N_{\gamma})$.
In the case $\rho>\rho_{0}$, the length of $\delta'$ is $\leq C(G)$
so the probability to follow $\delta$ is $\geq c(G)p_{\gamma}$.
We use lemma \ref{lem:arc} again to get that the probability of a
random walk starting from $y:=\delta(\len\delta)$ to hit \[
\alpha:=\partial B(y,\rho)\cap\{ z:d(z,\partial Q)>{\textstyle \frac{1}{2}}\rho\}\]
before $\partial H$ is $\geq C(G)/\rho$. Finally, lemma \ref{lem:DD}
shows that for any $z\in\alpha$ a random walk starting from $z$
has a probability $>c(G)$ to exit $B(v,r)$ before hitting $\partial H$.
To use lemma \ref{lem:DD} we need to assume that $\rho$ is large
enough, and this is the condition for $\rho_{0}$ which can now be
fixed. All three together give \ref{enu:prob}.

As for \ref{enu:cover}, it is easy to see that every $h\in\mathcal{H}$
can belong to only boundedly many $N_{\gamma}$ for which $\rho\leq\rho_{0}$.
Hence examine the case $\rho>\rho_{0}$ and let $h\in\mathcal{H}$.
If $h\in N_{\gamma}$ then $x\not\in H$ but after $y$ all points
of $h$ are in $H$ and the path between $x$ and $y$ is in $B(x,C(G))$.
Therefore if we define $e(h)$ as the last vertex in $h\setminus H$
we know that $x\in B(e(h)),C(G))$ and in particular has just $C(G)$
possibilities. Since $\gamma$ is simply the part of $h$ up to $x$
we see that it too has only $C(G)$ possibilities which shows \ref{enu:cover}
and the lemma. 
\end{proof}
\begin{lem}
\label{lem:IGamIC}Let $G$ be a three dimensional isotropic graph
and let $\epsilon>0$. Then there exist a $q=q(G)>0$ and a $\delta=\delta(\epsilon,G)>0$
such that the following holds: Let $v\in G$, let $r>C(\epsilon,G$)
and let $s\in[r,2r-\epsilon r]$. Let $\Gamma\subset B(v,s)$ be some
set such that \begin{equation}
\mathbb{P}^{v}(R[0,T_{v,4r}]\cap\Gamma\neq\emptyset)\leq\delta.\label{eq:IGamICassm}\end{equation}
Let $w\in\partial B(v,s)$ be admissible (see below). Then \begin{gather}
\mathbb{P}^{1,w}(\forall y\in B(w,\epsilon r),\,\mathbb{P}^{2,y}(\mathcal{I})>q)\,|\, R^{1}[0,T_{v,4r}^{1}]\cap\Gamma=\emptyset)>q,\label{eq:IGamIC}\\
\mathcal{I}:=\left\{ \cut(R^{1}[0,T_{v,4r}^{1}];T_{w,\epsilon r}^{1})\cap R^{2}[0,T_{v,4r}^{2}]\neq\emptyset\right\} .\nonumber \end{gather}

\end{lem}
We call $w$ admissible if there exists a path $\gamma\subset\overline{B(w,16C_{\ref{C:arc}})}$
starting from $w$ and ending outside $B(v,|v-w|+2C_{\ref{C:arc}})$
which does not intersect $\Gamma$ (the constant $16C_{\ref{C:arc}}$
will be used in lemma \ref{lem:L1L2} below to show that many admissible
points exist). 

In words, the lemma says that the fact that $\cut(R^{1})$ is hittable
does not change if one condition by not hitting $\Gamma$, even if
one starts very close to $\Gamma$ --- the only condition is that
$\Gamma$ is not very hittable ($<\delta$) from far away ($v$).
The fact that $\epsilon$ affects only $\delta$ but not $q$ will
play a significant role later on.

\begin{proof}
Let $\lambda=\lambda(G)$ be some parameter that will be fixed later.
Denote also $\mu:=\epsilon/4\lambda$ and $\rho=\mu r$. Let $\gamma$
be the path from the definition of admissibility of $w$ and assume
w.l.o.g.~that it is simple (say by taking $\LE$). We use lemma \ref{lem:RAH}
with the starting point being $w':=\gamma(\len\gamma)$; with $H'$
being the half space orthogonal to the segment $[v,w]$ such that
$w'\in\partial_{\textrm{cont}}H'$; and with the radius some $\rho'$
to be fixed later. (Note that the condition of lemma \ref{lem:RAH}
$d(\Gamma,H')>C_{\ref{C:arc}}$ will be fulfilled if $r$ is sufficiently
large). We get that if $\rho'>C(G)$ then \[
\mathbb{P}^{1,w'}(R^{1}(T_{w',\rho'}^{1})\in H'\,|\, R^{1}[0,T_{w',\rho'}^{1}]\cap\Gamma=\emptyset)\geq c(G).\]
Let $H$ be the translation of $H'$ such that $w\in\partial_{\textrm{cont}}H$.
On one side, the probability that a random walk $R^{1}$ starting
from $w$ will follow $\gamma$ until $w'$ is $\geq c(G)$. On the
other side, if $\rho'=\rho-C(G)$ for some $C(G)$ sufficiently large
then for any point $x\in\partial B(w',\rho')\cap H'$ there is a probability
$\geq c(G)$ that a random walk starting from $x$ will hit $\partial B(w,\rho)\cup\Gamma$
in $\partial B(w,\rho)\cap H$. All these allow us to drop the $'$
notations and we get\begin{equation}
\mathbb{P}^{1,w}(R^{1}(T_{w,\rho}^{1})\in H\,|\, R^{1}[0,T_{w,\rho}^{1}]\cap\Gamma=\emptyset)\geq c(G).\label{eq:HHc}\end{equation}
Denote this event by $\mathcal{H}$.

Next define $\mathcal{C}=\cut(R^{1}[T_{w,\rho}^{1},T_{v,4r}^{1}];T_{w,\epsilon r}^{1})$
and $A=B(w,\frac{1}{2}\epsilon r)\setminus\overline{B(w,\frac{1}{4}\epsilon r)}$.
The corollary to lemma \ref{lem:CARb0} shows that, for any $x\in\partial B(w,\rho)$
and some $c(G)$,\begin{eqnarray}
\lefteqn{\mathbb{P}^{1,x}(\forall y\in B(w,{\textstyle \frac{1}{8}}\epsilon r),\,\mathbb{P}^{2,y}(\mathcal{C}\cap R^{2}\cap A\neq\emptyset)>c(G))\geq}\nonumber \\
 &  & \geq\mathbb{P}^{1,x}(\forall y\in B(w,{\textstyle \frac{1}{8}}\epsilon r),\nonumber \\
 &  & \qquad\mathbb{P}^{2,y}(\cut(R^{1}[0,\infty[;T_{w,\epsilon r}^{1})\cap R^{2}[0,T_{w,\epsilon r}^{2}]\cap A\neq\emptyset)>c(G))\geq\nonumber \\
 &  & \geq c(G)\label{eq:yBwer8}\end{eqnarray}
where the notation $R^{i}$ (e.g.~$R^{2}$ above) stands for $R^{i}[0,T_{v,4r}^{i}]$,
and where the $T_{w,\rho}^{1}$ in the definition of $\mathcal{C}$
is considered to be $0$ when starting from $x$. Since I promised
to prove the lemma for any $y\in B(w,\epsilon r)$, just note that
for any such $y$ we have that with probability $>c(G)$ the walk
$R^{2}$ hits $B(w,\frac{1}{8}\epsilon r)$ and therefore (\ref{eq:yBwer8})
holds for the larger ball too, i.e.\begin{equation}
\mathbb{P}^{1,x}(\forall y\in B(w,\epsilon r),\,\mathbb{P}^{2,y}(\mathcal{C}\cap R^{2}\cap A\neq\emptyset)>c(G))\geq c(G).\label{eq:JJc}\end{equation}
Denote this event by $\mathcal{J}$.

Next we take into consideration $\lambda$. Using (\ref{eq:fEfRT})
with Green's function $G(\cdot,w;\linebreak[4]B(v,4r))$ shows that\begin{equation}
\mathbb{P}^{1,x}(R^{1}[T_{w,\epsilon r/4}^{1},T_{v,4r}^{1}]\cap B(w,\rho)\neq\emptyset)\leq\frac{C(G)}{\lambda}.\label{eq:KKClam}\end{equation}
Denote this event by $\mathcal{K}$.

The next step is saying, roughly, {}``if $\Gamma$ is not hittable,
then conditioning by not hitting $\Gamma$ has no effect''. Formally,
we assume that, for some $\nu=\nu(G)$ to be fixed later \begin{equation}
\mathbb{P}^{1,x}(R^{1}\cap\Gamma\neq\emptyset)\leq\nu\quad\forall x\in\partial B(w,\rho)\cap H.\label{eq:assmpni}\end{equation}
As we shall see later, this assumption will be satisfied with a proper
choice of $\delta$. For now, this allows us to preform the following
calculation, which will return us to a walk starting from $w$:\begin{align}
\lefteqn{\mathbb{P}^{1,w}(\mathcal{J}\setminus\mathcal{K}\,|\, R^{1}\cap\Gamma=\emptyset)\geq}\nonumber \\
 & \quad\geq\mathbb{P}^{1,w}((\mathcal{J}\setminus\mathcal{K})\cap\mathcal{H})\,|\, R^{1}\cap\Gamma=\emptyset)=\nonumber \\
 & \quad=\sum_{x\in\partial B(w,\rho)\cap H}\mathbb{P}^{1,w}((\mathcal{J}\setminus\mathcal{K})\cap\left\{ R^{1}(T_{w,\rho}^{1})=x\right\} \,|\, R^{1}\cap\Gamma=\emptyset)=\nonumber \\
 & \quad\!\stackrel{(*)}{=}\!\sum_{x\in\partial B(w,\rho)\cap H}\mathbb{P}^{1,w}(R^{1}(T_{w,\rho}^{1})=x\,|\, R^{1}[0,T_{w,\rho}^{1}]\cap\Gamma=\emptyset)\cdot\nonumber \\
 & \qquad\qquad\cdot\frac{\mathbb{P}^{1,x}\left((\mathcal{J}\setminus\mathcal{K})\cap\left\{ R^{1}\cap\Gamma=\emptyset\right\} \right)}{\mathbb{P}^{1,w}(R^{1}\cap\Gamma=\emptyset\,|\, R^{1}[0,T_{w,\rho}^{1}]\cap\Gamma=\emptyset)}\geq\nonumber \\
 & \quad\!\!\stackrel{(**)}{\geq}\left(c(G)-\frac{C(G)}{\lambda}-\nu\right)\sum_{x}\mathbb{P}^{1,w}(R^{1}(T_{w,\rho}^{1})=x\,|\, R^{1}[0,T_{w,\rho}^{1}]\cap\Gamma=\emptyset)=\nonumber \\
 & \quad=\left(c(G)-\frac{C(G)}{\lambda}-\nu\right)\mathbb{P}^{1,w}(\mathcal{H}\,|\, R^{1}[0,T_{w,\rho}^{1}]\cap\Gamma=\emptyset)\geq\nonumber \\
 & \quad\!\stackrel{(\ref{eq:HHc})}{\geq}\!\left(c(G)-\frac{C(G)}{\lambda}-\nu\right)c(G)\label{eq:long}\end{align}
where $(*)$ comes from the definition of conditioned probability;
and $(**)$ comes from the estimates (\ref{eq:JJc}) for $\mathcal{J}$
and (\ref{eq:KKClam}) for $\mathcal{K}$, from the assumption (\ref{eq:assmpni})
and from bounding the denominator by $1$. Picking $\lambda$ sufficiently
large and $\nu$ sufficiently small we get that the result of the
computation is positive and dependant on the isotropic structure constants
of $G$ only. 

Finally, notice that if $\mathcal{K}$ did not occur, i.e.~if $R^{1}[T_{w,\epsilon r/4}^{1},T_{v,4r}^{1}]$
doesn't return to the ball $B(w,\rho)$ then \[
\cut(R^{1}[T_{w,\rho}^{1},T_{v,4r}^{1}];T_{w,\epsilon r}^{1})\cap A\subset\cut(R^{1}[0,T_{v,4r}^{1}];T_{w,\epsilon r}^{1})\]
 so (\ref{eq:long}) gives us (\ref{eq:IGamIC}) with an appropriate
choice of $q$. Hence we need only justify (\ref{eq:assmpni}).

To see (\ref{eq:assmpni}) we use Harnack's inequality (lemma \ref{lem:Harnack-general})
for the family of domains (see figure \ref{cap:DzEz})%
\begin{figure}
\input{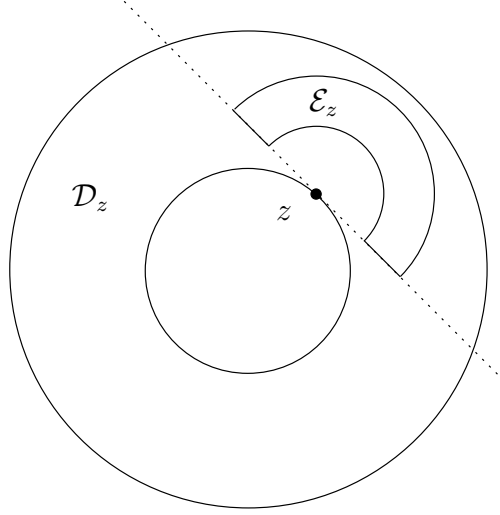}

\caption{\label{cap:DzEz}The domains $\mathcal{D}_{z}$ and $\mathcal{E}_{z}$.}
\end{figure}
 \begin{align*}
\mathcal{D}_{z} & =B(0,2)\setminus\overline{B(0,|z|)}\\
\mathcal{E}_{z} & =\left(B(z,2\mu)\setminus\overline{B(z,{\textstyle \frac{1}{2}}\mu)}\right)\cap\left\{ y:\langle y,z\rangle>|z|^{2}\right\} \end{align*}
defined for $z\in B(0,2-\frac{1}{2}\epsilon)\setminus B(0,1)$. The
function $f(x)=\mathbb{P}^{x}(R[0,T_{v,4r}]\cap\Gamma\neq\emptyset)$
is harmonic in $x$ in the domain $v+r\mathcal{D}_{z}$, $z=(w-v)/r$,
so we get for $r>C(\mu,G)$ that it is (up to constants depending
on $G,\mu$) independent of $x\in v+r\mathcal{E}_{z}$. Note that
for $r>C(\epsilon,G)$, $||z||\leq2-\epsilon-C(G)/r\leq2-\frac{1}{2}\epsilon$
and $\partial B(w,\rho)\cap H\subset v+r\mathcal{E}_{z}$.

We want to compare $f(x)$ with $f(v)$. We use (clause \ref{enu:lemDD-K}
of) lemma \ref{lem:DD} with $\mathcal{D}=B(0,2)$, $\mathcal{S}=B(0,\frac{1}{2})$
and $\mathcal{H}=\mathcal{E}_{z}$. We get\newc{c:RAHtmp}\newc{c:RAHtmp2}\begin{align}
f(v) & \stackrel{(*)}{\geq}\mathbb{P}^{v}(T(v+r\mathcal{E}_{z})<T_{v,4r})\mathbb{E}f(R(T(v+r\mathcal{E}_{z})))\geq\nonumber \\
 & \!\stackrel{(**)}{\geq}\! c(\mu,G)\mathbb{E}f(R(T(v+r\mathcal{E}_{z})))\geq\nonumber \\
 & \stackrel{(\dagger)}{\geq}c(\mu,G)f(x)\qquad\forall x\in\partial B(w,\rho)\cap H\label{eq:fvfy}\end{align}
where $(*)$ comes from the strong Markov property, $(**)$ comes
from lemma \ref{lem:DD} and $(\dagger)$ comes from Harnack's inequality.
This finishes the lemma --- choose $\delta(\epsilon,G)=\nu c(\mu,G)$
and the assumption (\ref{eq:IGamICassm}) will imply (\ref{eq:assmpni}).
\end{proof}

\subsection{Wiener's shell test}

For the next lemma we need to introduce a few notations. Let $v\in G$
and $r>1$ be some number and let \begin{equation}
A_{1}:=B(v,2r)\setminus\overline{B(v,r)}\quad A_{2}:=B(v,4r)\setminus\overline{B(v,{\textstyle \frac{1}{2}}r)}.\label{eq:defA12}\end{equation}
The notation $\overline{B(v,r)}$ relates here to closure in $G$,
not in $\mathbb{R}^{3}$. We shall denote $\partial A_{1}=\partial B(v,2r)\cup\partial B(v,r)$
(the general definition of boundary in $G$ might be a little smaller).
Both conventions apply to any annulus in this section.

If $\gamma_{1},\dotsc,\gamma_{n}$ are discontinuous paths (usually
we will consider $2$ or $3$) then we will consider $\LE(\gamma_{1}\cup\dotsb\cup\gamma_{n})$
as composed of $n$ pieces, each one {}``coming from some $\gamma_{i}$'',
and will denote them by $\LE_{i}(\gamma_{1}\cup\dotsb\cup\gamma_{n})$.
Formally, denote $t_{0}=1$, $t_{n}=\len\LE(\gamma_{1}\cup\dotsb\cup\gamma_{n})$
and $t_{i}$, $i=1,\dotsc,n-1$ to be the first $t$ such that $j_{t}>\len(\gamma_{1})+\dotsb+\len(\gamma_{i-1})$
($j_{t}$ from the definition of $\LE$, (\ref{eq:defLE})). Then\[
\LE_{i}(\gamma_{1}\cup\dotsb\cup\gamma_{n}):=\begin{cases}
\LE(\gamma_{1}\cup\dotsb\cup\gamma_{n})[t_{i-1},t_{i}-1] & t_{i}>t_{i-1}\\
\emptyset & \textrm{otherwise.}\end{cases}\]
Note that $\LE_{i}$ are simple and disjoint and that $\LE(\gamma_{1}\cup\dotsb\cup\gamma_{n})=\LE_{1}\cup\dotsb\cup\LE_{n}$.

\begin{lem}
\label{lem:L1L2}Let $G$ be a three dimensional isotropic graph and
let $\epsilon>0$ and $\eta>0$. Then there exists a $\delta=\delta(\epsilon,\eta,G)$
such that the following holds: Let $v\in G$, $r>C(\epsilon,\eta,G)$
and $s\in[r,2r-\eta r]$. Let $A_{1}$ and $A_{2}$ be as in (\ref{eq:defA12}).
Let $\gamma\subset A_{2}$ be a discontinuous path starting from $\partial B(v,\frac{1}{2}r)$
and let $R^{1}$ be a random walk starting from some point in $\partial A_{1}$
and stopped at $\partial A_{2}$. Let\[
L=\LE(\gamma\cup R^{1}),\quad L_{1}=\LE_{1}(\gamma\cup R^{1}),\quad L_{2}=\LE_{2}(\gamma\cup R^{1}).\]
Let $\mathcal{X}$ be the event that the following three events hold:
\begin{enumerate}
\item \label{enu:XX1}$L_{1}\subset B(v,s)$;
\item \label{enu:XX2}$L_{2}\not\subset B(v,s+\eta r)$;
\item \label{enu:XX3}$\mathbb{P}^{2,v}(R^{2}[0,T_{v,4r}^{2}]\cap L'\neq\emptyset)<\delta$
where $L'=L[1,t]$ and $t$ is the first time when $L$ hits $\partial B(v,s+\eta r)$
i.e.~$t:=\min\{ i:L(i)\not\in B(v,s+\eta r)\}$.
\end{enumerate}
Then\[
\mathbb{P}(\mathcal{X})<\epsilon.\]

\end{lem}
In words, the probability that $R^{1}$ extends $\gamma$ even by
a little ($\eta$) without being \linebreak[4]($\delta$-)hittable,
is small. To get a clearer geometric picture, think about $\gamma$
as the restriction of a continuous path to $A_{2}$, i.e.~as a sequence
of paths coming in and (except the last one) ending in $\partial A_{2}$;
and think about $R^{1}$ as starting from the end of $\gamma$.

We remark that $\gamma$ may be empty, in which case \ref{enu:XX1}
always holds, and in \ref{enu:XX2} one may replace $L_{2}$ with
$L$.

\begin{proof}
Let $q=q(G)$ be from lemma \ref{lem:IGamIC}. Let $K=K(\epsilon,G)>2$
be an integer such that \[
(1-q)^{K}<\epsilon.\]
 Now fix the parameter $\epsilon$ of lemma \ref{lem:IGamIC} to be
$\eta/2K$ and denote by $\lambda=\lambda(\epsilon,G)$ the result.

\noindent \begin{center}\begin{tabular}{|c|c|c|c|}
\hline 
Lemma \ref{lem:IGamIC}&
$\epsilon$&
$q$&
$\delta$\tabularnewline
\hline 
here&
$\eta/2K$&
$q$&
$\lambda$\tabularnewline
\hline
\end{tabular}\end{center}

Define $s_{i}=s+r\eta i/(K+2)$ for $i=1,\dots,K+1$. Let $j_{k}$
be as in the definition of loop-erasure (\ref{eq:defLE}) so that
$R^{1}[j_{k}+1,T^{1}(\partial A_{2})]$ is a random walk conditioned
not to hit $\LE(\gamma\cup R^{1}[0,j_{k}])$. Define \[
\tau_{i}':=\max\{ j_{k}\leq T^{1}(\partial A_{2}):\LE(\gamma\cup R^{1}[0,j_{k}])\subset B(v,s_{i})\}.\]
As in lemma \ref{lem:IGamIC}, let $j_{k}$ be admissible if there
exists a path $\delta\subset\overline{B(R^{1}(j_{k}),16C_{\ref{C:arc}})}$
from $R^{1}(j_{k})$ to $\partial B(v,|v-R^{1}(j_{k})|+2C_{\ref{C:arc}})$
which does not intersect $\LE(\gamma\cup R^{1}[0,j_{k}])$ and define
\[
\tau_{i}:=\min\{ j_{k}\geq\tau_{i}':j_{k}\textrm{ admissable}\}.\]
We note that if $L_{1}\subset B(v,s)$ and $L_{2}\not\subset B(v,s_{i}+C_{\ref{C:arc}})$
then $\tau_{i}$ is well defined. Indeed, let $t>\tau_{i}'$ ($\tau_{i}'$
is obviously well defined) be the first $j_{k}$ such that $R^{1}(t)\not\in B(v,s_{i}+C_{\ref{C:arc}})$,
and denote $s^{*}=|R^{1}(t)-v|-C_{\ref{C:arc}}$ so that $s^{*}>s_{i}$.
We use lemma \ref{lem:arc} with the radius being $8C_{\ref{C:arc}}$
and with $H$ being the half space tangent to $B(v,s^{*})$ and orthogonal
to the segment $[v,R^{1}(t)]$. We get that there exists a path $\delta'\subset B(R^{1}(t),8C_{\ref{C:arc}})\setminus H$
from $R^{1}(t)$ to some $x\in\partial B(R^{1}(t),8C_{\ref{C:arc}})\cap\{ x:d(x,\partial_{\textrm{cont}}H)>4C_{\ref{C:arc}}\}$.
If $r$ is sufficiently large then this implies $d(x,v)>s^{*}+3C_{\ref{C:arc}}$.
Let\[
q_{0}:=\max\{ q:\delta'(q)\in\LE(\gamma\cup R^{1}[0,t])\}\]
and let $t_{0}$ be the $j_{k}$ such that $R^{1}(t_{0})=\delta'(q_{0})$.
Clearly $\delta'(q_{0})\not\in B(v,s_{i})$ so $t_{0}>\tau_{i}'$.
Further, $\delta:=\delta'[q_{0},\len\delta']$ is a path from $R^{1}(t_{0})$
to $x$ not intersecting $\LE(\gamma\cup R^{1}[0,t_{0}])$. Since
$R^{1}(t_{0})\in B(v,s^{*}+C_{\ref{C:arc}})$ and $x\not\in B(v,s^{*}+3C_{\ref{C:arc}})$
and since $B(R^{1}(t),8C_{\ref{C:arc}})\subset B(R^{1}(t_{0}),16C_{\ref{C:arc}})$
we see that all the admissibility requirements are satisfied. Therefore
$t_{0}$ is admissible and hence $\tau_{0}$ is well defined and $\leq t_{0}$.

Fix some $i$ and assume that $\tau_{i}$ is well defined and that
$\tau_{i}<T^{1}(\partial A_{2})$. Examine $R^{1}$ after $\tau_{i}$.
The definition of $\tau_{i}$ considers only $\LE(R[0,\tau_{i}])$
therefore it only affects $R^{1}$ after $\tau_{i}$ by conditioning
it to not intersect $\Gamma_{i}:=\LE(\gamma\cup R^{1}[0,\tau_{i}]).$
Hence lemma \ref{lem:IGamIC} applies. We get\begin{equation}
\mathbb{P}(\mathcal{B}_{i}\,|\,\tau_{i},R^{1}[0,\tau_{i}])<1-q\label{eq:expcond}\end{equation}
where $\mathcal{B}_{i}$ is the event that\renewcommand{\theenumi}{(\alph{enumi})}
\begin{enumerate}
\item \label{enu:BB1}$\tau_{i}$ is well defined,
\item \label{enu:BB2}$\mathbb{P}^{2,v}(R^{2}[0,T_{v,4r}^{2}]\cap\Gamma_{i}\neq\emptyset)\leq\lambda$,
\item \label{enu:BB3}For some $y\in B(w_{i},\rho)$, $\rho:=\eta r/2K$,
\begin{equation}
\mathbb{P}^{2,y}(\cut(R^{1}[\tau_{i},T_{v,4r}^{1}];T_{i}^{*})\cap R^{2}[0,T_{v,4r}^{2}]\neq\emptyset)\leq q.\label{eq:tilTv4r}\end{equation}
where $T_{i}^{*}:=\min\{ t>\tau_{i}:R^{1}(t)\in\partial B(w_{i},\rho)\}$.
\end{enumerate}
\renewcommand{\theenumi}{(\roman{enumi})}The notation in (\ref{eq:expcond})
might deserve some explanation: we are conditioning here on the fact
that $\tau_{i}$ is well defined (so \ref{enu:BB1} is satisfied automatically),
on its value which gives some information on $R$ beyond $\tau_{i}$
and on the entire path from $0$ to $\tau_{i}$ (which gives $\Gamma_{i}$
and $w_{i}$).

It will be convenient to replace \ref{enu:BB1} with the stronger
\begin{enumerate}
\item [(a')]$\tau_{i+1}$is well defined
\end{enumerate}
and then replace (\ref{eq:tilTv4r}) with a slightly stronger condition:\begin{equation}
\mathbb{P}^{2,y}(\cut(R^{1}[\tau_{i},\tau_{i+1}];T_{i}^{*})\cap R^{2}[0,T_{v,4r}^{2}]\neq\emptyset)\leq q.\label{eq:TilTauip1}\end{equation}
This is possible since $T_{i}^{*}<\tau_{i+1}$ by our choice of $\rho$.
Denote $\mathcal{B}_{i}$ with \ref{enu:BB1} and (\ref{eq:tilTv4r})
replaced with (a') and (\ref{eq:TilTauip1}) by $\mathcal{B}_{i}'$
and get $\mathbb{P}(\mathcal{B}_{i}'\,|\,\tau_{i},R^{1}[0,\tau_{i}])<1-q$.
$\mathcal{B}_{i}'$ depends only on $\tau_{i+1}$ and $R[0,\tau_{i+1}]$
and hence we can write\[
\mathbb{P}(\mathcal{B}_{i}'\,|\,\mathcal{B}'_{0},\dotsc,\mathcal{B}_{i-1}')=\mathbb{EP}(\mathcal{B}_{i}'\,|\,\tau_{1},\dotsc,\tau_{i},R[0,\tau_{i}])=\mathbb{EP}(\mathcal{B}_{i}'\,|\,\tau_{i},R[0,\tau_{i}])<1-q\]
(the $\mathbb{E}$ signs above stand for conditional expectation with
respect to $\mathcal{B}_{i-1}',\dotsc,\mathcal{B}_{0}'$). Hence\[
\mathbb{P}\Big(\bigcap_{i=0}^{K-1}\mathcal{B}_{i}'\Big)<(1-q)^{K}<\epsilon.\]

The lemma will be finished when we show that for an appropriate choice
of $\delta$ we have $\mathcal{X}\subset\cap\mathcal{B}_{i}'$. As
explained above, conditions \ref{enu:XX1} and \ref{enu:XX2} in the
definition of $\mathcal{X}$, show that all $\tau_{i}$ are well defined.
Hence condition (a') in the definition of $\mathcal{B}_{i}'$ is satisfied
for all $i$. Setting $\delta<\lambda$ will ensure condition \ref{enu:BB2}
for all $i$, since $\Gamma_{i}\subset L'$ for all $i$. Finally,
lemma \ref{lem:DD} shows that for some $\nu(\epsilon,\eta,G)$, \[
\mathbb{P}^{v}(T(B(w_{i},\rho))<T_{v,4r})>\nu\quad\forall i\forall w_{i}\in\partial B(v,s_{i})\]
which shows that setting $\delta<q\nu$ ensures condition \ref{enu:BB3}
(we used here (\ref{eq:TilTauip1}) and $\cut(R^{1}[\tau_{i}+1,\tau_{i+1}];T_{i}^{*})\subset L'$).
Hence $\mathcal{X}\subset\cap\mathcal{B}_{i}'$ and the lemma is proved.
\end{proof}
\begin{lem}
\label{lem:wienersimp}Let $G$ be a three dimensional isotropic graph,
and let $K$ be some number. Then there exists a $\delta=\delta(K,G)>0$
such that the following holds. Let $v\in G$ and let $R^{1}$ be a
random walk starting from $v$. Let $r>s>1$. Let $\mathcal{B}$ be
the event that there exists some $T\geq0$ such that
\begin{enumerate}
\item $\LE(R^{1})\not\subset B(v,r)$ where $R^{1}:=R^{1}[0,T]$; and
\item $\mathbb{P}^{2,v}(R^{2}[0,T_{v,r}^{2}]\cap\LE(R^{1})\cap(B(v,r)\setminus\overline{B(v,s)})=\emptyset)>(s/r)^{\delta}$.
\end{enumerate}
Then $\mathbb{P}(\mathcal{B})<C(K,G)(s/r)^{K}$.
\end{lem}
It is not difficult to see that the probability that $\LE(R^{1})\subset B(v,s)$
is $\approx s/r$ (for $T\approx r^{2}$ which is the range that interests
us). Therefore in fact $\LE(R)$ has a rather big probability to be
{}``unhittable'' because it is small. The point about the lemma
is that if it isn't small, this probability is negligible.

\begin{proof}
We may assume w.l.o.g.~that both $r/s$ and $s$ are sufficiently
large (and the bound may depend on $K$). A Brownian motion starting
from a point in $\partial B(0,1)$ has probability $\frac{2}{3}$
to reach $\partial B(0,2)$ before $\partial B(0,\frac{1}{2})$. Hence
by lemma \ref{lem:sball}, if $\rho>C(G)$ then random walk starting
from $x\in\partial B(v,\rho)$ has probability $>\frac{7}{12}$ to
reach $\partial B(v,2\rho)$ before $\partial B(v,\frac{1}{2}\rho)$.
Hence if we define $\rho_{i}=s2^{i}$ and stopping times $\tau_{0}'=0$,\[
\tau_{j}'=\min\{ t>\tau_{j-1}':R^{1}(t)\in\partial B(v,\rho_{i'(j)})\textrm{ for }i'(j)\textrm{ s.t. }R^{1}(\tau_{j-1}')\not\in\partial B(v,\rho_{i'(j)})\}\]
then the process $i'(j)$ dominates a random walk on $\mathbb{N}$
with a drift to infinity. In particular, it follows that if we denote
\begin{equation}
n_{i}':=\#\{ j:i'(j)=i\}\quad N':=\sum_{i=0}^{I}n_{i}'\label{eq:defniNp}\end{equation}
then there exists a $\lambda=\lambda(K)$ such that\begin{equation}
\mathbb{P}(N'>\lambda I)\leq C2^{-2IK}.\label{eq:NlamK}\end{equation}
(\ref{eq:NlamK}) holds for any value of $I$, but we will define
$I:=\left\lfloor \log_{2}r/s\right\rfloor $ since we are interested
in what happens until $r$, and get $\mathbb{P}(N'>\lambda I)\leq C(s/r)^{K}$.

Next we take the $I$ annuli between $s$ and $r$ and denote the
even ones as {}``$A_{1}$-s'', i.e.~\[
A_{1,i}:=B(v,2\rho_{2i})\setminus\overline{B(v,\rho_{2i})}\quad A_{2,i}:=B(v,4\rho_{2i})\setminus\overline{B(v,{\textstyle \frac{1}{2}}\rho_{2i})}\]
and define stopping times $\tau_{j}$ {}``from one $A_{1}$ to the
next'', i.e. $\tau_{1}:=\tau_{1}'$ and\[
\tau_{j}=\min\{ t>\tau_{j-1}:R^{1}(t)\in\partial A_{1,i(j)}\textrm{ for some }i(j)\textrm{ such that }R^{1}(\tau_{j-1})\not\in\partial A_{1,i(j)}\}.\]
It is clear that the $\tau_{j}$ are a subsequence of the $\tau_{j}'$
hence if we define $n_{i}$ and $N=\sum_{i=0}^{\left\lfloor (I-1)/2\right\rfloor }n_{i}$
analogously to (\ref{eq:defniNp}) then the analog of (\ref{eq:NlamK})
will also hold. Define now $M:=\left\lceil 6\lambda\right\rceil $
and get that \[
\mathbb{P}(\#\{ i<{\textstyle \frac{1}{2}}I:n_{i}>M\}>{\textstyle \frac{1}{6}}I)\leq C(s/r)^{K}.\]

Next use lemma \ref{lem:L1L2} with \[
\eta=\frac{1}{M}\quad\epsilon=\frac{\mu}{M^{2}}\]
where $\mu=\mu(K)$ is some parameter which will be fixed later. Call
the $\delta$ of lemma \ref{lem:L1L2} $\nu$. Fix one $j$ and some
$m\in\{0,\dotsc,M-1\}$. Denote $r_{j}=\rho_{i(j)}$ and $s_{j,m}=r_{j}(1+m/M)$.
Define $\gamma_{i,j}:=\LE(R^{1}[0,\tau_{j}])\cap A_{2,i}$ and $\gamma_{j}:=\gamma_{i(j),j}$.
Define $\mathcal{X}_{j,m}$ the events that the following conditions
are all fulfilled:
\begin{enumerate}
\item \label{enu:YYj2}$L_{1,j}\subset B(v,s_{j,m})$ where\[
L_{k,j}:=\LE_{k}(\gamma_{j}\cup R^{1}[\tau_{j},\tau_{j+1}])\quad k=1,2;\]

\item \label{enu:YYj3}$L_{2,j}\not\subset B(v,s_{j,m+1})$;
\item \label{enu:YYj4}$\mathbb{P}^{2,v}(R^{2}[0,T_{v,4r_{j}}^{2}]\cap L_{j,m}'\neq\emptyset)<\nu$
where $L_{j,m}'=L_{1,j}\cup L_{2,j}[1,t]$ and $t$ is the first time
when $L_{2,j}$ hits $\partial B(v,s_{j,m+1})$ i.e.~$t:=\min\{ i:L_{2,j}(i)\not\in B(v,s_{j,m+1})\}$.
\end{enumerate}
This completes all the parameters of lemma \ref{lem:L1L2} (see table
\ref{cap:lemYYjm}) %
\begin{table}
\begin{tabular}{|c|c|c|c|c|c|c|}
\hline 
Lemma \ref{lem:L1L2}&
$\epsilon$&
$\eta$&
$\delta$&
$r$&
$s$&
$A_{1}$\tabularnewline
\hline
here&
$\frac{\mu}{M^{2}}$&
$\frac{1}{M}$&
$\nu$&
$r_{j}$&
$s_{j,m}$&
$A_{1,i(j)}$\tabularnewline
\hline
Lemma \ref{lem:L1L2}&
$A_{2}$&
$\gamma$&
$R^{1}$&
$L_{k}$&
$L'$&
\emph{$\mathcal{X}$}\tabularnewline
\hline
here&
$A_{2,i(j)}$&
$\gamma_{j}$&
$R^{1}[\tau_{j},\tau_{j+1}]$&
$L_{k,j}$&
$L_{j,m}'$&
$\mathcal{X}_{j,m}$\tabularnewline
\hline
\end{tabular}

\caption{\label{cap:lemYYjm}Use of lemma \ref{lem:L1L2} in lemma \ref{lem:wienersimp}.}
\end{table}
 and we may use it to get $\mathbb{P}(\mathcal{X}_{j,m}\,|\, R[0,\tau_{j}])<\mu/M^{2}$
and hence\[
\mathbb{P}\Big(\bigcup_{m}\mathcal{X}_{j,m}\,\Big|\, R[0,\tau_{j}]\Big)<\frac{\mu}{M}.\]
Denote this event by $\mathcal{X}_{j}$. The $\mathcal{X}_{j}$-s
are therefore dominated by a sequence of independent random variables
with probability $\mu/M$. A simple (and standard) calculation now
shows that if $\mu$ is taken sufficiently small one has \[
\mathbb{P}\Big(\#\{ j\leq\lambda I:\mathcal{X}_{j}\}>\frac{\lambda I}{M}\Big)\leq C(G)(s/r)^{K}.\]
Fix $\mu$ to satisfy this condition. Thus we may define $\mathcal{G}$
as the event that $N\leq\lambda I$ and $\mathcal{X}_{j}$ occurred
less than $\frac{\lambda}{M}I\leq\frac{1}{6}I$ times before $\lambda I$.
Our calculations show that $\mathbb{P}(\neg\mathcal{G})\leq C(s/r)^{K}$.
The lemma will be finished once we show that with an appropriate choice
of $\delta$ one has $\mathcal{B}\cap\mathcal{G}=\emptyset$.

Assume therefore from now on that both $\mathcal{B}$ and $\mathcal{G}$
occurred. There are $\left\lceil \frac{1}{2}I\right\rceil $ different
$A_{1,i}$-s, and under the assumption $\mathcal{G}$, in $\leq\frac{1}{6}I$
we have $n_{i}>M$ and in $\leq\frac{1}{6}I$ of them we have that
some $\mathcal{X}_{j}$ occurred when $R(\tau_{j})\in\partial A_{1,i}$.
Hence we get that at least $\frac{1}{6}I$ are {}``good'' in the
sense that they are visited $\leq M$ times and are {}``$\mathcal{X}_{j}$-free''.
Fix $i$ to be one such good index. We will now show that $A_{1,i}\cap\LE(R^{1})$
is {}``hittable'', and we shall show that by induction.

\begin{sublem}For any $j$ define \[
n_{i,j}:=\#\{ k<j:R^{1}(\tau_{k})\in A_{1,i}\},\]
$u_{i,j}:=\rho_{2i}(1+n_{i,j}/M)$ and $U_{i,j}:=B(v,u_{i,j})\setminus\overline{B(v,\rho_{2i})}$.
Let also $\gamma_{i,j}^{*}$ be $\gamma_{i,j}[0,t_{i,j}]$ where $t_{i,j}$
is the first time when $\gamma_{i,j}(t)\in\partial B(v,u_{i,j})$,
or empty if $\gamma_{i,j}\cap\partial B(v,u_{i,j})=\emptyset$. Then
for all $j$, either $\gamma_{i,j}^{*}=\emptyset$ or \[
\mathbb{P}^{2,w}\left(\gamma_{i,j}^{*}\cap U_{i,j}\cap R^{2}[0,T_{4\rho_{2i}}^{2}]=\emptyset\right)\geq\nu.\]

\end{sublem}

Note that $n_{i}\leq M$ (because $i$ is a good index) and hence
$u_{i,j}\leq2\rho_{2i}$ throughout the induction.
\begin{proof}
[Subproof]We use induction over $j$. If $R^{1}(\tau_{j})\not\in\partial A_{1,i}$,
then $R^{1}[\tau_{j},\tau_{j+1}]$ does not enter $A_{1,i}$ and hence
can only affect $\gamma_{i,j}^{*}\cap A_{1,i}$ by removing some components
from its end. In this case it must remove the last component of $\gamma_{i,j}^{*}\cap A_{1,i}$,
which is the component intersecting $\partial B(w,u_{i,j})$, and
all the components of $\gamma_{i,j}$ after $\gamma_{i,j}^{*}$. Hence
we get that $\gamma_{i,j+1}\cap\partial B(w,u_{i,j})=\emptyset$ (here
$u_{i,j}=u_{i,j+1}$) and the induction holds in this case. Therefore
we need only be interested in the case $R^{1}(\tau_{j})\in\partial A_{1,i}$.
We know $\mathcal{X}_{j}$ did not occur, and in particular $\mathcal{X}_{j,n_{i,j}}$
did not occur. Hence one of the three constituents of $\mathcal{X}_{j,n_{i,j}}$
must have failed. Let us review them in order.
\begin{enumerate}
\item If $L_{1,j}\not\subset B(v,u_{i,j})$ (note that $s_{j,n_{i,j}}=u_{i,j}$)
then

\begin{itemize}
\item $\gamma_{i,j}\cap\partial B(v,u_{i,j})\neq\emptyset$ and by the induction
hypothesis, $\gamma_{i,j}^{*}$ is hittable.
\item $R^{1}$ must have hit $\gamma_{i,j}$, if at all, after exiting $B(v,u_{i,j})$.
Therefore, if $\gamma_{i,j+1}^{*}\neq\emptyset$ then it contains
$\gamma_{i,j}^{*}$ so $\gamma_{i,j+1}^{*}\cap U_{i,j+1}$ is hittable
and the induction holds.
\end{itemize}
\item If $L_{1,j}\cup L_{2,j}\subset B(v,u_{i,j+1})$ then $\gamma_{i,j+1}\cap\partial B(v,u_{i,j+1})=\emptyset$
and the induction holds.
\item Finally, if \ref{enu:YYj2}-\ref{enu:YYj3} of the definition of $\mathcal{X}_{j}$
happened then we must have $\mathbb{P}^{2,w}(R^{2}[0,T_{w,4\rho_{2i}}^{2}]\cap L_{j,n_{i,j}}'\neq\emptyset)\geq\nu$
and then the induction holds because $L_{j,n_{i,j}}'=\gamma_{i,j+1}^{*}$.
\end{enumerate}
Hence the sublemma is proved.
\end{proof}
Assume now that $\mathcal{B}$ happened and let $T$ be the {}``bad''
time. Let $J$ be such that $\tau_{J}\leq T<\tau_{J+1}$. As in the
sublemma, the part of the walk on $[\tau_{j},T]$ can affect in an
adverse way only the part of the walk inside $A_{1,i(J)}$. Hence
the sublemma shows that for any good $i\neq i(J)$, either $\LE(R^{1})$
does not contain a crossing of $A_{1,i}$, or this crossing is hittable.
However, if $\mathcal{B}$ happened then $\LE(R^{1})$ crosses $B(v,r)\setminus\overline{B(v,s)}$
and hence all $A_{1,i}$-s so we get for all good $i\neq i(J)$ that\[
\mathbb{P}^{2,v}\left(\LE(R^{1})\cap A_{1,i}\cap R^{2}[0,T_{v,4\rho_{2i}}^{2}]\neq\emptyset\right)\geq\nu.\]
Harnack's inequality (lemma \ref{lem:Harnack-general}) shows that
\begin{equation}
\mathbb{P}^{2,w}\left(\LE(R^{1})\cap A_{1,i}\cap R^{2}[0,T_{v,4\rho_{2i}}^{2}]\neq\emptyset\right)>c\nu\quad\forall w\in B(v,{\textstyle \frac{1}{2}}\rho_{2i}).\label{eq:GBLARcnu}\end{equation}
 Hence let $i_{1}<i_{2}<\dotsc<i_{n}$ be good indices with $i_{k+1}-i_{k}\geq2$,
$i_{k}\neq I,i(J)$. We can take $n:=\left\lfloor I/12\right\rfloor -2$.
Then\begin{eqnarray*}
\lefteqn{\mathbb{P}^{2,v}\left(R^{2}[0,T_{v,r}^{2}]\cap\LE(R^{1})\cap\left(B(v,r)\setminus\overline{B(v,s)}\right)=\emptyset\right)\leq}\\
 &  & \leq\prod_{k=1}^{n-1}\mathbb{P}^{2,v}(R^{2}[T_{v,\frac{1}{2}\rho_{2i_{k}}}^{2},T_{v,4\rho_{2i_{k}}}^{2}]\cap\LE(R^{1})=\emptyset\,|\\
 &  & \qquad\qquad\qquad\qquad|\, R^{2}[0,T_{v,4\rho_{2i_{k-1}}}^{2}]\cap\LE(R^{1})=\emptyset)\leq\\
 &  & \leq\prod_{k=1}^{n-1}\max_{w\in B(v,\frac{1}{2}\rho_{2i_{k}})}\mathbb{P}^{2,w}(R^{2}[0,T_{v,4\rho_{2i_{k}}}^{2}]\cap\LE(R^{1})=\emptyset)<\\
 &  & \!\stackrel{(\ref{eq:GBLARcnu})}{<}\!\prod_{k=1}^{n-1}(1-c\nu)\leq C(1-c\nu)^{I/12}\leq C(s/r)^{(\log(1-c\nu))/12\log2}.\end{eqnarray*}
Therefore taking $\delta=(\log(1-c\nu))/24\log2$ and $r/s$ sufficiently
large we get a contradiction, so $\mathcal{G}\cap\mathcal{B}=\emptyset$,
and the lemma is proved.
\end{proof}

\subsection{Hittable sets}

The next step (lemmas \ref{lem:outgoing} and \ref{lem:incoming})
is to show that $\LE(R)$ is {}``hittable'' in a rather strong sense. 

\begin{defn*}
Let $\lambda>0$ and $1\leq\rho<\sigma$. We define $\mathcal{H}^{\textrm{out}}(\lambda,\rho,\sigma)$
to be the family of paths $\beta$ satisfying that for every $w$,
and every subpath $\gamma\subset\beta$ which is an outgoing crossing
of $B(w,\sigma)\setminus\overline{B(w,\rho)}$ i.e.~$\gamma(0)\in\partial B(w,\rho)$
and $\gamma(\len\gamma)\in\partial B(w,\sigma)$ one has\[
\mathbb{P}^{w}(R[0,T_{w,\sigma}]\cap\gamma)=\emptyset)\leq\lambda.\]
Define $\mathcal{H}^{\textrm{out}}(\lambda,\mu)=\cap_{\rho\geq1}\mathcal{H}^{\textrm{out}}(\lambda,\rho,\mu\rho)$.
\end{defn*}
Similarly we define $\mathcal{H}^{\textrm{in}}$ to be the set where
this holds for $\gamma$ incoming, i.e.~$\gamma(0)\in\partial B(w,\sigma)$
and $\gamma(\len\gamma)\in\partial B(w,\rho)$. We define $\mathcal{H}:=\mathcal{H}^{\textrm{out}}\cap\mathcal{H}^{\textrm{in}}$.
Note that these properties are hereditary: if $\beta\subset\gamma\in\mathcal{H}^{\textrm{out/in}}$
then $\beta\in\mathcal{H}^{\textrm{out/in}}$.

\begin{lem}
\label{lem:expescape}Let $G$ be a Euclidean net and let $v\in G$
and $r>1$. Then\[
\mathbb{P}^{v}(R[0,t]\subset B(v,r))\leq C(G)e^{-c(G)t/r^{2}}.\]

\end{lem}
\begin{proof}
This follows immediately by applying lemma \ref{lem:r2} repeatedly.
\end{proof}
\begin{lem}
\label{lem:outgoing}Let $G,\epsilon,v,r$ and $\mathcal{D}$ be as
in theorem \ref{thm:QL} and let $K>0$. Then there exists a $\delta=\delta(\epsilon,K,G)$
such that \[
\mathbb{P}^{1,v}(\exists T\leq T^{1}(\partial\mathcal{D}):\LE(R^{1}[0,T])\not\in\mathcal{H}^{\textrm{out}}(r^{-\delta},r^{\epsilon}))<C(\epsilon,K,G)r^{-K}.\]

\end{lem}
\begin{proof}
Clearly we may assume $r$ is sufficiently large (and the bound may
depend on $K$ and $\epsilon$). Since $\mathcal{H}^{\textrm{out}}(r^{-\delta},\rho,\sigma)$
is increasing in $\sigma$ and decreasing in $\rho$, it is enough
to show that for any fixed $\sigma=\frac{1}{2}\rho r^{\epsilon}$
we have \begin{equation}
\mathbb{P}^{1,v}(\exists T\leq T^{1}(\partial\mathcal{D}):\LE(R^{1}[0,T])\not\in\mathcal{H}^{\textrm{out}}(r^{-\delta},\rho,\sigma))\leq Cr^{-K-1}.\label{eq:fixrhosig}\end{equation}
Once (\ref{eq:fixrhosig}) is established, we can apply it for $\rho=1,2,4,\dotsc,\log r$
and bound the probability that (\ref{eq:fixrhosig}) holds for one
such $\rho$ by the sum, which is $\leq Cr^{-K-1}\log r\leq Cr^{-K}$
and the lemma would be proved.

Fix therefore one $\rho$ and $\sigma=\frac{1}{2}\rho r^{\epsilon}$.
Let $w\in\mathcal{D}$ be some point and let $A=B(w,\sigma)\setminus\overline{B(w,\rho)}$.
Let $x\in\overline{B(w,\rho)}$ be some point, and let $\rho':=\max2|x-w|,\rho_{0}$
for some constant $\rho_{0}(\epsilon,K,G)$ to be fixed later. Let
$A':=B(x,\frac{1}{2}\sigma)\setminus\overline{B(x,\rho')}$ so that
$A'\subset A$ for all $r$ sufficiently large. Let $t$ be some number.
Use lemma \ref{lem:wienersimp} with the parameters and notations
in table \ref{cap:lemWtW}. %
\begin{table}
\begin{tabular}{|c|c|c|c|c|c|c|c|}
\hline 
Lemma \ref{lem:wienersimp}&
$K$&
$\delta$&
$v$&
$T$&
$r$&
$s$&
\emph{$\mathcal{B}$}\tabularnewline
\hline 
Here&
$(K+10)/\epsilon$&
$\lambda$&
$x$&
$t$&
$\frac{1}{2}\sigma$&
$\rho'$&
$\mathcal{B}_{x}$\tabularnewline
\hline
\end{tabular}

\caption{\label{cap:lemWtW}Use of lemma \ref{lem:wienersimp} in lemma \ref{lem:outgoing}.}
\end{table}
We get for the event $\mathcal{B}_{x}$ that there exists some $t$
such that
\begin{enumerate}
\item $\LE(R^{1}[0,t])\not\subset B(x,\frac{1}{2}\sigma)$;
\item \label{enu:BBx2}$\mathbb{P}^{2,x}(R^{2}[0,T_{x,\sigma/2}^{2}]\cap\LE(R^{1}[0,t])\cap A'=\emptyset)>C(\epsilon,K,G)r^{-\lambda\epsilon}$
\end{enumerate}
that $\mathbb{P}^{1,x}(\mathcal{B}_{x})<C(\epsilon,K,G)r^{-K-10}$.
If $\rho_{0}$ is sufficiently large we can use Harnack's inequality
to change in \ref{enu:BBx2} the starting point of $R^{2}$ from $x$
to $w$ and pay only by increasing the constant.

Returning to a random walk starting from $v$, we define an event
$\mathcal{E}=\mathcal{E}(t_{1},t_{2},x)$ by
\begin{enumerate}
\item \label{enu:EExt1t2}$R^{1}(t_{1})=x$;
\item \label{enu:EExt1t22}$\LE(R^{1}[0,t_{1}])\cap R^{1}[t_{1}+1,t_{2}]=\emptyset$;
\item $\LE(R^{1}[t_{1}+1,t_{2}])\not\subset B(x,\sigma)$;
\item $\mathbb{P}^{2,w}(R^{2}[0,T_{w,\sigma}^{2}]\cap\LE(R^{1}[t_{1}+1,t_{2}])\cap A=\emptyset)>c(\epsilon,K,G)r^{-\lambda\epsilon}$.
\end{enumerate}
and get $\mathbb{P}^{1,v}(\bigcup_{t_{2}\geq t_{1}}\mathcal{E})\leq C(\epsilon,K,G)r^{-K-10}$.
Summing we get for any $U$,\begin{equation}
\mathbb{P}\Big(\bigcup_{(*)}\mathcal{E}(t_{1},t_{2},x)\Big)\leq Cr^{-K-4}U\label{eq:Et1t2x}\end{equation}
 where the union is over all $w\in\mathcal{D}$, $x\in\overline{B(w,\rho)}$,
$t_{1}\leq U$ and $t_{2}\geq t_{1}$. Now, lemma \ref{lem:expescape}
shows that \begin{align*}
\mathbb{P}(T(\partial\mathcal{D})>\mu r^{2}\log r) & \leq\mathbb{P}(T_{v,r}>\mu r^{2}\log r)\leq C(G)e^{-c(G)\mu\log r}\leq\\
\intertext{\textrm{and for $\mu=\mu(K,G)$ sufficiently large}} & \leq Cr^{-K-1}.\end{align*}
Therefore using $U=\mu r^{2}\log r$ in (\ref{eq:Et1t2x}) gives\[
\mathbb{P}\Big(\bigcup_{(*)}\mathcal{E}(t_{1},t_{2},x)\Big)\leq Cr^{-K-1}\]
where here the union is over all $w,x$ and $t_{1}\leq t_{2}\leq T^{1}(\partial\mathcal{D})$.
This finishes the lemma by taking $\delta=\epsilon\lambda/2$ and
$r$ sufficiently large, since if $\LE(R[0,T])\not\in\mathcal{H}^{\textrm{out}}$
occurred then one of the $\mathcal{E}(t_{1},t_{2},x)$ must have occurred.
Namely, let $\gamma$ be a subpath of $\LE(R[0,T])$ satisfying the
requirements from the definition of $\mathcal{H}^{\textrm{out}}$
and minimal (basically this means $\gamma\subset\overline{A}$) and
assume $\gamma=\LE(R[0,T])[l,m]$. Then we may take $t_{1}:=j_{l}$
($j$ from the definition of $\LE$, (\ref{eq:defLE})), $t_{2}:=j_{m}$
and $x=R^{1}(t_{1})$, and directly from the definitions, $\mathcal{E}(t_{1},t_{2},x)$
will be satisfied. Therefore (\ref{eq:fixrhosig}) holds for any appropriate
$\rho$ and $\sigma$ and the lemma is proved.
\end{proof}
The following lemma will be used only formally --- in effect the previous
lemma is enough. We include it here mainly for completeness.

\begin{lem}
\label{lem:incoming}Lemma \ref{lem:outgoing} holds with $\mathcal{H}^{\textrm{out}}$
replaced by $\mathcal{H}^{\textrm{in}}$.
\end{lem}
\begin{proof}
Let $x\in\overline{\mathcal{D}}$ satisfy that \begin{equation}
\mathbb{P}^{v}(R(T(\partial\mathcal{D}\cup\{ x\}))=x)\geq r^{-K-3}.\label{eq:xgood}\end{equation}
 The transience of $G$ shows that (compare to (\ref{eq:noWW2}))
that\[
\mathbb{P}^{v}(R(T(\partial\mathcal{D}\cup\{ v,x\}))=x)\geq c(G)r^{-K-3}.\]
Using the symmetry of random walk in the form (\ref{eq:symT}) we
get\begin{equation}
\mathbb{P}^{x}(R(T(\partial\mathcal{D}\cup\{ v,x\}))=v)\geq c(G)r^{-K-3}.\label{eq:xgoodx}\end{equation}
Next use lemma \ref{lem:outgoing} with the parameters in table \ref{cap:outgoing}
(if $x\in\partial\mathcal{D}$ then use the lemma with each neighbor
$x'$ of $x$ in $\mathcal{D}$ serving as the starting point instead
of $x$). %
\begin{table}
\begin{tabular}{|c|c|c|c|c|}
\hline 
Lemma \ref{lem:outgoing}&
$v$&
$\mathcal{D}$&
$K$&
$\delta$\tabularnewline
\hline 
Here&
$x$&
$\mathcal{D}$&
$2K+7$&
$\delta$\tabularnewline
\hline
\end{tabular}

\caption{\label{cap:outgoing}Use of lemma \ref{lem:outgoing} in lemma \ref{lem:incoming}.}
\end{table}
 For any time or stopping time $t$, define an event\begin{equation}
\mathcal{O}(t)=\{\LE(R[0,t])\in\mathcal{H}^{\textrm{out}}(r^{-\delta},r^{\epsilon})\}.\label{eq:defOt}\end{equation}
With this notation, the conclusion of lemma \ref{lem:outgoing} is
that\[
\mathbb{P}^{x}(\exists t\leq T(\partial\mathcal{D}):\mathcal{O}(t))\leq C(K,G)r^{-2K-7}.\]
In particular, $\mathbb{P}^{x}(\mathcal{O}(T(\partial\mathcal{D}\cup\{ v,x\})))\leq Cr^{-2K-7}$.
Hence\[
\mathbb{P}^{x}(\mathcal{O}(T(v))\,|\, R(T(\partial\mathcal{D}\cup\{ v,x\}))=v)\leq\frac{\mathbb{P}^{x}(\mathcal{O}(T(\partial\mathcal{D}\cup\{ v,x\})))}{\mathbb{P}^{x}(R(T(\partial\mathcal{D}\cup\{ v,x\})=v)}\stackrel{(\ref{eq:xgoodx})}{\leq}Cr^{-K-4}.\]
Now, loop-erased random walk conditioned on the end vertex is symmetric
(see e.g.~\cite[lemma 2]{K}) so,\[
\mathbb{P}^{v}(\mathcal{I}(T(x))\,|\, R(T(\partial\mathcal{D}\cup\{ v,x\}))=x)=\mathbb{P}^{x}(\mathcal{O}(T(v))\,|\, R(T(\partial\mathcal{D}\cup\{ v,x\}))=v)\]
 where $\mathcal{I}$ is $\mathcal{O}$ for the reversed path, or,
in other words, with $\mathcal{H}^{\textrm{out}}$ replaced by $\mathcal{H}^{\textrm{in}}$
in (\ref{eq:defOt}). Now, since $\LE(R[0,T(\partial\mathcal{D}\cup\{ x\})])$
does not depend on how many times we returned to $v$, we get\[
\mathbb{P}^{v}(\mathcal{I}(T(x))\,|\, R(T(\partial\mathcal{D}\cup\{ x\}))=x)=\mathbb{P}^{v}(\mathcal{I}(T(x))\,|\, R(T(\partial\mathcal{D}\cup\{ v,x\}))=x)\]
Lemma \ref{lem:omerB} now shows that, if $T_{n}=T_{n}(x)$ is the
time of the $n$-th return to $x$ before hitting $\partial\mathcal{D}$
then $\LE(R[0,T_{1}])\sim\LE(R[0,T_{n}])$ and therefore\begin{equation}
\mathbb{P}^{v}(\mathcal{I}(T_{n})\,|\, T_{n}<T(\partial\mathcal{D}))\leq Cr^{-K-4}\quad n=1,2,\dotsc\label{eq:ITn}\end{equation}

To finish the lemma we need only sum over $n$ and $x$. The transience
of $G$ shows that \[
\mathbb{P}^{v}(T_{n}>\lambda)\leq Ce^{-c\lambda}\]
and therefore for $\lambda=C(G)\log r$ with $C$ sufficiently large
we have that this probability is $\leq Cr^{-K-3}$. Therefore we get\begin{align*}
\mathbb{P}\Big(\bigcup_{n}\{ T_{n}<T(\partial\mathcal{D})\}\cap\mathcal{I}(T_{n})\Big) & \leq Cr^{-K-3}+\sum_{n=1}^{\lambda}\mathbb{P}(\mathcal{I}(T_{n})\,|\, T_{n}<T(\partial\mathcal{D}))\\
 & \leq Cr^{-K-3}+(C\log r)r^{-K-4}\leq Cr^{-K-3}.\end{align*}
Denoting by $\mathcal{X}$ the $x$-s satisfying (\ref{eq:xgood})
we can sum over $x$ and get\begin{eqnarray*}
\lefteqn{\mathbb{P}\Big(\bigcup_{n,x}\{ T_{n}(x)<T(\partial\mathcal{D})\}\cap\mathcal{I}(T_{n})\Big)\leq}\\
 & \qquad & \leq\sum_{x\in\mathcal{X}}\mathbb{P}\Big(\bigcup_{n}\{ T_{n}(x)<T(\partial\mathcal{D})\}\cap\mathcal{I}(T_{n})\Big)+\sum_{x\not\in\mathcal{X}}\mathbb{P}(T_{1}(x)<T(\partial\mathcal{D}))\\
 &  & \leq(\#\mathcal{X})\cdot Cr^{-K-3}+\#(\overline{\mathcal{D}}\setminus\mathcal{X})\cdot Cr^{-K-3}\leq Cr^{-K}.\end{eqnarray*}
Which finishes the lemma, since the event $\bigcup_{n,x}\{ T_{n}(x)<T(\partial\mathcal{D})\}\cap\mathcal{I}(T_{n})$
is exactly the desired event.
\end{proof}
The combination of lemmas \ref{lem:outgoing} and \ref{lem:incoming}
is that for an appropriate $\delta=\delta(\epsilon,K,G)$ we have\begin{equation}
\mathbb{P}(\exists t\leq T(\partial\mathcal{D}):\LE(R[0,t])\in\mathcal{H}(r^{-\delta},r^{\epsilon}))\leq Cr^{-K}.\label{eq:inout}\end{equation}

\subsection{Proof of theorem \refs{thm:QL} in three dimensions}

The proof works in three different scales, which we will denote by
$\rho\ll\sigma\ll\tau\ll r$. $\rho=r^{1-\epsilon}$ is the scale
of the {}``closeness'' of the two ends of the quasi-loop. $\tau$
is the scale of the diameter of the quasi-loop, so after all parameters
are fixed we shall fix some $\delta>0$ and then define $\tau:=r^{1-\delta}$.
$\sigma$ is an auxiliary scale --- rather than estimate for some
point $w\in\frac{1}{3}\rho\mathbb{Z}^{3}$ that the probability that
$w\in\mathcal{QL}(\rho,\tau,\LE(R[0,T(\partial\mathcal{D})]))$, we
shall show it simultaneously for all points in $B(w,\sigma)\cap\frac{1}{3}\rho\mathbb{Z}^{3}$. 

\begin{lem}
\label{lem:onegam}Let $G$ be a three dimensional isotropic graph,
and let $\epsilon>0$ and $\delta>0$ be given. Let $1\leq\rho\leq\sigma\leq\frac{1}{2}\tau\leq\frac{1}{2}r$
and $\sigma\geq2\rho(r^{\epsilon}+2)$. Let $w$ be some vertex and
let $b\subset B(w,\tau)$ be some set. Let $\gamma\in\mathcal{H}=\mathcal{H}(r^{-\delta},r^{\epsilon})$
be a path starting from some point in $\overline{B(w,\tau)}$ and
ending on $\partial B(w,\tau)$. Let $R$ be a random walk starting
from some $x\in\partial B(w,\sigma)$ and stopped at $\partial B(w,\tau)\cup b$;
and define \[
\Delta(\gamma):=\#\{ y\in{\textstyle \frac{1}{3}}\rho\mathbb{Z}^{3}\cap B(w,\sigma):d(y,\gamma)\leq\rho,\, d(y,R)\leq\rho\}\cdot\mathbf{1}_{\{ R\cap\gamma=\emptyset\}}\]
Then \begin{align}
\mathbb{E}^{x}\Delta(\gamma)\mathbf{1}_{\{ R\cap b=\emptyset\}} & \leq C(G)r^{-\delta}\log^{2}r,\label{eq:Delta1}\\
\mathbb{E}^{x}\Delta(\gamma)\mathbf{1}_{\{ R\cap b\neq\emptyset\}} & \leq C(G)\mathbb{P}^{x}(R\cap b\neq\emptyset)\log^{2}r.\label{eq:Delta2}\end{align}

\end{lem}
\begin{proof}
Denote the first by $\Delta_{1}(\gamma)$ and the second by $\Delta_{2}(\gamma)$.
Define stopping times $s_{1}\leq s_{1}^{*}\leq s_{2}\leq\dotsc\leq T(\partial B(w,\tau)\cup b)$
by \[
s_{1}:=\min\{ t:d(R(t),\gamma\cap B(w,\sigma+\rho))\leq2\rho\},\]
 and for $i\geq1$ \begin{align*}
s_{i}^{*} & :=\min\{ t\geq s_{i}:R(t)\in\partial B(R(s_{i}),8\rho)\}\\
s_{i} & :=\min\{ t\geq s_{i-1}^{*}:d(R(t),\gamma\cap B(w,\sigma+\rho))\leq2\rho\}\end{align*}
(for clarity we removed the conditions $s_{i},s_{i}^{*}\leq T(\partial B(w,\tau)\cup b)$
--- if $R$ hits $\partial B(w,\linebreak[1]\tau)\cup b$ after any
of them, consider the sequence stabilized at this point). Let $I$
be the number of $s_{i}$'s defined before the process is stopped
i.e.~$R(s_{I+1})\in\partial B(w,\tau)\cup b$. Lemma \ref{lem:Beurling}
gives for $i<I$ that the probability to intersect $\gamma$ between
$s_{i}$ and $s_{i}^{*}$ is $\geq c(G)/\log r$. This gives \begin{equation}
\mathbb{P}(I\geq i)\leq\left(1-\frac{c(G)}{\log r}\right)^{i-1}\quad\forall i.\label{eq:lexpologr}\end{equation}
 Further, for the time period between $s_{i}$ and $s_{I+1}$ the
definition of $\mathcal{H}$ gives that \[
\mathbb{P}^{x}(\left\{ R[s_{i},s_{I+1}]\cap\gamma=\emptyset\right\} \cap\left\{ R(s_{I+1})\not\in b\right\} \,|\, I\geq i)\leq r^{-\delta}.\]
Here is where we use the condition $\sigma\geq2\rho(r^{\epsilon}+2)$,
since $d(R(s_{i}),\partial B(w,\tau))\geq\tau-\sigma-3\rho\geq\sigma-3\rho$.
Finally, it is clear that $\Delta(\gamma)<CI$. Therefore\begin{align}
\mathbb{E}^{x}\Delta_{1}(\gamma) & =\sum_{i=0}^{\infty}\mathbb{E}^{x}(\Delta_{1}(\gamma)\cdot\mathbf{1}_{\{ I=i\}})\leq\sum_{i=0}^{\infty}Ci\cdot\mathbb{P}^{x}(\left\{ \Delta_{1}(\gamma)>0\right\} \cap\left\{ I\geq i\right\} )\leq\nonumber \\
 & \leq C\sum_{i=1}^{\infty}i\mathbb{P}^{x}(\left\{ R(s_{I+1})\not\in b\right\} \cap\{ R\cap\gamma=\emptyset\}\cap\{ I\geq i\})\leq\nonumber \\
 & \leq C\sum_{i=1}^{\infty}i\left(1-\frac{c(G)}{\log r}\right)^{i-1}\cdot Cr^{-\delta}\leq C(G)r^{-\delta}\log^{2}r.\label{eq:Xvg_Y}\end{align}

As for $\mathbb{E}^{x}\Delta_{2}(\gamma)$, we first notice that if
$b$ is empty then there is nothing to prove. If $b$ is not empty,
then the estimate of Green's function (lemma \ref{lem:a}) shows that
$\mathbb{P}(R\cap b\neq\emptyset)>C/r$. Therefore we may use (\ref{eq:lexpologr})
to show that If $\lambda=\lambda(G)$ is a sufficiently large constant
then\[
\mathbb{P}(I>\lambda\log^{2}r)\leq r^{-4}\]
Hence we get\begin{align*}
\mathbb{E}^{x} & \Delta_{2}(\gamma)\leq\mathbb{E}^{x}\Delta_{2}(\gamma)\mathbf{1}_{\{ I>\lambda\log^{2}r\}}+\mathbb{E}^{x}\Delta_{2}(\gamma)\mathbf{1}_{\{ I\leq\lambda\log^{2}r\}}\leq\\
 & \leq C\left(\frac{\sigma}{\rho}\right)^{3}r^{-4}+C\lambda(\log r)^{2}\mathbb{P}^{x}(R\cap b\neq\emptyset)\leq C\mathbb{P}^{x}(R\cap b\neq\emptyset)\log^{2}r.\qedhere\end{align*}

\end{proof}
Returning to the proof of theorem \ref{thm:QL}, we shall from this
point on assume that $2r^{1-\epsilon/2}\leq\sigma\leq r^{1-2\delta}$.
We fix some $w\in\mathcal{D}$ and some stopping time $T$ and define\[
\mathcal{X}(T):=\mathcal{X}_{w}(T):=\#(\mathcal{QL}(\rho,\tau,\LE(R[0,T]))\cap B(w,\sigma)\cap{\textstyle \frac{1}{3}}\rho\mathbb{Z}^{3})\]
so $\mathcal{X}:=\mathcal{X}(T(\partial\mathcal{D}))$ is what we
need to estimate. Define exit and entry stopping times by $T_{0}=0$
and\begin{align*}
T_{2i+1} & =\min\{ t\geq T_{2i}:R(t)\in\partial B(w,2\sigma)\cup\partial\mathcal{D}\}\\
T_{2i} & =\min\{ t\geq T_{2i-1}:R(t)\in\partial B(w,4\sigma)\cup\partial\mathcal{D}\}.\end{align*}
Let $I$ be the first $i$ such that $R(T_{i})\in\partial\mathcal{D}$.
We note that lemma \ref{lem:DDunbounded} shows that \[
\mathbb{P}(R(T_{2i+1})\in\partial\mathcal{D}\,|\, R[0,T_{2i}])\geq\mathbb{P}^{R(T_{2i})}(T_{w,2\sigma}=\infty)\geq c(G)\]
and hence we get that the probability that $I$ is large drops exponentially,
and hence if $\lambda=\lambda(G)$ is sufficiently large we get\begin{equation}
\mathbb{P}(I>\lambda\log r)\leq\frac{1}{r^{4}}.\label{eq:Ilamr3}\end{equation}
Denote $M:=\left\lfloor \lambda\log r\right\rfloor $. With this notation
we can write\begin{equation}
\mathbb{E}\mathcal{X}\leq C\left(\frac{\sigma}{\rho}\right)^{3}\mathbb{P}(I>M)+\sum_{i=1}^{M}\mathbb{E}(\mathcal{X}\cdot\mathbf{1}_{\{ I=i\}})\stackrel{(\ref{eq:Ilamr3})}{\leq}Cr^{-1}+\sum_{i=1}^{M}\mathbb{E}(\mathcal{X}\cdot\mathbf{1}_{\{ I=i\}}).\label{eq:X_sum_Xk}\end{equation}
In other words, we can ignore the first summand in (\ref{eq:X_sum_Xk}).
Let us therefore define the number of quasi-loops up to the $i$th
time $\mathcal{X}_{i}:=\mathcal{X}(T_{i})$. The process of loop-erasing
between $T_{2i}$ and $T_{2i+1}$ can only destroy quasi loops in
$B(w,\sigma)$ so we get \[
\mathcal{X}\cdot\mathbf{1}_{\{ I=2i+1\}}\leq\mathcal{X}_{2i}\cdot\mathbf{1}_{\{ I=2i+1\}}\leq\mathcal{X}_{2i}\cdot\mathbf{1}_{\{ I>2i\}}.\]
 With this in mind we define \[
\Delta_{i}:=(\mathcal{X}_{2i+2}-\mathcal{X}_{2i})\cdot\mathbf{1}_{\{ I>2i+2\}}.\]

\begin{lem}
\label{lem:Deltai}With the notations above and \begin{equation}
\eta:=\min\delta({\textstyle \frac{1}{2}}\epsilon,4,G),1\label{eq:defdelf}\end{equation}
 where $\delta(\cdot)$ is given by (\ref{eq:inout}),\[
\mathbb{E}\Delta_{i}\leq Cir^{-\eta}\log^{2}r.\]

\end{lem}
\begin{proof}
Examine some $y$,\[
y\in\mathcal{QL}(\rho,\tau,\LE(R[0,T_{2i+2}]))\setminus\mathcal{QL}(\rho,\tau,\LE(R[0,T_{2i}])).\]
 Directly from the definitions, we must have 
\begin{enumerate}
\item $y$ is $\rho$-near at least one component $\gamma$ of $\LE(R[0,T_{2i}])\cap B(w,4\sigma)$.
\item $R[T_{2i},T_{2i+1}]$ gets $\rho$-near $y$ and then fails to intersect
at least one of the segments $\gamma$ from (i), as well as $\partial\mathcal{D}\cap B(w,4\sigma)$. 
\end{enumerate}
In other words, the number of such $y$'s in $\frac{1}{3}\rho\mathbb{Z}^{3}\cap B(w,\sigma)$
with respect to a specific $\gamma$, can be bounded by $\Delta(\gamma)$,
$\Delta$ from lemma \ref{lem:onegam}, with the parameters in the
following table:

\noindent \begin{center}\begin{tabular}{|c|c|c|c|c|c|c|}
\hline 
Lemma \ref{lem:onegam}&
$\rho$&
$\sigma$&
$\tau$&
$\epsilon$&
$\delta$&
$b$\tabularnewline
\hline 
here&
$\rho$&
$2\sigma$&
$4\sigma$&
$\frac{1}{2}\epsilon$&
$\eta$&
$\partial\mathcal{D}\cap B(w,4\sigma)$\tabularnewline
\hline
\end{tabular}\end{center}

\noindent Note that (\ref{eq:inout}) shows that with probability
$\geq1-Cr^{-4}$ we have $\gamma\in\mathcal{H}=\mathcal{H}(r^{-\eta},\linebreak[1]r^{\epsilon/2})$
for all $\gamma\subset\LE(R[0,T_{2i+1}])$. 

With this in mind we denote by $\Gamma_{i}$ the collection of connected
components $\gamma$ of $\LE(R[0,T_{2i+1}])\cap B(w,4\sigma)$ satisfying
$\gamma\cap B(w,\sigma+\rho)\neq\emptyset$ and get \begin{align*}
\mathbb{E}\Delta_{i}\cdot\mathbf{1}_{\{\Gamma_{i}\subset\mathcal{H}\}} & \leq\max_{x\in\partial B(w,2\sigma)}\sum_{\gamma\in\Gamma_{i}}\mathbb{E}^{x}\Delta(\gamma)\cdot\mathbf{1}_{\{ R\cap\partial\mathcal{D}=\emptyset\}}\stackrel{(\ref{eq:Delta1})}{\leq}C(\#\Gamma_{i})r^{-\eta}\log^{2}r\\
\mathbb{E}\Delta_{i}\cdot\mathbf{1}_{\{\Gamma_{i}\not\subset\mathcal{H}\}} & \leq C(\sigma/\rho)^{3}\mathbb{P}(\Gamma_{i}\not\subset\mathcal{H})\leq Cr^{3}\cdot r^{-4}.\end{align*}
It easy to see that $\#\Gamma_{i}\leq i$, and the lemma is finished.
\end{proof}
Summing lemma \ref{lem:Deltai} up to $i$ we get\[
\mathbb{E}\mathcal{X}_{2i}\cdot\mathbf{1}_{\{ I>2i\}}\leq Cr^{-\eta}(\log r)^{2}i^{2}\]
($\eta$ being defined by (\ref{eq:defdelf})). Another summation,
up to $M$, will give us\begin{align}
\mathbb{E}\mathcal{X} & \leq Cr^{-\eta}\log^{5}r+\mathbb{E}\Delta'\label{eq:EX_not_boundary}\\
\Delta' & :=(\mathcal{X}_{I}-\mathcal{X}_{I-2})\cdot\mathbf{1}_{\{ I\textrm{ even}\}}.\nonumber \end{align}
 Thus we are left with the estimate of $\mathbb{E}\Delta'$, which
is the behavior near the boundary --- if $\partial\mathcal{D}\cap B(w,4\sigma)=\emptyset$
then of course $I$ is always odd and we get $\Delta'\equiv0$. It
is at this point that we utilize the difference between $\tau$ and
$\sigma$. 

\begin{lem}
\label{lem:Delpharm}Let $\omega=\omega(v,\mathcal{D})$ be the discrete
harmonic measure from $v$ i.e.~$\omega(A):=\mathbb{P}^{v}(R(T(\partial\mathcal{D}))\in A)$
for all $A\subset\partial\mathcal{D}$. Then\[
\mathbb{E}^{v}(\Delta')\leq C\frac{\sigma}{\tau}\omega(B(w,16\sigma))\log^{2}r.\]

\end{lem}
\begin{proof}
 Define stopping times $U_{0}\leq U_{1}\leq\dotsb$ using $U_{0}=0$
and\begin{align*}
U_{2j+1} & :=\min\{ t\geq U_{2j}:R(t)\in\partial B(w,8\sigma)\cup\partial\mathcal{D}\}\\
U_{2j} & :=\min\{ t\geq U_{2j-1}:R(t)\in\partial B(w,{\textstyle \frac{1}{2}}\tau)\cup\partial\mathcal{D}\}.\end{align*}
We assume here that $\sigma<\frac{1}{32}\tau$ which will hold if
$r$ is sufficiently large. Define $J$ to be the first $j$ such
that $R(U_{j})\in\partial\mathcal{D}$. Our first target is to connect
$J$ and $\omega(B(16,\sigma))$. Denote\[
p(x):=\mathbb{P}^{x}(T(\partial\mathcal{D}\cap B(w,4\sigma))<T_{w,\tau/2})\quad x\in\partial B(w,8\sigma).\]
We have\begin{align}
\omega(B(w,16\sigma)) & =\sum_{j=1}^{\infty}\mathbb{P}(\{ J=j\}\cap\{ R(U_{j})\in\partial\mathcal{D}\cap B(w,16\sigma)\})\geq\nonumber \\
 & \geq\mathbb{P}(\{ J=2\}\cap\{ R(U_{2})\in\partial\mathcal{D}\cap B(w,16\sigma)\})\geq\nonumber \\
 & \geq\mathbb{P}(J>1)\mathbb{E}\mathbb{P}^{R(U_{1})}(T(\partial\mathcal{D})<T_{w,16\sigma})\geq\nonumber \\
 & \geq\mathbb{P}(J>1)\min_{x\in\partial B(w,8\sigma)}\mathbb{P}^{x}(T(\partial\mathcal{D}\cap B(w,4\sigma))<T_{w,16\sigma})\geq\nonumber \\
 & \!\stackrel{(*)}{\geq}\! c(G)\mathbb{P}(J>1)\min p(x).\label{eq:harmBws4sig}\end{align}
where $(*)$ comes from using lemma \ref{lem:hitbsame} to change
from $T_{w,16\sigma}$ to $T_{w,\tau/2}$. $\min$ here and later
on also $\max$ always refers to the minimum over $\partial B(w,8\sigma)$.

The estimate of $\Delta'$ follows from the clause (\ref{eq:Delta2})
of lemma \ref{lem:onegam}, which we now use with the parameters as
follows:

\noindent \begin{center}\begin{tabular}{|c|c|c|c|c|c|c|}
\hline 
Lemma \ref{lem:onegam}&
$\rho$&
$\sigma$&
$\tau$&
$\epsilon$&
$\delta$&
$b$\tabularnewline
\hline 
here&
$\rho$&
8$\sigma$&
$\frac{1}{2}\tau$&
$\frac{1}{2}\epsilon$&
$\eta$&
$\partial\mathcal{D}\cap B(w,4\sigma)$\tabularnewline
\hline
\end{tabular}\end{center}

\noindent And an argumentation identical to that of the previous lemma
gives that, \[
\mathbb{E}(\Delta'\,|\, J>2j-1)\leq C(G)(\#\Gamma_{j}')\max p(x)\log^{2}r\]
where $\Gamma_{j}'$ is the collection of connected components $\gamma$
of $\LE(R[0,U_{2j-1}])\cap B(w,\frac{1}{2}\tau)$ satisfying $\gamma\cap B(w,\sigma+\rho)\neq\emptyset$.
Again, it is easy to see that $\#\Gamma_{j}'\leq j-1$ and in particular
$\#\Gamma_{1}'=0$. Hence we get\begin{equation}
\mathbb{E}(\Delta')=\sum_{j=1}^{\infty}\mathbb{E}\Delta'\cdot\mathbf{1}_{\{ J=2j\}}\leq C(G)\max p(x)\log^{2}r\sum_{j=2}^{\infty}(j-1)\mathbb{P}(J>2j-1)\label{eq:XXp}\end{equation}
The Green's function estimate (\ref{eq:atrans}) also shows that \begin{align*}
\mathbb{P}(J>2j+1\,|\, J>2j) & =\mathbb{EP}^{R(U_{2i})}(R(T(\partial B(w,8\sigma)\cup\partial\mathcal{D}))\in B(w,8\sigma))\leq\\
 & \leq\mathbb{EP}^{R(U_{2i})}(T_{w,8\sigma}<\infty)\leq C\frac{\sigma}{\tau}.\end{align*}
 Hence we have (for $r$ sufficiently large),\[
\sum_{j=1}^{\infty}j\mathbb{P}(J>2j+1)\leq\mathbb{P}(J>1)\sum_{j=1}^{\infty}j\left(C\frac{\sigma}{\tau}\right)^{j}\leq C\mathbb{P}(J>1)\frac{\sigma}{\tau}\]
and with (\ref{eq:harmBws4sig}), (\ref{eq:XXp}) and Harnack's inequality
(lemma \ref{lem:Harnack-general}) which shows that $\max p(x)\leq C(G)\min p(x)$,
we get\[
\mathbb{E}(\Delta')\leq C\frac{\sigma}{\tau}\omega(B(w,16\sigma))\log^{2}r.\qedhere\]

\end{proof}
Lemma \ref{lem:Delpharm} and (\ref{eq:EX_not_boundary}) give\begin{equation}
\mathbb{E}(\mathcal{X}_{w})\leq Cr^{-\eta}\log^{5}r+C\left(\frac{\sigma}{\tau}\right)\omega(B(w,16\sigma))\log^{2}r.\label{eq:XXfinal}\end{equation}
For any point in $z\in B(v,r)\cap\frac{1}{3}\sigma\mathbb{Z}^{3}$,
let $w_{z}$ be the point closest to $z$ in $G$. Clearly, if $r$
is sufficiently large, the balls $B(w_{z},\sigma)$ form a cover of
$\mathcal{D}$ and furthermore\[
\sum_{z}\omega(B(w_{z},16\sigma))\leq C\omega(\partial\mathcal{D})=C.\]
Hence we can sum (\ref{eq:XXfinal}) over $z$ and get\[
\mathbb{E}\QL(\rho,\tau,\LE(R[0,T(\partial\mathcal{D})]))\leq C\left(\frac{r}{\sigma}\right)^{3}r^{-\eta}\log^{5}r+C\left(\frac{\sigma}{\tau}\right)\log^{2}r.\]
Now is the time to pick $\delta$. We take $\delta=\min\frac{1}{11}\eta,\frac{1}{7}\epsilon$
and define $\tau=r^{1-\delta}$ and $\sigma=r^{1-3\delta}$. For $r$
sufficiently large the condition $\sigma>2r^{1-\epsilon/2}$ would
be fulfilled. We get\[
\mathbb{E}\QL(r^{1-\epsilon},r^{1-\delta},\LE(R[0,T(\partial\mathcal{D})]))\leq Cr^{-2\delta}\log^{5}r+Cr^{-2\delta}\log^{2}r\leq Cr^{-\delta}\]
and the theorem is proved.\qed

\subsection{Proof of theorem \refs{thm:QL} in two dimensions}

\begin{lem}
\label{lem:ann2}Let $G$ be a two dimensional Euclidean net. Let
$v\in G$, let $r>C(G)$ and let $\gamma$ be a path from $\partial B(v,r)$
to $\partial B(v,2r)$. Let $w\in B(v,r)$ and let $R$ be a random
walk starting from $w$ and stopped on $\partial B(v,2r)$. Then\[
\mathbb{P}(R\cap\gamma\neq\emptyset)>c(G).\]

\end{lem}
\begin{proof}
It is easy to see, applying lemma \ref{lem:DD}, say three times,
that there is a probability $>c(G)$ that $R$ does a loop around
the annulus $B(v,\frac{3}{2}r)\setminus\overline{B(v,r)}$. Two dimensional
geometry shows that in this case the linear extensions of $R$ and
$\gamma$ intersect (the linear extension of a path in $G$ is a path
in $\mathbb{R}^{2}$ which is composed of all points of $\gamma$
connected by linear segments). Since the length of edges in $G$ is
bounded, we get\[
\mathbb{P}(\exists t\leq T_{v,2r}:d(R(t),\gamma\cap B(v,{\textstyle \frac{3}{2}}r))<C(G))>c(G).\]
However, the first such $t$ is a stopping time, so we can consider
the walk after it as a regular random walk, and of course it has a
positive probability to hit $\gamma$.
\end{proof}
\begin{lem}
\label{lem:kesten}Let $G$ be a two dimensional Euclidean net. Let
$v\in G$ and let $s>r>1$. Let $\gamma$ be a path from $\partial B(v,r)$
to $\partial B(v,s)$. Let $R$ be a random walk starting from $v$.
Then\[
\mathbb{P}(R\cap\gamma=\emptyset)\leq C(G)\left(\frac{r}{s}\right)^{c(G)}.\]

\end{lem}
This follows immediately from the previous lemma and the Wiener shell
test.

The proof of theorem \ref{thm:QL} now proceeds exactly as in the
three dimensional case, with lemma \ref{lem:ann2} serving as a replacement
for lemma \ref{lem:Beurling} and lemma \ref{lem:kesten} serving
as a replacement to lemmas \ref{lem:outgoing} and \ref{lem:incoming}.
The only complication is lemma \ref{lem:Delpharm}, which no longer
holds as stated. It is necessary at this point to estimate $\mathbb{P}(J>2i+1\,|\, J>2i)$
by lemma \ref{lem:kesten} for the incoming walk, with the result
that one gets only $(\sigma/\tau)^{c}$ in the formulation of the
lemma. See \cite{K}, sublemmas 18.2 \& 18.3 for details.

\subsection{Continuity in the starting point}

We will need one additional corollary of the techniques of this chapter:

\begin{lem}
\label{lem:contsp}Let $G$ be an isotropic graph and let $v,w\in\mathcal{E}\subset\mathcal{D}\subset G$,
$\mathcal{D}$ finite. Denote\[
p^{x}:=\mathbb{P}^{x}(\LE(R[0,T(\partial\mathcal{D})])\subset\mathcal{E}),\quad x=v,w.\]
Then \[
|p^{v}-p^{w}|\leq C(G)\left(\frac{|v-w|}{d(v,\partial\mathcal{E})}\right)^{c(G)}.\]

\end{lem}
\begin{proof}
Denote $\mu=|v-w|/d(v,\partial\mathcal{E})$. We may assume w.l.o.g
that $\mu$ is sufficiently small. Let $H$ be the graph generated
by taking $\overline{\mathcal{D}}$ and identifying all the points
of $\partial\mathcal{D}$ (this process is often called {}``wiring
the boundary''). Let $T$ be the uniform spanning tree on $H$. Then
by Pemantle \cite{P91}, the distribution of the path in $T$ from
$x$ to $\partial\mathcal{D}$ is identical to the distribution of
a loop-erased random walk on \emph{}$H$ starting from $x$ and stopped
on $\partial\mathcal{D}$, which is identical to the distribution
on the graph $G$. Hence the two branches $\beta^{v}$ and $\beta^{w}$
of $T$ from $v$ and $w$ to $\partial\mathcal{D}$ is a coupling
of the two loop-erased walks. Denote by $\gamma^{x}$ the portion
of $\beta^{x}$ from $x$ until the unique intersection of $\beta^{v}$
and $\beta^{w}$; and by $p$ the probability that $\gamma^{v}\cup\gamma^{w}\subset\mathcal{E}$.
Then clearly, \[
|p^{v}-p^{w}|\leq1-p.\]
Now, by Wilson's algorithm \cite{W96}, $\gamma^{w}$ may be constructed
by first constructing $\beta^{v}$ and then taking a random walk starting
from $w$, stopping it when it first hits $\beta^{v}$ (possibly at
time $0$) and performing loop-erasure on the result. Hence we need
to show that\[
\mathbb{P}^{w}(T(\beta^{v})>T(\partial\mathcal{E}))\leq C(G)\mu^{c(G)}.\]
In three dimensions we use lemma \ref{lem:wienersimp}. Let $K=\frac{1}{2}$
and let $\delta$ be the $\delta(\frac{1}{2},G)$ of lemma \ref{lem:wienersimp}.
Let the $r$ and $s$ of lemma \ref{lem:wienersimp} be $d(w,\partial\mathcal{E})$
and $\max|v-w|,1$ respectively --- if $\mu$ is sufficiently small
we would get $r>s$. We get that (except for probability $C\mu^{1/2}$
in the walk starting from $v$)\[
\mathbb{P}^{2,w}(R^{2,w}[0,T_{w,r}^{2}]\cap\LE(R^{1}[0,T^{1}(\partial\mathcal{D})])\cap B(w,r)\setminus\overline{B(w,s)}=\emptyset)\leq C\mu^{\delta}.\]
Therefore $\mathbb{P}(\gamma^{w}\not\subset\mathcal{E})\leq C\mu^{\delta}+C\mu^{1/2}$.
The same holds for $\mathbb{P}(\gamma^{v}\not\subset\mathcal{E})$
and the three dimensional case is finished. The two dimensional case
follows similarly from lemma \ref{lem:kesten}.
\end{proof}
\begin{rem*}
The use of lemma \ref{lem:wienersimp} to estimate the probability
that a loop-erased random walk and a random walk starting from close
points will hit is somewhat an overkill. For example, in $\mathbb{Z}^{3}$,
if they start from the same point then the nice symmetry argument
of \cite{L99} can show that this probability is $>1-Cr^{-1/3}$.
Presumably, an equivalent argument would work in our case as well.
The arguments of \cite{AB99} should also give a usable estimate.
\end{rem*}

\section{\label{sec:Isotropic-gluing}Isotropic interpolation}

The purpose of this chapter is to compare random walks on two or more
graphs all of which are isotropic, with uniformly bounded structure
constants. We shall call such a collection $\mathcal{G}$ an \textbf{isotropic
family} and denote by $C(\mathcal{G})$, $c(\mathcal{G})$ etc.~constants
which depend only on the maximum of the isotropic structure constants
of all $G\in\mathcal{G}$.

\subsection{Hitting probabilities}

In this section we will compare probabilities by proving inequalities
of the sort $|p-q|\leq\epsilon\max p,q$. It will be convenient to
denote this by\[
p\stackrel{\epsilon}{\simeq}q.\]
When $\epsilon=Cr^{-c}$ for some constants $C(\mathcal{G})$ and
$c(\mathcal{G})$ we will usually omit it, and just write $p\simeq q$
($r$ will be clear from the context). Occasionally we will prove
instead that $|p-q|\leq\epsilon p$ or $|p-q|\leq\epsilon\min p,q$.
We will always assume $\epsilon\leq\frac{1}{2}$, and then they are
all equivalent up to constants. 

We will often use the following version of differentiation of product:
assume\[
p^{1}\stackrel{\alpha}{\simeq}p^{2}\quad q^{1}\stackrel{\beta}{\simeq}q^{2}\quad\alpha,\beta\leq\frac{1}{2}.\]
Then\begin{align}
|p^{1}q^{1}-p^{2}q^{2}| & \leq|p^{1}-p^{2}|q^{1}+p^{2}|q^{1}-q^{2}|\leq\alpha q^{1}\max p^{1},p^{2}+\beta p^{2}\max q^{1},q^{2}\nonumber \\
 & \leq\alpha q^{1}(1+2\alpha)p^{1}+\beta p^{2}(1+2\beta)q^{2}\leq\nonumber \\
 & \leq(\alpha+\beta+2(\alpha^{2}+\beta^{2}))\max p^{1}q^{1},p^{2}q^{2}\leq\label{eq:ddm}\\
 & \leq2(\alpha+\beta)\max p^{1}q^{1},p^{2}q^{2}\label{eq:ddm2}\end{align}
Another useful fact: if $p^{i}=\sum_{n}q_{n}^{i}$ then\begin{align}
|p^{1}-p^{2}| & \leq\sum_{n}|q_{n}^{1}-q_{n}^{2}|\leq\sum_{n}\beta\max q_{n}^{1},q_{n}^{2}\leq\sum_{n}\beta(q_{n}^{1}+q_{n}^{2})=\beta(p^{1}+p^{2})\nonumber \\
 & \leq2\beta\max p^{1},p^{2}.\label{eq:dds}\end{align}

\begin{lem}
\label{lem:nodir}Let $G$ be an isotropic graph and let $r_{1}>s>r_{2}>1$.
Let $v\in G$, $v^{1},v^{2}\in\partial B(v,s)$ and let $w\in A:=\partial B(v,r_{1})\cup\partial B(v,r_{2})$.
Let \[
p^{i}=\mathbb{P}^{v^{i}}(R(T(A))=w).\]
Then \begin{equation}
\left|p^{1}-p^{2}\right|\leq C(G)\left(\max\frac{s}{r_{1}},\frac{r_{2}}{s}\right)^{c(G)}\max p_{1},p_{2}.\label{eq:lemnodir}\end{equation}
Similarly if $A:=\partial B(v,r_{1})\cup\{ v\}$ then $|p_{1}-p_{2}|\leq C\left(\max\frac{s}{r_{1}},\frac{1}{s}\right)^{c(G)}\max p_{1},p_{2}$.
\end{lem}
\begin{proof}
Assume first that $r_{2}\geq r_{1}^{1-c_{\ref{c:RWclose}}/2}$ where
$c_{\ref{c:RWclose}}(G)$ is from lemma \ref{lem:rismall}, and assume
also that $r_{1}=2^{N}s=2^{2N}r_{2}$ for some integer $N$ (removing
both assumptions is easy, and we do it in the end of the lemma). Assume
also that $w\in\partial B(v,r_{2})$. Define \begin{align*}
a_{j} & :=s2^{-j}\quad j=0,1,\dotsc,N.\\
X_{j} & :=\partial B(v,a_{j})\cup\partial B(v,r_{1})\end{align*}
Define stopping times $T_{0}=0$ and \[
T_{j}:=\min\{ t\geq T_{j-1}:R(t)\in X_{j}\}.\]
Define probability measures $\mu_{j}^{i}$ on $\partial B(v,a_{j})$
by \[
\mu_{j}^{i}(x)=\mathbb{P}^{v^{i}}(R(T_{j})=x\,|\, R(T_{j})\in\partial B(v,a_{j})).\]

We proceed by examining how $\mu_{j}^{i}$ evolves with $j$. It will
be useful to use $L^{1}$ estimates during intermediate stages, so
define \[
\epsilon_{j}:=\sum_{x\in\partial B(v,a_{j})}|\mu_{j}^{1}(x)-\mu_{j}^{2}(x)|.\]
Clearly $\epsilon_{0}=2$. For $x\in\partial B(v,a_{j})$ and $y\in\partial B(v,a_{j+1})$
we define \[
\pi(x,y)=\mathbb{P}^{x}(R(T(X_{j+1}))=y\,|\, R(T(X_{j+1}))\in\partial B(0,a_{j+1}))\]
and then get\[
\mu_{j+1}^{i}(y)=\sum_{x\in\partial B(0,a_{j})}\mu_{j}^{i}(x)\pi(x,y).\]

Define now \begin{gather*}
A^{+}:=\{ x\in\partial B(0,a_{j}):\mu_{j}^{1}(x)\geq\mu_{j}^{2}(x)\},\quad A^{-}:=\partial B(0,a_{j})\setminus A^{+},\\
D^{\pm}(y):=\sum_{x\in A^{\pm}}|\mu_{j}^{1}(x)-\mu_{j}^{2}(x)|\pi(x,y).\end{gather*}
 Clearly, \begin{equation}
\sum_{y}D^{\pm}(y)=\sum_{x\in A^{\pm}}|\mu_{j}^{1}(x)-\mu_{j}^{2}(x)|\sum_{y}\pi(x,y)=\sum_{x\in A^{\pm}}|\mu_{j}^{1}(x)-\mu_{j}^{2}(x)|={\textstyle \frac{1}{2}}\epsilon_{j}.\label{eq:sumDpmeps}\end{equation}
Next, Harnack's inequality (lemma \ref{lem:Harnack-general}) shows
that $\pi(x,y)\approx\pi(x',y)$ for any $x,x'\in\partial B(0,a_{j})$
so $D^{+}(y)\approx D^{-}(y)$ ($\pi(x,y)$ is a quotient of two harmonic
functions and we use Harnack's inequality for both). This gives that\[
|D^{+}(y)-D^{-}(y)|\leq(1-c(G))(D^{+}(y)+D^{-}(y))\]
and hence\[
\epsilon_{j+1}=\sum_{y\in\partial B(v,a_{j+1})}|D^{+}(y)-D^{-}(y)|\leq(1-c)\sum_{y}D^{+}(y)+D^{-}(y)\stackrel{(\ref{eq:sumDpmeps})}{=}(1-c)\epsilon_{j}.\]
Therefore we get $\epsilon_{N-1}\leq2(1-c)^{N-1}=C(s/r_{1})^{c}$.
This establishes the $L^{1}$ estimate. We now use the last step (from
$N-1$ to $N$) to move to a uniform estimate. Return to our $w\in\partial B(v,r_{2})$.
We have \[
\mu_{N}^{i}(w)=\sum_{x\in\partial B(v,a_{N-1})}\mu_{N-1}^{i}(x)\pi(x,w)\geq\min_{x}\pi(x,w)\]
and on the other hand\begin{align*}
|\mu_{N}^{1}(w)-\mu_{N}^{2}(w)| & \leq\sum_{x\in\partial B(v,a_{N-1})}|\mu_{N-1}^{1}(x)-\mu_{N-1}^{2}(x)|\pi(x,w)\leq\\
 & \leq C(s/r_{1})^{c}\max_{x}\pi(x,w)\end{align*}
and by Harnack's inequality (lemma \ref{lem:Harnack-general}),\begin{equation}
|\mu_{N}^{1}(w)-\mu_{N}^{2}(w)|\leq C(s/r_{1})^{c}\max\mu_{N}^{1}(w),\mu_{N}^{2}(w).\label{eq:muxmuy}\end{equation}

The difference between $\mu_{N}^{i}(w)$ and $p^{i}$ is just that
$\mu_{N}^{i}(w)$ is conditioned, i.e.\[
p^{i}=\mu_{N}^{i}(w)\mathbb{P}^{v^{i}}(R(T_{N})\in\partial B(v,r_{2})).\]
To estimate the second term we use lemma \ref{lem:sball}. The probability
$q^{i}$ that Brownian motion starting from $v^{i}$ hits $\partial B(v,r_{2})$
before $\partial B(v,r_{1})$ is \begin{equation}
q^{i}=\frac{a(\left|v^{i}\right|)-a(r_{1})}{a(r_{2})-a(r_{1})}\quad a(t)=\begin{cases}
t^{2-d} & d\geq3\\
\log t & d=2\end{cases}\label{eq:qiprecis}\end{equation}
and in either case we get $|q^{1}-q^{2}|\leq C(G)s^{-1}\max q^{1},q^{2}$.
Therefore lemma \ref{lem:sball} (take e.g.\ the $\epsilon$ of lemma
\ref{lem:sball} to be $\frac{1}{2}$) gives that \begin{equation}
\mathbb{P}^{v^{1}}(R(T_{N})\in\partial B(0,r_{2}))\simeq\mathbb{P}^{v^{2}}(R(T_{N})\in\partial B(0,r_{2})).\label{eq:Pvir2same}\end{equation}
Together with (\ref{eq:muxmuy}) the lemma is proved in this case.

The case that $w\in\partial B(v,r_{1})$ is identical, with this time
defining $a_{j}=s2^{j}$ and $X_{j}=\partial B(v,a_{j})\cup\partial B(v,r_{2})$.
The argument about the exponential decrease of the $\epsilon_{j}$
works identically. Finally, (\ref{eq:qiprecis}) shows that $\mathbb{P}^{v^{i}}(R(T_{N})\in\partial B(v,r_{1}))>c$
and therefore (\ref{eq:Pvir2same}) implies an identical estimate
for $\mathbb{P}^{v^{i}}(R(T_{N})\in\partial B(v,r_{1}))$. Hence this
case is finished too.

Finally, assume one of the assumptions on $r_{1}$, $s$ and $r_{2}$
fails. If $r_{1}/s$ or $s/r_{2}<2$ the lemma holds trivially for
sufficiently large constants. Hence assume that both are $\geq2$.
It now follows easily that one can find $u_{i}$ such that $r_{1}>u_{1}>s>u_{2}>r_{2}$,
$\frac{u_{1}}{s}=\frac{s}{u_{2}}=2^{N}$ and $u_{2}\geq u_{1}^{1-c_{\ref{c:RWclose}}/2}$,
and further, \[
\frac{s}{u_{1}}\leq\left(\max\frac{s}{r_{1}},\frac{r_{2}}{s}\right)^{c(G)}.\]
We use the case already established and find that for any $x\in\partial B(v,u_{1})\cup\partial B(v,u_{2})$
one has\[
|p^{1}(x)-p^{2}(x)|\leq C(G)\left(\frac{s}{u_{1}}\right)^{c}\max p^{1}(x),p^{2}(x)\]
where we define $p^{i}(x):=\mathbb{P}^{v^{i}}(R(T(\partial B(v,u_{1})\cup\partial B(v,u_{2})))=x)$.
Define \[
\pi(x)=\mathbb{P}^{x}(R(T(\partial B(v,r_{1})\cup\partial B(v,r_{2})))=w)\]
and get\[
p^{i}=\sum_{x}p^{i}(x)\pi(x).\]
Therefore\begin{align*}
|p^{1}-p^{2}| & \leq\sum_{x}|p^{1}(x)-p^{2}(x)|\pi(x)\leq C(G)\left(\frac{s}{u_{1}}\right)^{c}\sum_{x}(p^{1}(x)+p^{2}(x))\pi(x)=\\
 & =C(G)\left(\frac{s}{u_{1}}\right)^{c}(p^{1}+p^{2})\leq C\left(\max\frac{s}{r_{1}},\frac{r_{2}}{s}\right)^{c}\max p^{1},p^{2}.\end{align*}
The case of a single point is proved identically.
\end{proof}
\begin{lem}
\label{lem:nodir2}Let $G,r_{1},s,v,v^{1},v^{2}$ be as in lemma \ref{lem:nodir},
and let $w\in\partial B(v,r_{1})$. Let $p^{i}=\mathbb{P}^{v^{i}}(R(T_{v,r_{1}})=w)$.
Then \[
|p^{1}-p^{2}|\leq C(G)\left(\frac{s}{r_{1}}\right)^{c(G)}\max p^{1},p^{2}.\]

\end{lem}
The proof is a simplified version of the proof of lemma \ref{lem:nodir}
and we shall omit it.

\begin{cor*}
The conclusion of lemma \ref{lem:nodir2} holds if $v^{1},v^{2}\in B(0,s)$
(and not on its boundary).
\end{cor*}
\begin{proof}
Apply lemma \ref{lem:nodir2} after the stopping time on $\partial B(0,s)$.
\end{proof}
\begin{lem}
\label{lem:ifhity}Let $\mathcal{G}=(G^{1},G^{2})$ be an isotropic
family, let $v\in\mathbb{R}^{d}$ and assume that for some $r>C(\mathcal{G})$
one has \[
G^{1}\cap(\overline{B(v,4r)}\setminus B(v,r))=G^{2}\cap(\overline{B(v,4r)}\setminus B(v,r)).\]
Let $w\in\overline{B(v,\frac{5}{2}r)}\setminus B(v,\frac{3}{2}r)$
and $x\in\overline{B(v,3r)}\setminus B(v,r)$ and let $R^{i}$ be
random walks on $G^{i}$ starting from $x$ (which is contained in
both $G^{i}$) and stopped on $\partial B(v,4r)$. Let \[
p^{i}=\mathbb{P}(w\in R^{i}[1,T_{v,4r}^{i}]).\]
Then $|p^{1}-p^{2}|\leq C(\mathcal{G})r^{-c(\mathcal{G})}\min\{ p^{1},p^{2}\}$.
\end{lem}
Here $G\cap(\overline{B(v,4r)}\setminus B(v,r))$ is the subgraph
of $G$ containing all vertices of $(\overline{B(v,4r)}\setminus B(v,r))$
and all edges between them; and the $=$ sign refers to equality of
graphs.

Notice that in the definition of $p^{i}$ we consider that $w$ was
hit only starting from the first step. Hence the lemma gives a non-trivial
estimate even if $x=w$ (a case that is actually important).

\begin{proof}
Let $X=\partial B(v,4r)\cup\{ w\}$. Let $s=r^{1-c_{\ref{c:RWclose}}(\mathcal{G})/2}$
where $c_{\ref{c:RWclose}}(\mathcal{G})=\min c_{\ref{c:RWclose}}(G_{1}),\linebreak[1]c_{\ref{c:RWclose}}(G_{2})$
and $c_{\ref{c:RWclose}}(G)$ is from lemma \ref{lem:rismall}. Define
stopping times as follows: $T_{0}^{i}=0$ and \begin{align*}
T_{2n+1}^{i} & :=\min\{ t\geq T_{2n}^{i}:R^{i}(t)\in\partial B(w,s)\cup X\}\\
T_{2n}^{i} & :=\min\{ t\geq T_{2n-1}^{i}:R^{i}(t)\in\partial B(w,{\textstyle \frac{1}{4}}r)\cup X\}.\end{align*}
Let $p_{n}^{i}=\mathbb{P}(R^{i}(T_{n}^{i})\not\in X).$ The core of
the lemma is showing that\begin{equation}
|p_{n}^{1}-p_{n}^{2}|\leq Knr^{-k}\max p_{n}^{1},p_{n}^{2}\label{eq:pn1pn2}\end{equation}
for any $n$ such that $Kn^{2}r^{-k}\leq\frac{1}{4}$. Here $K=K(\mathcal{G})$
and $k=k(\mathcal{G})$ are some sufficient constants (by which we
mean that $K$ is sufficiently large and $k>0$ sufficiently small)
which will be fixed later.

The proof of (\ref{eq:pn1pn2}) will be done by induction over $n$.
For $n=1$, (\ref{eq:pn1pn2}) follows from lemma \ref{lem:sball}
with $\epsilon=\frac{1}{8}$, if only $r$, $K$ and $k$ are sufficient.
Next, if $n\geq1$, define\[
\mu_{2n}^{i}(y):=\mathbb{P}(R^{i}(T_{2n}^{i})=y\,|\, T_{2n-1}^{i}\not\in X)\quad y\in\partial B(w,{\textstyle \frac{1}{4}}r)\cup\{ w\}\]
and then lemma \ref{lem:nodir} gives that the portion of the walk
between $T_{2n-1}$ and $T_{2n}$, which is the same since $G^{1}\cap B(v,\frac{1}{4}r)=G^{2}\cap B(v,\frac{1}{4}r)$,
erases most of the difference between the $\mu^{i}$-s and we have\begin{equation}
\mu_{2n}^{1}(y)\simeq\mu_{2n}^{2}(y).\label{eq:mu2n1mu2n2}\end{equation}
Denote $\alpha^{i}=\sum_{y\neq w}\mu_{2n}^{i}(y)$ so $\alpha^{1}\simeq\alpha^{2}$.
Now, $p_{2n}^{i}=p_{2n-1}^{i}\alpha^{i}$ so we can use (\ref{eq:ddm}),
note that $Cr^{-c}<\frac{1}{2}$ for $r$ sufficiently large and $K(2n-1)r^{-k}\leq\frac{2n-1}{4(2n)^{2}}<\frac{1}{8n}<\frac{1}{2}$,
and we get \begin{align}
\lefteqn{|p_{2n}^{1}-p_{2n}^{2}|\leq}\nonumber \\
 & \:\:\:\stackrel{(\ref{eq:ddm})}{\leq}\left(K(2n-1)r^{-k}+Cr^{-c}+2\left(\left(K(2n-1)r^{-k}\right)^{2}+\left(Cr^{-c}\right)^{2}\right)\right)\max p_{2n}^{1},p_{2n}^{2}\nonumber \\
 & \quad\leq\left(K\left(2n-{\textstyle \frac{1}{2}}\right)r^{-k}+Cr^{-c}\right)\max p_{2n}^{1},p_{2n}^{2}\stackrel{(*)}{\leq}K(2n)r^{-k}\max p_{2n}^{1},p_{2n}^{2}\label{eq:palpha1}\end{align}
where $(*)$ holds if only $K$ and $k$ are sufficient. Hence the
induction holds when moving from $2n-1$ to $2n$.

For the case of moving from $2n$ to $2n+1$, define, for any $y\in\partial B(w,\frac{1}{4}r$),
\[
\pi^{i}(y)=\mathbb{P}^{y}(R^{i}(T^{i}(\partial B(w,s)\cup X))\in\partial B(w,s)).\]
Then \[
p_{2n+1}^{i}=p_{2n}^{i}\sum_{y}\mu_{2n}^{i}(y)\pi^{i}(y).\]
 As in the previous part, we define $\alpha^{i}:=\sum\mu_{2n}^{i}(y)\pi^{i}(y)$.
Lemma \ref{lem:sball} shows that $\pi^{1}(y)\simeq\pi^{2}(y)$. Together
with (\ref{eq:mu2n1mu2n2}) we get, \[
|\alpha^{1}-\alpha^{2}|\stackrel{(\ref{eq:ddm2})}{\leq}\sum_{y}Cr^{-c}\max\mu_{2n}^{1}(y)\pi^{1}(y),\mu_{2n}^{2}(y)\pi(y)\stackrel{(\ref{eq:dds})}{\leq}Cr^{-c}\max\alpha^{1},\alpha^{2}\]
for $r$ sufficiently large. Thus we can repeat the calculations of
(\ref{eq:palpha1}) and get again that (\ref{eq:pn1pn2}) is preserved
if only $K$ and $k$ are sufficient. Hence the proof of (\ref{eq:pn1pn2})
is completed and we may fix the values of $K$ and $k$.

Next we need to ask how many $n$-s are actually relevant. Lemma \ref{lem:DD}
shows that \[
p_{2n+1}^{i}\leq(1-c(\mathcal{G}))p_{2n}^{i}\quad\forall n\geq1.\]
Hence the $p^{i}$ decrease exponentially. On the other hand, the
Green's function estimates (\ref{eq:arecur},\ref{eq:atrans}) shows
that \[
p^{i}>c(\mathcal{G})\begin{cases}
r^{2-d} & d\geq3\\
1/\log r & d=2\end{cases}.\]
Define therefore $N=C(\mathcal{G})\log r$ for some $C$ sufficiently
large and get\[
\sum_{N+1}^{\infty}p_{n}^{i}<r^{-1}p^{i}.\]
We can now calculate\[
p^{i}=\mathbb{P}(R^{i}(T_{1}^{i})=w)+\sum_{n=1}^{\infty}p_{2n-1}^{i}\mu_{2n}^{i}(w).\]
The first summand is non-zero only if $x\in B(w,s)$ and in this case
it is independent of $i$. Hence we may write, \begin{align*}
\lefteqn{|p^{1}-p^{2}|\leq}\\
 & \quad\stackrel{(\ref{eq:ddm2})}{\leq}\sum_{n=1}^{N}(2Knr^{-k}+Cr^{-c})\max p_{2n-1}^{1}\mu_{2n}^{1}(w),p_{2n-1}^{2}\mu_{2n}^{2}(w)+2r^{-1}\max p^{1},p^{2}\\
 & \quad\stackrel{(\ref{eq:dds})}{\leq}(4KNr^{-k}+Cr^{-c}+2r^{-1})\max p^{1},p^{2}\end{align*}
if only $r$ is sufficiently large such that we get $KN^{2}r^{-k}<\frac{1}{4}$.
This finishes the lemma.
\end{proof}
\begin{rem*}
In three dimensions this lemma may be simplified significantly, since
the probability to hit $w$ after $T_{2}$ is significantly smaller
than between $0$ and $T_{2}$, so there is no need for the induction.
\end{rem*}
\begin{lem}
\label{lem:iiexit}Let $\mathcal{G},r,v,w$ and $R^{i}$ be as in
lemma \ref{lem:ifhity} (perhaps with a different constant bounding
$r$ from below). Let $x\in\partial B(v,4r)$ and let \[
p^{i}:=\mathbb{P}^{w}(R^{i}(T_{v,4r}^{i})=x).\]
Then $p^{1}\simeq p^{2}$.
\end{lem}
\begin{proof}
The symmetry of random walk in the form (\ref{eq:symT}) shows that
\begin{equation}
p^{i}=\frac{\mathbb{P}^{x}(w\in R^{i}[1,T_{v,4r}^{i}])}{1-\mathbb{P}^{w}(w\in R^{i}[1,T_{v,4r}^{i}])}\nu\label{eq:pisym}\end{equation}
where the constant $\nu$ is the ratio of the degrees of $w$ and
$x$ and is independent of $i$. Denote the denominator by $1-a^{i}$.
Lemma \ref{lem:ifhity} shows that $a^{1}\simeq a^{2}$. The Green's
function estimates (\ref{eq:arecur},\ref{eq:atrans}) show that $1-a^{i}\geq c/\log r$
and therefore \begin{equation}
\left|\frac{1}{1-a^{1}}-\frac{1}{1-a^{2}}\right|\leq Cr^{-c}\log r\max\frac{1}{1-a^{1}},\frac{1}{1-a^{2}}\label{eq:1oa11oa2}\end{equation}
and we can drop the $\log$ factor from (\ref{eq:1oa11oa2}) and pay
in the constants only.

Next denote the nominator of (\ref{eq:pisym}) by $b^{i}$. For any
$y\in\partial B(v,3r)$ denote \[
\pi(y):=\mathbb{P}^{x}(R^{i}(T^{i}(\partial B(v,4r)\cup\partial B(v,3r)))=y)\quad\rho^{i}(y):=\mathbb{P}^{y}(w\in R^{i}[1,T_{v,4r}^{i}])\]
($\pi$ does not depend on $i$ because $G^{1}$ and $G^{2}$ are
identical on the relevant annulus). Hence we get $b^{i}=\sum_{y}\pi(y)\rho^{i}(y)$.
Lemma \ref{lem:ifhity} shows that $\rho^{1}(y)\simeq\rho^{2}(y)$
and therefore by (\ref{eq:dds}) we get $b^{1}\simeq b^{2}$. Together
with (\ref{eq:1oa11oa2}) the lemma is proved (again we use here (\ref{eq:ddm2})).
\end{proof}
\begin{cor*}
Let $\mathcal{G},r,v$ and $x$ be as in lemma \ref{lem:iiexit},
and let $v^{1},v^{2}\in\partial B(v,2r)$. Let\[
p^{i}:=\mathbb{P}^{v^{i}}(R^{i}(T_{v,4r}^{i})=x).\]
Then \[
|p^{1}-p^{2}|\leq C(\mathcal{G})\left(\frac{|v^{1}-v^{2}|}{r}\right)^{c(\mathcal{G})}\max p^{1},p^{2}.\]

\end{cor*}
\begin{proof}
If $|v^{1}-v^{2}|\leq\frac{1}{8}r$ then we can use lemma \ref{lem:nodir2}
to shows that the walk up to $B(v^{1},\frac{1}{4}r)$ erases the difference
between $v^{1}$ and $v^{2}$, and then use lemma \ref{lem:iiexit}
to show that the fact that the graphs are different has a small effect.
If $|v^{1}-v^{2}|>\frac{1}{8}r$ then for a $C$ sufficiently large,
there is nothing to prove.
\end{proof}
\begin{lem}
\label{lem:iiexit2}Lemma \ref{lem:iiexit} holds also when $w\in B(v,\frac{3}{2}r)$.
\end{lem}
Here $w$ is some point in $\mathbb{R}^{d}$ and the notation $\mathbb{P}^{w}$
refers to a random walk starting from the point of the relevant graph
closest to $w$.

\begin{proof}
Let $\{\Delta_{j}\}$ be a triangulation of $B(v,2r)$ by spherical
triangles such that\[
|\Delta_{j}|\geq cr^{-c_{\ref{c:striag}}/2}\quad\diam\Delta_{j}\leq Cr^{1-c_{\ref{c:striag}}/2d}\]
where $c_{\ref{c:striag}}$ is from lemma \ref{lem:striag} and $|\cdot|$
is the normalized volume. Let $\Delta_{j}^{*}$ be disjoint discrete
versions of the $\Delta_{j}$ covering $\partial B(v,2r)$. Let $w^{i}$
be the point of $G^{i}$ closest to $w$, let $p_{j}^{i}:=\mathbb{P}^{w^{i}}(R^{i}(T_{v,2r})\in\Delta_{j}^{*})$
and $q_{j}^{i}:=\mathbb{P}^{w^{i}}(W(S_{v,2r})\in\Delta_{j})$ the
Brownian motion analogs. Now, $q_{j}^{i}$ has a formula given from
the surface integral over the Poisson kernel \cite[II theorem 1.17]{B95}:\[
q_{j}^{i}=r^{d-2}\int_{\Delta_{j}}\frac{r^{2}-\left\Vert w^{i}\right\Vert ^{2}}{\left\Vert x-w^{i}\right\Vert ^{d}}\, dx.\]
which immediately shows, since $||w^{1}-w^{2}||\leq C(G)$ that $q_{j}^{1}\simeq q_{j}^{2}$.
Lemma \ref{lem:striag} shows that \[
|q_{j}^{i}-p_{j}^{i}|\leq Cr^{-c_{\ref{c:striag}}}\leq Cr^{-c_{\ref{c:striag}}/2}|\Delta_{j}|\leq Cr^{-c_{\ref{c:striag}}/2}q_{j}^{i}\]
So also $p_{j}^{1}\simeq p_{j}^{2}$.

Next, let $y^{1},y^{2}\in\Delta_{j}^{*}$. Denote $\pi^{i}(y)=\mathbb{P}^{y}(R^{i}(T_{v,4r}^{i})=x)$.
Because $\diam\Delta_{j}^{*}\leq Cr^{1-c}$ we can use the corollary
after lemma \ref{lem:iiexit} to get $\pi^{1}(x^{1})\simeq\pi^{2}(x^{2})$.
This finishes the lemma, since the probabilities to hit a given $\Delta_{j}^{*}$
are $Cr^{-c}$ similar, and the point in which you hit is unimportant
up to $Cr^{-c}$ error.
\end{proof}

\subsection{\label{sub:Definition-II}Definition}

Let $G^{1}$ and $G^{2}$ be two $d$-dimensional isotropic graphs,
and let $\alpha>0$. We say that $G^{1}$ and $G^{2}$ have an $\alpha$-\textbf{isotropic
interpolation} if the following holds. Let $L$ and $M$, $M\geq1$
integer, satisfy \begin{equation}
L>C(G^{i})\quad M\leq L^{\alpha}.\label{eq:iirMreq}\end{equation}
We assume $C(G^{i})\geq2$ always. Let $\xi\in\{1,2\}^{M^{d}}$ be
any configuration. Then there exist graphs $G(L,M,\xi)$ such that:

\begin{enumerate}
\item \label{enu:isinI}If all the coordinates of $\xi$ are $i$ then $G(L,M,\xi)\equiv G^{i}$. 
\item \label{enu:isinII}If $I:=[a_{1},b_{1}]\times\dotsb\times[a_{d},b_{d}]\subset\{0,\dotsc,M-1\}^{d}$
is some box in the configuration space , and if $\left.\xi_{1}\right|_{I}\equiv\left.\xi_{2}\right|_{I}$
then $G(L,M,\xi_{1})$ is equal to $G(L,M,\xi_{2})$ on a corresponding
box in $\mathbb{R}^{d}$ i.e.\begin{gather}
G(L,M,\xi_{1})\cap J=G(L,M,\xi_{2})\cap J,\nonumber \\
J:=[L(a_{1}+{\textstyle \frac{1}{3}}),L(b_{1}+{\textstyle \frac{2}{3}})]\times\dotsb\times[L(a_{d}+{\textstyle \frac{1}{3}}),L(b_{d}+{\textstyle \frac{2}{3}})].\label{eq:influ}\end{gather}

\item $G(L,M,\xi)$ is isotropic with the isotropic structure constants
bounded independently of $L$, $M$ and $\xi$.
\end{enumerate}
Notice that this definition gives special importance to the point
zero. This is just for convenience --- in practice, this fact will
have no significance.

The core of this paper is the proof of the following theorem.

\begin{thm}
\label{thm:interp}Let $G^{1}$ and $G^{2}$ be two $d$-dimensional
graphs with an $\alpha$-isotropic interpolation, $d=2,3$. Let $\mathcal{D}\subset[\frac{1}{4},\frac{3}{4}]^{d}$
be an open polyhedron and let $\mathcal{E}$ be some open set. Let
$a\in\mathcal{E}\cap\mathcal{D}$ be some point. Let $s>0$ be some
number, and let $R^{i}$ be random walks on $G^{i}$ starting from
$sa$ and stopped when hitting $\partial(s\mathcal{D})$. Then \[
\mathbb{P}(\LE(R^{i})\subset s(\mathcal{E}+B(0,Cs^{-c})))>\mathbb{P}(\LE(R^{3-i})\subset s\mathcal{E})-Cs^{-c}.\]

\end{thm}
Remember that {}``a random walk starting from $sa$'', means that
it starts from the point of $G^{i}$ closest to $sa$; that $\mathcal{E}+B(0,Cs^{-c})$
refers to the set of points within distance $Cs^{-c}$ from $\mathcal{E}$;
and that an {}``open polyhedron'' is any open set whose boundary
is made of non-degenerate linear polyhedra of dimension $d-1$, and
that we do not require that the boundary of the polyhedron be connected,
but we do not allow slits.

A comment is due on the use of constants here. They all depend on
$a$, $\mathcal{D}$ and on the graphs $G^{1}$ and $G^{2}$ (in fact
they don't really depend on $a$ but we will have no use for this
fact). Like in previous chapters, they depend only on $\alpha$ and
the global bound for the isotropic structure constants over all of
$G(L,M,\xi)$ and not on other properties of $G^{i}$. However, there
is no need to continue to point this fact out --- we only did so in
previous chapters in order to be able to analyze walks on $G(L,M,\xi)$
simultaneously, and we will not have families of isotropic interpolations
in the future.

\subsection{Proof of theorem \refs{thm:interp}}

The argumentation in this section is very similar to that of \cite[section 4]{K},
so we will be brief.

\begin{lem}
\label{lem:xizhit}Let $G^{i},a,s$ be as in theorem \ref{thm:interp}.
Assume $s\leq LM$ with $L$ and $M$ satisfying (\ref{eq:iirMreq})
and $L$ sufficiently large (in addition to the restriction of (\ref{eq:iirMreq})).
Let $\xi^{1}$,$\xi^{2}$ be two configurations which differ only
in one point $z$, so in particular $H:=G(L,M,\xi^{i})\cap([\frac{1}{3}L,s-\frac{1}{3}L]^{d}\setminus B(Lz,16L))$
does not depend on $i$. Let $\mathcal{B}\subset H$ be any subset
such that $a$ is in a finite component of $G(L,M,\xi^{i})\setminus\mathcal{B}$
(think about \emph{$\mathcal{B}$} as the boundary of some $\mathcal{D}$,
$a\in\mathcal{D}$). Let $R^{i}$ be random walks on $G(L,M,\xi^{i})$
starting from $sa$. Let $b\in\mathcal{B}$ and define\[
p^{i}=\mathbb{P}(R^{i}(T^{i}(\mathcal{B}))=b).\]
Then $|p^{1}-p^{2}|\leq C(\mathcal{G})L^{-c(\mathcal{G})}\max p^{1},p^{2}$.
\end{lem}
\begin{proof}
[Proof sketch]Define stopping times $T_{j}^{i}$ on $\partial B(Lz,3L)$
and $\partial B(Lz,12L)$ alternatively. The graphs $G(L,M,\xi^{i})$
are identical outside $B(Lz,3L)$ hence lemma \ref{lem:iiexit2} shows
that the transition probabilities, up to $CL^{-c}$ do not depend
on $i$. Further, the probability to reach $T_{j}^{i}$ without hitting
$\mathcal{B}$ drops like $e^{-cj}$ in three dimensions and like
$e^{-cj/\log M}$ in two dimensions. Harnack's inequality on $B(Lz,16L)$
shows that the probabilities to hit $b$ after $T_{j}$ conditioned
on not hitting it before are, up to a constant, independent of $j$.
Hence these $CL^{-c}$ errors do not accumulate to more than $CL^{-c}\log M\leq CL^{-c}$.
See \cite[lemma 16]{K} for a detailed argument.
\end{proof}
\begin{lem}
\label{lem:xizLE}Let $G^{i},a,s,L,M,\xi^{i},z,\mathcal{B}$ and $b$
be as in lemma \ref{lem:xizhit}. Let $R^{i}$ be random walks on
$G(L,M,\xi^{i})$ starting from $sa$ and conditioned to hit $\mathcal{B}$
in $b$. Let $\zeta^{i}$ be the segment of $\LE(R^{i})$ until first
hitting $\partial B(Lz,16L)$ or all of $\LE(R^{i})$ if $\LE(R^{i})\cap\partial B(Lz,16L)=\emptyset$.
Then\[
\sum_{\gamma}|\mathbb{P}(\zeta^{1}=\gamma)-\mathbb{P}(\zeta^{2}=\gamma)|\leq C(\mathcal{G})L^{-c(\mathcal{G})}.\]
where the sum is over all simple paths $\gamma$ from $sa$ to $\partial B(Lz,16L)\cup\partial(s\mathcal{D})$.
\end{lem}
\begin{proof}
[Proof sketch]The crucial point here is that $\zeta^{i}$ depends
only on what happens outside $B(Lz,16L)$ (quite unlike the other
portion of $\LE(R^{i})$). Therefore the same argument as in the previous
lemma works here. Conditioning on all the entry exit points from $\partial B(Lz,3L)$
and $\partial B(Lz,12L)$ the probabilities are identical, and in
average in $\gamma$ it is enough to consider $\log^{2}L$ such points.
See \cite[lemma 17]{K} for a detailed argument.
\end{proof}
The proof of theorem \ref{thm:interp} is also detailed in \cite{K}
where it is called the {}``main lemma''. Since this is a crucial
part of the argument, I prefer to bring it here in full.

Since the theorem is symmetric in $i$, fix it to be $1$. Let $\epsilon>0$
be some constant (depending on $\mathcal{G}$) which will be fixed
later. We define $L$ by \begin{equation}
34L=s^{1-\epsilon}\label{eq:defepst3}\end{equation}
and $M=\left\lceil s/L\right\rceil $. We assume $L$ and $M$ satisfy
the requirements (\ref{eq:iirMreq}) of the isotropic interpolations,
which will hold if $\epsilon<\min\frac{1}{2}\alpha,\frac{1}{2}$ and
$L$ is large enough. We will also assume $L$ is sufficiently large
such that lemmas \ref{lem:xizhit} and \ref{lem:xizLE} hold, and
also that all edges of all graphs $G(L,M,\xi)$ are shorter than $L$
and every ball of radius $L$ in $\mathbb{R}^{d}$ contains at least
one point from every $G(L,M,\xi$). All these requirements translate
to $s>C(\epsilon,\mathcal{G})$.

Let $\delta=\delta(\epsilon,\mathcal{G})$ be the quantity given by
theorem \ref{thm:QL} (page \pageref{thm:QL}) for the $\epsilon$
from (\ref{eq:defepst3}), for all the graphs $G(L,M,\xi)$ simultaneously.
In other words, we have, if $s>C(\epsilon,\mathcal{G})$, \[
\mathbb{E}^{sa}\QL(34L,s^{1-\delta},\LE(R[0,T(\partial\mathcal{D})]))\leq C(\mathcal{G})s^{-\delta}\]
which holds for $R$ a random walk on any $G(L,M,\xi)$. With this
$\delta$, define {}``bad'' subsets of $\{ x\in\{0,\dotsc,M-1\}^{d}:x/M\in\mathcal{D}\}$
as follows:\begin{align}
\Phi & :=\left\{ x:d(Lx,\partial_{\textrm{cont}}s\mathcal{E})\leq s^{1-\delta}+18L\right\} \label{eq:defPHI}\\
\Psi & :=\left\{ x:d(Lx,\partial_{\textrm{cont}}s\mathcal{D})\leq17L\right\} .\label{eq:defPSI}\\
\Theta & :=\left\{ x:d(Lx,sa)\leq17L\right\} \nonumber \end{align}
It should be noticed that any $x\not\in\Psi$ satisfies $d(Lx,\partial_{G(L,M,\xi)}(s\mathcal{D}))>17L$
and any $x\not\in\Phi$ satisfies $d(Lx,\partial_{G(L,M,\xi)}(s\mathcal{E}))>s^{1-\delta}+17L$,
both for any configuration $\xi$.

\begin{lem}
\label{lem:0_to_Y}With the definitions above, let $Y\subset\{0,\dotsc,M-1\}^{d}\setminus(\Phi\cup\Psi\cup\Theta)$.
Let $\xi^{1}$ and $\xi^{2}$ be two configurations such that $\xi^{i}|_{Y}\equiv i$
but $\xi^{1}|_{Y^{c}}\equiv\xi^{2}|_{Y^{c}}$. Let $R^{i}$ be random
walks on $G(L,M,\xi^{i})$ starting from $sa$ and stopped on $\partial\mathcal{D}$.
Let \[
p^{i}:=\mathbb{P}(\LE(R^{i})\subset s\mathcal{E}).\]
Then\[
|p^{1}-p^{2}|\leq C(\mathcal{G})s^{-\delta}\log\# Y+C(\mathcal{G})L^{-c(\mathcal{G})}\# Y.\]

\end{lem}
\begin{proof}
For every $0\leq k\leq\# Y$, let $D_{k}$ be a random subset of $Y$
of size $k$, let $\xi_{k}$ be the configuration which is 2 on $D_{k}$,
1 on $Y\setminus D_{k}$ and identical to $\xi^{1}$ outside $Y$.
Define $H_{k}:=G(L,M,\xi_{k})$. Let $a_{k}$ be the point of $H_{k}$
closest to $sa$. Let $S_{k}$ be a random walk on $H_{k}$ starting
from $a_{k}$ and stopped on $B_{k}:=\partial_{H_{k}}s\mathcal{D}$.
Let \[
p_{k}:=\mathbb{P}(\LE(S_{k})\subset s\mathcal{E})\]
where $\mathbb{P}$ here is over both the walk and the randomness
of the graph (notice that $p^{1}=p_{0}$ and $p^{2}=p_{\# Y}$). The
lemma will be proved once we show that

\begin{equation}
|p_{k}-p_{k+1}|\leq Cs^{-\delta}\left(\frac{1}{k+1}+\frac{1}{\# Y-k}\right)+CL^{-c}.\label{step1}\end{equation}

For this purpose, couple $D_{k}$ and $D_{k+1}$ such that $D_{k}\subset D_{k+1}$.
Let $\Delta_{k}\subset\Delta_{k+1}$ be subsets of $Y$ of sizes $k$
and $k+1$ and let $z=\Delta_{k+1}\setminus\Delta_{k}$. For most
of the rest of the lemma, we condition on the event (denote it by
$\mathcal{X}$) that $D_{k}=\Delta_{k}$ and $D_{k+1}=\Delta_{k+1}$.
Let $Z=B(Lz,16L)$. We construct $\LE(S_{k})$ as follows: 
\begin{itemize}
\item let $b_{k}$ be a random point on $B_{k}$ chosen with the hitting
probabilities of $S_{k}$. 
\item Let $\check{S}_{k}$ be a random walk from $a_{k}$ to $B_{k}$ conditioned
to hit $b_{k}$.
\item Let $\check{\gamma}_{k}$ be a random simple path from $a_{k}$ to
$\partial Z\cup\{ b_{k}\}$, which has the same distribution as the
segment of $\LE(\check{S}_{k})$ until $\partial Z$ (including the
first point in $\partial Z$), or all of $\LE(\check{S}_{k})$ if
$\LE(\check{S}_{k})\cap\partial Z=\emptyset$ (notice that $a_{k}\not\in Z$
from the requirement $Y\cap\Theta=\emptyset$).
\item Let $c_{k}$ be the point where $\check{\gamma}_{k}$ hits $\partial Z$
if it does.
\item Let $\check{T}_{k}$ be a random walk on $H_{k}$ starting from $b_{k}$
and conditioned to hit $B_{k}\cup\check{\gamma}_{k}$ in $c_{k}$,
or $\emptyset$ if $\check{\gamma}_{k}$ never hits $\partial Z$.
\item Let $\gamma_{k}=\check{\gamma}_{k}\cup\LE(\check{T}_{k})$.
\end{itemize}
An easy application of lemma \ref{lem:condLE_sym} (symmetry of conditioned
loop-erased random walk) shows that $\gamma_{k}\sim\LE(S_{k})$. Lemma
\ref{lem:xizhit} shows that, \begin{gather}
\sum_{b}|q_{k}^{1}-q_{k+1}^{1}|\leq CL^{-c}\label{bk}\\
q_{k}^{1}(b):=\mathbb{P}(b_{k}=b\,|\,\mathcal{X}).\nonumber \end{gather}
Next we use lemma \ref{lem:xizLE} for the random walk on $H_{k}$
starting from $a_{k}$, stopped on $B_{k}$, and conditioned to hit
$b_{k}$. The definition of $\Psi$ (\ref{eq:defPSI}) ensures the
condition $B_{k}\cap Z=\emptyset$ required by lemma \ref{lem:xizLE}.
This shows that \begin{gather}
\sum_{\gamma}|q_{k}^{2}-q_{k+1}^{2}|\leq CL^{-c}\qquad\forall b\in B_{k}\label{gamma_k}\\
q_{k}^{2}(b,\gamma):=\mathbb{P}(\check{\gamma}_{k}=\gamma\,|\,\mathcal{X},\, b_{k}=b).\nonumber \end{gather}
Thirdly, we again use lemma \ref{lem:xizLE}, this time for a random
walk starting from $b_{k}$, stopped on $B_{k}\cup\check{\gamma}_{k}$
and conditioned to hit $\check{c}_{k}$ to show that, when $\check{\gamma}_{k}'$
is the portion of $\LE(\check{T}_{k})$ up to $Z$,\begin{gather}
|q_{k}^{3}-q_{k+1}^{3}|\leq CL^{-c}\qquad\forall b,\gamma\label{gamma_k_prime}\\
q_{k}^{3}(b,\gamma):=\mathbb{P}(\check{\gamma}_{k}'\subset s\mathcal{E}\,|\,\mathcal{X},\, b_{k}=b,\,\check{\gamma}_{k}=\gamma).\nonumber \end{gather}
Summing (\ref{bk}), (\ref{gamma_k}) and (\ref{gamma_k_prime}) gives\begin{eqnarray}
\lefteqn{|\mathbb{P}((\check{\gamma}_{k}\cup\check{\gamma}_{k}')\subset s\mathcal{E}\,|\,\mathcal{X})-\mathbb{P}((\check{\gamma}_{k+1}\cup\check{\gamma}_{k+1}')\subset s\mathcal{E}\,|\,\mathcal{X})|=}\nonumber \\
 & \qquad & =\left|\sum_{b,\gamma\subset s\mathcal{E}}q_{k}^{1}q_{k}^{2}q_{k}^{3}-q_{k+1}^{1}q_{k+1}^{2}q_{k+1}^{3}\right|\leq\nonumber \\
 &  & \leq CL^{-c}+\left|\sum_{b,\gamma\subset s\mathcal{E}}(q_{k}^{1}q_{k}^{2}-q_{k+1}^{1}q_{k+1}^{2})q_{k}^{3}\right|\leq\dotsb\leq\nonumber \\
 &  & \leq CL^{-c}.\label{no_quasi_loop}\end{eqnarray}
 In other words, we have proved that the probabilities (for $k$ and
$k+1$) that both segments of $\LE(S_{k})$, leading up to $Z$ and
from $Z$ to $B_{k}$ to be in $s\mathcal{E}$ are close. Thus the
only case we have not covered is that $\check{\gamma}_{k}\cup\check{\gamma}_{k}'\subset s\mathcal{E}$
but $\LE(S_{k})\not\subset s\mathcal{E}$. But $Lz$ is far from $\partial s\mathcal{E}$
(because $Y\cap\Phi=\emptyset$ and the definition of $\Phi$ (\ref{eq:defPHI}))
so we get a quasi loop near $Lz$, namely\[
Lz\in\mathcal{QL}(17L,s^{1-\delta},\LE(S_{k})).\]
 Denote therefore\[
q_{k}^{4}(x):=\mathbb{P}(Lx\in\mathcal{QL}(17L,s^{1-\delta},\LE(S_{k}))\,|\,\mathcal{X})\]
and then write (\ref{no_quasi_loop}) as\[
|\mathbb{P}(\LE(S_{k})\subset s\mathcal{E}\,|\,\mathcal{X})-\mathbb{P}(\LE(S_{k+1})\subset s\mathcal{E}\,|\,\mathcal{X})|\leq CL^{-c}+q_{k}^{4}(z)+q_{k+1}^{4}(z).\]
We now integrate over $\mathcal{X}$. We get\begin{equation}
|p_{k}-p_{k+1}|\leq CL^{-c}+\mathbb{E}q_{k}^{4}(z)+\mathbb{E}q_{k+1}^{4}(z).\label{eq:quasi_loop_left}\end{equation}

The estimate of (\ref{eq:quasi_loop_left}) is where the random choice
of the sets $D_{k}$ plays its part. Theorem \ref{thm:QL} gives us
that \[
\sum_{x\in Y\setminus D_{k}}q_{k}^{4}(x)\leq\sum_{x\in Y}q_{k}^{4}(x)\leq C\mathbb{E}\QL(34L,s^{1-\delta},\LE(R[0,T(\partial\mathcal{D})]))\leq Cs^{-\delta}.\]
Now, if we think about the coupling of $D_{k}$ and $D_{k+1}$ as
{}``$D_{k+1}$ is the addition of a random $z$ to $D_{k}$'', then
$z$ is is obviously independent of the walk on $G(r,M,\xi_{k})$
so we have \[
\mathbb{E}q_{k}(z)\leq\frac{Cs^{-\delta}}{\#(Y\setminus D_{k})}=\frac{Cs^{-\delta}}{\# Y-k}.\]
For $q_{k+1}(z)$ we similarly think about $D_{k}$ as the removal
of a random $z$ from $D_{k+1}$ and get\[
\mathbb{E}q_{k+1}(z)\leq\frac{Cs^{-\delta}}{k+1}\]
and the lemma is proved. 
\end{proof}
Continuing the proof of the theorem, we first apply lemma \ref{lem:0_to_Y}
to the configurations $\xi^{1}\equiv1$ and \[
\xi^{2}=\begin{cases}
1 & \Phi\cup\Psi\cup\Theta\\
2 & \textrm{otherwise}\end{cases}.\]
Denoting by \[
p^{i}:=\mathbb{P}_{G(L,M,\xi^{i})}^{sa}(\LE(R)\subset s\mathcal{E})\]
($R$ here is $R[0,T(\partial\mathcal{D})]$ and will stay so for
a while). We get \begin{equation}
|p^{1}-p^{2}|\leq Cs^{-\delta}\log s+CL^{-c}M^{d}.\label{eq:step1}\end{equation}

Next, we wish to remove $\Phi$. For this purpose we define \begin{eqnarray*}
p^{3} & := & \mathbb{P}_{G(L,M,\xi^{2})}^{sa}(\LE(R)\subset s\mathcal{E}_{2})\\
\mathcal{E}_{2} & := & \mathcal{E}+B(0,2s^{-\delta}+37L/s)\end{eqnarray*}
and get  $p^{3}>p^{2}$. The definition of $\mathcal{E}_{2}$ gives
us that\[
d(Lx,\partial_{\textrm{cont}}s\mathcal{E}_{2})>s^{1-\delta}+18L\quad\forall x\in\Phi.\]
so we can use lemma \ref{lem:0_to_Y} with $\xi^{2}$, \[
\xi^{3}=\begin{cases}
1 & \Psi\cup\Theta\\
2 & \textrm{otherwise}\end{cases}\]
and the domain $\mathcal{E}_{2}$ to get \begin{equation}
|p^{3}-p^{4}|\leq Cs^{-\delta}\log s+CL^{-c}M^{d}.\label{eq:step2}\end{equation}
where $p^{4}=\mathbb{P}_{G(L,M,\xi^{3})}^{sa}(\LE(R)\subset s\mathcal{E}_{2})$.

The third step is to get rid of $\Theta$. For this purpose we define
\[
\mathcal{E}_{3}:=\mathcal{E}_{2}+B(0,s^{-\epsilon/2})\]
and $p^{5}=\mathbb{P}_{G(L,M,\xi^{3})}^{sa}(\LE(R)\subset s\mathcal{E}_{3})$
so that $p^{5}>p^{4}$. Find some point $b\in\mathcal{E}_{3}$ with
$|a-b|=35L/s$ (which can always be done if $s$ is sufficiently large)
and use lemma \ref{lem:contsp} to find that\begin{multline}
|p^{5}-p^{6}|\leq C\left(\frac{35L}{sd(a,\partial_{\textrm{cont}}(\mathcal{E}_{3}\cap\mathcal{D})}\right)^{c}\leq\\
\leq C\left(\frac{35L}{\min s^{1-\epsilon/2},sd(a,\partial_{\textrm{cont}}\mathcal{D})}\right)^{c}=C(\epsilon,a,\mathcal{D},\mathcal{G})s^{-c(\epsilon,a,\mathcal{D},G)}.\label{eq:step3a}\end{multline}
where $p^{6}=\mathbb{P}_{G(L,M,\xi^{3})}^{sb}(\LE(R)\subset s\mathcal{E}_{3})$.
We can now apply lemma \ref{lem:0_to_Y} with $\xi^{3}$, \[
\xi^{4}=\begin{cases}
1 & \Psi\\
2 & \textrm{otherwise,}\end{cases}\]
 the domain $\mathcal{E}_{3}$ and the point $b$ to get \begin{equation}
|p^{6}-p^{7}|\leq Cs^{-\delta}+CL^{-c}.\label{eq:step3b}\end{equation}
where $p^{7}=\mathbb{P}_{G(L,M,\xi^{4})}^{sb}(\LE(R)\subset s\mathcal{E}_{3})$.
We apply lemma \ref{lem:contsp} again to return to $a$: we get $|p^{7}-p^{8}|\leq Cs^{-c}$
where $p^{8}=\mathbb{P}_{G(L,M,\xi^{4})}^{sa}(\LE(R)\subset s\mathcal{E}_{3})$.

\label{page:begstep4}Finally, we need to get rid of $\Psi$. Let
$\mathcal{D}_{2}$ be a shrinking of $\mathcal{D}$ by $21L/s$ i.e.\[
\mathcal{D}_{2}:=\{ x:B(x,21L/s)\subset\mathcal{D}\}\]
so we have $L\Psi\cap s\mathcal{D}_{2}=\emptyset$. Let $\delta_{2}$
be the value given by theorem \ref{thm:QL} for $\epsilon/2$, $\epsilon$
from (\ref{eq:defepst3}), for all the graphs $G(L,M,\xi)$ simultaneously,
and define\[
\mathcal{E}_{4}:=\mathcal{E}_{3}+B(0,s^{-\delta_{2}}+s^{-\epsilon/2}).\]
Let \[
p^{9}=\mathbb{P}_{G(L,M,\xi^{4})}^{sa}(\LE(R[0,T(\partial\mathcal{D}_{2})])\subset s\mathcal{E}_{4})\]

\begin{lem}
\label{sublem:step4}\begin{equation}
p^{9}>p^{8}-Cs^{-c}\label{eq:step4}\end{equation}
where $C$ and $c$ may depend on $\mathcal{D}$ and $\epsilon$ in
addition to $\mathcal{G}$.
\end{lem}
\begin{proof}
$\mathbb{R}^{d}\setminus\mathcal{D}$ has a finite number of connected
components, $\{ Q_{i}\}$ which are all polyhedra. For each $Q_{i}$
we may use lemma \ref{lem:polyh} to get that, for every $1<r_{1}<r_{2}<s$,
every $\xi$ and every $v\in G(L,M,\xi)$, $d(v,sQ_{i})\leq r_{1}$,
\begin{equation}
\mathbb{P}(T_{v,r_{2}}<T(sQ_{i}))\leq C(\mathcal{G},Q_{i})\left(\frac{r_{1}}{r_{2}}\right)^{c(\mathcal{G},Q_{i})}.\label{eq:defKk}\end{equation}
 Let $K=\max_{i}C(\mathcal{G},Q_{i})$ and $k=\min_{i}c(\mathcal{G},Q_{i})$.

Now $p^{8}$ and $p^{9}$ measure walks on the same graph stopped
at $\partial s\mathcal{D}$ and $\partial s\mathcal{D}_{2}$ respectively.
Therefore we may couple these walks so that the first is a continuation
of the second. In other words, define $t_{1}>t_{2}$ be the stopping
times of $R$ on $\partial s\mathcal{D}$ and $\partial s\mathcal{D}_{2}$
(define $t_{2}=0$ if $a\not\in s\mathcal{D}_{2}$) then the question
reduces to an estimate of \begin{equation}
\mathbb{P}\left(\left\{ \LE(R[0,t_{2}])\not\subset s\mathcal{E}_{4}\right\} \cap\left\{ \LE(R[0,t_{1}])\subset s\mathcal{E}_{3}\right\} \right).\label{eq:St2_large_St1_small}\end{equation}
Now, the definition of $\mathcal{D}_{2}$ gives that the distance
of $R(t_{2})$ from the closest connected component of $\mathbb{R}^{d}\setminus(s\mathcal{D})$
is $\leq21L$. From the definition of $K$ and $k$ we get\begin{equation}
\mathbb{P}(R[t_{2},t_{1}]\text{ exits }B(R(t_{2}),{\scriptstyle \frac{1}{2}}s^{1-\epsilon/2}))\leq K\left(\frac{42L}{s^{1-\epsilon/2}}\right)^{k}=C(\mathcal{D},\mathcal{G})s^{-\epsilon c(\mathcal{D},\mathcal{G})}.\label{eq:r_large}\end{equation}
On the other hand, if $R[t_{2},t_{1}]\subset B(R(t_{2}),{\scriptstyle \frac{1}{2}}s^{1-\epsilon/2})$
and in addition the event of (\ref{eq:St2_large_St1_small}) holds
then we can conclude that $R(t_{2})\in\mathcal{QL}({\scriptstyle \frac{1}{2}}s^{1-\epsilon/2},s^{1-\delta_{2}},\LE(S_{3}[0,t_{2}]))$
which means that \[
\QL(s^{1-\epsilon/2},s^{1-\delta_{2}},\LE(S_{3}[0,t_{2}]))\geq1.\]
By theorem \ref{thm:QL} and Markov's inequality, the probability
for that is $\leq Cs^{-\delta_{2}}$. This ends the lemma.
\end{proof}
The definition of isotropic interpolation shows $\overline{G(L,M,\xi^{4})\cap\mathcal{D}_{2}}=\overline{G^{2}\cap\mathcal{D}_{2}}$.
Therefore we may define $R$ to be a random walk on $G^{2}$ starting
from $sa$ and get that $p^{9}=\mathbb{P}(\LE(R[0,T(\partial s\mathcal{D}_{2})])\subset s\mathcal{E}_{4})$.
We only need to return from $\mathcal{D}_{2}$ to $\mathcal{D}$ so
write\[
p^{10}:=\mathbb{P}(\LE(R[0,T(\partial s\mathcal{D})])\subset s\mathcal{E}_{5})\]
where\[
\mathcal{E}_{5}:=\mathcal{E}_{4}+B(0,s^{-\epsilon/2}).\]

\begin{lem}
With the definitions above\begin{equation}
p^{10}>p^{9}-Cs^{-c}.\label{eq:step5}\end{equation}
where again $C$ and $c$ may depend on $\mathcal{D}$ and $\epsilon$
in addition to $\mathcal{G}$.
\end{lem}
\begin{proof}
As in lemma \ref{sublem:step4}, it is enough to show that\[
\mathbb{P}(\LE(R[0,t_{1}])\not\subset s\mathcal{E}_{5},\,\LE(R[0,t_{2}])\subset s\mathcal{E}_{4})\leq Cs^{-c}.\]
with the same $t_{1}$ and $t_{2}$. Unlike in sublemma \ref{sublem:step4},
this requires no recourse to theorem \ref{thm:QL} but rather follows
directly from lemma \ref{lem:polyh} since this event implies that
$R[t_{2},t_{1}]\not\subset R(t_{2})+B(0,s^{1-\epsilon/2})$ whose
probability can be bounded by $K(21L/\linebreak[4]s^{1-\epsilon/2})^{k}=Cs^{-c}$
with the same $K$ and $k$ as in lemma \ref{sublem:step4}.\label{page:endstep5}
\end{proof}
Summing up (\ref{eq:step1}), (\ref{eq:step2}), (\ref{eq:step3a}),
(\ref{eq:step3b}), (\ref{eq:step4}) and (\ref{eq:step5}) we get\newc{c:Lpow}
for some $c_{\ref{c:Lpow}}(\mathcal{G})$, \begin{equation}
p^{10}>p^{1}-C(\epsilon,a,\mathcal{D},\mathcal{G})s^{-c(\epsilon,a,\mathcal{D},\mathcal{G})}-C(\mathcal{G})M^{d}L^{-c_{\ref{c:Lpow}}}.\label{just_choose_M}\end{equation}
The only thing left now is to choose $\epsilon$. For $\epsilon<c_{\ref{c:Lpow}}/2d$,
we will have $M^{d}L^{-c_{\ref{c:Lpow}}}\leq Cs^{-c_{\ref{c:Lpow}}/2}$.
This finishes the proof of theorem \ref{thm:interp}.\qed

\subsection{The limit process}

In this section we derive consequences from theorem \ref{thm:interp}:
we will prove the following theorems:

\begin{thm}
\label{thm:scallim}Let $G$ be a $d$-dimensional isotropic graph
with an isotropic interpolation to $2G$. Let $\mathcal{D}\subset\mathbb{R}^{d}$
be a polyhedron and let $a\in\mathcal{D}$. Let $\mathbb{P}_{n}$
be the distribution of the loop-erasure of a random walk on $2^{n}\mathcal{D}\cap G$
starting from $2^{n}a$ and stopped when hitting $\partial2^{n}\mathcal{D}$,
multiplied by $2^{-n}$. Then $\mathbb{P}_{n}$ converge in the space
$\mathcal{M}(\mathcal{H}(\overline{\mathcal{D}}))$.
\end{thm}
Recall that $\mathcal{H}(\mathcal{X})$ is the space of compact subsets
of $\mathcal{X}$ with the Hausdorff metric, and $\mathcal{M}(\mathcal{X})$
is the space of measures on $\mathcal{X}$ with the topology of weak
convergence. $2G$ is the graph gotten by stretching $G$ by $2$
uniformly and possibly multiplying all weights by a constant $\alpha$
(this last action does not change the process of course). In other
words, the theorem holds if for any $\alpha$ there is an isotropic
interpolation between $G$ and $2G$. Strangely enough, in 3 dimensions
$\alpha$ will usually not be 1.

The limit of $\mathbb{P}_{n}$ is called the {}``scaling limit''
(of the loop-erased random walk) and $G$ is said to {}``have a scaling
limit''.

\begin{thm}
\label{thm:univ}Let $G$ and $H$ be two $d$-dimensional isotropic
graphs with an isotropic interpolation, and assume $G$ has a scaling
limit. Then $H$ has a scaling limit and they are identical.
\end{thm}
We start with a proof of theorem \ref{thm:scallim}. In the following
we shall abuse notations by denoting, for any subset $A$ in a metric
space $X$ and any $\epsilon>0$,\begin{align*}
A+B(\epsilon) & :=\{ x\in X:d(x,A)<\epsilon\}\\
A+\overline{B}(\epsilon) & :=\{ x\in X:\exists a\in A,\, d(x,a)\leq\epsilon\}.\end{align*}
We shall also need the following notation: For a set $E\subset\overline{\mathcal{D}}$
relatively open, define subsets of $\mathcal{H}(\overline{\mathcal{D}})$,\begin{equation}
\mathcal{S}(E):=\{ K\subset E\},\quad\mathcal{I}(E):=\{ K\cap E\neq\emptyset\}.\label{eq:defSIE}\end{equation}

\begin{lem}
\label{lem:montrik}Let $G,a$ and $\mathcal{D}$ be as in theorem
\ref{thm:scallim} and let $\mathcal{E}_{1},\dotsc,\mathcal{E}_{k}\subset\mathcal{D}$
be open or closed. Then for almost every $\epsilon>0$ the limit\[
\lim_{n\to\infty}\mathbb{P}_{n}\Big(\mathcal{S}\Big(\bigcap_{i=1}^{k}\mathcal{E}_{i}+\overline{B(\epsilon)}\Big)\Big)\]
exists.
\end{lem}
Here and below {}``almost every'' can be replaced with {}``except
for a countable set''.

\begin{proof}
We may replace $\overline{B(\epsilon)}$ with $B(\epsilon)$: for
$\mathcal{E}_{i}$ open $\mathcal{E}_{i}+\overline{B(\epsilon)}=\mathcal{E}_{i}+B(\epsilon)$
and for $\mathcal{E}_{i}$ closed $\mathcal{E}_{i}+\overline{B(\epsilon)}=\overline{\mathcal{E}_{i}+B(\epsilon)}$
and there is only a countable number of $\epsilon$-s for which the
closing operation affects the problem at all. Further, we may also
assume all $\mathcal{E}_{i}$ are open, since replacing each with
$\mathcal{E}_{i}+B(\delta)$ and then taking $\delta$ to zero will
prove the general case. 

Denote by $\LE$ the expression $\LE(R[0,T(\partial2^{n}\mathcal{D})])$
when $n$ is assumed to be clear from the context. Denote\[
\overline{p}(\mathcal{E}):=\varlimsup_{n\to\infty}\mathbb{P}_{n}(\mathcal{S}(\mathcal{E}))\]
and similarly $\underline{p}(\mathcal{E})$. By theorem \ref{thm:interp}
we have\begin{align*}
\mathbb{P}_{G}^{2^{n}a}(\LE\subset2^{n}\cap\mathcal{E}_{i}) & =\mathbb{P}_{2G}^{2^{n+1}a}(\LE\subset2^{n+1}\cap\mathcal{E}_{i})<\\
 & <\mathbb{P}_{G}^{2^{n+1}a}\big(\LE\subset2^{n+1}((\cap\mathcal{E}_{i})+B(C2^{-nc}))\big)+C2^{-nc}\\
\intertext{\textrm{and inductively for any $m$,}} & <\mathbb{P}_{G}^{2^{n+m}a}\big(\LE\subset2^{n+m}((\cap\mathcal{E}_{i})+B(C(2^{-nc}+\dotsb\\
 & \qquad+2^{-(n+m)c})))\big)+C(2^{-nc}+\dotsb+2^{-(n+m)c})<\\
 & <\mathbb{P}_{G}^{2^{n+m}a}\big(\LE\subset2^{n+m}((\cap\mathcal{E}_{i})+B(C2^{-nc}/(1-2^{-c})))\big)\;+\\
 & \qquad+C2^{-nc}/(1-2^{-c})\leq\\
 & \leq\mathbb{P}_{G}^{2^{n+m}a}\big(\LE\subset2^{n+m}\cap(\mathcal{E}_{i}+B(C2^{-nc}))\big)+C2^{-nc}.\end{align*}
Taking $m$ to $\infty$ we get\[
\mathbb{P}_{G}^{2^{n}a}(\LE\subset2^{n}\cap\mathcal{E}_{i})\leq\underline{p}(\cap(\mathcal{E}_{i}+B(C2^{-nc})))+C2^{-nc}\]
and taking $n$ to $\infty$ gives\[
\overline{p}(\cap\mathcal{E}_{i})\leq\lim_{\epsilon\to0^{+}}\underline{p}(\cap(\mathcal{E}_{i}+B(\epsilon))).\]
Now, $\underline{p}(\cap(\mathcal{E}_{i}+B(\epsilon)))$ is a monotone
function hence it is continuous except at a countable number of points.
At each point $x$ of continuity we have\[
\overline{p}(\cap(\mathcal{E}_{i}+B(x)))\leq\lim_{\epsilon\to0^{+}}\underline{p}(\cap(\mathcal{E}_{i}+B(x+\epsilon)))=\underline{p}(\cap(\mathcal{E}_{i}+B(x)))\leq\overline{p}(\cap(\mathcal{E}_{i}+B(x)))\]
hence all the inequalities are equalities.
\end{proof}
\begin{rem*}
It is not very difficult to construct an example of an open set $\mathcal{E}$
(say with $G=\mathbb{Z}^{d}$, $\mathcal{D}=[\frac{1}{4},\frac{3}{4}]^{d}$
and $a=\vec{\frac{1}{2}}$) such that $\mathbb{P}_{n}(\mathcal{S}(\mathcal{E}))$
does not converge (construct $\mathcal{E}$ that it contains areas
$\mathcal{E}_{k}$ which become connected only for $n$ above some
$N_{k}$ and affect the probability significantly). 
\end{rem*}
\begin{lem}
\label{lem:OVOeps}Let $G,a$ and $\mathcal{D}$ be as in theorem
\ref{thm:scallim}. Let $\mathcal{O}\subset\mathcal{H}(\overline{\mathcal{D}})$
be an open set and let $\epsilon>0$. Then there exists an open set
$\mathcal{V}\subset\mathcal{H}(\overline{\mathcal{D}})$ with $\mathcal{O}\subset\mathcal{V}\subset\mathcal{O}+B(\epsilon)$
such that \begin{equation}
\lim_{n\to\infty}\mathbb{P}_{n}(\mathcal{V})\label{eq:lemOVOeps}\end{equation}
exists.
\end{lem}
\begin{proof}
Let $E\subset\overline{\mathcal{D}}$ be relatively open, let $v_{1},\dotsc,v_{k}\in\overline{\mathcal{D}}$
and let $a>b>0$. We define \[
\mathcal{Q}(E,v_{1},\dotsc,v_{k},a,b):=\mathcal{S}(E+B(a+b))\cap\bigcap_{i=1}^{k}\mathcal{I}(B(v_{i},a-b)).\]
Our main goal is to show that for $\mathcal{Q}_{1},\dotsc,\mathcal{Q}_{l}$
with a common $a$ and $b$,\[
\mathcal{Q}_{i}=\mathcal{Q}(E_{i},v_{1,i},\dotsc,v_{k_{i},i},a,b)\]
one has, for every $a$ and almost every $b$, that\begin{equation}
\lim_{n\to\infty}\mathbb{P}_{n}(\cap\mathcal{Q}_{i})\label{eq:capQi}\end{equation}
exists. Collect all the $v_{i,j}$-s into a single list, $v_{1},\dotsc,v_{m}$,
and let $I\subset\{1,\dotsc,m\}$ be arbitrary. We use lemma \ref{lem:montrik}
for $E_{1}+B(a),\dotsc,E_{l}+B(a)$ and for $\mathbb{R}^{d}\setminus B(v_{i},a)$
for all $i\in I$. We get that for almost every $0<b<a$, \begin{equation}
\lim_{n\to\infty}\mathbb{P}_{n}\Big(\mathcal{S}\Big(\Big(\bigcap_{i=1}^{k}E_{i}+B(a+b)\Big)\setminus\Big(\bigcup_{i\in I}B(v_{i},a-b)\Big)\Big)\Big)\label{eq:oneI}\end{equation}
exists. Now take some $b$ such that the limit (\ref{eq:oneI}) exists
for all $I$. Subtracting (\ref{eq:oneI}) for a given $I$ from (\ref{eq:oneI})
for $I=\emptyset$ gives that the limit\[
\lim_{n\to\infty}\mathbb{P}_{n}\Big(\mathcal{S}\Big(\bigcap_{i=1}^{k}E_{i}+B(a+b)\Big)\cap\bigcup_{i\in I}\mathcal{I}(B(v_{i},a-b))\Big)\]
exists. Since this holds for all $I$, we can use the inclusion-exclusion
principle to show that the limit\[
\lim_{n\to\infty}\mathbb{P}_{n}\Big(\mathcal{S}\Big(\bigcap_{i=1}^{k}E_{i}+B(a+b)\Big)\cap\bigcap_{i=1}^{m}\mathcal{I}(B(v_{i},a-b))\Big)\]
exists, which is equivalent to (\ref{eq:capQi}).

Proving the lemma is now easy. Take a finite set of $\left\{ K_{i}\right\} _{i=1}^{L}\in\mathcal{O}$
such that $B(K_{i},\frac{1}{4}\epsilon)$ cover $\mathcal{O}$ (this
is possible from the compactness of $\mathcal{H}$). For each $K_{i}$
one can take $E_{i}=K_{i}+B(\frac{1}{4}\epsilon)$ and $v_{1,i},\dotsc,v_{k_{i},i}$
to be a $\frac{1}{4}\epsilon$-net in $K_{i}$ and get that \[
B(K_{i},{\textstyle \frac{1}{4}}\epsilon)\subset\mathcal{Q}(E_{i},v_{1,i},\dotsc,v_{k_{i},i},{\textstyle \frac{1}{2}}\epsilon,b)\subset B(K_{i},\epsilon)\quad\forall0<b<{\textstyle \frac{1}{4}}\epsilon.\]
Denote this set by $\mathcal{Q}_{i}$. For every $I\subset\{1,\dotsc,L\}$
we use (\ref{eq:capQi}) to see that for almost every $0<b<\frac{1}{4}\epsilon$,
\[
\lim_{n\to\infty}\mathbb{P}_{n}\Big(\bigcap_{i\in I}\mathcal{Q}_{i}\Big)\]
exists. Take a $b$ such that the limit exists for all $I$. By the
inclusion-exclusion principle we get that\[
\lim_{n\to\infty}\mathbb{P}_{n}\Big(\bigcup_{i=1}^{L}\mathcal{Q}_{i}\Big)\]
exists. Define $\mathcal{V}=\cup\mathcal{Q}_{i}$ and the lemma is
finished.
\end{proof}

\begin{proof}
[Proof of theorem \ref{thm:scallim}]Let $f:\mathcal{H}(\overline{\mathcal{D}})\to[0,\infty)$
be a continuous function. Let $\epsilon>0$ and define\[
\mathcal{O}_{i}=f^{-1}[\epsilon i,\infty).\]
Then\[
\epsilon\sum_{i=1}^{\infty}\mathbf{1}_{\mathcal{O}_{i}}\leq f<\epsilon\Big(1+\sum_{i=1}^{\infty}\mathbf{1}_{\mathcal{O}_{i}}\Big).\]
By the compactness of $\mathcal{H}$ there exists some $\delta>0$
such that $\mathcal{O}_{i}+B(\delta)\subset\mathcal{O}_{i-1}$ for
all $i\geq1$ (note that $\mathcal{O}_{i}=\emptyset$ for $i$ sufficiently
large). For every $i\geq1$ such that $\mathcal{O}_{i}\neq\emptyset$
use lemma \ref{lem:OVOeps} to find a $\mathcal{O}_{i}\subset\mathcal{V}_{i}\subset\mathcal{O}_{i-1}$
such that the limit (\ref{eq:lemOVOeps}) exists (for larger $i$-s
define $\mathcal{V}_{i}=\emptyset$). We get\begin{align*}
\varlimsup_{n\to\infty}\int fd\mathbb{P}_{n} & \leq\varlimsup_{n\to\infty}\int\epsilon\Big(1+\sum_{i=1}^{\infty}\mathbf{1}_{\mathcal{V}_{i}}\Big)d\mathbb{P}_{n}=\varliminf_{n\to\infty}\int\epsilon\Big(1+\sum_{i=1}^{\infty}\mathbf{1}_{\mathcal{V}_{i}}\Big)d\mathbb{P}_{n}\\
 & \leq\varliminf_{n\to\infty}\int\epsilon\Big(1+\sum_{i=0}^{\infty}\mathbf{1}_{\mathcal{O}_{i}}\Big)d\mathbb{P}_{n}\leq2\epsilon+\varliminf_{n\to\infty}\int fd\mathbb{P}_{n}.\end{align*}
Since $\epsilon$ was arbitrary we see that the limit $\int fd\mathbb{P}_{n}$
exists for any positive $f$. Any function $f$ can be written as
$f^{+}-f^{-}$ hence we see that the limit exists for any continuous
$f$. This finishes the theorem since by compactness, if $\mathbb{P}_{n}$
does not converge it must have two subsequences converging to different
values, which is a contradiction.
\end{proof}
We now move to the proof of theorem \ref{thm:univ}.

\begin{lem}
\label{lem:SLtoeps}Let $G$ be a graph with a scaling limit and let
$a,\mathcal{D}$ and $\mathcal{E}_{1},\dotsc,\mathcal{E}_{k}$ be
as in lemma \ref{lem:montrik}. Then for almost every $\epsilon>0$,\[
\lim_{n\to\infty}\mathbb{P}_{n}(\mathcal{S}(\cap(\mathcal{E}_{i}+\overline{B(\epsilon)})))\]
exists.
\end{lem}
\begin{proof}
As in lemma \ref{lem:montrik} we may assume all $\mathcal{E}_{i}$
are open and replace $\overline{B}$ with $B$. Denote $\overline{p}(\mathcal{E})$
and $\underline{p}(\mathcal{E})$ as in lemma \ref{lem:montrik}.
If they are different on $\cap(\mathcal{E}_{i}+B(x))$ for an uncountable
number of $x$-s, then one of them would be a continuity point for
both $\overline{p}$ and $\underline{p}$. Hence we get an interval
$[\alpha,\beta]$ such that \[
\underline{p}(\cap(\mathcal{E}_{i}+B(x)))\leq A<B\leq\overline{p}(\cap(\mathcal{E}_{i}+B(x)))\quad\forall x\in[\alpha,\beta].\]
Define now\begin{align*}
\mu(K) & :=\min\{ x:K\subset\cap(\mathcal{E}_{i}+B(x))\}\quad K\in\mathcal{H}(\overline{\mathcal{D}})\\
f(K) & :=\begin{cases}
1 & \mu(K)\leq\alpha\\
0 & \mu(K)\geq\beta\\
\textrm{linear in }\mu(K) & \textrm{otherwise}.\end{cases}\end{align*}
It is easy to see that $f$ is continuous on $\mathcal{H}$. Therefore
by the definition of scaling limit, one must have that\[
\lim_{n\to\infty}\int fd\mathbb{P}_{n}\]
exists. However, \begin{align*}
\varliminf_{n\to\infty}\int fd\mathbb{P}_{n} & \leq\varliminf_{n\to\infty}\mathbb{P}_{n}(\mathcal{S}(\cap(\mathcal{E}_{i}+B(\beta))))\leq A\\
\varlimsup_{n\to\infty}\int fd\mathbb{P}_{n} & \geq\varlimsup_{n\to\infty}\mathbb{P}_{n}(\mathcal{S}(\cap(\mathcal{E}_{i}+B(\alpha))))\geq B\end{align*}
which is a contradiction.
\end{proof}
\begin{lem}
Let $G$ and $H$ be as in theorem \ref{thm:univ} and let $a,\mathcal{D},\mathcal{E}_{1},\dotsc,\mathcal{E}_{k}$
and $\epsilon$ be as in lemma \ref{lem:montrik}. Then for almost
all $\epsilon>0$,\[
\lim_{n\to\infty}\mathbb{P}_{G}^{2^{n}a}(2^{-n}\LE\subset\cap(\mathcal{E}_{i}+\overline{B(\epsilon)}))=\lim_{n\to\infty}\mathbb{P}_{H}^{2^{n}a}(2^{-n}\LE\subset\cap(\mathcal{E}_{i}+\overline{B(\epsilon)})).\]
In particular both limits exist.
\end{lem}
\begin{proof}
[Proof sketch]Denoting \[
\overline{p}_{G}(x):=\varlimsup_{n\to\infty}\mathbb{P}_{G}^{2^{n}a}(2^{-n}\LE\subset\cap(\mathcal{E}_{i}+B(x)))\]
a calculation analogous to that of lemma \ref{lem:montrik} gives
that\[
\underline{p}_{G}(x)\leq\lim_{\epsilon\to0^{+}}\underline{p}_{H}(x+\epsilon)\quad\overline{p}_{H}(x)\leq\lim_{\epsilon\to0^{+}}\overline{p}_{G}(x+\epsilon)\]
hence the lemma holds for every $x$ which is a continuity for all
$\overline{p}$ and $\underline{p}$ and for which $\underline{p}_{G}(x)=\overline{p}_{G}(x)$.
Lemma \ref{lem:SLtoeps} is used here to show that this last condition
holds almost everywhere.
\end{proof}
\begin{lem}
Let $G$ and $H$ be as in theorem \ref{thm:univ} and let $\mathcal{O}$
and $\epsilon$ be as in lemma \ref{lem:OVOeps}. Then there exists
an $\mathcal{O}\subset\mathcal{V}\subset\mathcal{O}+B(\epsilon)$
such that\[
\lim_{n\to\infty}\mathbb{P}_{G}^{2^{n}a}(2^{-n}\LE\in\mathcal{V})=\lim_{n\to\infty}\mathbb{P}_{H}^{2^{n}a}(2^{-n}\LE\in\mathcal{V}).\]

\end{lem}
The proof of this lemma is a complete analogue of the proof of lemma
\ref{lem:OVOeps}. Concluding theorem \ref{thm:univ} from it is similar
to the conclusion of theorem \ref{thm:scallim} from lemma \ref{lem:OVOeps}.
We shall omit both.

\section{\label{sec:Examples}Examples}

We start with a lemma on continuous functions

\begin{lem}
\label{lem:harmdiff}Let $f$ be harmonic on $B(0,1)$ and continuous
on $\overline{B(0,1)}$, and assume that on $\partial B(0,1)$ $f$
is $C^{2}$ and satisfies\[
||f||_{\infty}\leq1,\quad||f||_{1,\infty}\leq N\quad||f||_{2,\infty}\leq N^{2}\quad\forall i,j\]
for some $N\geq1$. Then $|\partial f/\partial x_{i}|\leq CN$ inside
$B(0,1)$.
\end{lem}
The notation $||f||_{k,p}$ stands for the Sobolev norm, which in
this case simply means $||f||_{1,\infty}=\max_{x,i}|(\partial f/\partial\theta_{i})(x)|$
and similarly for second derivatives.

\begin{proof}
The value of $f$ inside $B(0,1)$ is related to the value on the
boundary by the Poisson kernel,\[
f(x)=\int_{\partial B(0,1)}f(y)\frac{1-||x||^{2}}{||x-y||^{d}}\, dy\]
where $dy$ is the surface area measure on $\partial B(0,1)$ normalized
to be a probability measure. This immediately gives an estimate in
the ball $B(0,\frac{1}{2})$ since \[
\frac{\partial f}{\partial x_{i}}=\int_{\partial B(0,1)}f(y)\frac{\partial}{\partial x_{i}}\frac{1-||x||^{2}}{||x-y||^{d}}\, dy\leq C\int_{\partial B(0,1)}f(y)\, dy\leq C.\]
Hence we need to estimate the derivatives only near the boundary.

Let now $x\in B(0,1)\setminus B(0,\frac{1}{2})$, let $\theta$ be
some direction on the sphere and let $T_{\epsilon}$ be a rotation
by $\epsilon$ around the pole orthogonal to $x$ and $x+\epsilon\theta$
in the direction $\theta$. Then the rotational invariance of the
Poisson kernel allows to write\begin{align}
\frac{\partial f}{\partial\theta}(x) & =\lim_{\epsilon\to0}\int_{\partial B(0,1)}\frac{f(y)-f(T_{\epsilon}y)}{\epsilon}\cdot\frac{1-||x||^{2}}{||x-y||^{2}}\, dy\leq\nonumber \\
 & \leq\int_{\partial B(0,1)}N\frac{1-||x||^{2}}{||x-y||^{2}}\, dy=N.\label{eq:tangent}\end{align}
Since $||x||\geq\frac{1}{2}$ we get that the tangent derivatives
are $\leq2N$. Therefore the lemma will be proved once we get a similar
estimate for the radial derivative.

A calculation similar to (\ref{eq:tangent}) shows that $|\partial^{2}f/\partial x_{i}^{2}|\leq4N^{2}$
for $x_{i}$ a tangent direction. Hence by the harmonicity of $f$
we get $|\partial^{2}f/\partial r^{2}|\leq4(d-1)N^{2}$. Examine now
a point $x\in\partial B(0,1)$. Since $|f(x)|\leq1$ and $|f(x(1-1/N))|\leq1$
there must be a point $y$ on the interval $[x,x(1-1/N)]$ such that
$|(\partial f/\partial r)(y)|\leq2N.$ With the bound on $\partial^{2}f/\partial r^{2}$
this gives that $|(\partial f/\partial r)(x)|\leq CN$. This allows
to bound $\partial f/\partial r$ everywhere since\begin{align*}
\frac{\partial f}{\partial r}(x) & =\lim_{\epsilon\to0^{+}}\int_{\partial B(0,1)}\frac{f(y)-f(y-\epsilon/||x||)}{\epsilon}\cdot\frac{1-||x||^{2}}{||x-y||^{2}}\, dy=\\
 & =\frac{1}{||x||}\int_{\partial B(0,1)}\frac{\partial f}{\partial r}(y)\cdot\frac{1-||x||^{2}}{||x-y||^{2}}\, dy\leq\frac{CN}{||x||}\leq CN.\qedhere\end{align*}

\end{proof}
Our purpose at this point is to prove that $\mathbb{Z}^{d}$ and $2\mathbb{Z}^{d}$
have an isometric interpolation, which will allow to invoke theorem
\ref{thm:scallim}. We start with the case $d=3$ and for which we
weight $2\mathbb{Z}^{3}$ by $2$ i.e. assume all the edges have weight
$2$ (the statement of theorem \ref{thm:scallim}, page \pageref{thm:scallim},
allows us to do so). We shall show this for $\alpha=\frac{1}{9}$
(the $\alpha$ from the definition of isotropic interpolation, (\ref{eq:iirMreq})
on page \pageref{eq:iirMreq}). Therefore assume $L>6$ and $M\leq L^{1/9}$
($6$ is our choice for the $C$ from (\ref{eq:iirMreq})). Let $\xi\in\{1,2\}^{M^{3}}$.
We need to construct a graph $G=G(L,M,\xi)$. We do it as follows.

\begin{lyxlist}{00.00.0000}
\item [Vertices]For every $(x,y,z)\in\{0,\dotsc,M-1\}^{3}$ such that $\xi(x,y,z)=1$
we take every point of $\mathbb{Z}^{3}\cap[Lx,Lx+L)\times[Ly,Ly+L)\times[Lz,Lz+L)$
to be a vertex of $G$. We call such vertices {}``vertices of type
$1$''. If $\xi(x,y,z)=2$ we take every point of $2\mathbb{Z}^{3}\cap[Lx,Lx+L)\times[Ly,Ly+L)\times[Lz,Lz+L)$
to be a vertex of $G$ and call it a vertex of type $2$. Outside
$[0,LM)^{3}$ we choose between $\mathbb{Z}^{3}$ and $2\mathbb{Z}^{3}$
by a majority vote on $\xi$ (or in any other way that gives $1$
for $\xi\equiv1$ and $2$ for $\xi\equiv2$).
\item [Edges]Two vertices of type $1$ will have an edge if and only if
their distance is $1$, and in this case the edge will have weight
$1$. Two vertices of type $2$ will have an edge if and only if their
distance is $2$ and in this case the edge will have weight $2$.
If $v$ is of type $1$ and $w$ of type $2$ then we need to find
the vertex $x$ of type $1$ closest to $w$ (usually $|x-w|$ will
be either $1$ or $2$); and define the weight by\begin{equation}
\omega(v,w)=\frac{1}{|w-x|}\begin{cases}
1 & v=x\\
1/2 & |v-x|=1\\
1/4 & |v-x|=\sqrt{2}\\
0 & \textrm{otherwise}\end{cases}\label{eq:defomeclose}\end{equation}
($0$ here means no edge). See figure \ref{cap:type12}.%
\begin{figure}
\input{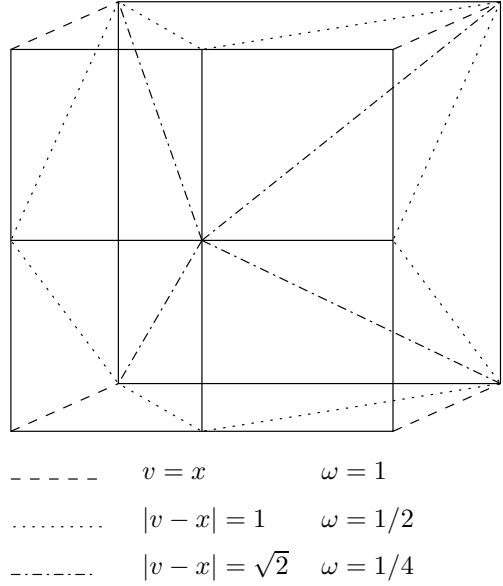}

\caption{\label{cap:type12}Lower-left square are vertices of type 1, upper-right
are vertices of type 2. The weights $\omega$ are under the assumption
that $|w-x|=1$; otherwise divide all weights by $2$.}
\end{figure}

\end{lyxlist}
Requirements \ref{enu:isinI} and \ref{enu:isinII} from the definition
of an isotropic interpolation are obvious. Hence the only thing we
need to show is that these graphs are isotropic. It is clear that
$G$ is always a Euclidean net so we need to verify the more delicate
condition on exit probabilities in the definition of an isotropic
graph (page \pageref{sub:DefinitionIsotropic}). This will follow
from a comparison of the continuous and discrete Laplacian. Let therefore
$f$ be a continuously harmonic function in a ball of radius $3$
around some vertex $v$. Let $\Delta$ be the discrete Laplacian.
Then

\begin{enumerate}
\item \label{enu:regpnt}If $d(v,X)>2$ for any square $X$ of the form
$\{ Ln\}\times[0,LM]\times[0,LM]$, $[0,LM]\times\{ Ln\}\times[0,LM]$
or $[0,LM]\times[0,LM]\times\{ Ln\}$ ($n\in\{0,\dotsc,M\}$), then
$(\Delta f)(v)\leq C||f||_{4,\infty}$, the Sobolev norm being on
the ball $B(v,3)$. In a ball of radius $r$ there are no more than
$C(r^{3}/L+r^{2})$ points not satisfying this condition.
\item \label{enu:plnpnt}If $d(v,X)>2$ for any segment $X$ of the form
$\{ Ln\}\times\{ Lm\}\times[0,LM]$, $\{ Ln\}\times[0,LM]\times\{ Lm\}$
or $[0,LM]\times\{ Ln\}\times\{ Lm\}$ ($n,m\in\{0,\dotsc,M\}$),
then $(\Delta f)(v)\leq C||f||_{2,\infty}$. In a ball of radius $r$
there are no more than $C(r^{3}/L^{2}+r)$ points not satisfying this
condition.
\item \label{enu:linepnt}For any $v$, $(\Delta f)(v)\leq C||f||_{1,\infty}$.
\end{enumerate}
Seeing \ref{enu:regpnt}-\ref{enu:linepnt} is not difficult. Write
a Taylor expansion for $f$ around $v$ of order $4$, $2$ or $1$
respectively, write\[
(\Delta f)(v)=\sum_{w\sim v}\omega(v,w)(f(w)-f(v))\]
and calculate and see that all terms except the error term vanish.
Note that case \ref{enu:plnpnt} is the one where all the fancy weights
in the {}``stitching'' between $\mathbb{Z}^{3}$ and $2\mathbb{Z}^{3}$
are needed, and also the reason we had to take $2\mathbb{Z}^{3}$
to have all the weights $2$. In fact \ref{enu:linepnt} holds in
any graph with bounded degree.

\begin{lem}
\label{lem:i2iiiiso}A $3$-dimensional Euclidean net satisfying \ref{enu:regpnt}-\ref{enu:linepnt}
is isotropic. 
\end{lem}
\begin{proof}
By the definition of an isotropic graph, we need to take some $v\in G$,
some $r>0$, and a spherical triangle $A\subset\partial_{\textrm{cont}}B(v,r)$.
We may assume that $r>1$ since otherwise (\ref{eq:defbrown}) holds
automatically for $K\geq2$. Let $A^{-}=\{ x\in A:d(x,\partial A)\geq r^{4/5}\}$
and $A^{+}:=\{ x\in\partial_{\textrm{cont}}B(v,r):d(x,A)\leq r^{4/5}\}$.
$A^{-}$ could be empty and $A^{+}$could be all of $\partial_{\textrm{cont}}B(v,r)$
but anyway we always have \[
\int\mathbf{1}_{A^{-}}+Cr^{-1/5}\geq\int\mathbf{1}_{A}\geq\int\mathbf{1}_{A^{+}}-Cr^{-1/5}\]
where $\int$ here is with respect to the surface area measure on
$\partial_{\textrm{cont}}B(v,r)$, normalized to have total area $1$.
The first step is to find two $C^{5}$ functions $f^{-}$ and $f^{+}$
on $\partial_{\textrm{cont}}B(v,r)$ such that $||f^{-}||_{k,\infty},||f^{+}||_{k,\infty}\leq Cr^{-4k/5}$
for $k=1,\dotsc,5$ and such that\begin{equation}
f^{-}\leq\mathbf{1}_{A^{-}}\leq\mathbf{1}_{A^{+}}\leq f^{+},\quad\int f^{-}+Cr^{-1/5}\geq\int\mathbf{1}_{A}\geq\int f^{+}-Cr^{-1/5}.\label{eq:f1Ag}\end{equation}
It is easy to construct such an $f^{\pm}$. For example, start with
a spherically symmetric function $\eta$ which is $1$ in a spherical
cap of radius $r^{4/5}$ and supported in a spherical cap of radius
$2r^{4/5}$, $||\eta||_{k,\infty}\leq Cr^{-4k/5}$. Cover the sphere
with a locally finite family of translations $\{ T_{1},\dotsc,T_{m}\}$
of the $r^{4/5}$ cap (so that $\sum_{j}T_{j}(\eta)\geq1$), and define
$\nu_{i}:=T_{i}(\eta)/\sum_{j}T_{j}(\eta)$ so that the $\nu_{i}$
form a division of unity. Define $f^{-}$ to be the sum of all $\nu_{i}$
supported inside $A^{-}$ and $f^{+}$ to be the sum of all $\nu_{i}$
such that $\supp\nu_{i}\cap A^{+}\neq\emptyset$. Verifying all the
properties of $f^{-}$ and $f^{+}$ is easy.

We wish to discretize (\ref{eq:f1Ag}). Let $A^{*}\subset\partial B(v,r)$
be a discrete version of $A$ as in the definition of isotropic graphs.
We want to find functions $F^{-}$ and $F^{+}$ satisfying (\ref{eq:f1Ag})
(with the integral being with respect to the discrete harmonic measure
of $\partial B(v,r)$ starting from $v$). We shall only show the
construction of $F^{-}$ --- the construction of $F^{+}$ is identical. 

Stretch $f^{-}$ to $\partial B(v,r+\lambda)$ where $\lambda$ is
some constant such that $\partial B(v,r)\subset B(v,r+\lambda)$ ---
for $\mathbb{Z}^{3}$ and $2\mathbb{Z}^{3}$ we can take $\lambda=3$.
Extend $f^{-}$ to a harmonic function on $B(v,r+\lambda)$ continuous
on $\overline{B(v,r+\lambda)}$. Call this extension $g^{1}$ and
notice that $||g^{1}||_{k,\infty}\leq Cr^{-4k/5}$ for $k=1,\dotsc,4$
(this follows from using lemma \ref{lem:harmdiff} for rescaled versions
of $g^{1}$ and its derivatives). Let $a(w,x):=-G(w,x;\overline{B(v,r)})$
be the discrete Green function. Define the following {}``correction''
for $g^{1}$:\[
D^{1}(w):=\sum_{x\in B(v,r)}(\Delta g^{1})(x)a(w,x).\]
Because $\Delta a(\cdot,x)$ is a delta function at $x$, we get that
$g^{2}:=g^{1}-D^{1}$ is discretely harmonic on $B(v,r)$. We also
note that $g^{2}\equiv g^{1}$ on $\partial B(v,r)$. What we need
is to estimate $D^{1}$ at $v$. Recall the estimate $a(v,w)\leq C|v-w|^{-1}$
from lemma \ref{lem:a} (\ref{eq:atrans}). Summing on points of type
\ref{enu:regpnt} we get that \begin{align*}
\sum_{x\textrm{ of type \ref{enu:regpnt}}}|(\Delta g^{1})(x)a(v,x)| & \leq Cr^{-16/5}\sum_{s=0}^{\left\lfloor \log_{2}r\right\rfloor }\sum_{2^{s}\leq|v-x|<2^{s+1}}a(v,x)\leq\\
 & \leq Cr^{-16/5}\sum_{s=0}^{\left\lfloor \log_{2}r\right\rfloor }4^{s}\leq Cr^{-6/5}.\end{align*}
Points of type \ref{enu:plnpnt} have a worse estimate, but are fewer.
In particular, a ball of radius $2^{s}$ around $v$ will contain
no more than $C\min(8^{s}/L+4^{s},L^{2}M^{3})$ such points. Assume
first that $r\leq LM$ which implies $L\geq r^{9/10}$. We get\begin{eqnarray*}
\lefteqn{{\sum_{x\textrm{ of type \ref{enu:plnpnt}}}|(\Delta g^{1})(x)a(v,x)|\leq Cr^{-8/5}\sum_{s=0}^{\left\lfloor \log_{2}r\right\rfloor }\sum_{2^{s}\leq|v-x|<2^{s+1}}a(v,x)\leq}}\\
 &  & \qquad\leq Cr^{-8/5}\sum_{s=0}^{\left\lfloor \log_{2}r\right\rfloor }2^{s}+\frac{4^{s}}{L}\leq C(r^{-3/5}+r^{2/5}/L)\leq Cr^{-1/2}.\end{eqnarray*}
If $r>LM$ (which implies $M\leq r^{1/10})$ a similar calculation
shows that \[
\sum\leq Cr^{-8/5}\Big(LM^{2}+\sum_{s=\left\lfloor \log_{2}LM\right\rfloor }^{\left\lfloor \log_{2}r\right\rfloor }L^{2}M^{3}2^{-s}\Big)\leq Cr^{-8/5}LM^{2}\leq Cr^{-1/2}.\]
Finally, for points of type \ref{enu:linepnt} we get that a ball
of radius $2^{s}$ will contain no more than $C\min(8^{s}/L^{2}+2^{s},LM^{3})$
and an identical calculation shows that\[
\sum_{w\textrm{ of type \ref{enu:linepnt}}}|(\Delta g^{1})(v)a(v,w)|\leq Cr^{-4/5}\log r+Cr^{-3/5}.\]
Summing these three terms we get $|D^{1}(v)|\leq Cr^{-1/2}$. 

The only reason not to use $g^{2}$ directly is that $g^{2}\not\leq\mathbf{1}_{A^{*}}$
on $\partial B(v,r)$. We do have, on $\partial B(v,r)$ that $g^{2}\equiv g^{1}\leq1$,
but we need to correct $g^{2}$ to be $0$ outside $A^{*}$. Let $w\in\partial B(v,r)\setminus A^{*}$
and let $w'=w\frac{r}{r+\lambda}$. Then\begin{align*}
g^{2}(w) & =g^{1}(w)=f^{-}(w')=r^{-1}\int_{\partial_{\textrm{cont}}B(v,r)}f^{-}(x)\frac{r^{2}-||wr'||^{2}}{||w'-x||^{3}}\, dx\leq\\
 & \leq r^{-1}\int_{A^{-}}\frac{r^{2}-||w'||^{2}}{||w'-x||^{3}}\, dx\leq C/d(w',A^{-})\leq Cr^{-4/5}.\end{align*}
Hence defining a second correction $D^{2}$ on $\partial B(v,r)$\[
D^{2}(w)=\begin{cases}
f^{2}(w) & w\not\in A^{*}\\
0 & \textrm{otherwise}\end{cases}\]
 and extending it to a discretely harmonic function on $B(v,r)$,
we get, from the discrete maximum principle that $D^{2}(v)\leq Cr^{-4/5}$.

Defining $F^{-}=g^{2}-D^{2}$ the lemma is now easy: By the definition
of the harmonic measure, $\int F^{-}=F^{-}(v)$ (discrete case) and
$\int f^{-}=f^{-}(v)$ (continuous case). $|F^{-}(v)-f^{-}(v)|\leq|D^{1}(v)|+|D^{2}(v)|\leq Cr^{-1/2}$.
Therefore\begin{align*}
p_{A^{*}} & =\int\mathbf{1}_{A^{*}}\geq\int F^{-}=F^{-}(v)\geq f^{-}(v)-Cr^{-1/2}=\int f^{-}-Cr^{-1/2}\geq\\
 & \!\!\stackrel{(\ref{eq:f1Ag})}{\geq}\int\mathbf{1}_{A}-Cr^{-1/5}=|A|-Cr^{-1/5}.\end{align*}
A similar calculation with the similarly defined $F^{+}$ will show
that $p_{A}\leq|A|+Cr^{-1/5}$ so (\ref{eq:defbrown}) is proved,
$G$ is isotropic and the lemma is proved.
\end{proof}
\begin{conclusion*}
Lemma \ref{lem:i2iiiiso}, the discussion before that and theorem
\ref{thm:scallim} together show that $\mathbb{Z}^{3}$ has a scaling
limit. Hence theorem \ref{thm:Z3scal} is proved in three dimensions.
\end{conclusion*}
For the case of $\mathbb{Z}^{2}$, we define $2\mathbb{Z}^{2}$ to
have all the weights $1$. We construct that graphs $G(L,M,\xi)$
equivalently, but with the weights defined as follows: If $v$ and
$w$ are of the same type then the edges between them are as in the
three dimensional cases and all the weights are $1$. If $v$ is of
type $1$ and $w$ is of type $2$ we again find the vertex $x$ of
type $1$ closest to $w$, and then define\[
\omega(v,w)=\frac{1}{|w-x|}\begin{cases}
1 & v=x\\
1/2 & |v-x|=1\\
0 & \textrm{otherwise}\end{cases}\]
(this is the equivalent of (\ref{eq:defomeclose})). A figure can
be found in \cite[figure 1, page 10]{K}. The equivalents of \ref{enu:regpnt}-\ref{enu:linepnt}
from page \pageref{enu:regpnt} are

\begin{enumerate}
\item If $d(v,X)>2$ for any segment $X$ of the form $\{ Ln\}\times[0,LM]$
or $[0,LM]\times\{ Ln\}$ ($n\in\{0,\dotsc,M\}$), then $(\Delta f)(v)\leq C||f||_{4,\infty}$.
In a ball of radius $r$ there are no more than $C(r^{2}/L+r)$ points
not satisfying this condition.
\item If $d(v,X)>2$ for any point $X$ of the form $\{ Ln\}\times\{ Lm\}$
($n,m\in\{0,\dotsc,\linebreak[4]M\}$), then $(\Delta f)(v)\leq C||f||_{2,\infty}$.
In a ball of radius $r$ there are no more than $C(r^{2}/L^{2}+1)$
points not satisfying this condition.
\item For any $v$, $(\Delta f)(v)\leq C||f||_{1,\infty}$.
\end{enumerate}
Which are easy to verify, and as in the three dimensional case, the
definition of $\omega$ at the stitches is used only for \ref{enu:plnpnt}.
The equivalent of lemma \ref{lem:i2iiiiso} is proved in the same
way (indeed, it is simpler as the construction of $f$ and $g$ is
easier; and since the estimate $a(v,w)\approx\log|v-w|$ means there
is no reason to divide into shells of size $2^{s}$ as in the three
dimensional case). This shows that $\mathbb{Z}^{2}$ has a scaling
limit, and concludes theorem \ref{thm:Z3scal}.

\begin{conjecture*}
Any non-trivial stitching of $\mathbb{Z}^{d}$ and $2\mathbb{Z}^{d}$
is an isotropic graph.
\end{conjecture*}
In other words, while condition \ref{enu:plnpnt} obviously does not
hold unless we insert weights in a manner similar to (\ref{eq:defomeclose}),
the conjecture states that the graph $G(L,M,\xi)$ would be isotropic
even if we, for example, just connect every vertex of type $2$ to
the nearest vertex of type $1$ and give the edge weight $18\frac{7}{10}$.
Be forewarned that the weighting of $2\mathbb{Z}^{d}$ by $2$ (in
the three dimensional case) is very much needed. A simple resistance
calculation would show that, for example for $\xi\equiv1$ on $[0,M/2]\times[0,M-1]^{2}$
and $2$ otherwise, the graph $G(L,M,\xi)$ can never be isotropic
(no matter what you put in the connecting layer) unless the edges
of length $2$ are weighted by $2$.

\subsection{\label{sub:Invariance}Invariance}

The definition of the scaling limit immediately implies that it is
invariant to multiplication by $2$, meaning that if $\mu(\mathcal{D},a)$
is the scaling limit of loop-erased random walk on $\mathbb{Z}^{3}\cap2^{n}\mathcal{D}$
starting from $2^{n}a$ then $\mu(2\mathcal{D},2a)=2\mu(\mathcal{D},a)$
where the notation {}``$2\mu$'' stands for a stretching of $\mu$
by $2$ in the natural way. 

In this section we give a few examples of additional invariances $\mu$
satisfies. The first example is multiplication by $3$. In other words,
we want to show that $\mu(3\mathcal{D},3a)=3\mu(\mathcal{D},a)$.
A moments reflection shows that this will follow if we show that $\mathbb{Z}^{3}$
and $3\mathbb{Z}^{3}$ have the same scaling limit, which would follow
by theorem \ref{thm:univ} if we show they have an isotropic interpolation.
We follow the same guidelines as in the previous section. Define $G(L,M,\xi)$
as $\mathbb{Z}^{3}$ and $3\mathbb{Z}^{3}$ with weight $3$ on the
internals of the cubes. In the stitches we let $x$ be the closest
vertex of type $1$ to $w$ and then define \begin{align*}
\omega(v,w) & =\frac{1}{|x-w|}\begin{cases}
1 & v=x\\
2/3 & |v-x|=1\textrm{ or }\sqrt{2}\\
1/3 & |v-x|=2\textrm{ or }\sqrt{8}.\end{cases}\end{align*}
As in the previous section, a calculation verifies \ref{enu:plnpnt}
which implies lemma \ref{lem:i2iiiiso}, isotropic interpolation and,
with theorem \ref{thm:univ}, the invariance of $\mu$.

Since (as is well known) the numbers $2^{k}3^{-n}$ are dense in $\mathbb{R}$,
this shows that $\mu$ is in fact invariant to a dense set of multiplications.
We shall now sketch a simple continuity argument which shows that
$\mu$ is in fact invariant to all multiplications. Let $\alpha>1$
and examine the situation of theorem \ref{thm:interp}, i.e.~we have
the graphs $\mathbb{Z}^{3}$ and $\alpha\mathbb{Z}^{3}$, some $a\in\mathcal{E}\cap\mathcal{D}$
and some $s>0$. We use theorem \ref{thm:interp} for $\mathbb{Z}^{3}$
and $2\mathbb{Z}^{3}$ repeatedly and rescale (as in the proof of
lemma \ref{lem:montrik}) to get that for any $k$,\[
\mathbb{P}_{\mathbb{Z}^{3}}^{2^{k}sa}(\LE(R)\subset2^{k}s(\mathcal{E}+B(0,Cs^{-c})))>\mathbb{P}_{\mathbb{Z}^{3}}^{sa}(\LE(R)\subset s\mathcal{E})-Cs^{-c}.\]
Next we use theorem \ref{thm:interp} for $\mathbb{Z}^{3}$ and $3\mathbb{Z}^{3}$
repeatedly and rescale in the opposite direction to get that, as long
as $2^{k}3^{-n}>1$,\[
\mathbb{P}_{\mathbb{Z}^{3}}^{2^{k}3^{-n}sa}(\LE(R)\subset2^{k}3^{-n}s(\mathcal{E}+Cs^{-c}B(0,1)))>\mathbb{P}_{\mathbb{Z}^{3}}^{sa}(\LE(R)\subset s\mathcal{E})-Cs^{-c}.\]
Notice that the various constants do not depend on $k$ and $n$ ---
in fact we can get as close as we want to $\alpha$, until $2^{k}3^{-n}\mathcal{D}\cap\mathbb{Z}^{3}=\alpha\mathcal{D}\cap\mathbb{Z}^{3}$
(and ditto for $\mathcal{E}$), with no price to pay. There might
be a problem that the point of $2^{k}3^{-n}\mathbb{Z}^{3}$ closest
to $a$ might be different by $C$ from the point of $\alpha\mathcal{D}\cap\mathbb{Z}^{3}$
closest to $a$, but we have lemma \ref{lem:contsp} to show us that
this affects the relevant probabilities by no more than $Cs^{-c}$
as well. Thus the conclusion of theorem \ref{thm:interp} holds for
$\mathbb{Z}^{3}$ and $\alpha\mathbb{Z}^{3}$ and therefore so does
theorem \ref{thm:univ}.

We next move to rotations. Examine the lattice $G$ in $\mathbb{R}^{3}$
spanned by the vectors \[
(4,3,0),(3,-4,0),(0,0,5).\]
Since the vectors are orthogonal and of equal length, condition \ref{enu:regpnt}
will continue to hold. By now the reader should have only technical
difficulties in producing a stitching of $G$ and, for example, $25\mathbb{Z}^{3}$
with the weights being $5$, which will satisfy condition \ref{enu:plnpnt}:
in fact our lattice is invariant to translations by $25\mathbb{Z}^{3}$
which reduces the verification of \ref{enu:plnpnt} to a small number
of cases. Therefore $G$ and $25\mathbb{Z}^{3}$ have an isotropic
interpolation and the same scaling limit. This shows that the scaling
limit is invariant to rotations by $\arctan\frac{3}{4}$ around the
$z$ axis. As is well known and not difficult to see, $\arctan\frac{3}{4}/\pi$
is irrational so the semigroup created by this rotation is dense,
and a similar continuity argument can be used to show that $\mu$
is invariant to any rotation around the $z$ axis. Since $\mu$ is
clearly invariant to a change of coordinates, and since any rotation
is a combinations of three rotations around the axes, we see that
$\mu$ is invariant to all rotations.

\end{document}